\documentclass[letterpaper,leqno,10pt,twoside]{amsart}
\usepackage{amsmath,amsthm,amsfonts,amssymb,euscript,mathrsfs}
\usepackage{mathtools}
\usepackage{tikz}
\usetikzlibrary{arrows.meta}
\usepackage[margin=1in]{geometry}
\usepackage[toc,page]{appendix}
\usepackage[shortlabels]{enumitem}
\usepackage{booktabs,array}
\usepackage{indentfirst}
\usepackage{float}

\allowdisplaybreaks

\numberwithin{equation}{section}

\usepackage[dvipsnames]{xcolor}

\usepackage{hyperref}
\hypersetup{colorlinks=true, pdfstartview=Fit,linkcolor=blue!70!black,citecolor=red!70!black, urlcolor=green!60!black}

\newtheorem{theorem}{Theorem}[section]
\newtheorem{lemma}[theorem]{Lemma}
\newtheorem{proposition}[theorem]{Proposition}
\newtheorem{corollary}[theorem]{Corollary}

\newtheorem*{restatement}{Proposition~\ref*{prop:auto-smoothing} (restated)}

\theoremstyle{definition}
\newtheorem{definition}{Definition}[section]
\theoremstyle{remark}
\newtheorem{remark}{Remark}[section]

\def\ue{e}%{\mathrm{e}}
\def\ui{\mathrm{i}}

\DeclarePairedDelimiter{\abs}{\lvert}{\rvert}
\DeclarePairedDelimiter{\babs}{\big\lvert}{\big\rvert}
\DeclarePairedDelimiter{\Babs}{\Big\lvert}{\Big\rvert}
\DeclarePairedDelimiter{\norm}{\lVert}{\rVert}
\DeclarePairedDelimiter{\bnorm}{\big\lVert}{\big\rVert}
\DeclarePairedDelimiter{\Bnorm}{\Big\lVert}{\Big\rVert}
\DeclarePairedDelimiterX{\pair}[2]{\langle}{\rangle}{#1,#2}
\DeclarePairedDelimiterX{\bpair}[2]{\big\langle}{\big\rangle}{#1,#2}
\DeclarePairedDelimiterX{\ip}[2]{(}{)}{#1,#2}
\newcommand{\R}{\mathbb R}
\newcommand{\cS}{\mathcal S}
\newcommand{\J}{\mathcal J}
\newcommand{\cQ}{\mathcal Q}
\newcommand{\cK}{\mathcal K}
\newcommand{\cM}{\mathcal M}
\newcommand{\cW}{\mathcal W}
\newcommand{\bL}{\mathbf L}
\newcommand{\bR}{\mathbf R}
\newcommand{\bG}{\mathbf G}
\newcommand{\mg}{\mathfrak g}
\newcommand{\Rf}{R^0}
\newcommand{\Lu}{L^2_u}
\newcommand{\loc}{\mathrm{loc}}
\DeclareMathOperator{\Dom}{Dom}
\newcommand{\sech}{\operatorname{sech}}
\DeclareMathOperator{\supp}{supp}
\newcommand{\dd}{\,\mathrm{d}}%{\,d}

\newcommand{\proofstep}[1]{\paragraph{\underline{#1}}}

\title[Unconditional Uniqueness for KdV and mKdV Green Solutions]{Unconditional Uniqueness of KdV and mKdV Green Solutions at Low Regularity}
\author{Zhimeng Ouyang}
\address{School of Mathematical Sciences, Peking University, Beijing 100871, China}
\email{ouyangzm@math.pku.edu.cn}

\subjclass[2020]{35Q53, 35B30, 35B65, 47A10}
\keywords{Korteweg--de Vries equation, modified Korteweg--de Vries equation, unconditional uniqueness, Green solution, low regularity, local smoothing, Gardner transform, Miura map}
\date{}

\begin{document}
\raggedbottom

\begingroup
\emergencystretch=1em
\begin{abstract}
We prove unconditional uniqueness for the Korteweg--de Vries (KdV) equation in $C_tH^{-1}_x(\mathbb R)$ and the real defocusing modified KdV equation in $C_tL^2_x(\mathbb R)$ within the class of Green solutions, 
an intrinsically well defined solution notion that encodes the evolution at these low regularities through the diagonal Green functions of the associated Lax operators.
% defined intrinsically at these low regularities by evolution identities for the diagonal Green functions of the associated Lax operators. 
A key contribution is automatic local smoothing: the Green identities imply $q\in L^2_{\mathrm{loc}}$ for KdV and $u_x,u^2\in L^2_{\mathrm{loc}}$ for mKdV in spacetime, with estimates uniform over spatial translates. The uniqueness results therefore identify every Green solution in these classes with the continuous extension of the corresponding Schwartz flow.
\end{abstract}
\endgroup

\maketitle

\setcounter{tocdepth}{1}%{2}
\makeatletter
\def\l@subsection{\@tocline{2}{0pt}{2.5pc}{5pc}{}}
\makeatother
\tableofcontents

% \newpage

%%%%%%%%%%%%%%%%%%%%%%%%%%%%%%%%%%%%%%%%%%%%%%%%%%%%%%%%%%%%%%%%%%%%%%%%%%%%%%%%%%
\section{Introduction}
%%%%%%%%%%%%%%%%%%%%%%%%%%%%%%%%%%%%%%%%%%%%%%%%%%%%%%%%%%%%%%%%%%%%%%%%%%%%%%%%%%

%%%%%%%%%%%%%%%%%%%%%%%%%%%%%%%%%%%%%%%%%%%%%%%%%%%%%%%%%%%%%%%%%%%%%%%%%%%%%%%%%%
\subsection{KdV Formulation}
%%%%%%%%%%%%%%%%%%%%%%%%%%%%%%%%%%%%%%%%%%%%%%%%%%%%%%%%%%%%%%%%%%%%%%%%%%%%%%%%%%

We study the Korteweg--de Vries (KdV) initial-value problem on the line. Fix $T>0$, set $I=[-T,T]$, and consider 

{
\begin{align} \tag{KdV}
 \label{eq:KdV} \begin{aligned}
 &q_t=-q_{xxx}+6qq_x=-q_{xxx}+3\partial_x(q^2), \qquad (t,x)\in I\times\R,\\
 &q(0,\cdot)=q_0\in H^{-1}(\R;\R).
 \end{aligned}
\end{align}
} The identity of \eqref{eq:KdV} has its usual classical meaning for smooth solutions. At the sharp Sobolev threshold $H^{-1}(\R)$ for global well-posedness \cite{Killip.Visan2019,Molinet2011}, the product $q^2$ need not have a canonical distributional meaning.
Consequently, \eqref{eq:KdV} cannot by itself be interpreted as a distributional evolution equation for an arbitrary curve $q\in C(I;H^{-1})$.

Killip--Vi\c{s}an \cite{Killip.Visan2019} constructed a global continuous extension of the Schwartz KdV flow to $H^{-1}(\R)$. 
Killip--Murphy--Vi\c{s}an \cite[Definition~7.3]{Killip.Murphy.Visan2020} subsequently introduced \emph{Green solutions}, trajectories satisfying an evolution identity for the reciprocal diagonal Green function of the Schr\"odinger Lax operator. 
For smooth solutions, this identity follows from the Lax formulation of KdV \cite[Proposition~3.1]{Killip.Visan2019}. Its terms remain well-defined for every curve in $C(I;H^{-1})$ at sufficiently large spectral parameters and depend continuously on the trajectory. The Green formulation thus gives an intrinsic solution notion at $H^{-1}$ that is compatible with classical solutions and stable under convergence in $C(I;H^{-1})$ when the defining identities hold on a common sufficiently large spectral range. This criterion allows us to ask whether every rough trajectory satisfying the identities is uniquely determined by its initial datum. We give the definition below and justify its terms in Section~\ref{sec:solution-notion}.

% Killip--Murphy--Vi\c{s}an \cite{Killip.Murphy.Visan2020} subsequently introduced an intrinsic solution notion that encodes the evolution through the diagonal Green function of a Schr\"odinger operator. We quote that notion here in the normalization used throughout this paper.
% 

For real-valued $q\in H^{-1}(\R)$, let $L_q=-\partial_x^2+q$ denote the self-adjoint, lower-semibounded Schr\"odinger operator associated with the quadratic form 
\begin{align}
 \label{eq:intro-form} \mathfrak l_q[\phi] =\int_\R \abs{\phi_x}^2\dd x+\pair*{q}{\abs{\phi}^2}, \qquad \Dom(\mathfrak l_q)=H^1(\R).
\end{align}
When discussing the Schr\"odinger operator $L_q$, we refer to its coefficient $q$ as the potential.

For every sufficiently large spectral parameter $\lambda>0$, the operator $L_q+\lambda^2$ is strictly positive. Let $G_\lambda(x,y;q)$ be the continuous integral kernel of its inverse and set
\begin{align}
 \label{eq:intro-resolvent} \begin{aligned}
 R_\lambda(q)&=(L_q+\lambda^2)^{-1},\qquad g_\lambda(x;q)=G_\lambda(x,x;q),\\
 h_\lambda(x;q)&=\frac{1}{2g_\lambda(x;q)},\qquad w_\lambda(x;q)=h_\lambda(x;q)-\lambda.
 \end{aligned}
\end{align}
For such $\lambda$, $g_\lambda-(2\lambda)^{-1}\in H^1$ and $h_\lambda-\lambda\in H^1$. The construction and mapping properties of these functions are reviewed in Section~\ref{sec:solution-notion}. Along a fixed trajectory $q$, we abbreviate $g_\lambda(q(t))$ by $g_\lambda(t)$. The same convention applies to $h_\lambda,w_\lambda$ and, once introduced, to the density and its mass. Adding the spatial argument, as in $g_\lambda(t,x)$, denotes evaluation at $x$.

{\begin{definition}[KdV Green solution]\label{def:green}
A real curve $q\in C(I;H^{-1}(\R;\R))$ is a \emph{Green solution to \eqref{eq:KdV}} if there exists $\lambda_*(q,I)\geq\lambda_0(q(I))$ 
% where $\lambda_0(q(I))$ is the admissibility threshold that Proposition~\ref{prop:static} attaches to the compact set $q(I)$, so that $h_\lambda(q(t))\in\lambda+H^1(\R)$ is defined for every $t\in I$, 
such that, for every $\lambda\geq\lambda_*(q,I)$ and $t\in I$, {
\begin{align}
 \label{eq:intro-green-identity} h_\lambda(q(t))-h_\lambda(q(0)) =\partial_x\int_0^t\J_\lambda(q(t'))\dd t' \quad\text{in }H^{-2}(\R),
\end{align}
} where the constant-normalized flux is
\begin{align}
 \label{eq:intro-green-flux} \J_\lambda(q) =\frac{q}{g_\lambda(q)} -2\lambda^2\left(\frac{1}{g_\lambda(q)}-2\lambda\right) \in H^{-1}(\R).
\end{align}
Here, $\lambda_0(q(I))$ is the threshold in Proposition~\ref{prop:static}, which ensures that $h_\lambda(q(t))\in\lambda+H^1(\R)$ is defined throughout $I$.
The product $q/g_\lambda(q)$ is defined as $2\lambda q+\big(g_\lambda(q)^{-1}-2\lambda\big)q$, using the continuous $H^1\times H^{-1}\to H^{-1}$ multiplication for the second term. The time integral in \eqref{eq:intro-green-identity} is a Bochner integral in $H^{-1}$, and the identity holds continuously in $t$ with values in $H^{-2}$. \end{definition}}

Definition~\ref{def:green} is equivalent to \cite[Definition~7.3 and equation~(7.17)]{Killip.Murphy.Visan2020}. Their raw flux differs from \eqref{eq:intro-green-flux} only by a spatial constant, which disappears after applying $\partial_x$.

The canonical Killip--Vi\c{s}an \cite{Killip.Visan2019} trajectories are known to be Green solutions. What remained open was whether the intrinsic Green identity alone determines the trajectory; Killip--Murphy--Vi\c{s}an \cite[p.~276]{Killip.Murphy.Visan2020} state their belief that it does and describe the question as a challenging open problem. 
% Here \emph{unconditional  uniqueness} means uniqueness among all real curves in $C(I;H^{-1})$ satisfying \eqref{eq:intro-green-identity} for every sufficiently large spectral parameter, without any approximation property, local smoothing estimate, or auxiliary resolution-space bound in the definition. Our main result answers this question.
Our main result proves this \emph{unconditional uniqueness}: uniqueness among all curves in Definition 1.1, without any approximation property or auxiliary spacetime bound imposed.

%%%%%%%%%%%%%%%%%%%%%%%%%%%%%%%%%%%%%%%%%%%%%%%%%%%%%%%%%%%%%%%%%%%%%%%%%%%%%%%%%%
\subsection{Main Results}
%%%%%%%%%%%%%%%%%%%%%%%%%%%%%%%%%%%%%%%%%%%%%%%%%%%%%%%%%%%%%%%%%%%%%%%%%%%%%%%%%%

{\begin{theorem}[Uniqueness for KdV in $H^{-1}$]\label{thm:main} Let $T>0$, $I=[-T,T]$, and let
\begin{align}
 \label{eq:main-a1} q_1,q_2\in C(I;H^{-1}(\R;\R))
\end{align}
be Green solutions in the sense of Definition~\ref{def:green}.

If $q_1(0)=q_2(0)$ in $H^{-1}$, then
\begin{align}
 \label{eq:main-a2} q_1(t)=q_2(t)\quad\text{in }H^{-1}(\R) \qquad\text{for every }t\in I.
\end{align}
\end{theorem}}

{\begin{corollary}
\label{cor:all-s}
For every $s\geq-1$, two Green solutions in $C(I;H^s(\R;\R))$ with the same initial datum coincide.
\end{corollary}}

{\begin{corollary}\label{cor:canonical}
Every $C(I;H^{-1}(\R;\R))$ Green solution equals the canonical Killip--Vi\c{s}an \cite{Killip.Visan2019} continuous extension of the Schwartz KdV flow with the same initial datum.
\end{corollary}}

Corollary~\ref{cor:all-s} concerns curves satisfying the Green identities. Applying it to a distributional $H^s$ solution requires verifying those identities. For $s\geq0$ it is already contained in the remark following Definition~7.3 of \cite{Killip.Murphy.Visan2020}: sending $\lambda\to\infty$ in \eqref{eq:intro-green-identity} shows that Green solutions in $C(I;L^2)$ are distributional solutions, which are unique by Zhou \cite{Zhou1997}. The new range is $-1\leq s<0$. To prove Corollary~\ref{cor:canonical}, we use the inclusion stated after \cite[Definition~7.3]{Killip.Murphy.Visan2020} and proved in Lemma~\ref{lem:canonical-is-green}: every canonical Killip--Vi\c{s}an \cite{Killip.Visan2019} flow is in the Green class. The converse follows from Theorem~\ref{thm:main}.

% Corollary~\ref{cor:all-s} concerns the Green class. For $s\geq0$, the limit $\lambda\to\infty$ in \eqref{eq:intro-green-identity} gives the distributional KdV equation, whose $L^2$ uniqueness is due to Zhou \cite{Zhou1997}, as noted after Definition~7.3 of \cite{Killip.Murphy.Visan2020}. The new range is $-1\leq s<0$. Corollary~\ref{cor:canonical} follows from Theorem~\ref{thm:main} and the inclusion of the canonical flow in the Green class, recorded in that reference and proved in Lemma~\ref{lem:canonical-is-green}.

%%%%%%%%%%%%%%%%%%%%%%%%%%%%%%%%%%%%%%%%%%%%%%%%%%%%%%%%%%%%%%%%%%%%%%%%%%%%%%%%%%
\subsection{mKdV Formulation and Theorem}%{mKdV Formulation and Result}
%%%%%%%%%%%%%%%%%%%%%%%%%%%%%%%%%%%%%%%%%%%%%%%%%%%%%%%%%%%%%%%%%%%%%%%%%%%%%%%%%%

We also study the real defocusing modified Korteweg--de Vries (mKdV) initial-value problem on the line, posed on the same time interval $I=[-T,T]$: 

{
\begin{align} \tag{mKdV}
 \label{eq:intro-mkdv} \begin{aligned}
 &u_t=-u_{xxx}+6u^2u_x=-u_{xxx}+2\partial_x(u^3), \qquad (t,x)\in I\times\R,\\
 &u(0,\cdot)=u_0\in L^2(\R;\R).
 \end{aligned}
\end{align}
} 
For $u\in C(I;L^2)$, the cube $u(t)^3\in L^{\frac23}$ need not be locally integrable pointwise in time, so the $L^2$ hypothesis alone does not give the cubic term a distributional meaning.
The class $C_tL^2_x$ thus plays for mKdV exactly the role that $C_tH^{-1}_x$ plays for KdV.
% As with \eqref{eq:KdV}, this identity has its usual classical meaning for smooth solutions. For an arbitrary curve $u\in C(I;L^2)$, however, the cubic term is not intrinsically defined: at a fixed time, the pointwise cube $\big(u(t)\big)^3$ belongs to $L^{\frac23}$ and need not be locally integrable, so the $L^2$ hypothesis alone does not give a distributional meaning to the cubic term. Consequently, \eqref{eq:intro-mkdv} cannot by itself be interpreted as a distributional evolution equation pointwise in time. The class $C_tL^2_x$ thus plays for mKdV exactly the role that $C_tH^{-1}_x$ plays for KdV.

% Harrop-Griffiths--Killip--Vi\c{s}an \cite{Harrop-Griffiths.Killip.Visan2020}  constructed the unique jointly continuous extension of the Schwartz defocusing mKdV flow to  $L^2(\R)$. 
% Harrop-Griffiths--Killip--Vi\c{s}an \cite[Theorem~1.1]{Harrop-Griffiths.Killip.Visan2020} proved that the Schwartz flow for real defocusing mKdV admits a unique jointly continuous extension to $H^s(\R)$ for every $s>-\frac12$. Forlano \cite[Theorem~1.1]{Forlano2022} subsequently gave a simpler proof in the range $0\leq s<\frac12$. We use the case $s=0$, which defines the continuous mKdV flow on $L^2(\R)$.
Harrop-Griffiths--Killip--Vi\c{s}an \cite[Theorem~1.1]{Harrop-Griffiths.Killip.Visan2020} constructed the unique jointly continuous extension of the Schwartz flow for complex-valued mKdV, in both the focusing and defocusing cases, to $H^s(\R)$ for every $s>-\frac12$. We use this result for real-valued initial data in the defocusing case at $s=0$.
Forlano--Killip--Vi\c{s}an \cite{Forlano.Killip.Visan2024} subsequently introduced an intrinsic solution notion that encodes the evolution through the diagonal Green functions of the mKdV Lax operator. We quote that notion below.
%here in the normalization used throughout this paper.

% For a real potential $u\in L^2(\R)$, the relevant matrix Lax operator is
% \begin{align}
%  \label{eq:intro-dirac} \bL_u =\begin{pmatrix}
%  -\partial_x&u\\
%  -u&\partial_x
%  \end{pmatrix}, \qquad \bR_\kappa(u)=(\bL_u+\kappa)^{-1}, \qquad \kappa\geq1,
% \end{align}
% where $\bL_u$ is anti-self-adjoint on its maximal domain \cite[Proposition~3.1]{Forlano.Killip.Visan2024}. Thus $\ui\bL_u$ is a self-adjoint operator of Dirac type. Since the spectrum of $\bL_u$ lies on the imaginary axis, $\bL_u+\kappa$ is boundedly invertible for every $\kappa\geq1$.

% The Schr\"odinger resolvent in \eqref{eq:intro-resolvent} is defined when $-\lambda^2$ lies outside the spectrum of $L_q$. For a general scalar potential, lower semiboundedness guarantees that $L_q+\lambda^2$ is strictly positive for all sufficiently large $\lambda$. The range can be larger for particular potentials. For the nonnegative Schr\"odinger operators associated with the Miura transforms in Section~\ref{sec:mkdv}, invertibility holds at every $\lambda>0$, as it does for the free multiplier $(\xi^2+\lambda^2)^{-1}$.

% Invertibility is in any case weaker than the quantitative hypothesis on the spectral parameter of Proposition~\ref{prop:static}, which may still require a large parameter. See Remark~\ref{rem:positivity-not-admissibility}.

For real-valued $u\in L^2(\R)$, the relevant Lax operator is no longer scalar. Let
\begin{align}
 \label{eq:intro-dirac} \bL_u =\begin{pmatrix}
 -\partial_x&u\\
 -u&\partial_x
 \end{pmatrix}, \qquad \bR_\kappa(u)=(\bL_u+\kappa)^{-1}, \qquad \kappa\geq1,
\end{align}
where $\bL_u$ is anti-self-adjoint on its maximal domain (cf. \cite[Proposition~3.1]{Forlano.Killip.Visan2024}). Its spectrum lies on the imaginary axis, so $\bL_u+\kappa$ is boundedly invertible for every $\kappa\geq1$, and $\ui\bL_u$ is a self-adjoint operator of Dirac type. For a general $q\in H^{-1}$, lower semiboundedness guarantees positivity of $L_q+\lambda^2$ for all sufficiently large $\lambda$. The Miura potentials considered below give non-negative Schr\"odinger operators, but the quantitative bounds of Proposition~\ref{prop:static} can still require a large parameter; see Remark~\ref{rem:positivity-not-admissibility}.

Write $\bG_\kappa(x,y;u)$ for the matrix integral kernel of $\bR_\kappa(u)$ and $\bG_{\kappa,0}$ for the free kernel. The kernel $\bG_\kappa$ inherits the diagonal jump of $\bG_{\kappa,0}$, but the difference $\bG_\kappa-\bG_{\kappa,0}$ is continuous and the free off-diagonal entries vanish. Hence the combinations
\begin{align}
 \gamma_\kappa&=[\bG_\kappa-\bG_{\kappa,0}]_{11}(x,x)+[\bG_\kappa-\bG_{\kappa,0}]_{22}(x,x), \label{eq:intro-dirac-fields}\\
 g_{-,\kappa}&=[\bG_\kappa]_{21}(x,x)-[\bG_\kappa]_{12}(x,x), \label{eq:intro-dirac-m}\\
 g_{+,\kappa}&=[\bG_\kappa]_{21}(x,x)+[\bG_\kappa]_{12}(x,x). \label{eq:intro-dirac-p}
\end{align}
are continuous functions of $x$. We follow the notation of \cite[equations~(3.19)--(3.21)]{Forlano.Killip.Visan2024}, with the spectral parameter displayed as a subscript. 
For every $\kappa\geq1$, these functions satisfy $\gamma_\kappa,g_{+,\kappa}\in H^1(\R)$ and, one derivative better, $g_{-,\kappa}\in H^2(\R)$, with locally Lipschitz dependence on $u\in L^2$ (cf. \cite[Proposition~3.7]{Forlano.Killip.Visan2024}). The construction and mapping properties of these functions are reviewed in Section~\ref{sec:mkdv}.

{\begin{definition}[mKdV Green solution]
\label{def:mkdv-green}
A real curve $u\in C(I;L^2(\R;\R))$ is a \emph{Green solution to \eqref{eq:intro-mkdv}} if, on writing $g_{-,\kappa}(t)=g_{-,\kappa}(\,\cdot\,;u(t))$ and $\gamma_\kappa(t)=\gamma_\kappa(\,\cdot\,;u(t))$, for every $\kappa\geq1$ and all $s,t\in I$, {
\begin{align}
 g_{-,\kappa}(t)-g_{-,\kappa}(s) &=\int_s^t\bigl[-g_{-,\kappa,xxx} +6u^2g_{-,\kappa,x}\bigr](t')\dd t', \label{eq:intro-mkdv-green-m}\\
 \gamma_\kappa(t)-\gamma_\kappa(s) &=\int_s^t\Bigl[-\gamma_{\kappa,xxx} +\partial_x\bigl\{6u^2(1+\gamma_\kappa)-12\kappa u g_{-,\kappa} -12\kappa^2\gamma_\kappa\bigr\}\Bigr](t')\dd t'. \label{eq:intro-mkdv-green-gamma}
\end{align}
} Every term here is meaningful at the stated regularity. Since $g_{-,\kappa,x}\in H^1\hookrightarrow L^\infty$ and $u^2\in L^1\hookrightarrow H^{-1}$, the first integrand is a continuous $H^{-1}$-valued function of time, while the second is continuous into $H^{-2}$. The two time integrals are Bochner integrals in $H^{-1}$ and $H^{-2}$, respectively. The identities hold continuously in $t$ in these spaces. Requiring them in $\mathcal D'(\R)$ is equivalent, since all terms already belong to the stated Sobolev spaces.
\end{definition}}

Definition~\ref{def:mkdv-green} is the unweighted case $b=0$ of \cite[Definition~3.11, equations~(3.75)--(3.76)]{Forlano.Killip.Visan2024}. It prescribes evolution identities for both $g_{-,\kappa}$ and $\gamma_\kappa$. 
% The extra regularity $g_{-,\kappa}\in H^2(\R)$ is needed for the chain rule in Section~\ref{sec:mkdv-chain-rule}. 
Both identities here and the extra derivative of $g_{-,\kappa}\in H^2$ enter the chain rule that transfers this evolution to the corresponding scalar reciprocals in Section~\ref{sec:mkdv-chain-rule}.
By comparison, the scalar identity \eqref{eq:intro-green-identity} governs $h_\lambda$, for which the definition gives only $h_\lambda-\lambda\in H^1$. 
% The proof uses both mKdV identities to obtain the evolution of the corresponding scalar reciprocal.

The canonical trajectories are again members of the intrinsic class: every curve of the canonical $L^2$ flow satisfies \eqref{eq:intro-mkdv-green-m}--\eqref{eq:intro-mkdv-green-gamma} on every compact interval; see Lemma~\ref{lem:canonical-mkdv-green}. Forlano--Killip--Vi\c{s}an \cite{Forlano.Killip.Visan2024} used this notion to construct measure-preserving trajectories for Gibbs-distributed data, but deterministic unconditional uniqueness within the class of Green solutions in $C(I;L^2)$ was not addressed there. 
% Here unconditional Green uniqueness means uniqueness among all real curves in $C(I;L^2)$ satisfying \eqref{eq:intro-mkdv-green-m}--\eqref{eq:intro-mkdv-green-gamma} for every $\kappa\geq1$, without any Miura-approximation property, local smoothing estimate, or auxiliary spacetime bound in the definition. Our second main result answers this question.
% Our second theorem proves (deterministic) unconditional uniqueness in its $C(I;L^2)$ class.
Our second theorem answers this question.

{\begin{theorem}[Uniqueness for mKdV in $L^2$]
\label{thm:mkdv-main}
Let $T>0$, $I=[-T,T]$, and let
\begin{align}
 \label{eq:mkdv-main-a1} u_1,u_2\in C(I;L^2(\R;\R))
\end{align}
be defocusing mKdV Green solutions in the sense of Definition~\ref{def:mkdv-green}.

If $u_1(0)=u_2(0)$ in $L^2$, then
\begin{align}
 \label{eq:mkdv-main-a2} u_1(t)=u_2(t)\quad\text{in }L^2(\R) \qquad\text{for every }t\in I.
\end{align}
Every such trajectory equals the canonical continuous extension of the Schwartz defocusing mKdV flow with the same datum.

\end{theorem}}

% Theorem~\ref{thm:mkdv-main} is proved in Section~\ref{sec:mkdv} through the two Miura transforms
% \begin{align}
%  \label{eq:intro-miura-partners} \begin{aligned}
%  \cM_\sigma&\colon L^2(\R;\R)\rightarrow H^{-1}(\R),\qquad \sigma\in\{+1,-1\},\\
%  q_\sigma&:=\cM_\sigma(u):=u^2+\sigma u_x,
%  \end{aligned}
% \end{align}
% where the derivative is distributional. Both images belong to $C(I;H^{-1})$ when $u\in C(I;L^2)$. Proposition~\ref{prop:mkdv-transfer} shows that both are KdV Green solutions. Applying Theorem~\ref{thm:main} for each sign then identifies $u$ through $q_+-q_-=2u_x$. The use of both signs avoids the non-injectivity of a single Miura map on $L^2$, which is discussed in Remark~\ref{rem:one-Miura-insufficient}.

% The global mapping properties needed for this argument hold in unweighted $L^2$. Extending the result to the weighted classes in \cite{Forlano.Killip.Visan2024} requires additional estimates, as explained in Remark~\ref{rem:mkdv-scope}. The defocusing sign is also needed: it makes the real transforms $u^2\pm u_x$ intertwine with \eqref{eq:KdV}. For the focusing equation the corresponding transforms are complex, outside the real self-adjoint theory used here. See \eqref{eq:classical-Miura-intertwining} and Remark~\ref{rem:defocusing-sign}.

% Theorem~\ref{thm:mkdv-main} is proved in Section~\ref{sec:mkdv}. Both solution notions also imply the spacetime regularity stated in Section~\ref{sec:intro-smoothing} below.

The link between the two theorems runs through the two Miura transforms
\begin{align}
 \label{eq:intro-miura-partners} \begin{aligned}
 \cM_\sigma&\colon L^2(\R;\R)\rightarrow H^{-1}(\R),\qquad \sigma\in\{+1,-1\},\\
 q_\sigma&:=\cM_\sigma(u):=u^2+\sigma u_x,
 \end{aligned}
\end{align}
where the derivative is distributional. Both images lie in $C(I;H^{-1})$ whenever $u\in C(I;L^2)$. We prove that the full pair \eqref{eq:intro-mkdv-green-m}--\eqref{eq:intro-mkdv-green-gamma} forces both $q_+$ and $q_-$ to be KdV Green solutions in the sense of Definition~\ref{def:green}. 
Theorem~\ref{thm:main} then applies once for each sign, and the difference $q_+-q_-=2u_x$ recovers $u$. We use both signs because that is what makes this recovery available: a single Miura map is not injective on $L^2(\R)$, so one image does not determine $u$; see Remark~\ref{rem:one-Miura-insufficient}. 
The rough transfer mechanism, which is the new input of this part of the paper, is outlined in Section~\ref{sec:idea} and carried out in Section~\ref{sec:mkdv}.

The global mapping properties needed for this argument hold in unweighted $L^2$. Extending the result to the weighted classes in \cite{Forlano.Killip.Visan2024} requires additional estimates, as explained in Remark~\ref{rem:mkdv-scope}. The defocusing sign is also needed: it makes the real transforms $u^2\pm u_x$ intertwine with \eqref{eq:KdV}. For the focusing equation the corresponding transforms are complex, outside the real self-adjoint theory used here. See \eqref{eq:classical-Miura-intertwining} and Remark~\ref{rem:defocusing-sign}.

% The solution notion in Theorem~\ref{thm:mkdv-main} is the established Dirac-resolvent notion of Forlano--Killip--Vi\c{s}an \cite[Definition~3.11]{Forlano.Killip.Visan2024}, restricted to unweighted $L^2(\R)$ trajectories. Unweighted $L^2(\R)$ control supplies the global mapping properties used to apply Theorem~\ref{thm:main} to the Miura transforms. Extending this argument to the weighted class treated probabilistically in \cite{Forlano.Killip.Visan2024} would require corresponding global estimates; see Remark~\ref{rem:mkdv-scope}. The defocusing sign is likewise essential. The expressions $u^2\pm u_x$ are real for every real $u$, whatever equation it solves. What the sign decides is whether these \emph{real} transforms intertwine with the real equation \eqref{eq:KdV}, as the identities \eqref{eq:classical-Miura-intertwining} record. For the focusing equation the corresponding transforms are complex, so that case requires a theory for complex Schr\"odinger potentials; see Remark~\ref{rem:defocusing-sign}.

%%%%%%%%%%%%%%%%%%%%%%%%%%%%%%%%%%%%%%%%%%%%%%%%%%%%%%%%%%%%%%%%%%%%%%%%%%%%%%%%%%
\subsection{Automatic Local Smoothing}
\label{sec:intro-smoothing}
%%%%%%%%%%%%%%%%%%%%%%%%%%%%%%%%%%%%%%%%%%%%%%%%%%%%%%%%%%%%%%%%%%%%%%%%%%%%%%%%%%

Definition~\ref{def:green} requires continuity in $H^{-1}$ and a single integrated identity at each sufficiently large spectral parameter. Definition~\ref{def:mkdv-green} requires continuity in $L^2$ and two integrated identities at every $\kappa\geq1$.
The parameter ranges reflect the different spectral properties of the Lax operators. For KdV, the lower bound of $L_q$ depends on $q$, and the quantitative Green-function estimates require a sufficiently large parameter. For real defocusing mKdV, $\bL_u$ is anti-self-adjoint, so $\bL_u+\kappa$ is invertible for every $\kappa>0$. We retain the range $\kappa\geq1$ from \cite[Definition~3.11]{Forlano.Killip.Visan2024}, although the uniqueness argument in Section~\ref{sec:mkdv} uses only sufficiently large $\kappa$.
These hypotheses imply the following Kato-type local smoothing estimates.
The KdV estimate supplies regularity needed for the Gardner transfer and uniqueness. The mKdV estimate follows from the same KdV result applied to both Miura transforms.

A representative of a curve $q$ in $L^2_{\loc}(J\times\R)$ is a function $\widetilde q$ in that space such that
\begin{align}
 \label{eq:representative-def} \pair*{q(t)}{\phi}=\int_\R\widetilde q(t,x)\phi(x)\dd x \qquad\text{for every }\phi\in C_c^\infty(\R),
\end{align}
for almost every $t\in J$, with one exceptional set for all tests $\phi$. Any two representatives agree almost everywhere, so their cell integrals coincide. We write $q$ for this representative when taking such integrals. For a single distribution $f$, an $L^2_{\loc}$ representative is a function realizing $f$ by the same test-function identity. 

{\begin{proposition}[Local smoothing for KdV]
\label{prop:auto-smoothing}
Let $q\in C(I;H^{-1}(\R;\R))$ be a KdV Green solution in the sense of Definition~\ref{def:green}. For every compact time interval $J\subset I$ there is a representative of $q$ in $L^2_{\loc}(J\times\R)$ such that
\begin{align}
 \label{eq:q-Kato-bound} \sup_{y\in\R}\int_J\int_{y-1}^{y+1}\abs{q(t,x)}^2\dd x\dd t<\infty
\end{align}
and
\begin{align}
 \label{eq:q-Kato-tail} \lim_{y\to\pm\infty} \int_J\int_{y-1}^{y+1}\abs{q(t,x)}^2\dd x\dd t=0.
\end{align}
\end{proposition}}

{\begin{proposition}[Local smoothing for mKdV]
\label{prop:mkdv-smoothing}
Let $u\in C(I;L^2(\R;\R))$ be a defocusing mKdV Green solution in the sense of Definition~\ref{def:mkdv-green}. For every compact interval $J\subset I$, the distributional derivative $u_x$ belongs to $L^2_{\loc}(J\times\R)$, $u$ belongs to $L^4_{\loc}(J\times\R)$, and
\begin{align}
 \label{eq:intro-mkdv-smoothing} \begin{aligned}
 &\sup_{y\in\R}\int_J\int_{y-1}^{y+1} \bigl(\abs{u_x}^2+\abs{u}^4\bigr)\dd x\dd t<\infty ,\\
 &\lim_{y\to\pm\infty}\int_J\int_{y-1}^{y+1} \bigl(\abs{u_x}^2+\abs{u}^4\bigr)\dd x\dd t=0 .
 \end{aligned}
\end{align}
\end{proposition}}

The two propositions are proved in Sections~\ref{sec:automatic-smoothing} and~\ref{sec:mkdv} respectively, and the second is deduced from the first: by Proposition~\ref{prop:mkdv-transfer} the two Miura transforms of $u$ are KdV Green solutions, Proposition~\ref{prop:auto-smoothing} applies to each of them, and the parallelogram identity \eqref{eq:Miura-parallelogram} converts the two conclusions into \eqref{eq:intro-mkdv-smoothing}. Both arguments apply to each Green trajectory directly.

Their roles in the paper are, however, opposite: 
% and a reader deciding what to look at should know which is which.

\textit{For KdV it is the decisive step.} Proposition~\ref{prop:auto-smoothing} is proved in Section~\ref{sec:automatic-smoothing} and supplies regularity needed in the proof of Theorem~\ref{thm:main}. The reason is that the two sides of the argument speak different languages. The hypothesis is spectral (an identity among Green functions) while the uniqueness theorem that is finally invoked, Lemma~\ref{lem:KKL-Gardner}, is posed in the Kato class $X_1$ and speaks about spacetime integrability. Proposition~\ref{prop:auto-smoothing} is the translation between the two, and two of the four gaps that Section~\ref{sec:Gardner} has to close are spent on transporting it: it is what puts the Gardner variable $r$ into the Kato conditions of $X_1$, locally in time, and what supplies the one further spatial derivative needed to identify the flux.

\textit{For mKdV it is an add-on.} Proposition~\ref{prop:mkdv-smoothing} is proved after Theorem~\ref{thm:mkdv-main} and records a further regularity consequence of the Green identities. The estimate follows from the Forlano--Killip--Vi\c{s}an \cite{Forlano.Killip.Visan2024} identities by applying the Killip--Murphy--Vi\c{s}an \cite{Killip.Murphy.Visan2020} formulation to both Miura transforms.

Proposition~\ref{prop:auto-smoothing} therefore does both jobs. In the first part it supplies the regularity needed for uniqueness. In the second part, applied to the two Miura transforms, it is what yields Proposition~\ref{prop:mkdv-smoothing}.

%%%%%%%%%%%%%%%%%%%%%%%%%%%%%%%%%%%%%%%%%%%%%%%%%%%%%%%%%%%%%%%%%%%%%%%%%%%%%%%%%%
\subsection{Literature Review}
\label{sec:literature}
%%%%%%%%%%%%%%%%%%%%%%%%%%%%%%%%%%%%%%%%%%%%%%%%%%%%%%%%%%%%%%%%%%%%%%%%%%%%%%%%%%

The low-regularity theory of \eqref{eq:KdV} and \eqref{eq:intro-mkdv} has been shaped by two distinct questions. The first concerns well-posedness: how rough an initial datum may be while still generating a solution through a continuous solution map. The second concerns uniqueness: at a given regularity, within how large a class of curves a solution is determined by its initial datum. 
For KdV the first question is resolved, and the answer is sharp: the solution map extends continuously to $H^{-1}(\R)$ and admits no continuous extension below it. For mKdV the problem is well-posed in $H^s(\R)$ for every $s>-\tfrac12$, and whether this threshold is optimal for the real defocusing equation considered here is not known. The present paper addresses the second question at two levels at which it had remained open: the KdV endpoint $H^{-1}(\R)$ and the mKdV Green class in $L^2(\R)$.
% Well-posedness constructs a continuous solution map, while uniqueness asks whether a specified class contains other trajectories with the same datum. 
We review these questions separately, with emphasis on the regularities of Theorems~\ref{thm:main} and~\ref{thm:mkdv-main}.

%%%%%%%%%%%%%%%%%%%%%%%%%%%%%%%%%%%%%%%%%%%%%%%%%%%%%%%%%%%%%%%%%%%%%%%%%%%%%%%%%%
\subsubsection*{Well-posedness of KdV}
%%%%%%%%%%%%%%%%%%%%%%%%%%%%%%%%%%%%%%%%%%%%%%%%%%%%%%%%%%%%%%%%%%%%%%%%%%%%%%%%%%
% \textit{Well-posedness of KdV.}
The classical theory uses energy estimates and compactness at high regularity. See Bona--Smith \cite{Bona.Smith1975}, Saut--Temam \cite{Saut.Temam1976}, and Kato \cite{Kato1983}. Kenig--Ponce--Vega \cite{Kenig.Ponce.Vega1991,Kenig.Ponce.Vega1993} used local smoothing and maximal function estimates to develop the contraction theory for generalized KdV. Bourgain \cite{Bourgain1993(=)} introduced Fourier restriction spaces and proved well-posedness for $s\geq0$, and Kenig--Ponce--Vega \cite{Kenig.Ponce.Vega1996} extended the local theory on $\R$ to $s>-\tfrac34$. At $s=-\tfrac34$, Christ--Colliander--Tao \cite[Theorem~6]{Christ.Colliander.Tao2003} constructed a locally Lipschitz extension of the smooth solution map, and Kishimoto \cite{Kishimoto2009} and Guo \cite{Guo2009} established local well-posedness. Colliander--Keel--Staffilani--Takaoka--Tao \cite{Colliander.Keel.Staffilani.Takaoka.Tao2003} proved global well-posedness for $s>-\tfrac34$, and Kishimoto \cite{Kishimoto2009} reached the endpoint. Guo \cite{Guo2009} also states the global endpoint result, with a reservation about the presented argument recorded in \cite[p.~652, footnote~2]{Kwon.Oh.Yoon2020}.

For $-1\leq s<-\tfrac34$, Christ--Colliander--Tao \cite[Theorem~7]{Christ.Colliander.Tao2003} exclude a uniformly continuous solution map on bounded sets for real KdV. The corresponding results for complex KdV and focusing mKdV are due to Kenig--Ponce--Vega \cite{Kenig.Ponce.Vega2001}. Complete integrability nevertheless gives continuous solution maps below this perturbative threshold. Kappeler--Topalov \cite{Kappeler.Topalov2006} proved global well-posedness in $H^{-1}(\mathbb T,\R)$. Killip--Vi\c{s}an--Zhang \cite{Killip.Visan.Zhang2018} obtained a priori $H^s$ bounds for $-1\leq s<1$, and Buckmaster--Koch \cite{Buckmaster.Koch2015} constructed $H^{-1}(\R)$ solutions as limits of smooth solutions. Killip--Vi\c{s}an \cite{Killip.Visan2019} then constructed the global, jointly continuous extension of the Schwartz flow to $H^{-1}(\R)$. Molinet \cite{Molinet2011,Molinet2012} showed that no continuous extension exists below $H^{-1}$, on either $\R$ or $\mathbb T$.

Related results treat the fifth-order KdV equation \cite{Bringmann.Killip.Visan2021}, the entire hierarchy \cite{Klaus.Koch.Liu2023}, perturbations of step-like backgrounds \cite{Laurens2022}, bounded distributional solutions \cite{Chapouto.Killip.Visan2024}, and invariant measures \cite{Killip.Murphy.Visan2020}.

%%%%%%%%%%%%%%%%%%%%%%%%%%%%%%%%%%%%%%%%%%%%%%%%%%%%%%%%%%%%%%%%%%%%%%%%%%%%%%%%%%
\subsubsection*{Well-posedness of mKdV}
%%%%%%%%%%%%%%%%%%%%%%%%%%%%%%%%%%%%%%%%%%%%%%%%%%%%%%%%%%%%%%%%%%%%%%%%%%%%%%%%%%
% \exstep{Well-posedness of mKdV.}
The Miura map \cite{Miura1968} connects the two equations. For mKdV, Kenig--Ponce--Vega \cite{Kenig.Ponce.Vega1993} proved local well-posedness for $s\geq\tfrac14$, and Colliander--Keel--Staffilani--Takaoka--Tao \cite{Colliander.Keel.Staffilani.Takaoka.Tao2003} proved global well-posedness for $s>\tfrac14$. The global endpoint is stated by Guo \cite{Guo2009}, while Kwon--Oh--Yoon \cite[p.~652, footnote~2]{Kwon.Oh.Yoon2020} attribute the established endpoint argument to Kishimoto \cite{Kishimoto2009}. Gr\"unrock--Vega \cite{Grunrock.Vega2009} proved local well-posedness in almost critical Fourier--Lebesgue spaces. For $-\tfrac14<s<\tfrac14$, Christ--Colliander--Tao \cite[Theorem~4]{Christ.Colliander.Tao2003} exclude a uniformly continuous solution map for the real defocusing equation.

Integrability extends the theory further. Kappeler--Topalov \cite{Kappeler.Topalov2005} proved global well-posedness in $L^2(\mathbb T,\R)$, shown to be sharp by Molinet \cite{Molinet2012}. On $\R$, Harrop-Griffiths--Killip--Vi\c{s}an \cite{Harrop-Griffiths.Killip.Visan2020} established well-posedness in $H^s$ for every $s>-\tfrac12$, for both mKdV and cubic NLS. This threshold is optimal for the complex-valued and focusing problems. Its optimality for real defocusing mKdV remains open in the cited work, as also noted by Forlano \cite{Forlano2022}. The space $L^2(\R;\R)$ in Theorem~\ref{thm:mkdv-main} lies within the constructed range. Forlano--Killip--Vi\c{s}an \cite{Forlano.Killip.Visan2024} introduced its Dirac-resolvent formulation in their work on invariant measures in infinite volume.

%%%%%%%%%%%%%%%%%%%%%%%%%%%%%%%%%%%%%%%%%%%%%%%%%%%%%%%%%%%%%%%%%%%%%%%%%%%%%%%%%%
\subsubsection*{Conditional and unconditional uniqueness}
%%%%%%%%%%%%%%%%%%%%%%%%%%%%%%%%%%%%%%%%%%%%%%%%%%%%%%%%%%%%%%%%%%%%%%%%%%%%%%%%%%
% \exstep{Conditional and unconditional uniqueness.}
Low-regularity constructions typically give uniqueness within an auxiliary resolution space, or characterize the solution as the continuous extension of the smooth flow. An arbitrary continuous Sobolev-valued trajectory need not satisfy either condition. Removing such an auxiliary hypothesis is the content of \emph{unconditional uniqueness}, a term introduced by Kato \cite{Kato1995} for NLS. The compactness construction of \cite{Buckmaster.Koch2015} and the a priori bounds of \cite{Killip.Visan.Zhang2018} do not themselves give uniqueness.

Zhou \cite{Zhou1997} proved unconditional uniqueness for KdV in $L^2(\R)$. Babin--Ilyin--Titi \cite{Babin.Ilyin.Titi2011} treated periodic KdV at $L^2$ regularity by normal forms. For mKdV, Molinet--Pilod--Vento \cite{Molinet.Pilod.Vento2018,Molinet.Pilod.Vento2019} used modified energies to prove unconditional well-posedness for $s>\tfrac13$ on $\R$ and $s\geq\tfrac13$ on $\mathbb T$. Kwon--Oh--Yoon \cite{Kwon.Oh.Yoon2020} reached $s>\tfrac14$ on $\R$ through iterated normal forms, and Kishimoto \cite{Kishimoto2019} developed an abstract version of that method. For the complex equation
\begin{align}\label{eq:complex-mkdv}
 \partial_t\varphi=-\partial_x^3\varphi\pm6\abs{\varphi}^2\partial_x\varphi ,
\end{align}
Killip--Ouyang--Vi\c{s}an--Wu \cite[Theorem~1.3]{AA027} proved unconditional uniqueness in $L^\infty_tH^1_x$ by spacetime estimates and a Gronwall argument, in their analysis of the Ablowitz--Ladik continuum limit. The cited normal form and modified energy results require a distributional interpretation of the nonlinearity and do not reach the Green classes considered here.

%%%%%%%%%%%%%%%%%%%%%%%%%%%%%%%%%%%%%%%%%%%%%%%%%%%%%%%%%%%%%%%%%%%%%%%%%%%%%%%%%%
\subsubsection*{Uniqueness at the endpoint regularity}
%%%%%%%%%%%%%%%%%%%%%%%%%%%%%%%%%%%%%%%%%%%%%%%%%%%%%%%%%%%%%%%%%%%%%%%%%%%%%%%%%%
% \exstep{Uniqueness at the endpoint regularity.}
At $H^{-1}$, a solution notion must specify how the undefined square in \eqref{eq:KdV} is replaced. The Green formulations of \cite{Killip.Murphy.Visan2020,Forlano.Killip.Visan2024} use resolvent identities. Klaus--Koch--Liu \cite[Theorem~2.24]{Klaus.Koch.Liu2023} instead prove weak uniqueness for the Gardner equation in a Kato class. Their theorem becomes applicable here after we derive the required spacetime bounds and Gardner evolution from the Green identities.

Chapouto--Killip--Vi\c{s}an \cite[Definition~1.1 and Theorem~1.2]{Chapouto.Killip.Visan2024} prove uniqueness for bounded distributional KdV solutions continuous in the weak-$*$ topology of $L^\infty$. Their proof uses relative Green functions, while their solution class requires boundedness and a distributionally defined square.

To the best of our knowledge, Theorem~\ref{thm:main} is the first uniqueness theorem for arbitrary Killip--Murphy--Vi\c{s}an Green solutions in $C_tH^{-1}_x$. Proposition~\ref{prop:auto-smoothing} supplies the spacetime bounds needed to connect that class to Gardner uniqueness. For mKdV, the new input is Proposition~\ref{prop:mkdv-transfer}, which uses both defining identities to put both Miura transforms in the KdV Green class. This requires the reciprocal evolution law, beyond the scalar diagonal-$g$ law in \cite[equation~(7.20)]{Forlano.Killip.Visan2024}.

%%%%%%%%%%%%%%%%%%%%%%%%%%%%%%%%%%%%%%%%%%%%%%%%%%%%%%%%%%%%%%%%%%%%%%%%%%%%%%%%%%
\subsubsection*{Significance and applications}
%%%%%%%%%%%%%%%%%%%%%%%%%%%%%%%%%%%%%%%%%%%%%%%%%%%%%%%%%%%%%%%%%%%%%%%%%%%%%%%%%%
% \exstep{Significance and potential applications.}
Together with the known canonical flows, our theorems give existence, uniqueness, and continuous dependence within the two Green classes. They characterize those flows by the identities in Definitions~\ref{def:green} and~\ref{def:mkdv-green}. Continuous dependence follows from the jointly continuous solution maps of \cite[Theorem~1.1]{Killip.Visan2019} and \cite[Theorem~1.1]{Harrop-Griffiths.Killip.Visan2020}.

The results also identify compactness limits. Suppose $\{q_n\}$ is precompact in $C(I;H^{-1})$, its initial data converge to $q_0$, and each subsequential limit is a Green solution. Corollary~\ref{cor:canonical} identifies every limit with the canonical trajectory from $q_0$, so the whole sequence converges.

The Green property can be checked through approximate identities. If $q_n\to q$ in $C(I;H^{-1})$ and \eqref{eq:intro-green-identity} holds with errors tending to zero in $C(I;H^{-2})$ for every $\lambda\geq\Lambda$, with $\Lambda$ independent of $n$, set $\mathcal Q=\bigcup_nq_n(I)\cup q(I)$. Uniform convergence and compactness of $I$ make this set compact. Proposition~\ref{prop:static}\,\textup{(ii)} and Lemma~\ref{lem:microscopic}\,\textup{(i)} give continuity of $h_\lambda$ and $\J_\lambda$ on $\mathcal Q$ for $\lambda\geq\lambda_0(\mathcal Q)$. Passing to the limit proves \eqref{eq:intro-green-identity} for $\lambda\geq\max\{\Lambda,\lambda_0(\mathcal Q)\}$. The same argument applies to mKdV with convergence in $C(I;L^2)$ and errors vanishing in the two Sobolev spaces of Definition~\ref{def:mkdv-green}, for every $\kappa\geq1$. The needed continuity is proved in Lemma~\ref{lem:canonical-mkdv-green}.

Continuum limits of integrable lattices provide a possible application. Killip--Ouyang--Vi\c{s}an--Wu \cite{AA027} establish the Ablowitz--Ladik limit to \eqref{eq:complex-mkdv} for $H^1$ data, and Liu--Koch \cite{Liu.Koch2026} establish the Toda-to-KdV limit for $H^1$ data with a rate in $C([-T,T];L^2)$. Our theorems can identify lower-regularity limits in $H^{-1}$ for KdV and real $L^2$ for defocusing mKdV, provided that the respective Green identities pass to the limit. Extensions of the smoothing and Miura arguments to higher KdV flows or to the nondecaying classes containing the invariant-measure trajectories of \cite{Killip.Murphy.Visan2020,Forlano.Killip.Visan2024} remain open.

\subsection{Idea and Strategy}
\label{sec:idea}
%%%%%%%%%%%%%%%%%%%%%%%%%%%%%%%%%%%%%%%%%%%%%%%%%%%%%%%%%%%%%%%%%%%%%%%%%%%%%%%%%%

The two theorems are proved in sequence, and the second uses the first. 
% Part~I, carried out in Sections~\ref{sec:solution-notion}--\ref{sec:Gardner}, proves the endpoint KdV statement. Part~II, carried out in Section~\ref{sec:mkdv}, proves that both Miura transforms of an arbitrary mKdV Green solution are KdV Green solutions and then applies the uniqueness theorem from Part~I. 
% The two parts share one principle: the intrinsic Green identity, imposed for a full family of spectral parameters, implies spacetime regularity/(local smoothing).
Sections~\ref{sec:solution-notion}--\ref{sec:Gardner} prove KdV uniqueness by deriving local smoothing and transferring the Green evolution to Gardner. Section~\ref{sec:mkdv} then transfers both Miura transforms of an mKdV Green solution to the KdV class and applies the first theorem.
The two parts share one principle: the intrinsic Green identity, imposed for a full family of spectral parameters, implies spacetime regularity.

%%%%%%%%%%%%%%%%%%%%%%%%%%%%%%%%%%%%%%%%%%%%%%%%%%%%%%%%%%%%%%%%%%%%%%%%%%%%%%%%%%
\subsubsection*{Part I: the endpoint KdV theorem}
%%%%%%%%%%%%%%%%%%%%%%%%%%%%%%%%%%%%%%%%%%%%%%%%%%%%%%%%%%%%%%%%%%%%%%%%%%%%%%%%%%

% The proof of the KdV theorem is a transfer argument with two main phases. The nonlinear term $q^2$ requires additional regularity at the $H^{-1}$ endpoint. We first derive that regularity from the Green identities for all sufficiently large spectral parameters. We then freeze one spectral parameter and use the reciprocal of the corresponding diagonal Green function to convert the KdV trajectory to a Gardner trajectory, where the uniqueness theorem of \cite{Klaus.Koch.Liu2023} becomes applicable after a globalization step.

% For every admissible spectral parameter $\lambda$, consider the functions $h_\lambda,w_\lambda$ from \eqref{eq:intro-resolvent}. At a fixed parameter, the recovery formula \eqref{eq:q-from-h} expresses $q$ in terms of $h_\lambda$ and its spatial derivatives. As $\lambda\to\infty$, $2\lambda w_\lambda$ converges to $q$. These two complementary roles lead to the following proof scheme.

The KdV argument uses the spectral parameter in two ways. The limit $\lambda\to\infty$ recovers $q$ from $2\lambda w_\lambda$ and yields its local smoothing. At one fixed parameter, the recovery formula \eqref{eq:q-from-h} converts the trajectory into a Gardner solution. The construction must supply the derivative bounds, weak equation, and global time interval required by the uniqueness theorem of \cite{Klaus.Koch.Liu2023}.
These two complementary roles lead to the following proof scheme.

\begin{enumerate}[
label=\textup{(\roman*)}, leftmargin=*, itemsep=.7em, before={\let\fullwidthdisplay\relax} ]

\item \emph{Green functions and a common spectral range}\ \ (Section~\ref{sec:solution-notion}). Proposition~\ref{prop:static} gives the mapping properties of $g_\lambda$ and its reciprocal. Together with \eqref{eq:intro-resolvent}, it gives the recovery identity \eqref{eq:q-from-h} for $q$. Definition~\ref{def:green} supplies the integrated evolution law for $h_\lambda$, equivalently \eqref{eq:integrated-w} for $w_\lambda$. When two Green solutions are being compared, Remark~\ref{rem:common-lambda} fixes one threshold above which every spectral parameter is admissible for both trajectories and at every time. Thus all later estimates can be carried out on a single common spectral range.

\item \emph{Large-parameter limits and coercivity}\ \ (Section~\ref{sec:alpha}). The density $\rho_\lambda$ and its mass $\alpha_\lambda$ are defined in \eqref{eq:def-rho}--\eqref{eq:def-alpha}. They measure the nonlinear part of $h_\lambda$. Lemma~\ref{lem:high-energy} proves that the density is nonnegative and integrable for sufficiently large $\lambda$. Starting from the estimate for small $H^{-1}$ potentials in \cite[Proposition~2.4]{Killip.Visan2019}, the scaling statement in Lemma~\ref{lem:reciprocal-density-scaling} yields Lemma~\ref{lem:high-energy}. In particular, uniformly on compact subsets of $H^{-1}$,
\begin{align}
 \label{eq:roadmap-high-energy} \rho_\lambda\geq0, \qquad \lambda\alpha_\lambda(q)\rightarrow0, \qquad 2\lambda w_\lambda(q)\rightarrow q \quad\text{in }H^{-1}.
\end{align}
On the other hand, inserting the recovery identity \eqref{eq:q-from-h} into the Green flux gives the decomposition in Lemma~\ref{lem:flux-decomposition}:
\begin{align}
 \label{eq:roadmap-flux-decomposition} \J_\lambda=-w_{\lambda,xx}+N_\lambda.
\end{align}
The remainder $N_\lambda$ is defined in \eqref{eq:J-decomposition}. The term $-w_{\lambda,xx}$ contributes derivatives of the weights in the later weighted identities. The coercive bound comes from $N_\lambda$: the $L^\infty$ estimate \eqref{eq:winfty-high} forces $\norm{w_\lambda}_\infty\leq\lambda/2$ for all sufficiently large $\lambda$, and part~\textup{(ii)} of Lemma~\ref{lem:flux-decomposition} then gives the almost-everywhere bound \eqref{eq:N-positive}, namely $\lambda N_\lambda\geq\abs{w_{\lambda,x}}^2+5\lambda^2\abs{w_\lambda}^2$.

\item \emph{Conservation of the density mass}\ \ (Section~\ref{sec:alpha}). Lemma~\ref{lem:two-parameter} compares the density evolutions at two large parameters $\kappa$ and $\lambda$ through the exact integrated identity \eqref{eq:S-evolution}. Applying Lemma~\ref{lem:zero-mass} gives the relation \eqref{eq:cross-alpha} between the corresponding mass increments. Sending $\lambda\to\infty$ by Lemma~\ref{lem:high-energy} proves Proposition~\ref{prop:alpha-conservation}.

\item \emph{Automatic Kato smoothing for $q$}\ \ (Section~\ref{sec:automatic-smoothing}). 
Pairing \eqref{eq:S-evolution} with the translated weight \eqref{eq:tanh-weight} and using \eqref{eq:N-positive} bounds $2\lambda w_\lambda$ in weighted spacetime $L^2$. Weak compactness and \eqref{eq:w-high} identify the limit with $q$. For the translated-cell tail, Lemma~\ref{lem:translation} removes the translated Sobolev pairings and conservation of $\alpha_\kappa$ cancels the limiting density term. These give both conclusions of Proposition~\ref{prop:auto-smoothing}.

% There the two-parameter identity is paired with the translated hyperbolic-tangent weight of \eqref{eq:tanh-weight}. Combining the resulting identity \eqref{eq:key-weighted-identity} with the coercive estimate \eqref{eq:N-positive} gives uniform weighted $L^2$ bounds for $2\lambda w_\lambda$. Weak compactness and the convergence \eqref{eq:w-high} then identify the limiting local $L^2$ function with the original potential $q$.

% The translated-cell tail requires one additional use of the integrable structure. After taking $\lambda\to\infty$, the translated Sobolev pairings vanish by Lemma~\ref{lem:translation}, while Proposition~\ref{prop:alpha-conservation} removes the possible nonzero constant arising from the limiting density term. The complete conclusion is Proposition~\ref{prop:auto-smoothing}, namely the uniform local bound \eqref{eq:q-Kato-bound} and the translated-cell tail \eqref{eq:q-Kato-tail}.

\item \emph{The Gardner variable and its regularity}\ \ (Sections~\ref{sec:Gardner-class}--\ref{sec:Gardner-smoothing}). At one fixed admissible parameter $\tau$, define $v$ and $r$ by \eqref{eq:v-r-def}. Proposition~\ref{prop:static} gives $v\in C_tH^1_x$ and $r\in C_tL^2_x$, and Lemma~\ref{lem:Miura} proves $q=r_x+2\tau r+r^2$. Squaring $r_x+r^2=q-2\tau r$ against a cutoff bounds $r_x$ and $r^2$ on each cell, because $2r_xr^2=\frac23\partial_x(r^3)$. Combined with Proposition~\ref{prop:auto-smoothing}, this gives the additional spatial derivative of $v$ in Lemma~\ref{lem:internal-H2}. Proposition~\ref{prop:r-X1-local} proves the derivative bound and translated-cell tail for $r$ required locally by the Gardner class of Definition~\ref{def:Gardner-solution}.

\item \emph{Gardner dynamics, globalization, and uniqueness}\ \ (Sections~\ref{sec:Gardner-dynamics}--\ref{sec:uniqueness}). Dividing the Green identity at $\tau$ gives \eqref{eq:v-integrated}. The regularity from Lemma~\ref{lem:internal-H2} permits spatial mollification in the logarithmic chain rule and flux calculation of Proposition~\ref{prop:green-Gardner}. Hence
\begin{align}
 \label{eq:novelty-Gardner-Miura} r_t=\partial_x\bigl(-r_{xx}+6\tau r^2+2r^3\bigr), \qquad q=r_x+2\tau r+r^2.
\end{align}
The uniqueness class $X_1$ is global in time. Lemma~\ref{lem:global-extension} extends each KdV competitor by the canonical flow of \cite{Killip.Visan2019}, with a single spectral threshold valid on every compact interval. Choose $\tau\geq K_{12}$, with the common threshold formed from the extended trajectories as in the proof of Theorem~\ref{thm:main}. Lemma~\ref{lem:Gardner-mass} gives the conserved norm $\norm{r(t)}_2^2=2\tau\alpha_\tau(q(t))$, so the global lift belongs to $X_1$. Lemma~\ref{lem:KKL-Gardner}, the uniqueness input from \cite{Klaus.Koch.Liu2023}, identifies the lifts. The Miura identity then gives $q_1=q_2$.
\end{enumerate}

The KdV uniqueness argument follows the chain
\begin{align}
 \label{eq:roadmap-implication-chain}
 &\text{Green identities for all sufficiently large }\lambda\\
 &\quad\Rightarrow \text{local smoothing and spatial tail decay for }q \notag\\
 &\quad\Rightarrow \text{weak Gardner dynamics and local }X_1\text{ bounds for }r \notag\\
 &\quad\Rightarrow \text{a global Gardner solution in }X_1\text{ after extension} \notag\\
 &\quad\Rightarrow \text{uniqueness of }r\text{ and hence of }q. \notag
\end{align}
\subsubsection*{Part II: the defocusing mKdV theorem}
%%%%%%%%%%%%%%%%%%%%%%%%%%%%%%%%%%%%%%%%%%%%%%%%%%%%%%%%%%%%%%%%%%%%%%%%%%%%%%%%%%

The mKdV result follows by applying the KdV uniqueness theorem to both Miura transforms. We must therefore show that these transforms of an arbitrary $C_tL^2_x$ Green solution belong to the class of Definition~\ref{def:green}. Definition~\ref{def:mkdv-green} prescribes evolution identities for $g_{-,\kappa}$ and $\gamma_\kappa$, whereas Definition~\ref{def:green} prescribes an evolution identity for $h_\lambda$. The scalar KdV formulation in \cite{Forlano.Killip.Visan2024} governs the diagonal Green function $g$. The transfer here requires an evolution identity for its reciprocal. 
% See Remark~\ref{rem:full-mkdv-law}. 
Passing between the two formulations by the formal rule $h_t=-g_t/(2g^2)$ is not permitted at this regularity, since it multiplies an $H^{-2}$ distribution by an $H^1$ function. The proof below uses the additional regularity of $g_{-,\kappa}$ to justify the required chain rule.

\begin{enumerate}[
label=\textup{(\roman*)}, leftmargin=*, itemsep=.7em, before={\let\fullwidthdisplay\relax} ]

\item \emph{Two Miura transforms and the operator factorization}\ \ (Section~\ref{sec:mkdv-static}). Section~\ref{sec:mkdv} works with both Miura transforms $q_\sigma$ of \eqref{eq:intro-miura-partners}. Lemma~\ref{lem:two-Miura-static} proves the form factorization \eqref{eq:Miura-form-factorization} and identifies the maximal first-order domain \eqref{eq:Miura-maximal-domain} with $H^1(\R)$. Hence $L_{q_\sigma}\geq0$ for every real $u\in L^2$, so that $L_{q_\sigma}+\kappa^2$ is boundedly invertible for every $\kappa\geq1$ and $h_{\sigma,\kappa}$ is defined there. Admissibility additionally requires the quantitative bounds \eqref{eq:g-quantitative} of Proposition~\ref{prop:static}, as the example in Remark~\ref{rem:positivity-not-admissibility} illustrates. Compactness of the two Miura trajectories supplies these bounds above the threshold $K_u$ of \eqref{eq:mkdv-common-threshold}. The Green identities are verified on this common spectral range. Both signs are used for the recovery: $\cM_+$ alone is not injective on $L^2(\R)$ by Remark~\ref{rem:one-Miura-insufficient}, whereas the pair determines the solution through \eqref{eq:two-Miura-recovery}.

% \item \emph{Expressing $h_{\sigma,\kappa}$ through $g_{-,\kappa}\in H^2$}\ \ (Section~\ref{sec:mkdv-transfer}). For the Miura potentials, write
% \begin{align}
%  \label{eq:mkdv-target-coordinate} \mg_{\sigma,\kappa}:=g_\kappa(\,\cdot\,;q_\sigma), \qquad h_{\sigma,\kappa}:=h_\kappa(\,\cdot\,;q_\sigma) =\bigl(2\mg_{\sigma,\kappa}\bigr)^{-1}.
% \end{align}
% Lemma~\ref{lem:mkdv-static-reciprocal} collects the required spatial identities and mapping properties. The factorized operators are those of \cite[Proposition~7.1]{Forlano.Killip.Visan2024}, so the diagonal formula \cite[Corollary~7.2, equation~(7.8)]{Forlano.Killip.Visan2024} applies, and combining it with the spatial identities \eqref{eq:AKNS-static-identities} gives the formula used to justify the chain rule,
% \begin{align}
%  \label{eq:intro-mkdv-reciprocal-bridge} \mg_{\sigma,\kappa} =\frac{1+\gamma_\kappa+\sigma g_{+,\kappa}}{2\kappa}, \qquad h_{\sigma,\kappa} =\frac{\kappa}{1+\gamma_\kappa+\sigma g_{+,\kappa}} =\kappa\,\frac{1+\gamma_\kappa-\sigma g_{+,\kappa}}{1-g_{-,\kappa}^2}, \qquad \sigma\in\{+1,-1\}.
% \end{align}
% The second expression is the decisive one. Its numerator lies in $1+H^1$, while the inversion has been moved onto a denominator built from $g_{-,\kappa}\in H^2$, and the pointwise bound $\abs{g_{-,\kappa}}<1$ upgrades to the uniform gap \eqref{eq:AKNS-uniform-gap} on a compact time interval. Thus the formula writes the reciprocal as a product of multipliers that act boundedly on the negative Sobolev spaces that contain the time derivatives of $g_{-,\kappa}$, $\gamma_\kappa$, and $g_{+,\kappa}$.

\item \emph{The reciprocal formula}\ \ (Section~\ref{sec:mkdv-transfer}). For the Miura potentials, write
\begin{align}
 \label{eq:mkdv-target-coordinate} \mg_{\sigma,\kappa}:=g_\kappa(\,\cdot\,;q_\sigma), \qquad h_{\sigma,\kappa}:=h_\kappa(\,\cdot\,;q_\sigma) =\bigl(2\mg_{\sigma,\kappa}\bigr)^{-1}.
\end{align}
The diagonal formula of \cite[Corollary~7.2]{Forlano.Killip.Visan2024} and the identities \eqref{eq:AKNS-static-identities} give
\begin{align}
 \label{eq:intro-mkdv-reciprocal-bridge} \mg_{\sigma,\kappa} =\frac{1+\gamma_\kappa+\sigma g_{+,\kappa}}{2\kappa}, \qquad h_{\sigma,\kappa} =\frac{\kappa}{1+\gamma_\kappa+\sigma g_{+,\kappa}} =\kappa\,\frac{1+\gamma_\kappa-\sigma g_{+,\kappa}}{1-g_{-,\kappa}^2}, \qquad \sigma\in\{+1,-1\}.
\end{align}
The last expression places the reciprocal factor on $1-g_{-,\kappa}^2$, with $g_{-,\kappa}\in H^2$. The uniform gap \eqref{eq:AKNS-uniform-gap} on compact time intervals then gives the multipliers needed for the time derivative. Lemma~\ref{lem:mkdv-static-reciprocal} records these properties.

\item \emph{An affine chain rule}\ \ (Sections~\ref{sec:mkdv-chain-rule} and~\ref{sec:mkdv-transfer}). Lemma~\ref{lem:affine-negative-chain} differentiates the reciprocal factor $(1-g_{-,\kappa}^2)^{-1}$ in $H^{-1}$ and its product with an affine $H^1$ factor in $H^{-2}$. Lemma~\ref{lem:mkdv-time-regularity} applies this rule to \eqref{eq:intro-mkdv-reciprocal-bridge}. The evolution of $g_{+,\kappa}$ follows from the $g_{-,\kappa}$ identity and $g_{-,\kappa,x}=-2\kappa g_{+,\kappa}$, while the second defining identity supplies the evolution of $\gamma_\kappa$. This yields \eqref{eq:weak-h-evolution}.

% \item \emph{An affine chain rule with negative-order time derivatives}\ \ (Sections~\ref{sec:mkdv-chain-rule} and~\ref{sec:mkdv-transfer}). Lemma~\ref{lem:affine-negative-chain} converts \eqref{eq:intro-mkdv-reciprocal-bridge} into a dynamical statement. Given $g_-\in C(I;H^2)$ whose increments are an $H^{-1}$-valued time integral, together with the gap \eqref{eq:negative-chain-gap}, it gives the reciprocal factor $F=(1-g_-^2)^{-1}$ with $\partial_tF=2g_-(1-g_-^2)^{-2}V$ in $H^{-1}$. Given moreover an affine factor $P\in1+C(I;H^1)$ whose increments are an $H^{-2}$-valued time integral, it gives the integrated product law \eqref{eq:negative-chain-product-integrated}. Lemma~\ref{lem:mkdv-time-regularity} carries this out, with $V=V_-$, $Z=V_\gamma-\sigma V_+$ and the expressions \eqref{eq:AKNS-vector-fields}--\eqref{eq:AKNS-vector-p}, and yields the integrated reciprocal evolution \eqref{eq:weak-h-evolution} for the expression \eqref{eq:AKNS-h-vector-field}. There the $g_{+,\kappa}$-equation is obtained from the $g_{-,\kappa}$-identity \eqref{eq:intro-mkdv-green-m} by one spatial derivative together with $g_{-,\kappa,x}=-2\kappa g_{+,\kappa}$, so both defining identities of Definition~\ref{def:mkdv-green} are used, and $\gamma_\kappa$ enters through $V_\gamma$.

\item \emph{Identification of the flux}\ \ (Section~\ref{sec:mkdv-identification}). Lemma~\ref{lem:mkdv-smooth-identification} identifies the resulting vector field with $\partial_x\J_\kappa(q_\sigma)$ for smooth compactly supported $u$, using the smooth evolution laws in \cite[equations~(3.72)--(3.74)]{Forlano.Killip.Visan2024} and the Miura intertwining \eqref{eq:classical-Miura-intertwining}. Lemmas~\ref{lem:mkdv-AKNS-continuity} and~\ref{lem:mkdv-map-continuity} prove continuity of both expressions from $L^2$ into $H^{-2}$. Density therefore extends this identity between maps of a fixed $u$ to all of $L^2$. Combining it with the chain rule proves Proposition~\ref{prop:mkdv-transfer}.

% \item \emph{Identification of the flux}\ \ (Section~\ref{sec:mkdv-identification}). It remains to identify the expression just obtained with $\partial_x\J_\kappa(q_\sigma)$. Both sides are maps of a fixed function $u$, so we prove their equality in two steps. Lemma~\ref{lem:mkdv-smooth-identification} proves the identity for compactly supported smooth $u$, using the evolution equations for $\gamma,g_-,g_+$ in \cite[equations~(3.72)--(3.74)]{Forlano.Killip.Visan2024}, the classical intertwining \eqref{eq:classical-Miura-intertwining}, and Lemma~\ref{lem:microscopic}\,\textup{(ii)}. The sign $\sigma=-1$ follows from $\sigma=+1$ applied to $-u$. Lemmas~\ref{lem:mkdv-AKNS-continuity} and~\ref{lem:mkdv-map-continuity} show that both maps are continuous from $L^2$ into $H^{-2}$, so the identity extends to every $u\in L^2$ by density. Together with Lemma~\ref{lem:mkdv-time-regularity}, this proves Proposition~\ref{prop:mkdv-transfer}: both $q_+$ and $q_-$ are KdV Green solutions in the sense of Definition~\ref{def:green}. The density argument approximates a fixed $L^2$ function by smooth functions inside an identity between two continuous maps.

\item \emph{Uniqueness and smoothing}\ \ (Section~\ref{sec:mkdv-uniqueness}). Theorem~\ref{thm:main} identifies the two Miura transforms of competing solutions. Their difference gives $\partial_x(u_1-u_2)=0$, hence $u_1=u_2$ in $L^2$. Lemma~\ref{lem:canonical-mkdv-green} then identifies this solution with the canonical flow. Independently, Proposition~\ref{prop:auto-smoothing} applies to both transforms, and \eqref{eq:Miura-parallelogram} yields Proposition~\ref{prop:mkdv-smoothing}.

% \item \emph{Uniqueness, identification, and automatic mKdV smoothing}\ \ (Section~\ref{sec:mkdv-uniqueness}). Theorem~\ref{thm:main} is now applied twice, once for each sign, to the two competitors, which share their Miura data at time zero. The difference of the two conclusions gives $\partial_x(u_1-u_2)=0$, and an $L^2$ distribution with vanishing derivative is zero. Lemma~\ref{lem:canonical-mkdv-green} shows that the canonical $L^2$ flow of \cite{Harrop-Griffiths.Killip.Visan2020} lies in the intrinsic class, so the same comparison identifies every mKdV Green solution with the canonical trajectory. The smoothing conclusion follows by applying Proposition~\ref{prop:auto-smoothing} to both Miura transforms, and the parallelogram identity \eqref{eq:Miura-parallelogram} converts their cell bounds and translated-cell tails into \eqref{eq:intro-mkdv-smoothing}. This is Proposition~\ref{prop:mkdv-smoothing}.
\end{enumerate}

The mKdV argument follows the chain
\begin{align}
 \label{eq:roadmap-mkdv-implication-chain}
 &\text{Both mKdV Green identities for every }\kappa\geq1\\
 &\quad\Rightarrow \text{integrated evolution of }h_{\sigma,\kappa}\text{ using }g_{-,\kappa}\in H^2 \notag\\
 &\quad\Rightarrow q_+\text{ and }q_-\text{ are KdV Green solutions} \notag\\
 &\quad\Rightarrow \text{KdV uniqueness and local smoothing for each }q_\sigma \notag\\
 &\quad\Rightarrow \text{uniqueness and local smoothing for }u. \notag
\end{align}

\subsection{Organization of the Paper}

Section~\ref{sec:preliminaries} collects the analytic preliminaries. Section~\ref{sec:solution-notion} constructs the Green functions and justifies the KdV Green formulation. Section~\ref{sec:alpha} proves the large-parameter estimates, conservation of $\alpha_\lambda$, and coercivity of the Green flux. Section~\ref{sec:automatic-smoothing} establishes automatic local smoothing and spatial tail estimates. Section~\ref{sec:Gardner} constructs the associated Gardner solutions and completes the KdV uniqueness proof. Section~\ref{sec:mkdv} uses the two Miura transforms to prove uniqueness and automatic local smoothing for real defocusing mKdV Green solutions.

\section{Preliminaries}\label{sec:preliminaries}
%%%%%%%%%%%%%%%%%%%%%%%%%%%%%%%%%%%%%%%%%%%%%%%%%%%%%%%%%%%%%%%%%%%%%%%%%%%%%%%%%%

%%%%%%%%%%%%%%%%%%%%%%%%%%%%%%%%%%%%%%%%%%%%%%%%%%%%%%%%%%%%%%%%%%%%%%%%%%%%%%%%%%
\subsection{Notation and Conventions}
%%%%%%%%%%%%%%%%%%%%%%%%%%%%%%%%%%%%%%%%%%%%%%%%%%%%%%%%%%%%%%%%%%%%%%%%%%%%%%%%%%

Unless a domain is displayed, spatial function spaces are over $\R$. The solutions $q,u,r$ are real-valued. Operator domains and quadratic-form tests may be complex-valued. Time integrals in Banach spaces are Bochner integrals. For a self-adjoint operator $A$ on $L^2(\R)$ we write $A>0$, and call $A$ strictly positive, when $A\geq c$ for some constant $c>0$, that is, when $\inf\sigma(A)>0$; such an operator is boundedly invertible. 

We use the unitary Fourier transform . For the scaling calculation we use the equivalent norm
\begin{align}
 \label{eq:Hminus-norm} \norm{f}_{H^{-1}}^2 :=\int_{\R}\frac{\abs{\widehat f(\xi)}^2}{4+\xi^2}\dd \xi.
\end{align}
Equivalent Sobolev norms are used interchangeably in estimates, with fixed comparison constants understood, while exact Fourier identities use the normalization stated locally.

We write $\pair*{F}{\phi}$ for the bilinear distributional pairing, which equals $\int_\R F\phi\dd x$ when integrable, and $\ip{f}{g}_{L^2}:=\int_\R f\overline g\dd x$ for the complex $L^2$ inner product, linear in the first argument. We abbreviate the spatial norm $\norm{f}_{L^p(\R)}$ by $\norm{f}_p$. Spacetime norm domains are displayed when needed.

For non-negative quantities $A,B$, we write
\begin{align}
 \label{eq:lesssim-asymp} A\lesssim B \ :\Longleftrightarrow\ A\leq CB, \qquad A\simeq B \ :\Longleftrightarrow\ A\lesssim B\ \text{ and }\ B\lesssim A, \ \text{ i.e. }\ C^{-1}B\leq A\leq CB,
\end{align}
for a finite constant $C>0$. Implicit constants may depend on fixed parameters, time intervals, or compact sets identified in the statement or argument. Uniformity in a varying spectral parameter or translation centre is stated where it is used.

%%%%%%%%%%%%%%%%%%%%%%%%%%%%%%%%%%%%%%%%%%%%%%%%%%%%%%%%%%%%%%%%%%%%%%%%%%%%%%%%%%
\subsection{Some Preliminary Lemmas}
%%%%%%%%%%%%%%%%%%%%%%%%%%%%%%%%%%%%%%%%%%%%%%%%%%%%%%%%%%%%%%%%%%%%%%%%%%%%%%%%%%

The following elementary facts are used repeatedly to interpret products of an $H^{-1}$ potential with $g_\lambda$ and its reciprocal, justify reciprocal multipliers, and regard integrable remainders as $H^{-1}$ distributions. The products with the potential are defined by Sobolev duality.

{\begin{lemma}\label{lem:products} In one dimension:
\begin{enumerate}[label=\textup{(\roman*)}, leftmargin=*, itemsep=.2em]
\item $H^1(\R)$ is a Banach algebra and embeds continuously into $C_0(\R)\subset L^\infty(\R)$.
{Here $C_0$ denotes the space of continuous functions that tend to zero as $x\to\pm\infty$.}
\item Multiplication by a given $H^1(\R)$ function is a continuous operator from $H^{-1}(\R)$ to $H^{-1}(\R)$.
\item $L^1(\R)$ embeds continuously into $H^{-1}(\R)$.
\end{enumerate}
% {Here $C_0$ denotes the space of continuous functions that tend to zero as $x\to\pm\infty$.}
\end{lemma}}

% {
% \begin{remark}
% Here $C_0$ denotes the space of continuous functions that tend to zero as $x\to\pm\infty$.
% \end{remark}
% }

\begin{proof}
The one-dimensional Sobolev inequality gives $H^1(\R)\hookrightarrow L^\infty(\R)$ and the standard product estimate
\begin{align}
 \label{eq:products-a1} \norm{fg}_{H^1} \lesssim \norm{f}_{H^1}\norm{g}_{H^1}.
\end{align}
Hence $H^1(\R)$ is a Banach algebra.

We also recall why the continuous representative of an $H^1$ function vanishes at spatial infinity. Every $f\in H^1(\R)$ has a locally absolutely continuous representative satisfying
\begin{align}
 \label{eq:products-a2} \abs{f(x)-f(y)} \leq \norm{f_x}_2\abs{x-y}^{\frac12}.
\end{align}
Hence this representative is uniformly continuous. If it did not tend to zero as $x\to+\infty$, there would exist $\varepsilon>0$ and a sequence $x_n\to+\infty$ such that $\abs{f(x_n)}\geq\varepsilon$. Uniform continuity would then give a fixed $\delta>0$ such that $\abs{f(x)}\geq\varepsilon/2$ whenever $\abs{x-x_n}<\delta$. After passing to a subsequence, these intervals are pairwise disjoint, which would imply
\begin{align}
 \label{eq:products-a3} \norm{f}_2^2 \geq \sum_n\int_{x_n-\delta}^{x_n+\delta}\abs{f(x)}^2\dd x =\infty,
\end{align}
a contradiction. The argument at $-\infty$ is identical. Therefore $f\in C_0(\R)$.

For the second assertion, let $a\in H^1$, $F\in H^{-1}$, and $\phi\in H^1$. Define multiplication by duality:
\begin{align}
 \label{eq:products-a4} \pair*{aF}{\phi}:=\pair*{F}{a\phi}.
\end{align}
The algebra estimate gives
\begin{align}
 \label{eq:products-a5} \abs{\pair*{aF}{\phi}} \lesssim \norm{F}_{H^{-1}}\norm{a\phi}_{H^1} \lesssim \norm{a}_{H^1}\norm{F}_{H^{-1}}\norm{\phi}_{H^1}.
\end{align}
Therefore $aF\in H^{-1}$ continuously.

Finally, for $f\in L^1$,
\begin{align}
 \label{eq:products-a6} \abs{\pair*{f}{\phi}} \leq \norm{f}_{L^1}\norm{\phi}_{L^\infty} \lesssim \norm{f}_{L^1}\norm{\phi}_{H^1},
\end{align}
which proves $L^1(\R)\hookrightarrow H^{-1}(\R)$.
\end{proof}

% We pause to recall the quadratic-form terminology used in the next lemma. 
We recall the quadratic-form terminology needed to construct $L_q$.
{All forms below act on the Hilbert space $L^2(\R)$. Their two-argument versions are sesquilinear, linear in the first argument, with $\mathfrak a[\phi]=\mathfrak a[\phi,\phi]$. A densely defined symmetric quadratic form $\mathfrak a$ is \emph{lower semibounded} if there is a constant $C\geq0$ such that
\begin{align}
 \label{eq:form-lower-bound-definition} \mathfrak a[\phi]\geq-C\norm{\phi}_2^2 \qquad\text{for every }\phi\in \Dom(\mathfrak a).
\end{align}
It is \emph{closed} if its domain is complete in the shifted form norm
\begin{align}
 \label{eq:closed-form-norm} \norm{\phi}_{\mathfrak a,C}^2 :=\mathfrak a[\phi]+(C+1)\norm{\phi}_2^2.
\end{align}
}
The free Dirichlet form
\begin{align}
 \label{eq:free-Dirichlet-form} \mathfrak a_0[\phi]=\norm{\phi_x}_2^2, \qquad \Dom(\mathfrak a_0)=H^1(\R),
\end{align}
is densely defined, non-negative, and closed.

We use the following special case of the Kato--Lions--Lax--Milgram--Nelson theorem, usually called the KLMN theorem:
{\begin{lemma}\label{lem:KLMN}
If $\mathfrak b$ is a symmetric form on $H^1(\R)$ and
\begin{align}
 \label{eq:KLMN-hypothesis} \abs{\mathfrak b[\phi]} \leq\vartheta\mathfrak a_0[\phi]+C_\vartheta\norm{\phi}_2^2 \qquad\text{for some }0\leq\vartheta<1,
\end{align}
then $\mathfrak a_0+\mathfrak b$, with domain $H^1(\R)$, is again closed and lower semibounded, and there is a unique self-adjoint operator $A$ on $L^2(\R)$ whose quadratic form it is, in the sense that
\begin{align}
 \label{eq:KLMN-representation} \Dom(A)\subset H^1(\R), \qquad (\mathfrak a_0+\mathfrak b)[\phi,\psi]=\ip{A\phi}{\psi}_{L^2} \quad\text{for }\phi\in\Dom(A),\ \psi\in H^1(\R).
\end{align}
\end{lemma}
}
This is the KLMN construction of \cite[Theorem~X.17]{Reed.Simon1975}, used for the Schr\"odinger operator in \cite[Proposition~2.1]{Killip.Visan2019}. 
The next lemma verifies its hypothesis for the potential form $\mathfrak b_q[\phi]:=\pair* q{\abs{\phi}^2}$.
% This quoted result is \cite[Theorem~X.17]{Reed.Simon1975}, in exactly the form in which Killip--Vi\c{s}an \cite[Proposition~2.1]{Killip.Visan2019} construct the Schr\"odinger operator.
% We verify its hypothesis for the potential form $\mathfrak b_q[\phi]=\pair* q{\abs{\phi}^2}$.

{\begin{lemma}\label{lem:form-bound} For every real-valued $q\in H^{-1}(\R)$ and every $\varepsilon>0$, there is a finite constant $C_{\varepsilon,q}$ such that
\begin{align}
 \label{eq:infinitesimal-form-bound} \abs[\big]{\pair* q{\abs{\phi}^2}} \leq \varepsilon\norm{\phi_x}_2^2 +C_{\varepsilon,q}\norm{\phi}_2^2, \qquad \phi\in H^1(\R).
\end{align}
Consequently, the quadratic form $\mathfrak l_q$ of \eqref{eq:intro-form} is closed and lower semibounded on $L^2(\R)$, with domain $H^1(\R)$.
\end{lemma}}

\begin{proof}
Fix $\varepsilon>0$. Since real-valued test functions are dense in the real subspace of $H^{-1}$, we may choose a real $q_0\in C_c^\infty(\R)$ such that $\norm{q-q_0}_{H^{-1}}\leq\delta$, where $\delta>0$ will be chosen below. By duality and the $H^1$ algebra estimate in Lemma~\ref{lem:products},
\begin{align}
 \label{eq:form-bound-a1} \abs[\big]{\pair*{q-q_0}{\abs{\phi}^2}} &\lesssim\norm{q-q_0}_{H^{-1}}\norm[\big]{\abs{\phi}^2}_{H^1}\leq C\delta\norm{\phi}_{H^1}^2 =C\delta\norm{\phi_x}_2^2+C\delta\norm{\phi}_2^2.
\end{align}
For the smooth part,
\begin{align}
 \label{eq:form-bound-a2} \abs[\big]{\pair*{q_0}{\abs{\phi}^2}} \leq\norm{q_0}_\infty\norm{\phi}_2^2.
\end{align}
Combining the two estimates yields
\begin{align}
 \label{eq:form-bound-before-delta} \abs[\big]{\pair* q{\abs{\phi}^2}} \leq C\delta\norm{\phi_x}_2^2 +\bigl(C\delta+\norm{q_0}_\infty\bigr)\norm{\phi}_2^2.
\end{align}
Choose $\delta$ so that $C\delta\leq\varepsilon$ and set
\begin{align}
 \label{eq:form-bound-a3} C_{\varepsilon,q}:=C\delta+\norm{q_0}_\infty.
\end{align}
This proves \eqref{eq:infinitesimal-form-bound}.

To obtain the final assertion, apply the estimate with, for example, $\varepsilon=\frac12$. It first gives the explicit lower bound
\begin{align}
 \label{eq:form-bound-a4} \mathfrak l_q[\phi] \geq\frac12\norm{\phi_x}_2^2-C_{\frac12,q}\norm{\phi}_2^2 \geq-C_{\frac12,q}\norm{\phi}_2^2.
\end{align}
Moreover, the perturbation $\mathfrak b_q[\phi]=\pair* q{\abs{\phi}^2}$ satisfies \eqref{eq:KLMN-hypothesis} with relative form bound $\frac12<1$. The quoted KLMN theorem in Lemma~\ref{lem:KLMN} therefore shows that $\mathfrak l_q$ is closed and lower semibounded on $L^2(\R)$, with form domain $H^1(\R)$. 
% This is also the construction used in \cite[Proposition~2.1]{Killip.Visan2019}.
\end{proof}

% {\begin{remark}
% 
% The reciprocal $g_\lambda^{-1}$ used later is a nonzero constant plus an $H^1$ function: $g_\lambda^{-1}=2\lambda+a_0$ with $a_0\in H^1$. We therefore record explicitly how functions in the affine space
The reciprocal $g_\lambda^{-1}$ used later has the form $g_\lambda^{-1}= 2\lambda+a_0$ with $a_0\in H^1$. To interpret its product with $q$, we describe how functions in the affine space
\begin{align}
 \label{eq:affine-Hone-algebra} \mathcal A:=\R+H^1(\R) =\big\{c+a_0:c\in\R,\ a_0\in H^1(\R)\big\}
\end{align}
act on $H^{-1}$.

\begin{enumerate}[
label=\textup{(\roman*)}, leftmargin=*, itemsep=.6em, before={\let\fullwidthdisplay\relax} ]
\item \emph{Multiplication.} Let $a=c+a_0\in\mathcal A$ and $F\in H^{-1}$. For $\phi\in H^1$, define
\begin{align}
 \label{eq:affine-multiplier-definition} \pair*{aF}{\phi}:=\pair* F{a\phi}.
\end{align}
Since $a\phi=c\phi+a_0\phi\in H^1$, Lemma~\ref{lem:products} gives
\begin{align}
 \label{eq:multipliers-a1} \norm{aF}_{H^{-1}} \lesssim\bigl(\abs{c}+\norm{a_0}_{H^1}\bigr)\norm{F}_{H^{-1}}.
\end{align}
Hence \eqref{eq:affine-multiplier-definition} is well-defined and is equivalent to writing $aF=cF+a_0F$.

\item \emph{Reciprocals.} Suppose in addition that $\abs{a(x)}\geq c_0>0$ for every $x$, where $c_0$ is a constant. Because $a_0(x)\to0$ as $\abs{x}\to\infty$, necessarily $c\neq0$. Moreover,
\begin{align}
 \label{eq:affine-reciprocal} a^{-1}=c^{-1}+b_0, \qquad b_0:=a^{-1}-c^{-1}=-\frac{a_0}{ca}.
\end{align}
The bounds
\begin{align}
 \label{eq:multipliers-a2} \norm{b_0}_2 &\leq \abs{c}^{-1}c_0^{-1}\norm{a_0}_2,\qquad\norm{(b_0)_x}_2 =\norm*{\frac{a_x}{a^2}}_2 \leq c_0^{-2}\norm{(a_0)_x}_2
\end{align}
show that $b_0\in H^1$. Hence a uniformly nonvanishing element of $\mathcal A$ has its reciprocal in $\mathcal A$.

\item \emph{Associativity and cancellation.} If $a,b\in\mathcal A$, then $ab\in\mathcal A$ by the $H^1$ algebra property. For $F\in H^{-1}$ and $\phi\in H^1$, the definition above gives
\begin{align}
 \label{eq:multipliers-a3} \pair*{a(bF)}{\phi} =\pair*{bF}{a\phi} =\pair* F{b(a\phi)} =\pair* F{(ab)\phi} =\pair*{(ab)F}{\phi}.
\end{align}
Therefore $a(bF)=(ab)F$ in $H^{-1}$. If $a$ is uniformly nonvanishing, part~\textup{(ii)} permits $b=a^{-1}$, and we obtain
\begin{align}
 \label{eq:multipliers-a4} a(a^{-1}F)=F.
\end{align}
\end{enumerate}
This convention is used whenever we write such expressions as $q/g_\lambda$, $h_\lambda^{-1}q$, or cancellations involving reciprocals in \eqref{eq:Drach} and the flux identities.
% 
% \end{remark}}

The spatial-tail argument in Proposition~\ref{prop:auto-smoothing} requires uniform decay when a translated $H^1$ weight is paired with a compact family in $H^{-1}$. For $y\in\R$, write
% The spatial-tail argument in Proposition~\ref{prop:auto-smoothing} pairs derivatives of a translated weight with $H^{-1}$ distributions.
% For $y\in\R$, write
\begin{align}
 \label{eq:translation-operator} (T_y\psi)(x):=\psi(x-y).
\end{align}

{\begin{lemma}
\label{lem:translation}
Let $\mathcal C\subset H^{-1}(\R)$ and $\mathcal D\subset H^1(\R)$ be compact. Then
\begin{align}
 \label{eq:translated-compact-pairing} \lim_{\abs{y}\to\infty} \sup_{F\in\mathcal C}\sup_{\psi\in\mathcal D} \abs[\big]{\pair* F{T_y\psi}}=0.
\end{align}
\end{lemma}}

\begin{proof}
We first prove weak decay for one fixed profile. For $\psi,\varphi\in H^1(\R)$, the usual $H^1$ inner product satisfies
\begin{align}
 \label{eq:translated-Hone-inner-product} (T_y\psi,\varphi)_{H^1}=\int_\R \ue^{-\ui y\xi}(1+\xi^2)\widehat\psi(\xi)\overline{\widehat\varphi(\xi)}\dd \xi.
\end{align}
The nonoscillatory factor belongs to $L^1_\xi$, because
\begin{align}
 \label{eq:translated-Hone-integrability} \int_\R\big(1+\xi^2\big)\abs{\widehat\psi(\xi)}\,\abs{\widehat\varphi(\xi)}\dd \xi\leq\norm{\psi}_{H^1}\norm{\varphi}_{H^1}.
\end{align}
The Riemann--Lebesgue lemma therefore gives
\begin{align}
 \label{eq:translated-Hone-weak-convergence} T_y\psi\rightharpoonup0 \quad\text{weakly in }H^1(\R) \quad\text{as }\abs{y}\to\infty.
\end{align}
As a consequence, for fixed $F\in H^{-1}$ and fixed $\psi\in H^1$, $\pair* F{T_y\psi}\to0$.

It remains to make this convergence uniform in both factors. Since the two sets are compact, they are bounded. Set
\begin{align}
 \label{eq:translated-compact-bounds} M_{-1}:=\sup_{F\in\mathcal C}\norm{F}_{H^{-1}}, \qquad M_1:=\sup_{\psi\in\mathcal D}\norm{\psi}_{H^1}.
\end{align}
Fix $\delta>0$. Choose finite $\delta$-nets $F_1,\dots,F_N\in\mathcal C$ and $\psi_1,\dots,\psi_M\in\mathcal D$ in $H^{-1}$ and $H^1$, respectively. For any $F\in\mathcal C$ and $\psi\in\mathcal D$, choose $j,k$ such that $\norm{F-F_j}_{H^{-1}}<\delta$ and $\norm{\psi-\psi_k}_{H^1}<\delta$.

Since $(4+\xi^2)^{1/2}\leq2(1+\xi^2)^{1/2}$, the norm convention \eqref{eq:Hminus-norm} gives $\abs{\pair* F\psi}\leq2\norm{F}_{H^{-1}}\norm{\psi}_{H^1}$. Translation invariance of the $H^1$ norm therefore yields
\begin{align}
 \label{eq:translated-two-net-estimate} \abs[\big]{\pair* F{T_y\psi}} &\leq \abs[\big]{\pair*{F_j}{T_y\psi_k}} +\abs[\big]{\pair*{F-F_j}{T_y\psi}} +\abs[\big]{\pair*{F_j}{T_y(\psi-\psi_k)}} \\
 &\leq \max_{\substack{1\leq j\leq N\\
 1\leq k\leq M}} \abs[\big]{\pair*{F_j}{T_y\psi_k}} +2\delta(M_1+M_{-1}).\notag
\end{align}
For each of the finitely many pairs $(j,k)$, the first term tends to zero by \eqref{eq:translated-Hone-weak-convergence}. Taking first $\abs{y}\to\infty$ in \eqref{eq:translated-two-net-estimate} and then $\delta\downarrow0$ proves \eqref{eq:translated-compact-pairing}.
\end{proof}

% {
% \begin{remark}
% In the notation of the proof of Proposition~\ref{prop:auto-smoothing}, in Section~\ref{sec:automatic-smoothing}, this lemma applies to $\phi_y'$, $\phi_y''$, and $\phi_y'''$, because they are translates of fixed $H^1$ functions. The weight $\phi_y$ of \eqref{eq:tanh-weight} lies outside $H^1(\R)$. Its pairing is treated there using conservation of the density mass and dominated convergence.
In Proposition~\ref{prop:auto-smoothing}, this lemma applies to $\phi_y'$, $\phi_y''$, and $\phi_y'''$, which are translates of fixed $H^1$ functions. The weight $\phi_y$ of \eqref{eq:tanh-weight} itself lies outside $H^1$. Its pairing with the density is treated there using mass conservation and dominated convergence.
% \end{remark}
% }

The following elementary observation uses cutoffs to show that an integrable derivative of an $H^{-1}$ distribution has zero integral. Applied to the two-parameter density identity, it yields an exact conservation law.

{\begin{lemma}\label{lem:zero-mass} If $f\in L^1(\R)$, $F\in H^{-1}(\R)$, and $f=\partial_xF$ in distributions, then $\int_\R f=0$.
\end{lemma}}

\begin{proof}
Choose $\chi\in C_c^\infty$ equal to one near the origin and put $\chi_R(x)=\chi(x/R)$. Dominated convergence gives $\int f\chi_R\to\int f$. On the other hand,
\begin{align}
 \label{eq:zero-mass-a1} \int f\chi_R=-\pair* F{\chi_R'}, \qquad \norm{\chi_R'}_{H^1}=O(R^{-\frac12}),
\end{align}
so the same integrals tend to zero.
\end{proof}

% {
% \begin{remark}
%     This seems to indicate $F\in H^{-1}\cap W^{1,1}$ and $f\in H^{-2}\cap L^1$.
% \end{remark}
% }

% \newpage

%%%%%%%%%%%%%%%%%%%%%%%%%%%%%%%%%%%%%%%%%%%%%%%%%%%%%%%%%%%%%%%%%%%%%%%%%%%%%%%%%%
\section{Green Solution Notion}\label{sec:solution-notion}
%%%%%%%%%%%%%%%%%%%%%%%%%%%%%%%%%%%%%%%%%%%%%%%%%%%%%%%%%%%%%%%%%%%%%%%%%%%%%%%%%%

This section justifies the Green-solution formulation stated in Definition~\ref{def:green}. We construct $g_\lambda$ and its reciprocal for rough potentials, verify that its flux is well-defined, and derive the defining identity for smooth KdV solutions. We then record the common spectral range and the integrated identities needed in the uniqueness proof.

The form $\mathfrak l_q$ and the operator $L_q$ are defined in \eqref{eq:intro-form}. Lemma~\ref{lem:products} makes the potential term meaningful, and Lemma~\ref{lem:form-bound} verifies the hypothesis of Lemma~\ref{lem:KLMN}. Thus the form determines a unique self-adjoint, lower-semibounded operator on $L^2(\R)$, as in \cite[Proposition~2.1]{Killip.Visan2019}.

Recall from \eqref{eq:intro-resolvent} that $R_\lambda(q)=(L_q+\lambda^2)^{-1}$ and $g_\lambda(x;q)=G_\lambda(x,x;q)$, where $G_\lambda$ is the resolvent kernel. Proposition~\ref{prop:static} below proves that this kernel is continuous for sufficiently large $\lambda$ and gives the mapping properties of its diagonal.

%%%%%%%%%%%%%%%%%%%%%%%%%%%%%%%%%%%%%%%%%%%%%%%%%%%%%%%%%%%%%%%%%%%%%%%%%%%%%%%%%%
\subsection{Identities of Green Functions}
%%%%%%%%%%%%%%%%%%%%%%%%%%%%%%%%%%%%%%%%%%%%%%%%%%%%%%%%%%%%%%%%%%%%%%%%%%%%%%%%%%

The next lemma is the form-level scaling that transfers the small-data resolvent theory at parameter $1$ of \cite{Killip.Visan2019} to arbitrary bounded subsets of $H^{-1}$ for sufficiently large $\lambda$. In words, it says that dilating the potential $q$ by $\varepsilon$ trades the spectral parameter $\lambda=\varepsilon^{-1}$ for the fixed parameter $1$: every statement about $q$ at parameter $\varepsilon^{-1}$ is equivalent to the corresponding statement about $q^{(\varepsilon)}$ at parameter $1$, and conversely. We first fix the two maps involved.

For $0<\varepsilon\leq1$, define the dilated potential $q^{(\varepsilon)}$ and the dilation operator $U_\varepsilon$ by
\begin{align}
 \label{eq:dilation-def} \bpair{q^{(\varepsilon)}}{\phi} :=\varepsilon\pair* q{\phi(\varepsilon^{-1}\cdot)} \quad(\phi\in H^1), \qquad (U_\varepsilon f)(x):=\varepsilon^{\frac12}f(\varepsilon x).
\end{align}
When $q$ is a function, the first formula is the pointwise rescaling $q^{(\varepsilon)}(x)=\varepsilon^2q(\varepsilon x)$.

{\begin{lemma}\label{lem:dilation}
Let $q\in H^{-1}(\R;\R)$, let $0<\varepsilon\leq1$, and let $q^{(\varepsilon)}$ and $U_\varepsilon$ be given by \eqref{eq:dilation-def}. Then the following hold.
\begin{enumerate}[
label=\textup{(\roman*)}, leftmargin=*, itemsep=.7em, before={\let\fullwidthdisplay\relax} ]
\item \emph{Forms and operators.} The dilated potential satisfies $q^{(\varepsilon)}\in H^{-1}(\R;\R)$ with
\begin{align}
 \label{eq:dilation-Hminus} \norm{q^{(\varepsilon)}}_{H^{-1}}^2 =\varepsilon^3\int_\R \frac{\abs{\widehat q(\eta)}^2}{4+\varepsilon^2\eta^2}\dd \eta \leq\varepsilon\norm{q}_{H^{-1}}^2,
\end{align}
the operator $U_\varepsilon$ is unitary on $L^2(\R)$ with $U_\varepsilon^{-1}=U_{\varepsilon^{-1}}$, and
\begin{align}
 \label{eq:operator-scaling} L_{q^{(\varepsilon)}}+1 =\varepsilon^2U_\varepsilon (L_q+\varepsilon^{-2})U_\varepsilon^{-1}
\end{align}
as an identity between self-adjoint operators associated with closed lower-semibounded forms on $H^1(\R)$.

\item \emph{Resolvents.} The two positivity conditions are equivalent:
\begin{align}
 \label{eq:dilation-positivity} L_q+\varepsilon^{-2}>0 \qquad\Longleftrightarrow\qquad L_{q^{(\varepsilon)}}+1>0.
\end{align}
When they hold, both operators are boundedly invertible and
\begin{align}
 \label{eq:resolvent-scaling} (L_{q^{(\varepsilon)}}+1)^{-1} =\varepsilon^{-2}U_\varepsilon (L_q+\varepsilon^{-2})^{-1}U_\varepsilon^{-1}.
\end{align}

\item \emph{Kernels and diagonals.} Assume the equivalent conditions in \eqref{eq:dilation-positivity}. Then $(L_q+\varepsilon^{-2})^{-1}$ admits a continuous integral kernel if and only if $(L_{q^{(\varepsilon)}}+1)^{-1}$ does. In that case the two kernels are related by
\begin{align}
 \label{eq:G-scaling} G_1(x,y;q^{(\varepsilon)}) =\varepsilon^{-1}G_{\frac1\varepsilon}(\varepsilon x,\varepsilon y;q), \qquad x,y\in\R,
\end{align}
and their diagonals by
\begin{align}
 \label{eq:g-scaling} g_1(x;q^{(\varepsilon)}) =\varepsilon^{-1}g_{\frac{1}{\varepsilon}}(\varepsilon x;q), \qquad x\in\R.
\end{align}
\end{enumerate}
\end{lemma}}

Part~\textup{(i)} holds for every real $q\in H^{-1}$ and $0<\varepsilon\leq1$. Under the positivity conditions of part~\textup{(ii)}, the resolvent and kernel identities transfer in both directions. The proof of Proposition~\ref{prop:static} starts from the small-ball hypotheses at parameter $1$ and concludes at parameter $\lambda=\varepsilon^{-1}$. Lemma~\ref{lem:reciprocal-density-scaling} uses the reverse direction.

\begin{proof}
\textup{(i)} Since $\widehat{q^{(\varepsilon)}}(\xi)=\varepsilon\widehat q(\varepsilon^{-1}\xi)$, the substitution $\xi=\varepsilon\eta$ in \eqref{eq:Hminus-norm} gives the equality in \eqref{eq:dilation-Hminus}. The inequality follows from $\varepsilon\leq1$, because then $\varepsilon^2(4+\eta^2)\leq4+\varepsilon^2\eta^2$ for every $\eta$. Unitarity of $U_\varepsilon$ and the identity $U_\varepsilon^{-1}=U_{\varepsilon^{-1}}$ are changes of variables.

For \eqref{eq:operator-scaling}, let $\phi\in H^1$ and put $\psi=U_\varepsilon^{-1}\phi$. A change of variables in the kinetic term, together with \eqref{eq:dilation-def}, gives
\begin{align}
 \label{eq:dilation-form-identity} \mathfrak l_{q^{(\varepsilon)}}[\phi]+\norm{\phi}_2^2 =\varepsilon^2\left\{ \mathfrak l_q[\psi] +\varepsilon^{-2}\norm{\psi}_2^2\right\}.
\end{align}
Both sides are closed and lower semibounded on $H^1(\R)$ by Lemma~\ref{lem:form-bound}, and $\phi\mapsto\psi$ is a bijection of $H^1(\R)$ onto itself. The uniqueness of the self-adjoint operator associated with the form, stated in Lemma~\ref{lem:KLMN}, therefore gives \eqref{eq:operator-scaling}.

\textup{(ii)} Conjugation by a unitary operator and multiplication by the positive constant $\varepsilon^2$ preserve strict positivity, and both operations are invertible. Hence \eqref{eq:operator-scaling} gives \eqref{eq:dilation-positivity}, and inverting \eqref{eq:operator-scaling} gives \eqref{eq:resolvent-scaling}.

\textup{(iii)} If a bounded operator $A$ has integral kernel $K_A$, then $U_\varepsilon AU_\varepsilon^{-1}$ has integral kernel $\varepsilon K_A(\varepsilon x,\varepsilon y)$, which is continuous precisely when $K_A$ is. The same computation applied to $U_\varepsilon^{-1}(\,\cdot\,)U_\varepsilon$ gives the reverse implication. Applying this to $A=(L_q+\varepsilon^{-2})^{-1}$, whose kernel is $G_{\varepsilon^{-1}}(\cdot,\cdot;q)$ by \eqref{eq:intro-resolvent} with $\lambda=\varepsilon^{-1}$, and using \eqref{eq:resolvent-scaling}, gives \eqref{eq:G-scaling}. Setting $y=x$ gives \eqref{eq:g-scaling}.
\end{proof}

The construction of $g_\lambda$ and its reciprocal, together with their continuous dependence on $q$, follows from \cite[Propositions~2.1 and~2.2]{Killip.Visan2019} after the scaling in Lemma~\ref{lem:dilation}. The recovery formula is equation~(2.25) of \cite[Lemma~2.14]{Killip.Murphy.Visan2020}. Since that lemma is stated for a strictly admissible complex spectral parameter $k$ with $\operatorname{Im}(k^2)$ large, the next lemma gives the approximation argument needed at a real positive parameter.

The next lemma deduces the recovery identity from regularity and convergence assumptions on $g_\lambda$ and its reciprocal. Proposition~\ref{prop:static} supplies these assumptions: positivity follows from part~\textup{(i)}, and the $H^1$ convergence follows from the continuity in part~\textup{(ii)} applied to $q_n\to q$. Parts~\textup{(i)} and~\textup{(ii)} are proved from \cite[Proposition~2.2]{Killip.Visan2019} and Lemma~\ref{lem:dilation} alone. Only part~\textup{(iii)} uses the lemma below.

{\begin{lemma}
\label{lem:Drach-density}
Fix $\lambda>0$. Let $q\in H^{-1}(\R;\R)$, and assume that there exists a sequence $q_n\in C_c^\infty(\R;\R)$ such that $q_n\to q$ in $H^{-1}$. Suppose that $L_{q_n}+\lambda^2$ and $L_q+\lambda^2$ are strictly positive, and denote the corresponding diagonal Green functions by $g_n$ and $g$. Assume further that
\begin{align}
 \label{eq:Drach-density-coordinate-convergence} g-\frac1{2\lambda}\in H^1(\R), \qquad g_n-\frac1{2\lambda} \rightarrow g-\frac1{2\lambda} \quad\text{in }H^1(\R), \qquad g(x)>0\quad(x\in\R).
\end{align}
Then
\begin{align}
 \label{eq:Drach-density-conclusion} q=\partial_x\!\left(\frac{g_x}{2g}\right) +\left(\frac{g_x}{2g}\right)^2 +\frac1{4g^2}-\lambda^2 \quad\text{in }H^{-1}(\R).
\end{align}
\end{lemma}}

{
\begin{remark}
The approximating sequence in Lemma~\ref{lem:Drach-density} is used to pass the recovery identity from smooth potentials to $q$. In Proposition~\ref{prop:static}, parts~\textup{(i)} and~\textup{(ii)} supply the positivity and convergence hypotheses needed to apply this lemma in part~\textup{(iii)}.
\end{remark}
}

\begin{proof}
We first prove the identity for a smooth compactly supported potential by expressing the Green function through two solutions of \eqref{eq:Weyl-ode}, one decaying at each end of the line. We then pass to the limit along the sequence $q_n$ using \eqref{eq:Drach-density-coordinate-convergence}.

\proofstep{Step 1: the smooth case}
Fix $n$ and abbreviate $p=q_n$, $g=g_n$. Choose $R>0$ such that $\supp p\subset[-R,R]$. Consider the ordinary differential equation
\begin{align}
 \label{eq:Weyl-ode} -\psi''+p\psi=-\lambda^2\psi \qquad\text{on }\R.
\end{align}
Define $\psi_+(x)=\ue^{-\lambda x}$ for $x\geq R$ and $\psi_-(x)=\ue^{\lambda x}$ for $x\leq-R$. The existence and uniqueness theorem for linear ordinary differential equations extends each function uniquely to a solution on $\R$. Thus $\psi_+$ is square integrable near $+\infty$, and $\psi_-$ is square integrable near $-\infty$.

They are linearly independent. Otherwise a nonzero common multiple would be square integrable on all of $\R$ and would be an eigenfunction of $L_p$ at $-\lambda^2$, contradicting $L_p+\lambda^2>0$. Define their Wronskian by
\begin{align}
 \label{eq:Weyl-wronskian} W:=\psi_-'\psi_+-\psi_-\psi_+'.
\end{align}
The equation gives $W_x=\psi_-^{\prime\prime}\psi_+-\psi_-\psi_+^{\prime\prime}=0$, and linear independence gives $W\neq0$. The resolvent kernel is
\begin{align}
 \label{eq:Weyl-green} G_\lambda(x,y;p) =\frac{\psi_-(\min\{x,y\})\,\psi_+(\max\{x,y\})}{W}, \qquad\text{so }\ g=\frac{\psi_+\psi_-}{W}.
\end{align}
To verify the kernel formula, fix $y$. Its right-hand side solves \eqref{eq:Weyl-ode} for $x\neq y$, is continuous at $x=y$, and its first $x$-derivative has jump $-1$ there. Hence applying $L_p+\lambda^2$ gives $\delta_y$. The prescribed decay at both ends identifies the inverse kernel. For $p\equiv0$, $\psi_\pm=\ue^{\mp\lambda x}$, $W=2\lambda$, and $g=\frac1{2\lambda}$, which fixes the normalization.

Differentiating \eqref{eq:Weyl-green} twice and using \eqref{eq:Weyl-ode} in the form $\psi_\pm''=(p+\lambda^2)\psi_\pm$ gives
\begin{align}
 \label{eq:Drach-density-a1} g_x=\frac{\psi_+'\psi_-+\psi_+\psi_-'}{W}, \qquad g_{xx}=\frac{2(p+\lambda^2)\psi_+\psi_-+2\psi_+'\psi_-'}{W}.
\end{align}
Hence
\begin{align}
 \label{eq:Drach-density-a2} 2gg_{xx}-g_x^2 =\frac{4(p+\lambda^2)\psi_+^2\psi_-^2 -\bigl(\psi_+'\psi_--\psi_+\psi_-'\bigr)^2}{W^2} =4(p+\lambda^2)g^2-1,
\end{align}
where the last step used $\psi_+'\psi_--\psi_+\psi_-'=-W$. Restoring the subscript gives
\begin{align}
 \label{eq:Drach-density-smooth-first-integral} 2g_ng_{n,xx}-g_{n,x}^2 -4(q_n+\lambda^2)g_n^2=-1.
\end{align}
To justify division by $g_n^2$, write $\mathfrak a_p[z]:=\mathfrak l_p[z]+\lambda^2\norm{z}_2^2$. Strict positivity gives $\mathfrak a_p[z]\geq c\norm{z}_2^2$ for some $c>0$. With $M:=\norm{(p+\lambda^2)_-}_\infty$, we have $\norm{z}_{H^1}^2\leq\bigl(1+(M+1)/c\bigr)\mathfrak a_p[z]$. Thus the form is bounded and coercive on $H^1$, and Lax--Milgram extends the resolvent boundedly from $H^{-1}$ onto $H^1$, agreeing with its $L^2$ realization. For the real nonzero function $z=R_\lambda(p)\delta_x$, this gives $g_n(x)=\pair*{\delta_x}{z}=\mathfrak a_p[z]>0$. Dividing \eqref{eq:Drach-density-smooth-first-integral} by $4g_n^2$ and regrouping therefore gives the equivalent form
\begin{align}
 \label{eq:Drach-density-smooth-formula} q_n=\partial_x\!\left(\frac{g_{n,x}}{2g_n}\right) +\left(\frac{g_{n,x}}{2g_n}\right)^2 +\frac1{4g_n^2}-\lambda^2.
\end{align}

\proofstep{Step 2: passage to the limit}
We first note that the pointwise positivity of $g$ self-improves to a uniform bound. Indeed, $g-\frac1{2\lambda}\in H^1(\R)\hookrightarrow C_0(\R)$ by Lemma~\ref{lem:products}, so $g$ is continuous and $g(x)\to\frac1{2\lambda}>0$ as $\abs{x}\to\infty$. Being positive on the remaining compact set, it satisfies
\begin{align}
 \label{eq:Drach-density-uniform-positivity} \inf_{x\in\R}g(x)>0.
\end{align}

We pass to the limit term by term. By \eqref{eq:Drach-density-coordinate-convergence} and $H^1(\R)\hookrightarrow L^\infty(\R)$, $g_n\to g$ uniformly. Hence, by \eqref{eq:Drach-density-uniform-positivity}, for all sufficiently large $n$,
\begin{align}
 \label{eq:Drach-density-uniform-gap} \inf_{x\in\R}g_n(x) \geq\frac12\inf_{x\in\R}g(x)>0.
\end{align}
The identities
\begin{align}
 g_n^{-1}-g^{-1} &=\frac{g-g_n}{g_ng},\label{eq:Drach-density-inverse-difference}\\
 \partial_x(g_n^{-1}-g^{-1}) &=-\frac{g_{n,x}-g_x}{g_n^2} -g_x\left(\frac1{g_n^2}-\frac1{g^2}\right) \label{eq:Drach-density-inverse-derivative}
\end{align}
then show that
\begin{align}
 \label{eq:Drach-density-inverse-convergence} g_n^{-1}-2\lambda \rightarrow g^{-1}-2\lambda \quad\text{in }H^1(\R).
\end{align}
Indeed, the first terms on the right of \eqref{eq:Drach-density-inverse-difference} and \eqref{eq:Drach-density-inverse-derivative} converge in $L^2$, while the last coefficient in \eqref{eq:Drach-density-inverse-derivative} converges in $L^\infty$.

Set
\begin{align}
 \label{eq:Drach-density-a3} \ell_n:=\frac{g_{n,x}}{2g_n}, \qquad \ell:=\frac{g_x}{2g}.
\end{align}
Using \eqref{eq:Drach-density-uniform-gap}, \eqref{eq:Drach-density-coordinate-convergence}, and \eqref{eq:Drach-density-inverse-convergence}, we obtain
\begin{align}
 \norm{\ell_n-\ell}_2 &\leq \frac12\norm{g_n^{-1}}_\infty \norm{g_{n,x}-g_x}_2 +\frac12\norm{g_x}_2 \norm{g_n^{-1}-g^{-1}}_\infty \rightarrow0.\label{eq:Drach-density-log-convergence}
\end{align}
Consequently,
\begin{align}
 \partial_x\ell_n&\rightarrow\partial_x\ell \qquad\text{in }H^{-1},\label{eq:Drach-density-derivative-convergence}\\
 \ell_n^2&\rightarrow\ell^2 \qquad\text{in }L^1\hookrightarrow H^{-1}, \label{eq:Drach-density-square-convergence}
\end{align}
where the second assertion follows from
\begin{align}
 \label{eq:Drach-density-a4} \norm{\ell_n^2-\ell^2}_1 \leq(\norm{\ell_n}_2+\norm{\ell}_2)\norm{\ell_n-\ell}_2.
\end{align}
Finally, if $b_n=g_n^{-1}-2\lambda$ and $b=g^{-1}-2\lambda$, then $b_n\to b$ in $H^1$ and
\begin{align}
 \label{eq:Drach-density-constant-cancellation} \frac1{4g_n^2}-\lambda^2 =\lambda b_n+\frac14b_n^2 \,\rightarrow\, \lambda b+\frac14b^2 =\frac1{4g^2}-\lambda^2 \quad\text{in }H^1,
\end{align}
by the $H^1$ algebra property. Passing to the limit in \eqref{eq:Drach-density-smooth-formula}, using $q_n\to q$ in $H^{-1}$, proves \eqref{eq:Drach-density-conclusion}.
\end{proof}

{\begin{proposition}\label{prop:static} Let $\cQ\subset H^{-1}(\R;\R)$ be compact, and put $M_\cQ:=\sup_{q\in\cQ}\norm{q}_{H^{-1}}$. There is an absolute constant $\delta\in(0,\frac12]$ such that, with
\begin{align}
 \label{eq:static-threshold} \lambda_0(\cQ):=\max\left\{1,\ \frac{4M_\cQ^2}{\delta^2}\right\},
\end{align}
the following four statements hold for every $q\in\cQ$ and every $\lambda\geq\lambda_0(\cQ)$.
\begin{enumerate}[
label=\textup{(\roman*)}, leftmargin=*, itemsep=.7em, before={\let\fullwidthdisplay\relax} ]
\item \emph{Positivity and continuity of the kernel.} The operator $L_q+\lambda^2$ is strictly positive, its inverse admits a continuous integral kernel $G_\lambda(\cdot,\cdot;q)$, and the diagonal \eqref{eq:intro-resolvent} satisfies $g_\lambda(x;q)>0$ for every $x\in\R$.

\item \emph{Regularity and continuity.} The functions $g_\lambda$ and $g_\lambda^{-1}$ are $H^1$ perturbations of their free values,
\begin{align}
 \label{eq:coordinate-reg} g_\lambda(q)-\frac1{2\lambda}\in H^1(\R), \qquad \frac1{g_\lambda(q)}-2\lambda\in H^1(\R),
\end{align}
and, at each fixed $\lambda\geq\lambda_0(\cQ)$, the two maps
\begin{align}
 \label{eq:coordinate-maps} p\mapsto g_\lambda(p)-\frac1{2\lambda}, \qquad p\mapsto \frac1{g_\lambda(p)}-2\lambda
\end{align}
are continuous from an open neighborhood of $\cQ$ in $H^{-1}(\R;\R)$ into $H^1(\R)$.

\item \emph{Recovery.} The identity
\begin{align}
 \label{eq:Drach} q=\partial_x\!\left(\frac{g_{\lambda,x}}{2g_\lambda}\right) +\left(\frac{g_{\lambda,x}}{2g_\lambda}\right)^2 +\frac1{4g_\lambda^2}-\lambda^2
\end{align}
holds in $H^{-1}(\R)$. In particular the diagonal Green function determines the potential: if $q_1,q_2\in\cQ$ satisfy $g_\lambda(q_1)=g_\lambda(q_2)$ for one such $\lambda$, then $q_1=q_2$.

\item \emph{Quantitative bounds.} The diagonal Green function is comparable to its free value,
\begin{align}
 \label{eq:g-quantitative} \frac1{4\lambda}\leq g_\lambda(x;q) \leq\frac3{4\lambda} \qquad(x\in\R),
\end{align}
and, at each fixed $\lambda\geq\lambda_0(\cQ)$, the two $H^1$ norms in \eqref{eq:coordinate-reg} are bounded uniformly for $q\in\cQ$.
\end{enumerate}
More generally, fix $\lambda\geq\lambda_0(\cQ)$. At this $\lambda$, the assertions of \textup{(i)}--\textup{(iii)} and the bound \eqref{eq:g-quantitative} hold for every $p$ in the open ball $\{p\in H^{-1}(\R;\R):\norm{p}_{H^{-1}}<\delta\lambda^{\frac12}\}$, which contains $\cQ$, and the maps \eqref{eq:coordinate-maps} are continuous on that ball.
\end{proposition}}

{
We call a parameter $\lambda>0$ \emph{admissible for $\cQ$} when $\lambda\geq\lambda_0(\cQ)$, that is, when the four conclusions above are available at that parameter. } A single potential $q$ and the image $q(I)$ of a curve $q\in C(I;H^{-1})$ are compact, so the same word applies to them. The phrases ``admissible for $q$'' and ``admissible for the trajectory'' mean $\lambda\geq\lambda_0(\{q\})$ and $\lambda\geq\lambda_0(q(I))$, respectively. By \eqref{eq:static-threshold} the threshold depends on $\cQ$ only through $M_\cQ$ and is nondecreasing in it, so that $\lambda_0(\cQ)\leq\lambda_0(\cQ')$ whenever $M_\cQ\leq M_{\cQ'}$, in particular whenever $\cQ\subset\cQ'$. On the small ball of \cite[Proposition~2.2]{Killip.Visan2019}, the maps \eqref{eq:coordinate-maps} are real-analytic diffeomorphisms onto their images. We use their continuity below.

{\begin{remark}[Positivity is weaker than admissibility]
\label{rem:positivity-not-admissibility}
Admissibility means the threshold condition just stated and is strictly stronger than invertibility of $L_q+\lambda^2$. This distinction applies in Section~\ref{sec:mkdv}, where the Schr\"odinger operators associated with the Miura transforms satisfy $L_{q_\sigma}\geq0$. The spectral point associated with \eqref{eq:intro-resolvent} is $-\kappa^2$: for a nonnegative operator,
\begin{align}
 \label{eq:negative-spectral-point} L_Q\geq0,\quad\kappa>0 \qquad\Longrightarrow\qquad L_Q+\kappa^2\geq\kappa^2>0, \qquad \norm[\big]{(L_Q+\kappa^2)^{-1}}_{L^2\to L^2}\leq\kappa^{-2},
\end{align}
so that $-\kappa^2$ lies in the resolvent set of $L_Q$. The point $+\kappa^2$ can lie in the spectrum. For $Q=0$ the spectrum is $[0,\infty)$.

For $a>0$ let
\begin{align}
 \label{eq:delta-counterexample-field} u_a(x):=\begin{cases}
 0, & x<0,\\[2pt]
 \dfrac{a}{1+ax}, & x>0 .
 \end{cases}
\end{align}
Then $u_a\in L^2(\R)$ with $\norm{u_a}_2^2=a$, and a direct computation gives, in the sense of distributions,
\begin{align}
 \label{eq:delta-counterexample-potential} \cM_+(u_a)=u_a^2+\partial_xu_a=a\delta_0 .
\end{align}
The operator $L_{a\delta_0}=-\partial_x^2+a\delta_0$ is nonnegative, so $-\lambda^2$ belongs to its resolvent set for every $\lambda>0$ and $L_{a\delta_0}+\lambda^2$ is boundedly invertible. Inserting the rank-one interaction into the free resolvent gives the diagonal value at the origin,
\begin{align}
 \label{eq:delta-counterexample-diagonal} g_\lambda(0;a\delta_0)=\frac1{2\lambda+a} .
\end{align}
At $\lambda=1$ and $a>2$ this is smaller than $\frac1{4\lambda}$, so the lower bound in \eqref{eq:g-quantitative} fails. Thus positivity of $L_q+\lambda^2$ does not imply admissibility, including for Miura transforms of $L^2$ functions. Proposition~\ref{prop:static} supplies admissibility above a threshold determined by a bound on the $H^{-1}$ norms of the potentials.
\end{remark}}

\begin{proof}
All four assertions are transported from the fixed-parameter small-ball theory of \cite[Proposition~2.2]{Killip.Visan2019} by the dilation of Lemma~\ref{lem:dilation}. We first record what is imported, then fix the threshold, and then prove \textup{(i)}--\textup{(iv)} in the order \textup{(ii)}, \textup{(iv)}, \textup{(i)}, \textup{(iii)}.

\proofstep{Step 1: the imported small-ball statement}
By \cite[Propositions~2.1 and~2.2]{Killip.Visan2019} at the parameter $1$, there is $\delta\in(0,\frac12]$ such that for every real $p\in H^{-1}(\R)$ with $\norm{p}_{H^{-1}}<\delta$ the operator $L_p+1$ is strictly positive:
\begin{align}
 \label{eq:small-ball-coercive} L_p+1\geq(1-\delta)\bigl(-\partial_x^2+1\bigr)
\end{align}
in the sense of quadratic forms on $H^1(\R)$, because the key estimate (2.4) of that reference gives
\begin{align}
 \label{eq:small-ball-key-estimate} \norm{\sqrt{R_0}\,p\,\sqrt{R_0}}_{L^2\to L^2}\leq\norm{p}_{H^{-1}}<\delta
\end{align}
in the norm \eqref{eq:Hminus-norm}, its inverse has a continuous integral kernel, and the two maps
\begin{align}
 \label{eq:small-ball-maps} p\mapsto g_1(p)-\frac12, \qquad p\mapsto\frac1{g_1(p)}-2
\end{align}
are continuous from the ball $\{\norm{p}_{H^{-1}}<\delta\}$ into $H^1(\R)$, and vanish at $p=0$. Since $H^1(\R)\hookrightarrow L^\infty(\R)$, we may shrink $\delta$ once and for all so that it does not exceed the radius $\delta_0$ of \cite[Proposition~3.1]{Killip.Visan2019}, used in Lemma~\ref{lem:microscopic}, and so that, in addition,
\begin{align}
 \label{eq:small-ball-sup} \Bnorm{g_1(p)-\frac12}_\infty\leq\frac14, \qquad\text{hence}\qquad \frac14\leq g_1(x;p)\leq\frac34 \quad(x\in\R),
\end{align}
throughout that ball. This fixes the absolute constant $\delta$ in \eqref{eq:static-threshold}.

\proofstep{Step 2: the threshold and a neighborhood of $\cQ$}
Write $\varepsilon=\lambda^{-1}$, and note that $M_\cQ$ is finite because $\cQ$ is compact. By part~\textup{(i)} of Lemma~\ref{lem:dilation}, $\norm{p^{(\varepsilon)}}_{H^{-1}}\leq\varepsilon^{\frac12}\norm{p}_{H^{-1}}$ for every $p$. Consequently the open set
\begin{align}
 \label{eq:static-neighborhood} \mathcal O_\lambda:=\bigl\{p\in H^{-1}(\R;\R): \norm{p}_{H^{-1}}<\delta\lambda^{\frac12}\bigr\}
\end{align}
is mapped by $p\mapsto p^{(\varepsilon)}$ into the small ball of Step~1. The threshold \eqref{eq:static-threshold} gives $\delta\lambda^{\frac12}\geq2M_\cQ$, and hence $\cQ\subset \mathcal O_\lambda$, for every $\lambda\geq\lambda_0(\cQ)$. Thus $\mathcal O_\lambda$ is an open neighborhood of $\cQ$ in $H^{-1}(\R;\R)$, and every $p\in \mathcal O_\lambda$ has its dilate in the small ball.

\proofstep{Step 3: proof of \textup{(ii)}}
Fix $\lambda\geq\lambda_0(\cQ)$ and $p\in \mathcal O_\lambda$. Step~2 and Step~1 give $L_{p^{(\varepsilon)}}+1>0$, so parts~\textup{(ii)} and~\textup{(iii)} of Lemma~\ref{lem:dilation} apply and yield $L_p+\lambda^2>0$ together with the diagonal identity \eqref{eq:g-scaling}. Conjugating \eqref{eq:small-ball-coercive} by \eqref{eq:operator-scaling} gives the quantitative form $L_p+\lambda^2\geq(1-\delta)(-\partial_x^2+\lambda^2)$. Substituting $\varepsilon=\lambda^{-1}$ in \eqref{eq:g-scaling} and taking reciprocals gives the two scaling identities
\begin{align}
 \label{eq:g-coordinate-scaling} \begin{aligned}
 g_\lambda(y;p)-\frac1{2\lambda} &=\varepsilon\left[ g_1(y/\varepsilon;p^{(\varepsilon)})-\frac12\right],\\
 \frac1{g_\lambda(y;p)}-2\lambda &=\varepsilon^{-1}\left[ \frac1{g_1(y/\varepsilon;p^{(\varepsilon)})}-2\right].
 \end{aligned}
\end{align}
The right-hand sides are dilates of the $H^1$ functions supplied by Step~1, and dilation preserves $H^1$. This proves the memberships \eqref{eq:coordinate-reg}. Moreover each map in \eqref{eq:coordinate-maps} is the composition of the bounded linear map $p\mapsto p^{(\varepsilon)}$, one of the continuous maps \eqref{eq:small-ball-maps}, and the fixed dilation appearing in \eqref{eq:g-coordinate-scaling}. All three are continuous, so the composition is continuous from $\mathcal O_\lambda$ into $H^1(\R)$. This proves \textup{(ii)}.

\proofstep{Step 4: proof of \textup{(iv)}, and of the positivity in \textup{(i)}}
Apply \eqref{eq:small-ball-sup} to $p^{(\varepsilon)}$ and insert it into the first identity in \eqref{eq:g-coordinate-scaling}, or equivalently into $g_\lambda(y;p)=\varepsilon g_1(y/\varepsilon;p^{(\varepsilon)})$. Since $\varepsilon=\lambda^{-1}$, this gives exactly \eqref{eq:g-quantitative}. In particular $g_\lambda(x;q)\geq\frac1{4\lambda}>0$ for every $x$, which is the last assertion of \textup{(i)}. The strict positivity of $L_q+\lambda^2$ and the continuity of the kernel were obtained in Step~3 from Lemma~\ref{lem:dilation}, so \textup{(i)} is complete.

For the uniform $H^1$ bounds, note that $p\mapsto p^{(\varepsilon)}$ is continuous on $H^{-1}$, so $\{q^{(\varepsilon)}:q\in\cQ\}$ is a compact subset of the small ball. Its image under either map in \eqref{eq:small-ball-maps} is then a compact, hence bounded, subset of $H^1(\R)$. Finally, at fixed $\varepsilon$, the two dilations in \eqref{eq:g-coordinate-scaling} are bounded operators on $H^1(\R)$. Explicitly, for $F\in H^1$ and $0<\varepsilon\leq1$,
\begin{align}
 \label{eq:dilation-H1-bounds} \norm[\big]{\varepsilon F(\cdot/\varepsilon)}_{H^1}^2 =\varepsilon^3\norm{F}_2^2+\varepsilon\norm{F_x}_2^2 \leq\varepsilon\norm{F}_{H^1}^2, \qquad \norm[\big]{\varepsilon^{-1}F(\cdot/\varepsilon)}_{H^1}^2 \leq\varepsilon^{-3}\norm{F}_{H^1}^2 .
\end{align}
Composing gives the asserted bounds, uniform in $q\in\cQ$ at each fixed $\lambda$. This proves \textup{(iv)}.

\proofstep{Step 5: proof of \textup{(iii)}}
Fix $q\in\cQ$ and $\lambda\geq\lambda_0(\cQ)$. By Steps~2 and~3, the open set $\mathcal O_\lambda$ of \eqref{eq:static-neighborhood} contains $q$ and supplies the hypotheses of Lemma~\ref{lem:Drach-density}: positivity of $L_p+\lambda^2$ for every $p\in \mathcal O_\lambda$, so that any sequence $q_n\in C_c^\infty(\R;\R)$ with $q_n\to q$ in $H^{-1}$ lies in $\mathcal O_\lambda$ after discarding finitely many terms, and continuity of $p\mapsto g_\lambda(p)-\frac1{2\lambda}$ from $\mathcal O_\lambda$ into $H^1(\R)$, which gives the convergence in \eqref{eq:Drach-density-coordinate-convergence}. Step~4 gives the pointwise positivity of $g_\lambda(q)$. Lemma~\ref{lem:Drach-density} therefore applies and yields \eqref{eq:Drach}.

Finally, if $q_1,q_2\in\cQ$ have $g_\lambda(q_1)=g_\lambda(q_2)$ for one admissible $\lambda$, then the right-hand side of \eqref{eq:Drach} is the same for both, so $q_1=q_2$.

Steps~3 and~5, and Step~4 apart from the uniform $H^1$ bounds, use about the potential only that it lies in $\mathcal O_\lambda$. This proves the last assertion of the proposition.
\end{proof}

{
\begin{remark}
The identity \eqref{eq:Drach} is an identity in $H^{-1}(\R)$, and each of its terms is meaningful by the properties just established. Forming the reciprocal $g_\lambda^{-1}$ is legitimate by \eqref{eq:affine-reciprocal}--\eqref{eq:multipliers-a2}, whose hypothesis of a uniform positive lower bound is supplied by \eqref{eq:g-quantitative}. That bound also gives $\norm{g_\lambda^{-1}}_\infty\leq4\lambda$. Since $g_{\lambda,x}\in L^2(\R)$ by \eqref{eq:coordinate-reg}, it follows that
\begin{align}
 \label{eq:Drach-log-derivative} \frac{g_{\lambda,x}}{2g_\lambda}\in L^2(\R).
\end{align}
Hence the first term of \eqref{eq:Drach} is the distributional derivative of an $L^2$ function, so it lies in $H^{-1}(\R)$, and the second term is the square of an $L^2$ function, so it lies in $L^1(\R)$, which embeds into $H^{-1}(\R)$ by part~\textup{(iii)} of Lemma~\ref{lem:products}. For the last two terms, factor
\begin{align}
 \label{eq:Drach-constant-term} \frac1{4g_\lambda^2}-\lambda^2 =\frac14\left(\frac1{g_\lambda}-2\lambda\right) \left(\frac1{g_\lambda}+2\lambda\right) \in H^1(\R)\subset H^{-1}(\R).
\end{align}
The first factor lies in $H^1(\R)$ by \eqref{eq:coordinate-reg}, and the second is $4\lambda+(g_\lambda^{-1}-2\lambda)$, an element of the affine algebra $\mathcal A=\R+H^1(\R)$ of \eqref{eq:affine-Hone-algebra}. The product of an $H^1$ function with an element of $\mathcal A$ lies again in $H^1(\R)$, by the algebra property in part~\textup{(i)} of Lemma~\ref{lem:products} and the decomposition in \eqref{eq:affine-Hone-algebra}. Each of $\frac1{4g_\lambda^2}$ and $\lambda^2$ has the nonzero limit $\lambda^2$ at spatial infinity. The cancellation of this common constant gives the $H^1$ membership in \eqref{eq:Drach-constant-term}.
\end{remark}
}

%%%%%%%%%%%%%%%%%%%%%%%%%%%%%%%%%%%%%%%%%%%%%%%%%%%%%%%%%%%%%%%%%%%%%%%%%%%%%%%%%%
\subsection{Green Solution}
%%%%%%%%%%%%%%%%%%%%%%%%%%%%%%%%%%%%%%%%%%%%%%%%%%%%%%%%%%%%%%%%%%%%%%%%%%%%%%%%%%

Recall the functions $h_\lambda,w_\lambda$ from \eqref{eq:intro-resolvent}:
\begin{align}
 \label{eq:reciprocal-coordinates-recalled} h_\lambda=(2g_\lambda)^{-1}, \qquad w_\lambda=h_\lambda-\lambda.
\end{align}
Proposition~\ref{prop:static} gives $w_\lambda\in H^1$. Substituting $g_\lambda=(2h_\lambda)^{-1}$ into \eqref{eq:Drach} expresses $q$ in terms of $h_\lambda$:
\begin{align}
 \label{eq:q-from-h} q=-\frac{h_{\lambda,xx}}{2h_\lambda} +\frac{3h_{\lambda,x}^2}{4h_\lambda^2} +h_\lambda^2-\lambda^2 \quad\text{in }H^{-1}.
\end{align}

The next lemma explains the choice of evolution law in Definition~\ref{def:green}. For smooth KdV solutions, $h_\lambda$ satisfies a local conservation law. Its flux contains no product $q^2$, and each term remains meaningful at $H^{-1}$ regularity.

{\begin{lemma}\label{lem:microscopic}
Let $\lambda$ be admissible in the sense of Proposition~\ref{prop:static}.
\begin{enumerate}[
label=\textup{(\roman*)}, leftmargin=*, itemsep=.7em, before={\let\fullwidthdisplay\relax} ]
\item \emph{The flux is well-defined and continuous.} Let $\cQ\subset H^{-1}(\R;\R)$ be compact, let $\lambda\geq\lambda_0(\cQ)$, and let $\mathcal O_\lambda$ be the open neighborhood of $\cQ$ on which Proposition~\ref{prop:static}\,\textup{(ii)} provides continuity. For $q\in\mathcal O_\lambda$, the flux \eqref{eq:intro-green-flux} satisfies
\begin{align}
 \label{eq:green-flux} \J_\lambda(q)=2\lambda q+2w_\lambda(q)q-4\lambda^2w_\lambda(q).
\end{align}
Then $\J_\lambda(q)\in H^{-1}(\R)$, and $q\mapsto\J_\lambda(q)$ is continuous from $\mathcal O_\lambda$ into $H^{-1}(\R)$. In particular, if $q\in C(I;H^{-1})$ and $\lambda\geq\lambda_0(q(I))$, then $t\mapsto\J_\lambda(q(t))$ lies in $C(I;H^{-1})$, so $\int_0^t\J_\lambda(q(t'))\dd t'$ is a Bochner integral in $H^{-1}(\R)$ and its spatial derivative lies in $H^{-2}(\R)$.

\item \emph{Smooth solutions obey the law.} If $q$ is a Schwartz solution of \eqref{eq:KdV}, then
\begin{align}
 \label{eq:microscopic-h} \partial_th_\lambda(q)=\partial_x\J_\lambda(q),
\end{align}
and consequently, for every $t\in I$,
\begin{align}
 \label{eq:green-identity-smooth} h_\lambda(q(t))-h_\lambda(q(0)) =\partial_x\int_0^t\J_\lambda(q(t'))\dd t' .
\end{align}
\end{enumerate}
\end{lemma}}

\begin{proof}
\textup{(i)} By Proposition~\ref{prop:static}\,\textup{(ii)} $g_\lambda(q)^{-1}=2\lambda+2w_\lambda(q)$ with $w_\lambda(q)\in H^1(\R)$. Substitution into \eqref{eq:intro-green-flux} gives \eqref{eq:green-flux}. Each of its three terms lies in $H^{-1}(\R)$: the first because $q$ does, the second by the multiplication $H^1\times H^{-1}\to H^{-1}$ of part~\textup{(ii)} of Lemma~\ref{lem:products}, and the third because $w_\lambda(q)\in H^1(\R)\subset H^{-1}(\R)$. Continuity in $q$ follows from the same two estimates together with the continuity of $q\mapsto w_\lambda(q)$ into $H^1(\R)$, which is Proposition~\ref{prop:static}\,\textup{(ii)}. The final assertions are then standard properties of Bochner integrals of continuous Banach-valued maps.

\textup{(ii)} Let $q$ be a Schwartz solution of \eqref{eq:KdV}. The diagonal Green function obeys the local conservation law
\begin{align}
 \label{eq:microscopic-g} \partial_tg_\lambda =\partial_x\bigl( 2g_{\lambda,xx}-6qg_\lambda-12\lambda^2g_\lambda\bigr),
\end{align}
which is \cite[Proposition~3.1, equation~(3.5)]{Killip.Visan2019} rewritten with the aid of the first integral \eqref{eq:Drach-density-smooth-first-integral}, and is recorded as equation~(7.20) of \cite{Forlano.Killip.Visan2024}. That proposition is stated for Schwartz solutions with data in the ball $B_{\delta_0}$ and $\lambda\geq1$. Since \eqref{eq:KdV} is invariant under $q\mapsto\varepsilon^2q(\varepsilon^3t,\varepsilon x)$, since \eqref{eq:g-scaling} scales the diagonal Green function accordingly, and since the radius $\delta$ of Step~1 in the proof of Proposition~\ref{prop:static} was chosen not to exceed $\delta_0$, the law holds for every Schwartz solution at every $\lambda$ admissible for its trajectory. Equation~(3.6) of \cite{Killip.Visan2019} states \eqref{eq:microscopic-h} with the flux normalized as in \eqref{eq:intro-green-flux}. We derive it from \eqref{eq:microscopic-g} for use again in Section~\ref{sec:mkdv}. Eliminating $g_\lambda$ in favour of $h_\lambda=(2g_\lambda)^{-1}$ by the recovery identity \eqref{eq:Drach} turns the flux in \eqref{eq:microscopic-g} into
\begin{align}
 \label{eq:microscopic-g-flux-in-h} 2g_{\lambda,xx}-6qg_\lambda-12\lambda^2g_\lambda =\frac{h_{\lambda,xx}}{2h_\lambda^2} -\frac{h_{\lambda,x}^2}{4h_\lambda^3} -3h_\lambda-\frac{3\lambda^2}{h_\lambda} =:\mathcal B_\lambda .
\end{align}
Since $h_\lambda=(2g_\lambda)^{-1}$, we have $\partial_th_\lambda=-2h_\lambda^2\,\partial_tg_\lambda=-2h_\lambda^2\,\partial_x\mathcal B_\lambda$. Using \eqref{eq:Drach} a second time to express $q$ through $h_\lambda$ and accounting for the spatial dependence of the factor $-2h_\lambda^2$, we obtain
\begin{align}
 \label{eq:microscopic-exact} -2h_\lambda^2\,\partial_x\mathcal B_\lambda =\partial_x\bigl(2h_\lambda q-4\lambda^2(h_\lambda-\lambda)\bigr) =\partial_x\J_\lambda(q),
\end{align}
where the last equality follows from \eqref{eq:intro-green-flux}, since $2h_\lambda=g_\lambda^{-1}$ and $h_\lambda-\lambda=w_\lambda$. This proves \eqref{eq:microscopic-h}, and integrating it from $0$ to $t$ proves \eqref{eq:green-identity-smooth}.
\end{proof}

We can now return to Definition~\ref{def:green} with all its terms justified. Recall that a real curve $q\in C(I;H^{-1})$ is a Green solution precisely when there is a common threshold $\lambda_*(q,I)$ such that, for every $\lambda\geq\lambda_*(q,I)$ and every $t\in I$,
\begin{align}
 \label{eq:green-identity-recalled} h_\lambda(q(t))-h_\lambda(q(0)) =\partial_x\int_0^t\J_\lambda(q(t'))\dd t' \quad\text{in }H^{-2}(\R).
\end{align}
This is the identity \eqref{eq:intro-green-identity}, recalled here after constructing $g_\lambda$ and its reciprocal. The flux is given by \eqref{eq:green-flux}, and the threshold is chosen so that $g_\lambda$ and its reciprocal are defined along the whole trajectory.

Subtracting \eqref{eq:green-identity-recalled} at two times gives, for $s,t\in I$,
\begin{align}
 \label{eq:integrated-w} w_\lambda(t)-w_\lambda(s) =\partial_x\int_s^t\J_\lambda(q(t'))\dd t' \quad\text{in }H^{-2}.
\end{align}

Lemma~\ref{lem:microscopic}\,\textup{(i)} makes the time integral a Bochner integral in $H^{-1}$. The product $w_\lambda q$ is defined by the $H^1\times H^{-1}\to H^{-1}$ multiplication of Lemma~\ref{lem:products}. Thus every term in \eqref{eq:green-identity-recalled} is meaningful under the assumed $C(I;H^{-1})$ regularity.

Lemma~\ref{lem:microscopic}\,\textup{(ii)} proves that every Schwartz KdV solution satisfies this law. Definition~\ref{def:green} imposes the same identity directly on a rough curve. The quantifier over all sufficiently large $\lambda$ will be used in the two-parameter comparison of Lemma~\ref{lem:two-parameter}.

The canonical trajectories of \cite{Killip.Visan2019} belong to the class so defined. This inclusion is stated in the paragraph following equation~(7.17) of \cite{Killip.Murphy.Visan2020}. Since Corollary~\ref{cor:canonical} rests on it, and its argument is used again in Lemma~\ref{lem:global-extension}, we record the short proof.

{\begin{lemma}[Canonical trajectories are Green solutions]\label{lem:canonical-is-green}
Let $\Phi_{\mathrm{KdV}}\colon\R\times H^{-1}(\R;\R)\to H^{-1}(\R;\R)$ be the jointly continuous extension of the Schwartz KdV flow of \cite[Theorem~1.1]{Killip.Visan2019}. For every $q_0\in H^{-1}(\R;\R)$ and every $T>0$, the curve $q(t):=\Phi_{\mathrm{KdV}}(t,q_0)$, $t\in I$, is a Green solution in the sense of Definition~\ref{def:green}.
\end{lemma}}

\begin{proof}
Choose $q_{0,n}\in\cS(\R;\R)$ with $q_{0,n}\to q_0$ in $H^{-1}$, and let $q_n(t):=\Phi_{\mathrm{KdV}}(t,q_{0,n})$ be the corresponding Schwartz solutions of \eqref{eq:KdV}. Joint continuity of $\Phi_{\mathrm{KdV}}$ on the compact set $I\times(\{q_{0,n}:n\geq1\}\cup\{q_0\})$ gives $q_n\to q$ in $C(I;H^{-1})$. Hence
\begin{align}
 \label{eq:canonical-is-green-compact} \cQ:=q(I)\cup\bigcup_{n\geq1}q_n(I)
\end{align}
is compact, being the image of the compact space $I\times(\mathbb N\cup\{\infty\})$ under a continuous map. Fix $\lambda\geq\lambda_0(\cQ)$. Since $\lambda$ is admissible for every $q_n(I)\subset\cQ$, Lemma~\ref{lem:microscopic}\,\textup{(ii)} gives \eqref{eq:green-identity-smooth} for each $q_n$. As $n\to\infty$, the left side converges to $h_\lambda(q(t))-h_\lambda(q(0))$ in $H^1(\R)$ by Proposition~\ref{prop:static}\,\textup{(ii)}. On the right, $\J_\lambda$ is continuous on the compact set $\cQ$ by Lemma~\ref{lem:microscopic}\,\textup{(i)}, hence uniformly continuous there, so $\J_\lambda(q_n(t'))\to\J_\lambda(q(t'))$ in $H^{-1}(\R)$ uniformly in $t'\in I$, and the Bochner integrals converge in $H^{-1}(\R)$. Applying $\partial_x$ shows that \eqref{eq:green-identity-recalled} holds for $q$ at every $\lambda\geq\lambda_0(\cQ)$, so $q$ is a Green solution with $\lambda_*(q,I):=\lambda_0(\cQ)$.
\end{proof}

{
\begin{remark}[Why an intrinsic notion is needed]\label{rem:why-green}
For a general $q\in H^{-1}(\R)$, the Sobolev hypothesis supplies no canonical product $q^2$, so it does not define the right-hand side of \eqref{eq:KdV} as a distribution. Definition~\ref{def:green} gives a solution criterion at this regularity through quantities defined by Lemma~\ref{lem:microscopic}. The same lemma proves that smooth KdV solutions satisfy this criterion.

The solutions of \cite{Killip.Visan2019} are constructed by extending the classical solution map. One solves \eqref{eq:KdV} for Schwartz data and proves that the solution map extends uniquely and \emph{continuously} to $H^{-1}$. The corresponding classical solutions therefore converge in $C(I;H^{-1})$ whenever their data converge in $H^{-1}$. The flow at a general datum $q_0\in H^{-1}$ is then defined as this limit. Uniform continuity on bounded sets fails already for $s<-\frac34$ by \cite[p.~251]{Killip.Visan2019}. On a fixed compact set of data and times, continuity is uniform, as used below. This construction determines a canonical trajectory for each datum. Definition~\ref{def:green} supplies a criterion for comparing this trajectory with an arbitrary curve in $C(I;H^{-1})$.

To check the criterion, at each time one forms the diagonal Green function of $L_{q(t)}+\lambda^2$, determined by $q(t)$ through \eqref{eq:intro-resolvent}, and verifies \eqref{eq:green-identity-recalled}. The definition imposes these identities on curves of arbitrary size in $C(I;H^{-1})$. Theorem~\ref{thm:main} proves uniqueness in this class. Lemma~\ref{lem:canonical-is-green} proves that the canonical trajectories belong to it, as also recorded after Definition~7.3 of \cite{Killip.Murphy.Visan2020}. Hence Corollary~\ref{cor:canonical} identifies every Green solution with its canonical trajectory.

The identity is required at all sufficiently large $\lambda$. This supplies the two-parameter comparison of Lemma~\ref{lem:two-parameter} and, through it, the conservation law and smoothing of Sections~\ref{sec:alpha}--\ref{sec:automatic-smoothing}. Its integrated form is meaningful under the assumed continuity in time.
\end{remark}
}

{
\begin{remark}[Comparison with the normalization of \cite{Killip.Murphy.Visan2020}]
\label{rem:KMV-normalization}
In the notation of \cite[Section~7]{Killip.Murphy.Visan2020}, $\kappa=\lambda$, and the identity there reads
\begin{align}
 \label{eq:KMV-original-green-identity} h_\lambda(q(t))-h_\lambda(q(0)) =\partial_x\int_0^t \left[ \frac{q(t')}{g_\lambda(q(t'))} -\frac{2\lambda^2}{g_\lambda(q(t'))} \right]\dd t',
\end{align}
interpreted there as an integrated identity in $C(I;H^{-2}(\R))$. Their raw flux satisfies
\begin{align}
 \label{eq:KMV-normalization-a1} \frac{q}{g_\lambda}-\frac{2\lambda^2}{g_\lambda} =\J_\lambda(q)-4\lambda^3,
\end{align}
and the discrepancy is the spatial constant $-4\lambda^3$, which disappears under $\partial_x$. Hence \eqref{eq:KMV-original-green-identity} and \eqref{eq:green-identity-recalled} define the same class of curves.

The normalized flux \eqref{eq:intro-green-flux} lies in $H^{-1}(\R)$. To identify the normalization, expand the two terms in the original flux. By Proposition~\ref{prop:static}\,\textup{(ii)}, $g_\lambda^{-1}=2\lambda+2w_\lambda$ with $w_\lambda\in H^1(\R)$, so
\begin{align}
 \label{eq:raw-flux-splitting} \frac{q}{g_\lambda}=2\lambda q+2w_\lambda q\in H^{-1}(\R), \qquad -\frac{2\lambda^2}{g_\lambda}=-4\lambda^3-4\lambda^2w_\lambda .
\end{align}
In \eqref{eq:raw-flux-splitting}, the second term contributes the constant $-4\lambda^3$, so the raw flux lies in $-4\lambda^3+H^{-1}(\R)$. Adding $4\lambda^3$ gives \eqref{eq:intro-green-flux} and preserves the spatial derivative. The free case is an immediate check: at $q=0$, $g_\lambda\equiv(2\lambda)^{-1}$ and $w_\lambda\equiv0$, so the raw flux equals $-4\lambda^3$ while $\J_\lambda(0)=0$.
\end{remark}
}

{
\begin{remark}[Common spectral range]\label{rem:common-lambda}
Three separate largeness requirements are imposed on a spectral parameter in this paper. The first is the threshold $\lambda_*$ of Definition~\ref{def:green}, above which the Green identity is asserted. The second is the admissibility threshold $\lambda_0$ of Proposition~\ref{prop:static}, above which $g_\lambda$ and its reciprocal exist with the quantitative bounds \eqref{eq:g-quantitative}. The third is the threshold $\Lambda$ of Lemma~\ref{lem:high-energy} below, above which the density is nonnegative and integrable. That lemma applies to any compact set of potentials. A fourth requirement, $\tau\geq2$, is imposed in Section~\ref{sec:Gardner-class}, where the Gardner parameter is chosen so as to meet every parameter restriction of the result of \cite{Klaus.Koch.Liu2023} quoted in Lemma~\ref{lem:KKL-Gardner}. It is convenient to collect them once.

For a single Green solution $q\in C(I;H^{-1})$ we set
\begin{align}
 \label{eq:single-K} K(q):=\max\bigl\{2,\ \lambda_*(q,I),\ \lambda_0(q(I)),\ \Lambda(q(I))\bigr\} .
\end{align}
This is legitimate because the trajectory
\begin{align}
 \label{eq:single-compact-set} q(I):=\{q(t):t\in I\}\subset H^{-1}(\R;\R)
\end{align}
is compact, being the continuous image of the compact interval $I$, and compactness is exactly the hypothesis of Proposition~\ref{prop:static} and of Lemma~\ref{lem:high-energy}. This is the step at which the time continuity demanded of a Green solution is converted into a spectral threshold valid along the entire curve, and it is what turns hypotheses about a single potential into hypotheses uniform in $t$. Only four numbers are combined, so the maximum is finite.

Whenever two Green solutions are compared, one parameter must serve both. In the proof of Theorem~\ref{thm:main} it must moreover serve them on every interval $[-L,L]$ with $L\geq T$, after both have been extended beyond $I$ by Lemma~\ref{lem:global-extension}. This is possible because the thresholds \eqref{eq:static-threshold} and \eqref{eq:high-energy-threshold} depend on a compact set only through the supremum of the $H^{-1}$ norm over it, and are nondecreasing in that supremum. Accordingly, let $q_1$ and $q_2$ be Green solutions on intervals $I_1=[-T_1,T_1]$ and $I_2=[-T_2,T_2]$, all of whose values have $H^{-1}$ norm at most $R$, and fix, once and for all,
\begin{align}
 \label{eq:common-K} K_{12}:=\max\left\{2,\ \lambda_*(q_1,I_1),\ \lambda_*(q_2,I_2),\ \frac{4R^2}{\delta_1^2}\right\},
\end{align}
with $\delta_1$ the absolute constant of Lemma~\ref{lem:high-energy}. Since $\delta_1\leq\delta$, the thresholds \eqref{eq:static-threshold} and \eqref{eq:high-energy-threshold} of the compact sets $q_j(I_j)$ are at most $\max\{1,4R^2/\delta_1^2\}$, and therefore
\begin{align}
 \label{eq:K-monotone} K(q_j)\leq K_{12}\qquad(j=1,2),
\end{align}
where $K(q_j)$ is formed as in \eqref{eq:single-K} on the interval $I_j$. Consequently every statement proved below for a single trajectory under the hypothesis $\tau\geq K(q)$ is available for both competitors at the one value $\tau\geq K_{12}$, and all later estimates are carried out on this single range. For two Green solutions on the same interval $I$, one may take for $R$ the largest $H^{-1}$ norm on the compact set $q_1(I)\cup q_2(I)$.

The bound $K_{12}\geq2$ ensures that the Gardner parameter $\tau\geq K_{12}$ chosen in Step~1 of the proof of Theorem~\ref{thm:main} satisfies the restriction $\tau\geq2$ of Section~\ref{sec:Gardner-class}.

Four symbols of this kind occur in the paper, and they are kept notationally distinct because they come from different results and are enlarged at different places. In order of generality: $\Lambda(\cQ)$ is the threshold of Lemma~\ref{lem:high-energy}, attached to an arbitrary compact subset of $H^{-1}$. $K(q)$ of \eqref{eq:single-K} is attached to one KdV Green solution. $K_{12}$ of \eqref{eq:common-K} is attached to a pair of competitors. $K_u$ of \eqref{eq:mkdv-common-threshold} is attached to an mKdV solution through its two Miura transforms, and is used only in Section~\ref{sec:mkdv} and in the Part~II outline of the introduction. Within the proof of Proposition~\ref{prop:auto-smoothing} a fifth symbol $K_\sharp\geq K(q)$ records one further enlargement, made there and used nowhere else. The calligraphic letter $\cK$ denotes a compact set or a spatial cell.
\end{remark}
}

% Subtracting \eqref{eq:green-identity-recalled} at two times gives, for $s,t\in I$,
% \begin{align}
%  \label{eq:integrated-w} w_\lambda(t)-w_\lambda(s) =\partial_x\int_s^t\J_\lambda(q(t'))\dd t' \quad\text{in }H^{-2}.
% \end{align}

% \newpage

%%%%%%%%%%%%%%%%%%%%%%%%%%%%%%%%%%%%%%%%%%%%%%%%%%%%%%%%%%%%%%%%%%%%%%%%%%%%%%%%%%
\section{Density and Conservation of \texorpdfstring{$\alpha_\lambda$}{alpha}}
\label{sec:alpha}
%%%%%%%%%%%%%%%%%%%%%%%%%%%%%%%%%%%%%%%%%%%%%%%%%%%%%%%%%%%%%%%%%%%%%%%%%%%%%%%%%%

The potential $q$ and the function $h_\lambda$ define the density $\rho_\lambda$ in \eqref{eq:def-rho}. This section proves that the mass $\alpha_\lambda$ of \eqref{eq:def-alpha} is conserved along every Green solution. The conservation law is Proposition~\ref{prop:alpha-conservation}. It is the one piece of \emph{global} information that the argument extracts from the Green identity \eqref{eq:green-identity-recalled}, and it is what will later remove a stray constant from the translated-cell tail estimate in Section~\ref{sec:automatic-smoothing}. Through Lemma~\ref{lem:Gardner-mass} it also bounds the Gardner variable of Section~\ref{sec:Gardner} in $L^2$ uniformly in time. Everything before it is preparation, and the section closes with the one structural fact about the Green flux, Lemma~\ref{lem:flux-decomposition}, that Section~\ref{sec:automatic-smoothing} needs.

One might expect the conservation law to follow at once: derive a microscopic evolution for the density, integrate it in $x$, and read off that its mass is constant. Two things stand in the way, however, and they organize the section.

\textit{Scaling the small-data estimates.} The nonnegativity and integrability of the density are known for small $H^{-1}$ data at the parameter $\lambda=1$. We extend these properties to arbitrary data in a compact set at sufficiently large $\lambda$. This gap is purely one of scaling: dilating a potential trades its size for the spectral parameter. Lemma~\ref{lem:reciprocal-density-scaling} records how the density transforms, and Lemma~\ref{lem:high-energy} converts the small-data input into the large-data, large-$\lambda$ statement actually used.

\textit{Time increments.} The flux integral in \eqref{eq:green-identity-recalled} is $C^1$ in $t$ with values in $H^{-1}$, so $h_\lambda$ satisfies $\partial_th_\lambda=\partial_x\J_\lambda(q(t))$ in $C(I;H^{-2})$. Differentiating $\alpha_\lambda(q(t))$ and integrating its evolution over the line would require two further regularity statements.

First, the linear term $2\lambda\Rf_\lambda q$ in \eqref{eq:def-rho} requires a time evolution for $q$. Recovering it by differentiating \eqref{eq:q-from-h} would involve the product $2h_\lambda\,\partial_th_\lambda$, with $h_\lambda\in\lambda+H^1(\R)$ and $\partial_th_\lambda\in H^{-2}$. Multiplication by an $H^1$ function on $H^{-2}$ is not justified by the available bounds. Second, integrating a distributional density evolution over the whole line requires control of its pairing with the constant function $1$. Lemma~\ref{lem:zero-mass} supplies this control for an $L^1$ function that is the spatial derivative of an $H^{-1}$ distribution. Applying it to $\partial_t\rho_\lambda$ would require $\partial_t\rho_\lambda\in L^1$, whereas the density estimates give $\rho_\lambda(t)\in L^1$ at each fixed time.

Both difficulties disappear once the derivative is replaced by a \emph{difference} in time, which is in any case the form in which the Green identity is stated. The increment $\rho_\lambda(q(t))-\rho_\lambda(q(s))$ does lie in $L^1$, being a difference of two $L^1$ functions, and \eqref{eq:integrated-w} already exhibits the increment of $w_\lambda$ as a spatial derivative of an $H^{-1}$ distribution. The remaining linear term in \eqref{eq:def-rho} has multiplier $2\lambda(\xi^2+4\lambda^2)^{-1}$, whose value at $\xi=0$ is $(2\lambda)^{-1}$. Consequently $2\lambda \Rf_\lambda(q(t)-q(s))$ need not have a spatial primitive in the class $H^{-1}(\R)$ required by Lemma~\ref{lem:zero-mass}. A second spectral parameter supplies the necessary cancellation. In the weighted combination $\kappa\rho_\kappa-\lambda\rho_\lambda$ the two linear terms merge into $2(\kappa^2\Rf_\kappa-\lambda^2\Rf_\lambda)$, whose symbol has a \emph{double} zero at $\xi=0$ and which is consequently an exact second derivative, as \eqref{eq:multiplier-identity} records. Lemma~\ref{lem:two-parameter} assembles this into an identity to which Lemma~\ref{lem:zero-mass} applies. Integrating that identity gives equality of the mass increments at the two parameters. Letting one parameter tend to infinity then shows that both increments vanish.

{
\begin{remark}[Comparison with the argument of \cite{Killip.Visan2019}]
\label{rem:why-not-KV}
With $\kappa=\lambda$, the functional $\alpha_\lambda$ defined in \eqref{eq:def-rho}--\eqref{eq:def-alpha} is exactly the $\alpha(\kappa;q)$ of \cite[equations~(2.18)--(2.19)]{Killip.Visan2019}, identified there with $-\log\det{}_2\bigl(1+\sqrt{R_0}\,q\sqrt{R_0}\bigr)$, essentially the logarithm of the transmission coefficient. Here the source notation $R_0$ means $(-\partial_x^2+\lambda^2)^{-1}=R_\lambda(0)$, by \cite[equation~(2.2)]{Killip.Visan2019}. The auxiliary operator $\Rf_\lambda$ in \eqref{eq:R0} uses the doubled parameter $2\lambda$. Their conservation argument uses time differentiation and approximation by Schwartz solutions.

The argument there proceeds in two steps. For a Schwartz solution, \cite[Proposition~3.1]{Killip.Visan2019} derives the microscopic density law (3.7); integrating it over the line gives conservation (3.8). The restriction to Schwartz data is then removed by approximation: the $H^{-1}$ trajectory of \cite{Killip.Visan2019} is by construction a limit of Schwartz solutions \cite[Corollary~5.2]{Killip.Visan2019}. At each fixed sufficiently large $\kappa$, Lemma~\ref{lem:high-energy} makes $q\mapsto\alpha_\kappa(q)$ continuous on the common compact family of approximating trajectories, as in Lemma~\ref{lem:canonical-is-green}, so the conservation law passes to the limit.

For a Green solution, Definition~\ref{def:green} gives the time evolution of $h_\lambda$ in $H^{-2}$. The functional derivative $(2\lambda)^{-1}-g_\lambda$ of \cite[equation~(2.34)]{Killip.Visan2019}, which lies in $H^1(\R)$ by \eqref{eq:coordinate-reg}, would require corresponding control of $\partial_tq$ to justify the pairing $\pair*{\partial_tq}{\delta\alpha/\delta q}$. Approximation of the whole trajectory by Schwartz solutions is an additional hypothesis. Definition~\ref{def:green} specifies a curve in $C(I;H^{-1})$ satisfying \eqref{eq:green-identity-recalled}, and Theorem~\ref{thm:main} identifies it with the canonical trajectory. We therefore prove conservation from \eqref{eq:integrated-w} through the two-parameter identity of Lemma~\ref{lem:two-parameter}, using the estimates of \cite[Proposition~2.4]{Killip.Visan2019} for a fixed potential and parameter in the form \eqref{eq:KV-alpha}.
\end{remark}
}

%%%%%%%%%%%%%%%%%%%%%%%%%%%%%%%%%%%%%%%%%%%%%%%%%%%%%%%%%%%%%%%%%%%%%%%%%%%%%%%%%%
\subsection{Properties and Scaling}
%%%%%%%%%%%%%%%%%%%%%%%%%%%%%%%%%%%%%%%%%%%%%%%%%%%%%%%%%%%%%%%%%%%%%%%%%%%%%%%%%%

Throughout this section $\lambda$ is admissible for the potential or trajectory under discussion, in the sense fixed after Proposition~\ref{prop:static}. To separate the linear part of $w_\lambda$ from its nonlinear remainder, introduce the free operator
\begin{align}
 \label{eq:R0} \Rf_\lambda=(4\lambda^2-\partial_x^2)^{-1},
\end{align}
The auxiliary free operator \eqref{eq:R0} uses the doubled spectral parameter. For $f\in L^2(\R)$,
\begin{align}
 \label{eq:R0-kernel} \Rf_\lambda=R_{2\lambda}(0),\qquad \big(\Rf_\lambda f\big)(x)=\int_\R\frac{\ue^{-2\lambda\abs{x-y}}}{4\lambda}f(y)\dd y.
\end{align}
Hence $2\lambda \Rf_\lambda$ has kernel $\frac12\ue^{-2\lambda\abs{x-y}}$, which matches the normalization of the density \eqref{eq:def-rho} in \cite[Proposition~2.4]{Killip.Visan2019}. The operator $4\lambda^2\Rf_\lambda$ is the approximate identity with multiplier $4\lambda^2/(\xi^2+4\lambda^2)$. The free diagonal of $R_\lambda(0)$ is $(2\lambda)^{-1}$, the reference value for $g_\lambda$ in Proposition~\ref{prop:static}.

The density and its mass are
{
\begin{align}
 \rho_\lambda(q) &:=\lambda-h_\lambda(q)+2\lambda \Rf_\lambda q =-w_\lambda(q)+2\lambda \Rf_\lambda q, \label{eq:def-rho}\\
 \alpha_\lambda(q)&:=\int_\R\rho_\lambda(x;q)\dd x . \label{eq:def-alpha}
\end{align}
}
The mass is defined whenever $\rho_\lambda\in L^1(\R)$. The next two lemmas establish this condition uniformly on compact sets of potentials for sufficiently large $\lambda$.

Admissibility of $\lambda$ already gives $\rho_\lambda\in H^1(\R)$, since $w_\lambda\in H^1$ by Proposition~\ref{prop:static}\,\textup{(ii)} and $\Rf_\lambda\colon H^{-1}\to H^1$. Lemma~\ref{lem:high-energy} supplies the additional properties $\rho_\lambda\geq0$ and $\rho_\lambda\in L^1$ above its threshold. We use \eqref{eq:def-alpha} for parameters in that range. The sign matters as much as the finiteness: it is nonnegativity that identifies $\alpha_\lambda$ with $\norm{\rho_\lambda}_{1}$ and thereby converts a bound on the number $\alpha_\lambda$ into an $H^{-1}$ bound on the function $\rho_\lambda$, both in the proof of Lemma~\ref{lem:high-energy} and again in Step~3 of Proposition~\ref{prop:auto-smoothing}.

Two estimates below are stated with the comparison symbol $\simeq$ of \eqref{eq:lesssim-asymp}. In both of them the implicit constant is absolute: it depends neither on the spectral parameter nor on the potential within the relevant compact set. That uniformity is the entire content of the estimates, and it is what makes the limits of Lemma~\ref{lem:high-energy} uniform.

The following scaling lemma is the transfer mechanism announced above. Lemma~\ref{lem:dilation} already scales the Green function. For the density \eqref{eq:def-rho}, we also calculate the scaling of its free term $2\lambda \Rf_\lambda q$. The lemma says that the two terms scale compatibly, so that the whole density behaves like a single object of weight one.

{\begin{lemma}
\label{lem:reciprocal-density-scaling}
Let $q\in H^{-1}(\R;\R)$, let $0<\varepsilon\leq1$, and let $q^{(\varepsilon)}$ be the dilated potential \eqref{eq:dilation-def}. Assume that $\varepsilon^{-1}$ is admissible for $q$ in the sense of Proposition~\ref{prop:static}, and that $1$ is admissible for $q^{(\varepsilon)}$. Then, for every $x\in\R$,
\begin{align}
 w_1(x;q^{(\varepsilon)}) &=\varepsilon w_{\frac{1}{\varepsilon}}(\varepsilon x;q), \label{eq:w-scaling}\\
 \rho_1(x;q^{(\varepsilon)}) &=\varepsilon\rho_{\frac{1}{\varepsilon}}(\varepsilon x;q). \label{eq:rho-scaling}
\end{align}
Consequently $\alpha_1(q^{(\varepsilon)})=\alpha_\lambda(q)$ with $\lambda=\varepsilon^{-1}$, whenever either side is defined.
\end{lemma}}

By \eqref{eq:dilation-positivity} the two positivity requirements in the lemma are equivalent, and by \eqref{eq:resolvent-scaling} and \eqref{eq:g-scaling} the corresponding resolvents and kernels determine one another. Only the two \emph{thresholds} are chosen separately, and in the application below both hold because $q^{(1/\lambda)}$ lies in the small ball of \cite[Proposition~2.2]{Killip.Visan2019}.

\begin{proof}
Taking reciprocals in \eqref{eq:g-scaling} gives $h_1(x;q^{(\varepsilon)})=\varepsilon h_{\frac{1}{\varepsilon}}(\varepsilon x;q)$. By \eqref{eq:intro-resolvent} the constant subtracted on the left to form $w_1$ is $1$, while on the right $\varepsilon w_{\frac{1}{\varepsilon}}=\varepsilon h_{\frac{1}{\varepsilon}}-\varepsilon\cdot\varepsilon^{-1}=\varepsilon h_{\frac{1}{\varepsilon}}-1$. The same constant is subtracted on both sides, which proves \eqref{eq:w-scaling}.

For the free term, the kernel \eqref{eq:R0-kernel} and the substitution $y\mapsto\varepsilon y$ in \eqref{eq:dilation-def} give
\begin{align}
 \label{eq:reciprocal-density-scaling-a1} \big(\Rf_1q^{(\varepsilon)}\big)(x) =\big(\Rf_{\frac{1}{\varepsilon}}q\big)(\varepsilon x),
\end{align}
The factor $2\lambda$ in the definition of the density gives this term the same scaling as $w_\lambda$. The linear part of \eqref{eq:def-rho} at parameter $1$ is $2\Rf_1q^{(\varepsilon)}$, while $\varepsilon$ times the linear part at parameter $\varepsilon^{-1}$ is $\varepsilon\cdot2\varepsilon^{-1}\Rf_{\frac{1}{\varepsilon}}q=2\Rf_{\frac{1}{\varepsilon}}q$. The explicit factor $2\lambda$ in \eqref{eq:def-rho} thus supplies exactly the missing power of $\varepsilon$. Combining this with \eqref{eq:w-scaling} proves \eqref{eq:rho-scaling}, and integrating \eqref{eq:rho-scaling} in $x$ proves the last assertion.
\end{proof}

We now import the small-data input. For $\norm{q}_{H^{-1}}$ small, Killip--Vi\c{s}an \cite[Proposition~2.4]{Killip.Visan2019} show that $q\mapsto\rho_1(x;q)$ is nonnegative and convex for each fixed $x$ (this is stated explicitly in the sentence following equation~(2.18) there), that $\rho_1\in L^1\cap H^1$, and that
\begin{align}
 \label{eq:KV-alpha} \alpha_1(q)\simeq\norm{q}_{H^{-1}}^2 ,
\end{align}
with an absolute constant, uniformly over that small ball. In writing \eqref{eq:KV-alpha} we have used that the right-hand side of the estimate in \cite{Killip.Visan2019} is precisely the norm \eqref{eq:Hminus-norm}. The same proposition and its proof, immediately before equations~(2.31)--(2.33) there, also give real analyticity of $q\mapsto\rho_1(q)$ into $L^1\cap H^1$. The cited result supplies the pointwise positivity used below.

The next lemma is the version of these statements for arbitrary data and sufficiently large $\lambda$, and it provides the quantitative density input from \cite{Killip.Visan2019} used below. Its three limits play three different and independent roles later: \eqref{eq:alpha-high} gives the vanishing mass limit needed in Proposition~\ref{prop:alpha-conservation}. Equation~\eqref{eq:w-high} says that $2\lambda w_\lambda$ converges to $q$, and is what identifies the weak limit in Proposition~\ref{prop:auto-smoothing}. Equation~\eqref{eq:winfty-high} is the pointwise smallness that makes the Green flux coercive in Lemma~\ref{lem:flux-decomposition}\,\textup{(ii)}.

{\begin{lemma}
\label{lem:high-energy}
Let $\cQ\subset H^{-1}(\R;\R)$ be compact, and let $M_\cQ$ and $\delta$ be as in Proposition~\ref{prop:static}. There is an absolute constant $\delta_1\in(0,\delta]$ such that, with
\begin{align}
 \label{eq:high-energy-threshold} \Lambda(\cQ):=\max\left\{1,\ \frac{4M_\cQ^2}{\delta_1^2}\right\},
\end{align}
for every $q\in\cQ$ and every $\lambda\geq\Lambda(\cQ)$ the parameter $\lambda$ is admissible for $\cQ$, the density \eqref{eq:def-rho} satisfies
\begin{align}
 \label{eq:rho-nonneg} \rho_\lambda(x;q)\geq0\quad(x\in\R), \qquad \rho_\lambda(q)\in L^1(\R)\cap H^1(\R),
\end{align}
and, at each fixed such $\lambda$, the map $q\mapsto\rho_\lambda(q)$ is continuous from $\cQ$ into $L^1\cap H^1$.

Moreover, with constants uniform for $q\in\cQ$,
\begin{align}
 \label{eq:alpha-large} \alpha_\lambda(q) \simeq\frac1\lambda\int_\R \frac{\abs{\widehat q(\xi)}^2}{\xi^2+4\lambda^2}\dd \xi \qquad(\lambda\geq\Lambda(\cQ)),
\end{align}
and the following three limits hold as $\lambda\to\infty$, uniformly for $q\in\cQ$:
\begin{align}
 \lambda\alpha_\lambda(q)&\rightarrow0, \label{eq:alpha-high}\\
 2\lambda w_\lambda(q)&\rightarrow q \quad\text{in }H^{-1}, \label{eq:w-high}\\
 \lambda^{-1}\bnorm{w_\lambda(q)}_{L^\infty}&\rightarrow0. \label{eq:winfty-high}
\end{align}
Finally, with an absolute constant $C$,
\begin{align}
 \label{eq:high-energy-bounds} \lambda\alpha_\lambda(q)+\bnorm{\lambda w_\lambda(q)}_{H^{-1}} \leq C\bigl(\norm{q}_{H^{-1}}+\norm{q}_{H^{-1}}^2\bigr) \qquad(q\in\cQ,\ \lambda\geq\Lambda(\cQ)).
\end{align}
\end{lemma}}

\begin{proof}
Every assertion is reduced to the unit parameter by the scaling of Lemma~\ref{lem:reciprocal-density-scaling} and then read off from \cite{Killip.Visan2019} at small data. The steps below carry this out one assertion at a time, in the order in which they are stated. Two of them also use the compactness argument stated first to upgrade pointwise limits to uniform limits.

\proofstep{A uniformity principle} Two of the steps below convert a limit that holds for each fixed potential into one that is uniform over $\cQ$, and it is worth isolating the mechanism once. Let $(T_j)$ be a family of bounded operators on $H^{-1}$, indexed by a parameter $j$ tending to a limit, such that
\begin{align}
 \label{eq:uniformity-hypotheses} M:=\sup_{j}\norm{T_j}_{H^{-1}\to H^{-1}}<\infty \qquad\text{and}\qquad T_jf\rightarrow0\ \text{ in }H^{-1}\ \text{ for each fixed }f .
\end{align}
Then $T_jq\to0$ uniformly for $q\in\cQ$. Indeed, let $\varepsilon>0$. Compactness of $\cQ$ supplies finitely many $f_1,\dots,f_n\in\cQ$ whose $\varepsilon/(2M+2)$-balls cover $\cQ$, and the second hypothesis in \eqref{eq:uniformity-hypotheses} supplies an index beyond which $\norm{T_jf_i}_{H^{-1}}<\varepsilon/2$ simultaneously for the finitely many $i$. Given $q\in\cQ$, choose $i$ with $\norm{q-f_i}_{H^{-1}}<\varepsilon/(2M+2)$. Then
\begin{align}
 \label{eq:high-energy-a1} \norm{T_jq}_{H^{-1}} \leq\norm{T_jf_i}_{H^{-1}}+M\norm{q-f_i}_{H^{-1}} <\frac\varepsilon2+\frac{M\varepsilon}{2M+2}<\varepsilon ,
\end{align}
and the bound is independent of $q$. This argument gives uniform convergence on $\cQ$ from strong convergence and the uniform operator bound.

\proofstep{Reduction to unit parameter} Let $0<\varepsilon\leq1$. By part~\textup{(i)} of Lemma~\ref{lem:dilation},
\begin{align}
 \label{eq:scaled-Hminus} \norm{q^{(\varepsilon)}}_{H^{-1}}^2 =\varepsilon^3\int_\R\frac{\abs{\widehat q(\eta)}^2}{4+\varepsilon^2\eta^2}\dd \eta \leq\varepsilon\norm{q}_{H^{-1}}^2 .
\end{align}
Let $\delta_1\in(0,\delta]$ be the smaller of the constant $\delta$ of Proposition~\ref{prop:static} and the radius of the small ball of \cite[Proposition~2.4]{Killip.Visan2019}, and define $\Lambda(\cQ)$ by \eqref{eq:high-energy-threshold}. Then $\Lambda(\cQ)\geq\lambda_0(\cQ)$ by \eqref{eq:static-threshold}, so every $\lambda\geq\Lambda(\cQ)$ is admissible for $\cQ$, and the inequality in \eqref{eq:scaled-Hminus} gives $\norm{q^{(1/\lambda)}}_{H^{-1}}\leq\lambda^{-\frac12}M_\cQ\leq\delta_1/2$ for $q\in\cQ$. Thus $q^{(1/\lambda)}$ lies in the small balls of \cite[Propositions~2.2 and~2.4]{Killip.Visan2019}, and the parameter $1$ is admissible for it. Fix $\lambda\geq\Lambda(\cQ)$ and put $\varepsilon=1/\lambda$.

\proofstep{Sign, regularity and continuity} At the unit parameter the three assertions are quoted: $\rho_1(q^{(1/\lambda)})$ is nonnegative and lies in $L^1\cap H^1$ by \cite[Proposition~2.4]{Killip.Visan2019}, and depends continuously on its argument there, since the proof of that proposition records that $q\mapsto\rho_1(q)$ is real analytic from the small ball into $L^1\cap H^1$. What has to be checked is that the scaling \eqref{eq:rho-scaling} carries each of them back to $\rho_\lambda(q)$.

It does, and for a different reason in each case. Read at $\varepsilon=1/\lambda$, that identity says
\begin{align}
 \label{eq:rho-scaling-inverted} \rho_\lambda(x;q) =\lambda\,\rho_1\bigl(\lambda x;q^{(1/\lambda)}\bigr) ,
\end{align}
so $\rho_\lambda(\cdot;q)$ is a positive multiple of a dilate of $\rho_1(\cdot;q^{(1/\lambda)})$. Nonnegativity is preserved by any such operation, being pointwise. Membership in $L^1\cap H^1$ is preserved because dilation is a bounded bijection of each of those spaces, the factor $\lambda$ merely rescaling the norms. This proves \eqref{eq:rho-nonneg}. Continuity at fixed $\lambda$ is the composition of three continuous maps, $q\mapsto q^{(1/\lambda)}$, the quoted map $\rho_1$, and the dilation just described. The first of these is bounded on $H^{-1}$ by \eqref{eq:scaled-Hminus}. This proves the assertions of the first paragraph of the lemma. The three limits are asymptotic statements about $\lambda\to\infty$ and are proved in the steps that follow. The closing bound \eqref{eq:high-energy-bounds}, which is uniform in $\lambda$, is proved at the very end.

\proofstep{The asymptotic \eqref{eq:alpha-large}} By the last assertion of Lemma~\ref{lem:reciprocal-density-scaling} and by \eqref{eq:KV-alpha} applied to $q^{(1/\lambda)}$,
\begin{align}
 \label{eq:high-energy-a2} \alpha_\lambda(q)=\alpha_1\bigl(q^{(1/\lambda)}\bigr) \simeq\norm[\big]{q^{(1/\lambda)}}_{H^{-1}}^2 ,
\end{align}
and the equality in \eqref{eq:scaled-Hminus} with $\varepsilon=1/\lambda$ evaluates the right side as $\lambda^{-1}\int\abs{\widehat q}^2(\xi^2+4\lambda^2)^{-1}\dd \xi$. This is \eqref{eq:alpha-large}, with the absolute constant of \eqref{eq:KV-alpha}.

\proofstep{Proof of \eqref{eq:alpha-high}} By \eqref{eq:alpha-large} it suffices to show that
$\int_\R\abs{\widehat q(\xi)}^2(\xi^2+4\lambda^2)^{-1}\dd \xi\to0$ uniformly on $\cQ$. The integral splits at a frequency $R$ that must be chosen before $\lambda$ is sent to infinity, and the point of the splitting is that the two regions are controlled for opposite reasons: the high frequencies are small because $\cQ$ is compact, uniformly in $\lambda$, while the low frequencies are small because the multiplier itself is, once $R$ is fixed.

For the high frequencies, let $P_{>R}$ denote the Fourier projection onto $\{\abs{\xi}>R\}$, so that $\norm{P_{>R}q}_{H^{-1}}^2=\int_{\abs{\xi}>R}\abs{\widehat q(\xi)}^2(4+\xi^2)^{-1}\dd \xi$ in the normalization \eqref{eq:Hminus-norm}. The family $(P_{>R})_{R\geq1}$ is bounded by $1$ on $H^{-1}$ and satisfies $P_{>R}f\to0$ for each fixed $f$, by dominated convergence, so it meets the hypotheses \eqref{eq:uniformity-hypotheses}. The uniformity principle therefore supplies, for each $\varepsilon>0$, a radius $R$ with
\begin{align}
 \label{eq:high-energy-a3} \sup_{q\in\cQ}\int_{\abs{\xi}>R} \frac{\abs{\widehat q(\xi)}^2}{4+\xi^2}\dd \xi<\varepsilon .
\end{align}
For $\lambda\geq1$, $(\xi^2+4\lambda^2)^{-1}\leq(4+\xi^2)^{-1}$, so the region $\abs{\xi}>R$ contributes at most $\varepsilon$, uniformly in $\lambda\geq1$. On $\abs{\xi}\leq R$ the multiplier is bounded by $(4\lambda^2)^{-1}$, whence
\begin{align}
 \label{eq:high-energy-a4} \sup_{q\in\cQ}\int_{\abs{\xi}\leq R} \frac{\abs{\widehat q(\xi)}^2}{\xi^2+4\lambda^2}\dd \xi \leq\frac{4+R^2}{4\lambda^2}\, \sup_{q\in\cQ}\norm{q}_{H^{-1}}^2 .
\end{align}
Choosing $R$ first and then letting $\lambda\to\infty$ proves \eqref{eq:alpha-high}.

\proofstep{Proof of \eqref{eq:w-high}} Multiplying \eqref{eq:def-rho} by $-2\lambda$ gives the exact identity
\begin{align}
 \label{eq:high-energy-a5} 2\lambda w_\lambda=4\lambda^2\Rf_\lambda q-2\lambda\rho_\lambda .
\end{align}
For the first term, put $T_\lambda:=4\lambda^2\Rf_\lambda-\mathrm{Id}$, a Fourier multiplier with symbol
\begin{align}
 \label{eq:high-energy-a6} \frac{4\lambda^2}{\xi^2+4\lambda^2}-1=-\frac{\xi^2}{\xi^2+4\lambda^2} ,
\end{align}
which is bounded by $1$ in modulus and tends to $0$ pointwise in $\xi$ as $\lambda\to\infty$. Hence $\norm{T_\lambda}_{H^{-1}\to H^{-1}}\leq1$, and $T_\lambda f\to0$ in $H^{-1}$ for each fixed $f$ by dominated convergence. The hypotheses \eqref{eq:uniformity-hypotheses} hold, and the uniformity principle gives $4\lambda^2\Rf_\lambda q\to q$ in $H^{-1}$ uniformly on $\cQ$. For the second term, nonnegativity \eqref{eq:rho-nonneg} identifies $\norm{\rho_\lambda}_1$ with $\alpha_\lambda$, and $L^1\hookrightarrow H^{-1}$ (part~\textup{(iii)} of Lemma~\ref{lem:products}) gives
\begin{align}
 \label{eq:high-energy-a7} \norm{2\lambda\rho_\lambda}_{H^{-1}} \lesssim\lambda\norm{\rho_\lambda}_1=\lambda\alpha_\lambda\rightarrow0
\end{align}
by \eqref{eq:alpha-high}. This proves \eqref{eq:w-high}.

\proofstep{Proof of \eqref{eq:winfty-high}} By \eqref{eq:w-scaling} with $\varepsilon=1/\lambda$,
$\lambda^{-1}\norm{w_\lambda(q)}_\infty=\norm{w_1(q^{(1/\lambda)})}_\infty$. The map $q\mapsto w_1(q)$ is continuous from a neighbourhood of the origin in $H^{-1}$ into $H^1$ by Proposition~\ref{prop:static}\,\textup{(ii)}, and it vanishes at $q=0$: the kernel of the free resolvent $(-\partial_x^2+1)^{-1}$ has diagonal $g_1(x;0)=\tfrac12$, whence $h_1(\cdot;0)=1$ and $w_1(\cdot;0)=0$ by \eqref{eq:intro-resolvent}. Consequently $\norm{w_1(q)}_{H^1}\to0$ as $q\to0$ in $H^{-1}$. By \eqref{eq:scaled-Hminus}, $\norm{q^{(1/\lambda)}}_{H^{-1}}\leq\lambda^{-\frac12}M_\cQ$, so the argument $q^{(1/\lambda)}$ lies in that neighbourhood and tends to $0$ there, uniformly on $\cQ$ and at a rate depending on $\cQ$ only through $M_\cQ$. The embedding $H^1\hookrightarrow L^\infty$ of Lemma~\ref{lem:products}\,\textup{(i)} then gives \eqref{eq:winfty-high}.

\proofstep{Proof of \eqref{eq:high-energy-bounds}} For $\lambda\geq1$ the multiplier in \eqref{eq:alpha-large} satisfies $(\xi^2+4\lambda^2)^{-1}\leq(4+\xi^2)^{-1}$, so \eqref{eq:alpha-large} gives $\lambda\alpha_\lambda(q)\lesssim\norm{q}_{H^{-1}}^2$, with the absolute constant of \eqref{eq:KV-alpha}. In \eqref{eq:high-energy-a5}, the operator $4\lambda^2\Rf_\lambda$ has norm at most one on $H^{-1}$, its symbol being bounded by one, and \eqref{eq:high-energy-a7} bounds the second term by a constant multiple of $\lambda\alpha_\lambda(q)$. Hence $\norm{\lambda w_\lambda(q)}_{H^{-1}}\lesssim\norm{q}_{H^{-1}}+\norm{q}_{H^{-1}}^2$, which proves \eqref{eq:high-energy-bounds}.
\end{proof}

The following lemma expresses the increment of an integrable combination of densities as a spatial derivative of an $H^{-1}$ distribution, as required by Lemma~\ref{lem:zero-mass}. Identity \eqref{eq:S-evolution} holds for each pair of parameters in the common spectral range. The weights $\kappa$ and $\lambda$ combine the terms of \eqref{eq:def-rho} that are linear in $q$ into a second derivative by \eqref{eq:multiplier-identity}.

{\begin{lemma}\label{lem:two-parameter}
Let $q\in C(I;H^{-1}(\R;\R))$ be a Green solution in the sense of Definition~\ref{def:green}, and let $\kappa,\lambda$ be two parameters that are simultaneously admissible for the trajectory and at least $\lambda_*(q,I)$, so that \eqref{eq:integrated-w} holds at both. Such a common range exists by Remark~\ref{rem:common-lambda}. Assume also that $\kappa,\lambda\geq\Lambda(q(I))$ from Lemma~\ref{lem:high-energy}, so that both densities are nonnegative and integrable.

Set
\begin{align}
 \label{eq:two-parameter-a1} S_{\kappa\lambda}:=\kappa\rho_\kappa-\lambda\rho_\lambda, \qquad B_{\kappa\lambda}:=\tfrac12(\Rf_\kappa-\Rf_\lambda).
\end{align}
Then $B_{\kappa\lambda}$ maps $H^{-1}$ boundedly into $H^{1}$, the two free operators satisfy the exact identity
\begin{align}
 \label{eq:multiplier-identity} 2\bigl(\kappa^2\Rf_\kappa-\lambda^2\Rf_\lambda\bigr) =\partial_x^2B_{\kappa\lambda},
\end{align}
and, for all $s,t\in I$,
\begin{align}
 \label{eq:S-evolution} S_{\kappa\lambda}(t)-S_{\kappa\lambda}(s) =\partial_x\bigg[ -\kappa\int_s^t\J_\kappa(q(t'))\dd t' +\lambda\int_s^t\J_\lambda(q(t'))\dd t' +\partial_xB_{\kappa\lambda}\bigl(q(t)-q(s)\bigr)\bigg].
\end{align}
The identity \eqref{eq:S-evolution} holds in tempered distributions. Its left side belongs to $L^1(\R)$ and its bracketed primitive belongs to $H^{-1}(\R)$, which verifies the hypotheses of Lemma~\ref{lem:zero-mass}. The flux integrals are Bochner integrals in $H^{-1}(\R)$ by part~\textup{(i)} of Lemma~\ref{lem:microscopic}.
\end{lemma}}

\begin{proof}
Both sides of \eqref{eq:multiplier-identity} are Fourier multipliers, equal respectively to
\begin{align}
 \label{eq:two-parameter-a2} 2\left(\frac{\kappa^2}{\xi^2+4\kappa^2}-\frac{\lambda^2}{\xi^2+4\lambda^2}\right) \quad\text{and}\quad -\frac{\xi^2}{2} \left(\frac{1}{\xi^2+4\kappa^2}-\frac{1}{\xi^2+4\lambda^2}\right),
\end{align}
and both reduce to $2\xi^2(\kappa^2-\lambda^2)\bigl[(\xi^2+4\kappa^2)(\xi^2+4\lambda^2)\bigr]^{-1}$. The same computation identifies the multiplier of $B_{\kappa\lambda}$ as $2(\lambda^2-\kappa^2)\bigl[(\xi^2+4\kappa^2)(\xi^2+4\lambda^2)\bigr]^{-1}$, which is $O\bigl((1+\xi^2)^{-2}\bigr)$. A Fourier multiplier with symbol $\Theta$ maps $H^{-1}$ boundedly into $H^1$ as soon as $(1+\xi^2)\Theta(\xi)$ is bounded, since then $(1+\xi^2)\abs{\Theta}^2(4+\xi^2)$ is bounded. Hence $B_{\kappa\lambda}\colon H^{-1}\to H^1$.

Substituting \eqref{eq:def-rho} into the definition of $S_{\kappa\lambda}$ and then using \eqref{eq:multiplier-identity} gives the exact spatial identity
\begin{align}
 \label{eq:S-expanded} S_{\kappa\lambda} =-\kappa w_\kappa+\lambda w_\lambda +2\bigl(\kappa^2\Rf_\kappa-\lambda^2\Rf_\lambda\bigr)q =-\kappa w_\kappa+\lambda w_\lambda +\partial_x^2B_{\kappa\lambda}q .
\end{align}
The weights are forced, and it is worth seeing why. A general linear combination $a\rho_\kappa+b\rho_\lambda$ has, by \eqref{eq:def-rho}, the linear-in-$q$ part $2\bigl(a\kappa \Rf_\kappa+b\lambda \Rf_\lambda\bigr)$, whose symbol
\begin{align}
 \label{eq:general-weights-symbol} 2\left(\frac{a\kappa}{\xi^2+4\kappa^2}+\frac{b\lambda}{\xi^2+4\lambda^2}\right) \qquad\text{takes the value}\qquad \frac12\left(\frac a\kappa+\frac b\lambda\right) \quad\text{at }\xi=0 .
\end{align}
An operator of the form $\partial_x^2(\,\cdot\,)$ has a symbol vanishing to second order at the origin, so the combination can be a second derivative only if $a/\kappa+b/\lambda=0$, that is, only if $(a,b)$ is proportional to $(\kappa,-\lambda)$. This determines $S_{\kappa\lambda}$ up to an irrelevant scalar. Conversely, that necessary condition is already sufficient: the symbol in \eqref{eq:general-weights-symbol} is an even function of $\xi$, so its first-order term at the origin vanishes automatically, and killing the value at $\xi=0$ therefore produces a full double zero. The required mapping property follows from the explicit factorization \eqref{eq:multiplier-identity}, whose factor $B_{\kappa\lambda}$, defined in \eqref{eq:two-parameter-a1}, is the operator $\frac12(\Rf_\kappa-\Rf_\lambda)$, bounded from $H^{-1}$ to $H^1$. For $b=0$ and $a\neq0$, the value $a/(2\kappa)\neq0$ gives the obstruction at a single parameter described at the start of the section. Now subtract \eqref{eq:S-expanded} at the times $s$ and $t$, and use the integrated Green identity \eqref{eq:integrated-w} at the parameters $\kappa$ and $\lambda$ to replace $w_\kappa(t)-w_\kappa(s)$ and $w_\lambda(t)-w_\lambda(s)$:
\begin{align}
 \label{eq:two-parameter-a3} S_{\kappa\lambda}(t)-S_{\kappa\lambda}(s) =-\kappa\partial_x\int_s^t\J_\kappa(q(t'))\dd t' +\lambda\partial_x\int_s^t\J_\lambda(q(t'))\dd t' +\partial_x\bigl[\partial_xB_{\kappa\lambda}(q(t)-q(s))\bigr].
\end{align}
Factoring out the common derivative gives \eqref{eq:S-evolution}. As for the regularity, the left side lies in $L^1$ because both densities do, by \eqref{eq:rho-nonneg}. In the bracket, the two flux integrals lie in $H^{-1}$ by part~\textup{(i)} of Lemma~\ref{lem:microscopic}, and $\partial_xB_{\kappa\lambda}(q(t)-q(s))\in L^2\subset H^{-1}$ because $B_{\kappa\lambda}\colon H^{-1}\to H^1$.
\end{proof}

%%%%%%%%%%%%%%%%%%%%%%%%%%%%%%%%%%%%%%%%%%%%%%%%%%%%%%%%%%%%%%%%%%%%%%%%%%%%%%%%%%
\subsection{Conservation Law and Coercivity}
%%%%%%%%%%%%%%%%%%%%%%%%%%%%%%%%%%%%%%%%%%%%%%%%%%%%%%%%%%%%%%%%%%%%%%%%%%%%%%%%%%

{\begin{proposition}\label{prop:alpha-conservation}
Let $q\in C(I;H^{-1})$ be a Green solution, and let $K(q)$ be the threshold \eqref{eq:single-K}. Then, for every $\kappa\geq K(q)$,
\begin{align}
 \label{eq:alpha-conservation-a2} \alpha_\kappa(q(t))=\alpha_\kappa(q(s)) \qquad(s,t\in I).
\end{align}
\end{proposition}}

\begin{proof}
Fix $\kappa\geq K(q)$ and $s,t\in I$, and let $\lambda\geq K(q)$ be a second parameter. Lemma~\ref{lem:two-parameter} places the left side of \eqref{eq:S-evolution} in $L^1$ and its bracketed primitive in $H^{-1}$. Applying Lemma~\ref{lem:zero-mass} directly gives $\int_\R\big[S_{\kappa\lambda}(t)-S_{\kappa\lambda}(s)\big]=0$.

Since $\int_\R S_{\kappa\lambda}=\kappa\alpha_\kappa-\lambda\alpha_\lambda$ by \eqref{eq:def-alpha}, this reads
\begin{align}
 \label{eq:cross-alpha} \kappa\bigl[\alpha_\kappa(q(t))-\alpha_\kappa(q(s))\bigr] =\lambda\bigl[\alpha_\lambda(q(t))-\alpha_\lambda(q(s))\bigr].
\end{align}
Identity \eqref{eq:cross-alpha} is already the substance of the proposition, and it is worth reading it on its own terms: it says that the product
\begin{align}
 \label{eq:alpha-conservation-a4} \lambda\bigl[\alpha_\lambda(q(t))-\alpha_\lambda(q(s))\bigr]
\end{align}
takes the same value at every $\lambda\geq K(q)$.

It remains only to evaluate it, and the natural place to do so is at $\lambda=\infty$, where Lemma~\ref{lem:high-energy} shows that both terms separately vanish. Accordingly, fix the left side of \eqref{eq:cross-alpha}, which is independent of $\lambda$, and let $\lambda\to\infty$ on the right. By nonnegativity \eqref{eq:rho-nonneg} and the uniform limit \eqref{eq:alpha-high} of Lemma~\ref{lem:high-energy}, applied to the compact set $q(I)$,
\begin{align}
 \label{eq:alpha-conservation-a5} \abs[\big]{\lambda[\alpha_\lambda(q(t))-\alpha_\lambda(q(s))]} \leq\lambda\alpha_\lambda(q(t))+\lambda\alpha_\lambda(q(s)) \rightarrow0 .
\end{align}
Since the left side of \eqref{eq:cross-alpha} is independent of $\lambda$, it must equal this limit, namely zero. As $\kappa>0$, the assertion follows.
\end{proof}

It remains to record the structural fact about the Green flux that Section~\ref{sec:automatic-smoothing} will use. Written as in \eqref{eq:intro-green-flux}, the flux $\J_\lambda$ contains a term linear in $q$, multiplied by a function of $h_\lambda$, with a sign that may vary. Substituting the recovery identity \eqref{eq:q-from-h}, which expresses $q$ through $h_\lambda$ alone, converts it into an expression entirely in $w_\lambda$: a linear term $-w_{\lambda,xx}$, plus a remainder $N_\lambda$ that is nonlinear.

Pairing the linear term with the weight leaves $-\pair*{w_\lambda}{\phi_y'''}$ by \eqref{eq:wxx-transfer}. This term is retained in \eqref{eq:key-weighted-identity} and controlled uniformly by \eqref{eq:w-high}. For the nonlinear remainder, the cubic $2w_\lambda^3$ changes sign and is dominated by the quadratic once $\abs{w_\lambda}\leq\lambda/2$, in the form
\begin{align}
 \label{eq:cubic-domination} 6\lambda w_\lambda^2+2w_\lambda^3\geq5\lambda w_\lambda^2 \qquad\text{whenever }\abs{w_\lambda}\leq\frac\lambda2 ,
\end{align}
which is what part~\textup{(ii)} below uses.

That domination is the coercivity referred to in the introduction, and part~\textup{(ii)} makes it quantitative: the bound \eqref{eq:N-positive} controls both $\abs{w_{\lambda,x}}^2$ and $\lambda^2\abs{w_\lambda}^2$ by $\lambda N_\lambda$. Since $2\lambda w_\lambda\to q$ by \eqref{eq:w-high}, the second of these is a bound for $q$. This is exactly the mechanism that produces Kato smoothing in Proposition~\ref{prop:auto-smoothing}. Part~\textup{(ii)} needs the pointwise smallness $\norm{w_\lambda}_\infty\leq\lambda/2$, which is why the lemma is placed after Lemma~\ref{lem:high-energy}.

{\begin{lemma}\label{lem:flux-decomposition}
Let $q\in H^{-1}(\R;\R)$ and let $\lambda$ be admissible for $q$.
\begin{enumerate}[
label=\textup{(\roman*)}, leftmargin=*, itemsep=.7em, before={\let\fullwidthdisplay\relax} ]
\item \emph{Decomposition.} The Green flux \eqref{eq:green-flux} decomposes as
\begin{align}
 \label{eq:J-decomposition} \J_\lambda(q)=-w_{\lambda,xx}+N_\lambda, \qquad N_\lambda :=\frac{3w_{\lambda,x}^2}{2h_\lambda} +6\lambda w_\lambda^2+2w_\lambda^3 ,
\end{align}
where $-w_{\lambda,xx}\in H^{-1}(\R)$ and $N_\lambda\in L^1(\R)\subset H^{-1}(\R)$.

\item \emph{Coercivity.} If in addition $\norm{w_\lambda}_{L^\infty}\leq\lambda/2$, then $\lambda/2\leq h_\lambda\leq3\lambda/2$ and, for almost every $x\in\R$,
\begin{align}
 \label{eq:N-positive} \lambda N_\lambda =\frac{3\lambda}{2h_\lambda}\abs{w_{\lambda,x}}^2 +2\lambda w_\lambda^2\bigl(3\lambda+w_\lambda\bigr) \geq\abs{w_{\lambda,x}}^2+5\lambda^2\abs{w_\lambda}^2 .
\end{align}
In particular $N_\lambda\geq0$.

\item \emph{Along a trajectory.} If $q\in C(I;H^{-1})$ and $\lambda$ is admissible for $q(I)$, then
\begin{align}
 \label{eq:flux-decomposition-a1} \inf_{(t,x)\in I\times\R}h_\lambda(t,x)\geq\tfrac23\lambda>0,
\end{align}
and $t\mapsto N_\lambda(t)$ is continuous from $I$ into $L^1(\R)$. Moreover the hypothesis of \textup{(ii)}, and hence \eqref{eq:N-positive}, holds at every $t\in I$ once $\lambda$ exceeds a threshold depending only on $q(I)$.
\end{enumerate}
\end{lemma}}

\begin{proof}
\textup{(i)} All the operations below are legitimate in $H^{-1}$: by Proposition~\ref{prop:static}\,\textup{(ii)} both $h_\lambda=\lambda+w_\lambda$ and $h_\lambda^{-1}=2g_\lambda$ are constants plus $H^1$ functions, $h_\lambda$ is bounded away from zero by \eqref{eq:g-quantitative}, and \eqref{eq:affine-multiplier-definition}--\eqref{eq:multipliers-a4} permit the resulting multiplications and reciprocal cancellations. Regrouping the expansion \eqref{eq:green-flux} by $2h_\lambda=2\lambda+2w_\lambda$ gives
\begin{align}
 \label{eq:J-2hq} \J_\lambda=2h_\lambda q-4\lambda^2w_\lambda .
\end{align}
Now insert the recovery identity \eqref{eq:q-from-h} for $q$. Multiplying it by $2h_\lambda$ clears every denominator but one, and, since $h_{\lambda,x}=w_{\lambda,x}$ and $h_{\lambda,xx}=w_{\lambda,xx}$,
\begin{align}
 \label{eq:2hq-expanded} 2h_\lambda q =-w_{\lambda,xx} +\frac{3w_{\lambda,x}^2}{2h_\lambda} +2h_\lambda^3-2\lambda^2h_\lambda .
\end{align}
Subtracting $4\lambda^2w_\lambda=4\lambda^2(h_\lambda-\lambda)$ as \eqref{eq:J-2hq} requires, and then using $h_\lambda=\lambda+w_\lambda$ in the resulting polynomial, gives
\begin{align}
 \label{eq:flux-decomposition-a2} \J_\lambda &=-w_{\lambda,xx} +\frac{3w_{\lambda,x}^2}{2h_\lambda} +2h_\lambda^3-6\lambda^2h_\lambda+4\lambda^3\\
 &=-w_{\lambda,xx} +\frac{3w_{\lambda,x}^2}{2h_\lambda} +6\lambda w_\lambda^2+2w_\lambda^3 ,\notag
\end{align}
the last step being the expansion of $2h_\lambda^3-6\lambda^2h_\lambda+4\lambda^3=2(h_\lambda-\lambda)^2(h_\lambda+2\lambda)$ at $h_\lambda=\lambda+w_\lambda$. Since $w_\lambda\in H^1$ and $h_\lambda^{-1}\in L^\infty$, each of the three terms of $N_\lambda$ lies in $L^1$:
\begin{align}
 \label{eq:flux-decomposition-a3} \Bnorm{\frac{w_{\lambda,x}^2}{h_\lambda}}_1 \leq\norm{h_\lambda^{-1}}_\infty\norm{w_{\lambda,x}}_2^2, \qquad \norm{w_\lambda^2}_1=\norm{w_\lambda}_2^2, \qquad \norm{w_\lambda^3}_1\leq\norm{w_\lambda}_\infty\norm{w_\lambda}_2^2 .
\end{align}
Finally $w_\lambda\in H^1$ gives $w_{\lambda,xx}\in H^{-1}$, and $L^1\subset H^{-1}$ is part~\textup{(iii)} of Lemma~\ref{lem:products}.

\textup{(ii)} From $h_\lambda=\lambda+w_\lambda$ and $\abs{w_\lambda}\leq\lambda/2$ we get $\lambda/2\leq h_\lambda\leq3\lambda/2$. The displayed equality in \eqref{eq:N-positive} is \eqref{eq:J-decomposition} multiplied by $\lambda$, after writing $6\lambda w_\lambda^2+2w_\lambda^3=2w_\lambda^2(3\lambda+w_\lambda)$. For the inequality, the upper bound $h_\lambda\leq3\lambda/2$ gives $3\lambda/(2h_\lambda)\geq1$, and $w_\lambda\geq-\lambda/2$ gives $3\lambda+w_\lambda\geq5\lambda/2$, so that $2\lambda w_\lambda^2(3\lambda+w_\lambda)\geq5\lambda^2w_\lambda^2$.

\textup{(iii)} The image $q(I)$ is compact, so \eqref{eq:g-quantitative} applies uniformly and gives $h_\lambda=(2g_\lambda)^{-1}\geq\frac23\lambda$ on $I\times\R$. Consequently
\begin{align}
 \label{eq:flux-decomposition-a4} \norm{h_\lambda(t)^{-1}-h_\lambda(s)^{-1}}_\infty \leq\frac{\norm{h_\lambda(t)-h_\lambda(s)}_\infty} {\bigl(\inf_{I\times\R}h_\lambda\bigr)^2},
\end{align}
so $t\mapsto h_\lambda(t)^{-1}$ is continuous into $L^\infty$. Together with $w_\lambda\in C(I;H^1)$, which is Proposition~\ref{prop:static}\,\textup{(ii)}, and elementary estimates such as $\norm{a_x^2-b_x^2}_1\leq\norm{a_x-b_x}_2(\norm{a_x}_2+\norm{b_x}_2)$, this makes each of the three terms of $N_\lambda$ continuous into $L^1$. The last assertion is \eqref{eq:winfty-high} of Lemma~\ref{lem:high-energy} applied to the compact set $q(I)$: it gives $\lambda^{-1}\norm{w_\lambda}_\infty\to0$ uniformly on $q(I)$, so in particular $\norm{w_\lambda}_\infty\leq\lambda/2$ for all large $\lambda$, uniformly in $t\in I$.
\end{proof}

% \newpage

%%%%%%%%%%%%%%%%%%%%%%%%%%%%%%%%%%%%%%%%%%%%%%%%%%%%%%%%%%%%%%%%%%%%%%%%%%%%%%%%%%
\section{Local Smoothing}
\label{sec:automatic-smoothing}
%%%%%%%%%%%%%%%%%%%%%%%%%%%%%%%%%%%%%%%%%%%%%%%%%%%%%%%%%%%%%%%%%%%%%%%%%%%%%%%%%%

% We use the translated weight
% \begin{align}
%  \label{eq:tanh-weight} \phi_y(x)=\tanh(x-y), \qquad \phi_y'(x)=\sech^2(x-y).
% \end{align}

This section derives local smoothing from Definition~\ref{def:green}. A Green solution is a continuous $H^{-1}$-valued curve satisfying \eqref{eq:green-identity-recalled} at every large parameter. Proposition~\ref{prop:auto-smoothing} shows that these identities imply that $q$ is square integrable in spacetime on every translated spatial cell of radius one, with a bound uniform in the translation and a tail vanishing at spatial infinity.
The conclusion is a regularity upgrade from a curve of distributions to the representative specified in \eqref{eq:representative-def}.

% {
% \begin{remark}
% The conclusion is a regularity upgrade from a curve of distributions to the representative specified in \eqref{eq:representative-def}. Step~4 constructs it as a weak $L^2_{\loc}$ limit of $2\lambda w_\lambda$. Remark~\ref{rem:measurability} supplies a common exceptional set of times, independent of the cell and its centre.
% \end{remark}
% }

Since the argument is somewhat long, we describe it first.
We use the translated weight
\begin{align}
 \label{eq:tanh-weight} \phi_y(x)=\tanh(x-y), \qquad \phi_y'(x)=\sech^2(x-y).
\end{align}

\textit{The mechanism.} Everything rests on the coercivity of Lemma~\ref{lem:flux-decomposition}\,\textup{(ii)}. The two-parameter identity \eqref{eq:S-evolution} is paired with the translated profile $\phi_y$ itself, through a compact cutoff, and the outer spatial derivative is then transferred onto the test function. The resulting weight is $\phi_y'>0$. The result puts the quantity $\lambda\iint\phi_y'N_\lambda$ on one side, and on the other four terms that will be bounded uniformly both in $\lambda$ and in the translation centre $y$. By \eqref{eq:N-positive} the left side dominates $\iint\phi_y'\abs{2\lambda w_\lambda}^2$, and by \eqref{eq:w-high} $2\lambda w_\lambda$ converges to $q$. A uniform bound on $2\lambda w_\lambda$ is thereby converted into a uniform weighted $L^2$ bound for $q$, which is \eqref{eq:q-Kato-bound}.

\textit{Where the conservation law enters.} The tail \eqref{eq:q-Kato-tail} is the same computation carried one step further. After $\lambda\to\infty$ the four terms become pairings of \emph{fixed} objects against translates of a fixed profile. Three of them vanish as $y\to\pm\infty$ by Lemma~\ref{lem:translation}, because an $H^{-1}$ distribution paired with a translated $H^1$ function tends to zero. The fourth is different: it pairs the $L^1$ function $\kappa[\rho_\kappa(q(b))-\rho_\kappa(q(a))]$ with $\phi_y$ itself, and $\phi_y(x)\to\mp1$ for each fixed $x$ as $y\to\pm\infty$. That pairing therefore tends to $\mp\int_\R\kappa[\rho_\kappa(q(b))-\rho_\kappa(q(a))]$, which is exactly $\mp\kappa[\alpha_\kappa(q(b))-\alpha_\kappa(q(a))]$, and it vanishes precisely because $\alpha_\kappa$ is conserved. This is the one place in this section where Proposition~\ref{prop:alpha-conservation} is used. It is used once more, through Lemma~\ref{lem:Gardner-mass}, to bound the Gardner variable of Section~\ref{sec:Gardner} in $L^2$ uniformly in time.

\textit{The weight.} These two requirements already suggest a bounded, nonconstant, nondecreasing test function. Coercivity asks for $\phi_y'\geq0$, and the mass term asks that $\phi_y$ tend to different limits at $\pm\infty$. The second is forced by the first, because
\begin{align}
 \label{eq:weight-not-compactly-supported} \phi\in C_c^\infty(\R) \quad\Longrightarrow\quad \int_\R\phi'(x)\dd x=0 ,
\end{align}
so a compactly supported weight cannot have a nonnegative derivative unless that derivative vanishes identically. We therefore choose a bounded increasing function with distinct limits. Its derivative may be compactly supported: if $\chi\in C_c^\infty(\R)$ is nonnegative, positive on $[-1,1]$, and normalized by $\int_\R\chi=2$, then
\begin{align}
 \label{eq:weight-alternative} \phi(x):=-1+\int_{-\infty}^x\chi(z)\,dz
\end{align}
satisfies $\phi'=\chi$ and $\phi(\pm\infty)=\pm1$, and its translates carry the same coercivity and endpoint-mass mechanism. The choice \eqref{eq:tanh-weight} is convenient because $\phi_y'=\sech^2(\cdot-y)$ decays exponentially, so that it and its derivatives are translates of fixed $H^1$ functions, while $\phi_y$ itself is bounded, pairs with an $L^1$ density, and tends to $\mp1$ at $\mp\infty$, which is what exposes the conserved mass.

The translated family $\{\phi_y\}_{y\in\R}$ gives a bound uniform in the centre through
\begin{align}
 \label{eq:weight-cell-lower-bound} \inf_{\abs{x-y}\leq1}\phi_y'(x)=\sech^2(1)>0 \qquad\text{for every }y\in\R,
\end{align}
so that for every measurable $f$,
\begin{align}
 \label{eq:weight-cell-comparison} \int_{y-1}^{y+1}\abs{f(x)}^2\dd x \leq\cosh^2(1)\int_\R\phi_y'(x)\abs{f(x)}^2\dd x ,
\end{align}
with a constant independent of $y$: the weight moves with the cell, and every norm of a derivative profile appearing in the estimates is a translate of a fixed one. The need to translate the weight with the cell is illustrated by the following example. If $f\in C_c^\infty((-\frac12,\frac12))$ with $\norm{f}_2=1$ and $f_R:=f(\cdot-R)$, then the cell mass of $f_R$ at the centre $R$ equals one for every $R$, whereas $\int_\R\sech^2(x)\abs{f_R(x)}^2\dd x\leq\sech^2\bigl(R-\tfrac12\bigr)\to0$ as $R\to+\infty$.

{\begin{remark}[Why the estimate needs two spectral parameters]
\label{rem:why-two-parameters}
Local smoothing is obtained in \cite{Killip.Visan2019} from the dynamics of the single density $\rho_\kappa$. For the Green solutions of Definition~\ref{def:green}, a difference of two densities makes the required weighted pairings finite.

The weight $\phi_y$ is bounded and approaches nonzero limits at spatial infinity, so $\phi_y\notin H^1(\R)$. Its pairing with an $L^1$ function is well-defined by integration. Its derivative $\phi_y'\in H^1(\R)$ also admits the Sobolev pairing with an $H^{-1}$ distribution. These are the two pairings used below.

The Green identity governs $w_\lambda$. Its right side can be paired with the weight after integration by parts. In the form \eqref{eq:integrated-w} it already exhibits the increment of $w_\lambda$ as $\partial_x$ of an $H^{-1}$ distribution, so one derivative moves onto the weight and leaves an $H^{-1}$--$H^1$ pairing against $\phi_y'$. The left side has the regularity $w_\lambda\in H^1$. The operator $\Rf_\lambda$ maps $H^{-1}(\R)$ onto $H^1(\R)$, and $H^1(\R)\not\subset L^1(\R)$, so \eqref{eq:def-rho} requires additional control of the free term to establish integrability of $w_\lambda$. Thus the pairing $\pair*{w_\lambda(b)-w_\lambda(a)}{\phi_y}$ need not converge.

For the density, by \eqref{eq:rho-nonneg}, $\rho_\lambda\geq0$ with $\norm{\rho_\lambda}_1=\alpha_\lambda$, so the increment $\rho_\lambda(b)-\rho_\lambda(a)$ is integrable and pairs with $\phi_y$ without difficulty. Subtracting \eqref{eq:def-rho} at the two times and using \eqref{eq:integrated-w} gives
\begin{align}
 \label{eq:single-parameter-obstruction} \rho_\lambda(b)-\rho_\lambda(a) =\partial_x\bigg[-\int_a^b\J_\lambda(q(t'))\dd t'\bigg] +2\lambda \Rf_\lambda\bigl(q(b)-q(a)\bigr),
\end{align}
where the second term has regularity $H^1$. The symbol $2\lambda(\xi^2+4\lambda^2)^{-1}$ of $2\lambda \Rf_\lambda$ equals $(2\lambda)^{-1}\neq0$ at $\xi=0$, which prevents factoring it as a derivative with the required $H^{-1}$ primitive. Pairing this term with $\phi_y$ would again require additional integrability.

The combination $S_{\kappa\lambda}$ of \eqref{eq:two-parameter-a1} removes both defects at once. Its increment is a difference of $L^1$ functions, so the left side pairs with $\phi_y$. The two obstructing terms of \eqref{eq:single-parameter-obstruction}, taken at $\kappa$ and at $\lambda$, fuse by the multiplier identity \eqref{eq:multiplier-identity} into the single exact second derivative $\partial_x^2B_{\kappa\lambda}$, which is therefore admitted under the outer derivative of \eqref{eq:S-evolution}. One of its two derivatives moves onto the weight, and what survives is the pairing of the two $L^2$ functions $\partial_xB_{\kappa\lambda}(q(b)-q(a))$ and $\phi_y'$. This is the third term of \eqref{eq:key-weighted-identity}, and it is finite for that reason alone.

For a Schwartz solution, $w_\lambda$ and $\Rf_\lambda q$ decay rapidly, every pairing above converges, and the single-parameter argument of \cite{Killip.Visan2019} applies. That estimate passes to limits of Schwartz solutions by continuity of the solution map, as discussed in Remark~\ref{rem:why-green}. For an arbitrary Green solution, we use the two-parameter identity of Lemma~\ref{lem:two-parameter} and spatial cutoffs of the test function. The same identity supplies an $H^{-1}$ primitive for Lemma~\ref{lem:zero-mass} in Proposition~\ref{prop:alpha-conservation} and finite weighted pairings in the present proof.
\end{remark}}

% \exstep{Comparison with \cite{Killip.Visan2019}.} Local smoothing is available there because the solutions are limits of Schwartz solutions: one proves an a~priori bound along smooth trajectories, where every manipulation is classical, and passes it to the limit. As explained in Remark~\ref{rem:why-not-KV}, no such route exists for an arbitrary Green solution, and the estimate must instead be produced from the identity itself. The proof below accordingly never approximates the trajectory in time; the only approximations used are in the spatial variable, and only to legitimate a test function.

\textit{Structure of the proof.} Step~1 legitimates the non-compactly-supported weight and produces the exact weighted identity \eqref{eq:key-weighted-identity}. Step~2 bounds the left side of that identity from below by a weighted $L^2$ norm, thereby reducing the proposition to an upper bound for its right side. Step~3 bounds the four terms on the right of \eqref{eq:key-weighted-identity} uniformly in $\lambda$ and $y$. Step~4 extracts a weak limit of $2\lambda w_\lambda$, identifies it with $q$, and proves \eqref{eq:q-Kato-bound}. Step~5 takes $\lambda\to\infty$ first and $y\to\pm\infty$ afterwards, and proves \eqref{eq:q-Kato-tail}.

{\begin{restatement}
Let $q\in C(I;H^{-1}(\R;\R))$ be an arbitrary Green solution. For every compact time interval $J\subset I$ there is a representative of $q$ in $L^2_{\loc}(J\times\R)$, in the sense of \eqref{eq:representative-def}, satisfying the uniform cell bound \eqref{eq:q-Kato-bound} and the translated-cell tail \eqref{eq:q-Kato-tail}. Thus an arbitrary Green solution satisfies, locally in time, the uniformly local Kato smoothing bound and the translated-cell tail condition.
\end{restatement}}

\begin{proof}[Proof of Proposition~\ref{prop:auto-smoothing}]
Write $J=[a,b]$. For every $\lambda\geq K(q)$, Remark~\ref{rem:common-lambda} supplies the bounds for $g_\lambda$ and its reciprocal, the nonnegative integrable density, and the Green identity uniformly along $q(I)$.

We make one further enlargement of the parameter range. In Step~2 a value
\begin{align}
 \label{eq:sharp-threshold} K_\sharp\geq K(q)
\end{align}
is fixed, above which the pointwise smallness \eqref{eq:winfty-half} of $w_\lambda/\lambda$ holds as well. It is the \emph{enlargement} beyond $K(q)$ that depends on the trajectory alone, through $\sup_{t\in I}\norm{q(t)}_{H^{-1}}$. $K_\sharp$ itself still carries the Green-law parameter $\lambda_*(q,I)$, inherited from $K(q)$ through \eqref{eq:single-K}. All parameters below are taken $\geq K_\sharp$ from the outset. The symbol $K_\sharp$ is local to this proof.

Two parameters are drawn from this range, and their roles are different. A \emph{base} parameter $\kappa\geq K_\sharp$ is fixed for the whole proof. Holding it fixed is what keeps the conservation law usable: $\rho_\kappa$ remains one fixed nonnegative $L^1$ function, so the endpoint difference $\kappa[\rho_\kappa(q(b))-\rho_\kappa(q(a))]$ that survives the limit in Step~5 is an $L^1$ function whose total integral vanishes, by Proposition~\ref{prop:alpha-conservation}. A \emph{second} parameter $\lambda\geq K_\sharp$ is then sent to infinity. It is what makes the flux remainder coercive, through \eqref{eq:winfty-high} and \eqref{eq:N-positive}, and what identifies the limit of $2\lambda w_\lambda$ with $q$, through \eqref{eq:w-high}. Compare Remark~\ref{rem:full-spectral-family}.

\proofstep{Step 1}
This step turns the distributional identity \eqref{eq:S-evolution} into a scalar identity in which the coercive quantity stands alone. There are three tasks. First, admit $\phi_y$ as a test function, although it is neither compactly supported nor an element of $H^1(\R)$. Second, record the resulting weak formulation. Third, split the flux at the parameter $\lambda$ so that the term to be bounded from below is isolated. Throughout the step $y$ is fixed.

\emph{The cutoff.} Choose $\chi\in C_c^\infty$ equal to one near zero and put $\chi_R(x)=\chi(x/R)$, so that $\chi_R\phi_y\in C_c^\infty(\R)$ is a legitimate test function at every $R$. For each fixed $y$,
\begin{align}
 \label{eq:cutoff-convergence} (\chi_R\phi_y)'\rightarrow\phi_y' \quad\text{in }H^1.
\end{align}
Indeed,
\begin{align}
 \label{eq:auto-smoothing-a1} \begin{aligned}
 (\chi_R\phi_y)'-\phi_y' &=(\chi_R-1)\phi_y'+\chi_R'\phi_y,\\
 \partial_x\bigl[(\chi_R\phi_y)'-\phi_y'\bigr] &=(\chi_R-1)\phi_y''+2\chi_R'\phi_y'+\chi_R''\phi_y.
 \end{aligned}
\end{align}
For fixed $y$, the terms containing $\chi_R-1$ converge to zero in $L^2$ by the exponential decay of $\phi_y'$ and $\phi_y''$. The remaining terms satisfy
\begin{align}
 \label{eq:auto-smoothing-a2} \norm{\chi_R'\phi_y}_2=O(R^{-\frac12}),\qquad \norm{\chi_R'\phi_y'}_2=O(R^{-1}),\qquad \norm{\chi_R''\phi_y}_2=O(R^{-\frac32}).
\end{align}
This proves \eqref{eq:cutoff-convergence}.

Pair \eqref{eq:S-evolution}, taken at $s=a$ and $t=b$, with the admissible test function $\chi_R\phi_y$ and move the outer derivative onto it,
\begin{align}
 \label{eq:outer-derivative-transfer} &\bpair{S_{\kappa\lambda}(b)-S_{\kappa\lambda}(a)}{\chi_R\phi_y}=-\pair[\Big]{-\kappa\int_a^b\J_\kappa(q(t'))\dd t' +\lambda\int_a^b\J_\lambda(q(t'))\dd t' +\partial_xB_{\kappa\lambda}(q(b)-q(a))}{(\chi_R\phi_y)'} .
\end{align}
The two sides now pass to the limit $R\to\infty$ for different reasons, and this is exactly the asymmetry recorded in Lemma~\ref{lem:two-parameter}. On the left, $S_{\kappa\lambda}(b)-S_{\kappa\lambda}(a)\in L^1(\R)$ while $\chi_R\phi_y\to\phi_y$ pointwise with $\abs{\chi_R\phi_y}\leq\norm{\chi}_\infty$, so dominated convergence applies. On the right, the bracket lies in $H^{-1}(\R)$ and is therefore continuous against the $H^1$ convergence \eqref{eq:cutoff-convergence}.

\emph{The weak formulation.} Letting $R\to\infty$ in \eqref{eq:outer-derivative-transfer} and interchanging each time integral with the pairing, which is legitimate because a bounded linear functional commutes with a Bochner integral, gives
\begin{align}
 \label{eq:weak-form-phi} \bpair{S_{\kappa\lambda}(b)-S_{\kappa\lambda}(a)}{\phi_y} ={}&\kappa\int_a^b\pair*{\J_\kappa(q(t'))}{\phi_y'}\dd t' -\lambda\int_a^b\pair*{\J_\lambda(q(t'))}{\phi_y'}\dd t'\\
 &-\pair*{\partial_xB_{\kappa\lambda}(q(b)-q(a))}{\phi_y'} .\notag
\end{align}
Every pairing here is finite, and for three different reasons: the two flux integrals are between $H^{-1}(\R)$ and $\phi_y'\in H^1(\R)$. The left side is between $L^1(\R)$ and $\phi_y\in L^\infty(\R)$. $\partial_xB_{\kappa\lambda}(q(b)-q(a))\in L^2(\R)$, since $B_{\kappa\lambda}$ maps $H^{-1}$ into $H^1$, is paired with $\phi_y'\in L^2(\R)$. Identity \eqref{eq:weak-form-phi} is the weak form of the Green identity that the rest of the proof exploits.

\emph{The decomposition.} We split the flux at $\lambda$ by \eqref{eq:J-decomposition} to isolate the coercive remainder,
\begin{align}
 \label{eq:flux-split-paired} \lambda\int_a^b\pair*{\J_\lambda(q(t'))}{\phi_y'}\dd t' =\lambda\int_a^b\pair*{-w_{\lambda,xx}(t')}{\phi_y'}\dd t' +\lambda\int_a^b\int_\R\phi_y'N_\lambda\dd x\dd t' .
\end{align}
The last term is an ordinary double integral because $N_\lambda(t)\in L^1(\R)$ by \eqref{eq:J-decomposition} and $\phi_y'\in L^\infty(\R)$. It is the only term of the whole identity that will be bounded from below. In the remaining term both derivatives are moved onto the weight,
\begin{align}
 \label{eq:wxx-transfer} \pair*{-w_{\lambda,xx}}{\phi_y'} =-\pair*{w_\lambda}{\phi_y'''},
\end{align}
which is legitimate because $w_\lambda\in H^1$ by Proposition~\ref{prop:static}\,\textup{(ii)} and $\phi_y'''\in H^1$. The derivatives of $\phi_y'$ decay exponentially, which gives $\phi_y'''\in H^1$. The undifferentiated weight $\phi_y$ has nonzero limits at spatial infinity, so its use as a test function was justified by the cutoff argument above.

\emph{The rearrangement.} Substituting \eqref{eq:flux-split-paired} and \eqref{eq:wxx-transfer} into \eqref{eq:weak-form-phi} and solving for the term isolated above gives the exact identity
{
\begin{align}
 \label{eq:key-weighted-identity} \lambda\int_a^b\int_\R\phi_y'N_\lambda\dd x\dd t' &={}\kappa\int_a^b\pair*{\J_\kappa(q(t'))}{\phi_y'}\dd t' -\bpair{S_{\kappa\lambda}(b)-S_{\kappa\lambda}(a)}{\phi_y}\\
 &\quad-\pair*{\partial_xB_{\kappa\lambda}(q(b)-q(a))}{\phi_y'} +\lambda\int_a^b\pair*{w_\lambda(t')}{\phi_y'''}\dd t'.\notag
\end{align}
}
This is the form in which the identity will be used. Its left side is the single quantity to which the coercivity of Step~2 applies. Each of the four terms on its right pairs a translate of a fixed profile with an object whose relevant norm is bounded uniformly in $\lambda$ and in $y$. Step~3 bounds the four terms uniformly in $\lambda$, and Step~5 takes their limits.

\proofstep{Step 2}
This step records the positivity of the flux remainder, pointwise almost everywhere. It is where the smallness for sufficiently large $\lambda$ becomes a coercive bound.

By \eqref{eq:winfty-high} of Lemma~\ref{lem:high-energy}, applied to the compact trajectory $q(I)$, the ratio $\lambda^{-1}\norm{w_\lambda(t)}_{L^\infty}$ tends to zero as $\lambda\to\infty$, uniformly in $t\in I$. In particular
\begin{align}
 \label{eq:winfty-half} \norm{w_\lambda(t)}_{L^\infty}\leq\frac\lambda2 \qquad(t\in I)
\end{align}
holds once $\lambda$ is large enough. This is the enlargement announced at the start of the proof: fix $K_\sharp$ as in \eqref{eq:sharp-threshold}, so large that \eqref{eq:winfty-half} is valid for every $\lambda\geq K_\sharp$. By the proof of \eqref{eq:winfty-high}, the enlargement depends on $q(I)$ only through $\sup_{t\in I}\norm{q(t)}_{H^{-1}}$, hence neither on $J$, nor on $y$, nor on $\kappa$, and it is made once and for all.

Condition \eqref{eq:winfty-half} is exactly the hypothesis of part~\textup{(ii)} of Lemma~\ref{lem:flux-decomposition}. That lemma therefore gives $\lambda/2\leq h_\lambda\leq3\lambda/2$ on $J\times\R$, together with the pointwise bound \eqref{eq:N-positive},
\begin{align}
 \label{eq:auto-smoothing-a3} \lambda N_\lambda\geq\abs{w_{\lambda,x}}^2+5\lambda^2\abs{w_\lambda}^2\geq0 ,
\end{align}
valid at almost every point of $J\times\R$ and for every $\lambda\geq K_\sharp$. Both sides involve $w_{\lambda,x}$, which is defined only up to a null set.

It is worth putting this at once into the integrated form in which it is used. The weight $\phi_y'$ is nonnegative, so multiplying \eqref{eq:auto-smoothing-a3} by it and integrating over $J\times\R$ preserves the inequality:
{
\begin{align}
 \label{eq:coercivity-integrated} \int_J\int_\R\phi_y' \bigl(\abs{w_{\lambda,x}}^2+5\lambda^2\abs{w_\lambda}^2\bigr)\dd x\dd t \leq\lambda\int_J\int_\R\phi_y'N_\lambda\dd x\dd t \qquad(\lambda\geq K_\sharp,\ y\in\R).
\end{align}
}
Both sides are nonnegative, and the right side is exactly the left side of the identity \eqref{eq:key-weighted-identity}. The two displays together therefore reduce the proposition to a single task. Whatever upper bound can be proved for the four terms on the right of \eqref{eq:key-weighted-identity} is transmitted by \eqref{eq:key-weighted-identity} to $\lambda\iint\phi_y'N_\lambda$, and then by \eqref{eq:coercivity-integrated} to a weighted $L^2$ norm of $w_{\lambda,x}$ and of $\lambda w_\lambda$ simultaneously. Producing that upper bound, uniformly in $\lambda$ and in $y$, is the whole content of Step~3.

\proofstep{Step 3}
This step supplies the upper bound to which Step~2 reduced the proposition: it bounds each of the four terms on the right of \eqref{eq:key-weighted-identity} uniformly in \emph{both} $\lambda$ and $y$. Uniformity in $\lambda$ is what allows the limit $\lambda\to\infty$ to be taken in Step~4, and uniformity in $y$ gives the uniformly local bound.

Throughout this step $H^1$ carries the Fourier weight $4+\xi^2$, so that $\abs{\pair*{F}{\psi}}\leq\norm{F}_{H^{-1}}\norm{\psi}_{H^1}$ holds with constant one for the norm \eqref{eq:Hminus-norm}. Uniformity in $y$ is the same in all four terms and may be disposed of once. Every weight that occurs is a translate of a fixed profile, $\phi_y^{(k)}=\phi_0^{(k)}(\cdot-y)$ for $k\geq1$, and $\norm{\phi_y}_\infty=\norm{\phi_0}_\infty=1$. Consequently every norm of a weight appearing below is independent of $y$. Uniformity in $\lambda$ has a different source in each term and is discussed as it arises. Write $\Delta q=q(b)-q(a)$ and $\Delta S_{\kappa\lambda}=S_{\kappa\lambda}(b)-S_{\kappa\lambda}(a)$, and take the four terms in the order in which they appear in \eqref{eq:key-weighted-identity}.

\emph{The flux at the fixed parameter.} Since $\phi_y'\in H^1$ and $\J_\kappa(q(t))\in H^{-1}$ by part~\textup{(i)} of Lemma~\ref{lem:microscopic},
\begin{align}
 \label{eq:term-flux-kappa} \Babs{\kappa\int_a^b\pair*{\J_\kappa(q(t'))}{\phi_y'}\dd t'} \leq\kappa\abs{J}\sup_{t\in J}\bnorm{\J_\kappa(q(t))}_{H^{-1}} \norm{\phi_0'}_{H^1}.
\end{align}
The supremum is finite because the same part of that lemma makes $t\mapsto\J_\kappa(q(t))$ continuous into $H^{-1}$, so its image over the compact interval $J$ is compact. The bound is independent of $\lambda$ because $\kappa$ is fixed.

\emph{The density endpoint.} Here the pairing is between $L^1$ and $L^\infty$, so
\begin{align}
 \label{eq:term-density} \abs[\big]{\pair*{\Delta S_{\kappa\lambda}}{\phi_y}} \leq\norm{\Delta S_{\kappa\lambda}}_1\norm{\phi_0}_\infty =\norm{\Delta S_{\kappa\lambda}}_1 ,
\end{align}
and it remains to bound the $L^1$ norm uniformly in $\lambda$. By \eqref{eq:rho-nonneg} each density is nonnegative, so its $L^1$ norm is its mass, $\norm{\rho_\mu}_1=\alpha_\mu$. Applying the triangle inequality to $S_{\kappa\lambda}=\kappa\rho_\kappa-\lambda\rho_\lambda$ at the two endpoints gives
\begin{align}
 \label{eq:S-L1-bound} \norm{\Delta S_{\kappa\lambda}}_1 \leq2\kappa\sup_I\alpha_\kappa+2\lambda\sup_I\alpha_\lambda ,
\end{align}
whose right side is bounded uniformly for $\lambda\geq K_\sharp$: the first term is fixed with $\kappa$, and $\lambda\alpha_\lambda$ is bounded on $q(I)\times[K_\sharp,\infty)$ by the closing bound \eqref{eq:high-energy-bounds} of Lemma~\ref{lem:high-energy}. For fixed $\kappa$, boundedness follows from the continuity of $q\mapsto\rho_\kappa(q)$ into $L^1$ on the compact trajectory. Conservation will give the vanishing of the density endpoint term in Step~5.

\emph{The multiplier endpoint.} Both factors are square integrable, since $B_{\kappa\lambda}$ maps $H^{-1}$ into $H^1$ and $\phi_y'\in L^2$, so by Cauchy--Schwarz and $\norm{\partial_xf}_2\leq\norm{f}_{H^1}$,
\begin{align}
 \label{eq:term-multiplier} \abs[\big]{\pair*{\partial_xB_{\kappa\lambda}\Delta q}{\phi_y'}} \leq\norm{B_{\kappa\lambda}\Delta q}_{H^1}\norm{\phi_0'}_2 \leq\frac12\norm{\Delta q}_{H^{-1}}\norm{\phi_0'}_2 .
\end{align}
The last constant is absolute, and in particular independent of $\lambda$. Indeed, with the norm \eqref{eq:Hminus-norm} on $H^{-1}$ and the equivalent Fourier weight $4+\xi^2$ on $H^1$, a Fourier multiplier $\Theta$ obeys $\norm{\Theta}_{H^{-1}\to H^1}\leq\sup_\xi(4+\xi^2)\abs{\Theta(\xi)}$, while the multiplier of $B_{\kappa\lambda}$ computed in the proof of Lemma~\ref{lem:two-parameter} satisfies
\begin{align}
 \label{eq:B-symbol-uniform} (4+\xi^2)\,\frac{2\abs[\big]{\lambda^2-\kappa^2}} {(\xi^2+4\kappa^2)(\xi^2+4\lambda^2)} \leq\frac{4+\xi^2}{2\bigl(\xi^2+4\min\{\kappa,\lambda\}^2\bigr)} \leq\frac12 \qquad(\kappa,\lambda\geq1).
\end{align}
For the first inequality, write $M=\max\{\kappa,\lambda\}$ and $m=\min\{\kappa,\lambda\}$. Then $\abs{\lambda^2-\kappa^2}\leq M^2$, while the factor of the denominator built from $M$ is at least $4M^2$ and the factor built from $m$ is at least $\xi^2+4m^2$. The second inequality is $m\geq1$. The estimate applies to both orderings of $\kappa$ and $\lambda$, as required when $\kappa$ is fixed and $\lambda$ ranges over $[K_\sharp,\infty)$.

\emph{The term containing $w_\lambda$.} Again an $H^{-1}$--$H^1$ pairing, now with the third derivative of the weight:
\begin{align}
 \label{eq:term-coordinate} \Babs{\lambda\int_a^b\pair*{w_\lambda(t')}{\phi_y'''}\dd t'} \leq\abs{J}\sup_{t\in J}\bnorm{\lambda w_\lambda(t)}_{H^{-1}} \norm{\phi_0'''}_{H^1}.
\end{align}
The factor $\lambda$ has been combined with $w_\lambda$, and $\norm{\lambda w_\lambda}_{H^{-1}}$ is bounded on $q(I)\times[K_\sharp,\infty)$ by \eqref{eq:high-energy-bounds}. This is the second of the two places where that uniform bound is consumed.

Every right side in \eqref{eq:term-flux-kappa}, \eqref{eq:term-density}, \eqref{eq:term-multiplier} and \eqref{eq:term-coordinate} is finite, independent of $y$, and bounded uniformly in $\lambda$. Denote by $C_J$ a constant dominating their sum (it depends on $J$, on the fixed parameter $\kappa$ and on the trajectory $q(I)$, but neither on $\lambda\geq K_\sharp$ nor on $y$), so that, with $K_\sharp$ the threshold \eqref{eq:sharp-threshold} fixed in Step~2, \eqref{eq:key-weighted-identity} becomes
\begin{align}
 \label{eq:coercive-upper} \lambda\int_J\int_\R\phi_y'N_\lambda\dd x\dd t\leq C_J \qquad(\lambda\geq K_\sharp,\ y\in\R).
\end{align}
This is the upper bound promised in Step~2. Chaining it with the coercive bound \eqref{eq:coercivity-integrated}, and discarding the factor $5$, gives
{
\begin{align}
 \label{eq:w-weighted-bound} \sup_{\lambda\geq K_\sharp}\sup_{y\in\R} \int_J\int_\R\phi_y' \bigl(\abs{w_{\lambda,x}}^2+\lambda^2\abs{w_\lambda}^2\bigr) \dd x\dd t\leq C_J.
\end{align}
}

\proofstep{Step 4}
This step does two things, and they are worth separating. Part~\textup{(a)} constructs the representative of $q$ whose existence the proposition asserts. Part~\textup{(b)} transfers the bound \eqref{eq:w-weighted-bound} from $2\lambda w_\lambda$ to that representative, uniformly in the translation centre, and so proves \eqref{eq:q-Kato-bound}. Both parts combine the same two ingredients, one qualitative and one quantitative, and both are recorded first.

\emph{The approximation.} Put $\widetilde q_\lambda:=2\lambda w_\lambda$. The limit \eqref{eq:w-high} of Lemma~\ref{lem:high-energy} asserts that $2\lambda w_\lambda(p)\to p$ in $H^{-1}$ as $\lambda\to\infty$, uniformly for $p$ in a compact subset of $H^{-1}(\R;\R)$. Applied to the trajectory $q(I)$, which is compact and contains $q(t)$ for every $t\in I$, it gives a limit that is uniform in time,
\begin{align}
 \label{eq:u-to-q-Hminus} \sup_{t\in J}\,\bnorm{\widetilde q_\lambda(t)-q(t)}_{H^{-1}}\rightarrow0 , \qquad\text{that is,}\qquad \widetilde q_\lambda\rightarrow q \quad\text{in }C(J;H^{-1}).
\end{align}
The convergence \eqref{eq:u-to-q-Hminus} identifies $q$ as the weak limit of the approximants bounded by \eqref{eq:w-weighted-bound}.

% {
% \begin{remark}
% The convergence \eqref{eq:w-high} identifies the weak limit in this step. The recovery identity \eqref{eq:Drach-density-conclusion} enters through the coercive flux decomposition of Lemma~\ref{lem:flux-decomposition}.
% \end{remark}
% }

\emph{An identification principle.} Let $\omega\in C^\infty(\R)$ be nonnegative with $\omega$ and $\omega'$ bounded, so that $\omega\varphi\in C_c^\infty((a,b)\times\R)$ whenever $\varphi$ is, and $\norm{\omega\varphi(t)}_{H^1}\leq C_\omega\norm{\varphi(t)}_{H^1}$. Write $\pair*{F}{\varphi}_{t,x}:=\int_J\pair*{F(t)}{\varphi(t,\cdot)}\dd t$ for $F\in C(J;H^{-1})$ and $\varphi\in C_c^\infty((a,b)\times\R)$; this is the pairing of $F$, regarded as a distribution on $(a,b)\times\R$, with $\varphi$, and it equals $\iint F\varphi\dd x\dd t$ when $F\in L^2_{\loc}(J\times\R)$. For such a weight, \eqref{eq:u-to-q-Hminus} gives
\begin{align}
 \label{eq:identification-principle} \abs[\big]{\pair*{\omega\widetilde q_\lambda-\omega q}{\varphi}_{t,x}} =\abs[\big]{\pair*{\widetilde q_\lambda-q}{\omega\varphi}_{t,x}} \leq C_\omega\int_J\norm{\widetilde q_\lambda(t)-q(t)}_{H^{-1}} \norm{\varphi(t)}_{H^1}\dd t \rightarrow0 .
\end{align}
Hence $\omega\widetilde q_\lambda\to\omega q$ in $\mathcal D'((a,b)\times\R)$. Since a weak $L^2$ limit is in particular a distributional one, and a distribution has at most one value on each test function, it follows that
\begin{align}
 \label{eq:identification-conclusion} \omega\widetilde q_{\lambda_n}\rightharpoonup F \ \text{ weakly in }L^2_{\loc} \qquad\Longrightarrow\qquad F=\omega q \ \text{ almost everywhere.}
\end{align}
Part~\textup{(a)} uses this with $\omega\equiv1$ and part~\textup{(b)} with $\omega=\sqrt{\phi_y'}=\sech(\cdot-y)$, both of which satisfy the hypothesis.

\emph{Part \textup{(a)}: construction of the representative.} Choose any sequence $\lambda_n\to\infty$. On every cylinder
\begin{align}
 \label{eq:auto-smoothing-a4} \Omega_N=J\times(-N,N),\qquad N\in\mathbb N,
\end{align}
the family $\{\widetilde q_{\lambda_n}\}$ is bounded in $L^2(\Omega_N)$. Indeed $(-N,N)$ is covered by the $2N+1$ cells centered at the integers $j$ with $\abs{j}\leq N$, so the cell comparison \eqref{eq:weight-cell-comparison} and the bound \eqref{eq:w-weighted-bound}, the latter rewritten through $\abs{\widetilde q_\lambda}^2=4\lambda^2\abs{w_\lambda}^2$, give
\begin{align}
 \label{eq:qtilde-cell-bound} \int_{\Omega_N}\abs{\widetilde q_{\lambda}}^2\dd x\dd t \leq\sum_{\abs{j}\leq N}\int_J\int_{j-1}^{j+1}\abs{\widetilde q_{\lambda}}^2\dd x\dd t \leq(2N+1)\cosh^2(1)\,4C_J
\end{align}
for every $\lambda\geq K_\sharp$, with a bound depending on $N$ but not on $\lambda$.

Extract a weakly convergent subsequence on $\Omega_1$, then a further one on $\Omega_2$, and so on, and take the diagonal subsequence. Two limits obtained on nested cylinders agree on the smaller one, because weak convergence in $L^2(\Omega_{N+1})$ implies weak convergence in $L^2(\Omega_N)$. This compatibility produces a single $\widetilde q\in L^2_{\loc}(J\times\R)$ with $\widetilde q_{\lambda_n}\rightharpoonup\widetilde q$ in $L^2(\Omega_N)$ for every $N$. Applying \eqref{eq:identification-conclusion} with $\omega\equiv1$ identifies it:
\begin{align}
 \label{eq:qtilde-equals-q} \widetilde q=q \qquad\text{almost everywhere on }J\times\R .
\end{align}
This is the step that converts bounds for $2\lambda w_\lambda$ into bounds for $q$. By \eqref{eq:qtilde-equals-q}, together with the countability argument of Remark~\ref{rem:measurability}, which produces one exceptional set of times for all test functions and all cells, the function $\widetilde q$ is a representative of $q$ in the sense of \eqref{eq:representative-def}. Any two such representatives agree almost everywhere, so the representative is independent of the integer cover used in the extraction. The existence assertion of the proposition is thereby proved, and $q$ denotes this function from here on.

\emph{Part \textup{(b)}: the uniform cell bound.} Fix $y\in\R$. Rewriting \eqref{eq:w-weighted-bound} in terms of $\widetilde q_\lambda$ gives
\begin{align}
 \label{eq:qtilde-weighted-bound} \int_J\int_\R\phi_y'\abs{\widetilde q_\lambda}^2\dd x\dd t =4\lambda^2\int_J\int_\R\phi_y'\abs{w_\lambda}^2\dd x\dd t \leq4C_J ,
\end{align}
with $C_J$ independent of $\lambda$ and of $y$. Hence $\{\sqrt{\phi_y'}\,\widetilde q_\lambda\}_{\lambda\geq K_\sharp}$ is bounded in $L^2(J\times\R)$, and along a sequence $\lambda\to\infty$ realizing the lower limit in \eqref{eq:qtilde-weighted-bound} we may extract a weak limit. By \eqref{eq:identification-conclusion}, now with $\omega=\sech(\cdot-y)$, that limit is $\sqrt{\phi_y'}\,q$:
\begin{align}
 \label{eq:weighted-weak-limit} \sqrt{\phi_y'}\,\widetilde q_\lambda \rightharpoonup\sqrt{\phi_y'}\,q \qquad\text{weakly in }L^2(J\times\R).
\end{align}
Weak lower semicontinuity of the $L^2$ norm, together with \eqref{eq:weighted-weak-limit} and \eqref{eq:qtilde-weighted-bound}, gives
\begin{align}
 \label{eq:q-weighted-bound} \int_J\int_\R\phi_y'\abs{q}^2\dd x\dd t \leq\liminf_{\lambda\to\infty} \int_J\int_\R\phi_y'\abs{\widetilde q_\lambda}^2\dd x\dd t \leq4C_J ,
\end{align}
and the bound is uniform in $y$ because $C_J$ is. The cell comparison \eqref{eq:weight-cell-comparison}, whose constant is also independent of $y$, converts \eqref{eq:q-weighted-bound} into \eqref{eq:q-Kato-bound} and completes the first half of the proposition.

\proofstep{Step 5}
Step~3 bounded the right side of \eqref{eq:key-weighted-identity} uniformly in $y$. This step shows that it in fact \emph{tends to zero} as the translation centre leaves every bounded set, and hence that the weighted cell integral of $q$ does too.

\emph{Reduction to the four terms.} Write
\begin{align}
 \label{eq:E-lambda-def} E_\lambda(y):={}& \Babs{\kappa\int_a^b\pair*{\J_\kappa(q(t'))}{\phi_y'}\dd t'} +\abs[\big]{\bpair{S_{\kappa\lambda}(b)-S_{\kappa\lambda}(a)}{\phi_y}} \\
 &+\abs[\big]{\pair*{\partial_xB_{\kappa\lambda}(q(b)-q(a))}{\phi_y'}} +\Babs{\lambda\int_a^b\pair*{w_\lambda(t')}{\phi_y'''}\dd t'}\notag
\end{align}
for the sum of the absolute values of the four terms on the right of \eqref{eq:key-weighted-identity}. Discarding the gradient term in \eqref{eq:N-positive} and using $\abs{\widetilde q_\lambda}^2=4\lambda^2\abs{w_\lambda}^2$ gives the almost-everywhere bound $\lambda N_\lambda\geq\frac54\abs{\widetilde q_\lambda}^2$, so that, after integrating against $\phi_y'\geq0$ and comparing with \eqref{eq:key-weighted-identity},
\begin{align}
 \label{eq:auto-smoothing-a7} \frac54\int_J\int_\R\phi_y'\abs{\widetilde q_\lambda}^2\dd x\dd t \leq\lambda\int_J\int_\R\phi_y'N_\lambda\dd x\dd t \leq E_\lambda(y) .
\end{align}
Part~\textup{(b)} of Step~4, applied at this fixed $y$ along a sequence realizing the lower limit, therefore gives
\begin{align}
 \label{eq:Xi-E} \int_J\int_\R\phi_y'\abs{q}^2\dd x\dd t\leq\frac45\liminf_{\lambda\to\infty}E_\lambda(y),
\end{align}
and it suffices to show that the right side tends to zero as $y\to\pm\infty$.

Each of the four terms of \eqref{eq:E-lambda-def} is treated on its own, and always in the same order: first $\lambda\to\infty$ at a \emph{fixed} $y$, which produces a limit no longer involving $\lambda$. Then we take $y\to\pm\infty$ in that limit. Two devices recur. In the two time-integrated terms the Bochner integral may be taken outside the pairing,
\begin{align}
 \label{eq:bochner-outside} \int_a^b\pair*{G(t')}{\psi}\dd t'=\pair[\Big]{\int_a^bG(t')\dd t'}{\psi},
\end{align}
so that the outer limit is applied to a single fixed element of $H^{-1}(\R)$. And Lemma~\ref{lem:translation} disposes of any pairing against $\phi_y'$, $\phi_y''$ or $\phi_y'''$, these being translates of fixed $H^1$ functions by \eqref{eq:tanh-weight}. The term with the undifferentiated weight $\phi_y\notin H^1(\R)$ is treated using conservation of the density mass. The four terms are taken in the order of \eqref{eq:E-lambda-def}.

\emph{The flux at the fixed parameter.} This term is independent of $\lambda$. By \eqref{eq:bochner-outside} it may be written as
\begin{align}
 \label{eq:tail-flux-kappa} \kappa\int_a^b\pair*{\J_\kappa(q(t'))}{\phi_y'}\dd t' =\kappa\pair[\Big]{\int_a^b\J_\kappa(q(t'))\dd t'}{\phi_y'},
\end{align}
a pairing of one fixed element of $H^{-1}(\R)$ with a translate of $\phi_0'$. Lemma~\ref{lem:translation} sends it to zero as $y\to\pm\infty$.

\emph{The density endpoint.} The two parameters behave differently here, and the $\lambda$ part disappears in the inner limit. By the nonnegativity \eqref{eq:rho-nonneg} and the limit \eqref{eq:alpha-high},
\begin{align}
 \label{eq:rho-lambda-endpoint} \bnorm{\lambda\big[\rho_\lambda(q(b))-\rho_\lambda(q(a))\big]}_1 \leq\lambda\alpha_\lambda(q(b))+\lambda\alpha_\lambda(q(a)) \rightarrow0,
\end{align}
and $\norm{\phi_y}_\infty=1$, so this part of the pairing tends to zero, and does so uniformly in $y$. Only the fixed-$\kappa$ part survives:
\begin{align}
 \label{eq:density-endpoint-limit} \lim_{\lambda\to\infty} \bpair{S_{\kappa\lambda}(b)-S_{\kappa\lambda}(a)}{\phi_y} =\pair*{f_\kappa}{\phi_y}, \qquad f_\kappa:=\kappa\bigl[\rho_\kappa(q(b))-\rho_\kappa(q(a))\bigr]\in L^1(\R).
\end{align}
Since $\phi_y(x)\to\mp1$ as $y\to\pm\infty$, the limit of this pairing is determined by the mass of $f_\kappa$. Proposition~\ref{prop:alpha-conservation} gives $\int_\R f_\kappa=0$, so that, as $y\to+\infty$,
\begin{align}
 \label{eq:auto-smoothing-a10} \pair*{f_\kappa}{\phi_y} =\int_\R f_\kappa(x)\bigl(\phi_y(x)+1\bigr)\dd x\rightarrow0 ;
\end{align}
dominated convergence applies because $\abs{\phi_y+1}\leq2$ and $\phi_y+1\to0$ pointwise. For $y\to-\infty$ use $\phi_y-1$ in place of $\phi_y+1$. The equality in \eqref{eq:auto-smoothing-a10} uses $\int_\R f_\kappa=0$. This is the only use made of Proposition~\ref{prop:alpha-conservation} in this proof.

\emph{The multiplier endpoint.} Since $\kappa$ is fixed and $B_{\kappa\lambda}=\frac12(\Rf_\kappa-\Rf_\lambda)$ by Lemma~\ref{lem:two-parameter}, the inner limit is settled by showing
\begin{align}
 \label{eq:Rlambda-strong} \Rf_\lambda\rightarrow0 \quad\text{strongly from }H^{-1}\text{ to }H^1,
\end{align}
whence $B_{\kappa\lambda}\to\frac12\Rf_\kappa$ strongly. To prove \eqref{eq:Rlambda-strong}, equip $H^1$ with the equivalent Fourier weight $4+\xi^2$ and $H^{-1}$ with the norm \eqref{eq:Hminus-norm}. Then, for $f\in H^{-1}$,
\begin{align}
 \label{eq:auto-smoothing-a8} \norm{\Rf_\lambda f}_{H^1}^2 =\int_\R\frac{4+\xi^2}{(\xi^2+4\lambda^2)^2}\abs{\widehat f(\xi)}^2\dd \xi =\int_\R \Theta_\lambda(\xi)\, \frac{\abs{\widehat f(\xi)}^2}{4+\xi^2}\dd \xi , \qquad \Theta_\lambda(\xi):=\frac{(4+\xi^2)^2}{(\xi^2+4\lambda^2)^2} .
\end{align}
For $\lambda\geq1$ we have $4+\xi^2\leq\xi^2+4\lambda^2$, so $0\leq\Theta_\lambda\leq1$ for every $\xi$, while $\Theta_\lambda(\xi)\to0$ for each fixed $\xi$ (the supremum $\sup_\xi\Theta_\lambda=1$ also shows that the operator norm remains one). Dominated convergence gives \eqref{eq:Rlambda-strong}, and the uniformity principle of Lemma~\ref{lem:high-energy} upgrades it to convergence uniform on compact subsets of $H^{-1}$, in particular on the endpoint data $\{q(a),q(b)\}$. The inner limit of this term is therefore $\frac12\babs{\pair*{\partial_x\Rf_\kappa\Delta q}{\phi_y'}}$, and one integration by parts puts it in the reach of Lemma~\ref{lem:translation}:
\begin{align}
 \label{eq:auto-smoothing-a9} \pair*{\partial_x\Rf_\kappa\Delta q}{\phi_y'} =-\pair*{\Rf_\kappa\Delta q}{\phi_y''} \rightarrow0 \qquad(y\to\pm\infty),
\end{align}
since $\Rf_\kappa\Delta q\in H^1(\R)$ is a fixed element and $\phi_y''$ is a translate of $\phi_0''$.

\emph{The term containing $w_\lambda$.} By \eqref{eq:u-to-q-Hminus}, $\widetilde q_\lambda=2\lambda w_\lambda\to q$ in $C(J;H^{-1})$, so $\lambda w_\lambda\to q/2$ in that space, and pairing against the fixed $H^1$ function $\phi_y'''$ converges accordingly. With \eqref{eq:bochner-outside}, the inner limit is
\begin{align}
 \label{eq:coordinate-inner-limit} \lim_{\lambda\to\infty} \lambda\int_a^b\pair*{w_\lambda(t')}{\phi_y'''}\dd t' =\frac12\int_a^b\pair*{q(t')}{\phi_y'''}\dd t' =\frac12\pair[\Big]{\int_a^bq(t')\dd t'}{\phi_y'''} .
\end{align}
Once again one fixed element of $H^{-1}(\R)$ is paired with a translate of a fixed $H^1$ function, so Lemma~\ref{lem:translation} sends it to zero as $y\to\pm\infty$.

\emph{Conclusion.} Collecting the four inner limits gives, for each fixed $y$,
\begin{align}
 \label{eq:E-infinity} \lim_{\lambda\to\infty}E_\lambda(y)=E_\infty(y):={}& \Babs{\kappa\pair[\Big]{\int_a^b\J_\kappa(q(t'))\dd t'}{\phi_y'}} +\abs[\big]{\pair*{f_\kappa}{\phi_y}} \\
 &+\frac12\abs[\big]{\pair*{\partial_x\Rf_\kappa\Delta q}{\phi_y'}} +\frac12\Babs{\pair[\Big]{\int_a^bq(t')\dd t'}{\phi_y'''}} ,\notag
\end{align}
and each of the four terms was shown above to vanish as $y\to\pm\infty$. Hence $E_\infty(y)\to0$, and \eqref{eq:Xi-E} gives
\begin{align}
 \label{eq:q-sech-tail} \int_J\int_\R\sech^2(x-y)\abs{q(t,x)}^2\dd x\dd t \rightarrow0\quad(y\to\pm\infty).
\end{align}
Applying \eqref{eq:weight-cell-comparison} once more, again with a $y$-independent constant, proves \eqref{eq:q-Kato-tail}.
\end{proof}

{\begin{remark}[The exceptional set of times]\label{rem:measurability}
Step~4 identifies $\widetilde q$ with $q$ as a spacetime distribution. To obtain the representative in \eqref{eq:representative-def}, we need a full-measure set of times on which $\widetilde q(t,\cdot)$ realizes $q(t)$ against every $\phi\in C_c^\infty(\R)$. Fubini initially gives such a set for each fixed $\phi$. Since $C_c^\infty(\R)$ is uncountable, a countable dense family is needed to obtain a common set. The same construction will cover all translation centres in \eqref{eq:q-Kato-bound}, \eqref{eq:q-Kato-tail} and the estimates of Section~\ref{sec:Gardner}.

We obtain a common set by countability. For each $N\in\mathbb N$ choose a countable subset of $C_c^\infty((-N,N))$ dense in $H^1_0((-N,N))$, and intersect, over that family and over $N$, the full-measure sets of times on which the local $L^2$ slice exists and on which its pairing with the test function agrees with the $H^{-1}$ pairing. This is a countable intersection of full-measure sets, hence itself of full measure. On it, density upgrades the identity from the chosen tests to every compactly supported spatial test, both sides being continuous in $\phi\in H^1_0((-N,N))$: the left by $q(t)\in H^{-1}$, the right by $\widetilde q(t,\cdot)\in L^2((-N,N))$. Fubini applied to the countable family of integer-centred cells gives the same common set for all the later cell estimates. Every other bounded cell is covered by finitely many of these. In particular the exceptional set of times is independent both of the cell and of the translation centre $y$.

At every time the curve has its $C_tH^{-1}$ trace from Definition~\ref{def:green}. The $L^2_{\loc}$ slices are defined on the full-measure set constructed above. Every subsequent use of the representative is made at almost every $t$ and always on this one common set: the cutoff estimate in Lemma~\ref{lem:internal-H2}, the Kato conditions for $r$ in Proposition~\ref{prop:r-X1-local}, and the identification of the distributional multiplier product with the pointwise one in Proposition~\ref{prop:green-Gardner}. Where a centre $y$ occurs, in the first two, every centre is treated separately and a supremum or a limit in $y$ is taken only afterwards, of the resulting time integrals. The common exceptional set lets every slicewise identity below hold at one time for all cells at once.
\end{remark}}

{\begin{remark}
\label{rem:full-spectral-family}
The base parameter $\kappa$ retains the exact conservation information, while the second parameter $\lambda\to\infty$ produces the coercive approximation $2\lambda w_\lambda\to q$. The argument uses the defining Green identity throughout this range of spectral parameters.
\end{remark}}

% \newpage

%%%%%%%%%%%%%%%%%%%%%%%%%%%%%%%%%%%%%%%%%%%%%%%%%%%%%%%%%%%%%%%%%%%%%%%%%%%%%%%%%%
\section{From Green Solutions to Gardner Solutions}
\label{sec:Gardner}
%%%%%%%%%%%%%%%%%%%%%%%%%%%%%%%%%%%%%%%%%%%%%%%%%%%%%%%%%%%%%%%%%%%%%%%%%%%%%%%%%%

This section constructs Gardner solutions from Green solutions and applies Gardner uniqueness.

Klaus--Koch--Liu \cite{Klaus.Koch.Liu2023} prove that, for a fixed parameter $\tau$, the first $\tau$-Gardner equation \eqref{eq:Gardner} possesses a weak-solution theory in a Kato class $X_1$, and within that class a solution is determined uniquely by its $L^2$ datum. Section~\ref{sec:Gardner-class} records the solution class (Definition~\ref{def:Gardner-solution}) and the uniqueness theorem itself (Lemma~\ref{lem:KKL-Gardner}).

By the end of Section~\ref{sec:automatic-smoothing}, an arbitrary Green solution $q$ satisfies the following three properties.
\begin{enumerate}[
label=\textup{(\arabic*)}, leftmargin=*, itemsep=.4em, before={\let\fullwidthdisplay\relax} ]
\item $q\in C(I;H^{-1}(\R;\R))$ and the integrated identity \eqref{eq:green-identity-recalled} holds at every $\lambda\geq\lambda_*(q,I)$. This is Definition~\ref{def:green}.
\item $\alpha_\lambda(q(t))$ is independent of $t\in I$ for every $\lambda$ above the threshold \eqref{eq:single-K}, by Proposition~\ref{prop:alpha-conservation}.
\item For every compact $J\subset I$, the curve $q$ has a representative in $L^2_{\loc}(J\times\R)$ obeying the uniform cell bound \eqref{eq:q-Kato-bound} and the translated-cell tail \eqref{eq:q-Kato-tail}, by Proposition~\ref{prop:auto-smoothing}.
\end{enumerate}
We use \textup{(1)}--\textup{(3)} and, outside $I$, the canonical flow of \cite{Killip.Visan2019} to construct a global Gardner solution $r$ of class $X_1$. The construction and verification proceed in four stages.

\textit{The variable \textup{(Section~\ref{sec:Gardner-lift})}.} The Gardner theory governs an $L^2$-valued function $r$, while $q$ takes values in $H^{-1}$. They are related by the shifted Miura map \eqref{eq:Miura-map}, which by Remark~\ref{rem:shifted-Miura-injective} is injective on $L^2(\R;\R)$ for $\tau>0$. We construct $r$ by an explicit formula in $q$ and the fixed parameter $\tau$. This formula must give $r\in C(I;L^2)$ and satisfy $M_\tau(r)=q$. Its dependence on $q$ ensures that equal initial data give equal initial values of $r$. The recovery identity \eqref{eq:q-from-h} of Section~\ref{sec:solution-notion} supplies such a formula, and property~\textup{(1)} gives its basic regularity.

\textit{The regularity \textup{(Section~\ref{sec:Gardner-smoothing})}.} Membership in $X_1$ requires the smoothing and tail conditions \eqref{eq:X1-bounds} and \eqref{eq:X1-tail} on $r_x$, whereas $r$ is a priori only $C(I;L^2)$. Since $r_x=\tau v_x-\frac12\mu_x$ with $\mu=v_x/(1+v)$, what is needed is one spatial derivative of $v$ more than property~\textup{(1)} provides, uniformly over unit cells and integrated in time, together with its translated-cell tail. The Miura identity, written as $r_x+r^2=q-2\tau r$ and squared against a cutoff, bounds $r_x$ and $r^2$ on each cell: property~\textup{(3)} controls the term in $q$, the cross term $2r_xr^2$ is a derivative and integrates by parts, and the tails of $q$ and $r$ give the tail of $r_x$. The identity $v_{xx}=(1+v)(\mu_x+\mu^2)$, with $\mu=2(\tau v-r)$, then returns the derivative to $v$. This step uses the local smoothing and tail estimates of Section~\ref{sec:automatic-smoothing}.

\textit{The equation \textup{(Sections~\ref{sec:Gardner-bootstrap} and \ref{sec:Gardner-dynamics})}.} To prove that $r$ solves \eqref{eq:Gardner}, we divide property~\textup{(1)} by $\tau$ to obtain a conservation law for $v$. We then convert its flux $\Psi$ into the Gardner flux through $\tau\Psi-\Psi_x/(2(1+v))$, as in \eqref{eq:rough-flux-identity}. That conversion rests on an algebraic identity between \emph{smooth} functions, and transporting it to $v$ at the available regularity by mollification costs one spatial derivative more than \textup{(1)} provides. The derivative estimate in Section~\ref{sec:Gardner-bootstrap} therefore also prepares this flux identification, which Section~\ref{sec:Gardner-dynamics} uses to prove the weak Gardner equation.

\textit{The time interval \textup{(Section~\ref{sec:uniqueness})}.} The class $X_1$ is global in time, while $I$ is compact. The remedy is to extend $q$ beyond $I$ by the canonical flow before the lift is taken. The extension is a Green solution on every interval $[-L,L]$ with $L\geq T$, with one threshold for all of them, so that a single parameter $\tau$ serves all times, and property~\textup{(2)} at that parameter bounds the lift in $L^2$ uniformly in time. Gardner uniqueness then gives Theorem~\ref{thm:main}.

Lemma~\ref{lem:KKL-Gardner} applies after the global extension. Property~\textup{(2)}, the conservation of $\alpha_\lambda$, is used in Section~\ref{sec:automatic-smoothing} to prove property~\textup{(3)} and again at the fixed parameter $\tau$ to bound the global lift.

%%%%%%%%%%%%%%%%%%%%%%%%%%%%%%%%%%%%%%%%%%%%%%%%%%%%%%%%%%%%%%%%%%%%%%%%%%%%%%%%%%
\subsection{The Gardner Equation and Its Uniqueness Class}
\label{sec:Gardner-class}
%%%%%%%%%%%%%%%%%%%%%%%%%%%%%%%%%%%%%%%%%%%%%%%%%%%%%%%%%%%%%%%%%%%%%%%%%%%%%%%%%%

This subsection records, in the normalization used here, the solution class and the uniqueness theorem of Klaus--Koch--Liu \cite{Klaus.Koch.Liu2023}. It fixes the target of the transfer. The numbering cited refers to version~2, revised 5 May 2026.

Throughout, $\tau\geq2$ is a fixed parameter, chosen in Section~\ref{sec:Gardner-lift} to be admissible for the trajectory. It corresponds to the fixed Gardner parameter $\tau_0$ in \cite{Klaus.Koch.Liu2023}. The Gardner theory in \cite[Sections~2.5--2.9]{Klaus.Koch.Liu2023} assumes a positive fixed parameter. The Kato-smoothing input \cite[Lemma~6.2]{Klaus.Koch.Liu2023} is stated for parameter at least one. Hence $\tau\geq2$ satisfies every parameter restriction used in \cite[Theorem~2.24]{Klaus.Koch.Liu2023}.

For smooth, rapidly decaying $r$, define the functional
\begin{align}
 \label{eq:Gardner-Hamiltonian} H^{\mathrm{Gardner}}_1(r,\tau) =\frac12\int_\R\bigl(r_x^2+r^4+4\tau r^3\bigr)\dd x.
\end{align}
The first $\tau$-Gardner equation is
\begin{align}
 \label{eq:Gardner} r_t=\partial_x\frac{\delta H^{\mathrm{Gardner}}_1}{\delta r} =\partial_x\bigl(-r_{xx}+6\tau r^2+2r^3\bigr).
\end{align}
Here $\delta H^{\mathrm{Gardner}}_1/\delta r$ is the variational derivative. For $\phi\in C_c^\infty(\R)$, differentiating $H^{\mathrm{Gardner}}_1(r+\varepsilon\phi,\tau)$ at $\varepsilon=0$ and integrating by parts gives the integral of $(-r_{xx}+6\tau r^2+2r^3)\phi$.

Equation~\eqref{eq:Gardner} is the case $N=1$ of \cite[Section~2.3, equation~(2.24)]{Klaus.Koch.Liu2023}. The functional \eqref{eq:Gardner-Hamiltonian} need not be finite for the Gardner solutions considered below. Its expression on smooth functions identifies the normalization of the equation.

The shifted Miura map links the Gardner variable to the KdV potential:
\begin{align}
 \label{eq:Miura-map} M_\tau\colon L^2(\R;\R)\rightarrow H^{-1}(\R), \qquad M_\tau(r):=r_x+2\tau r+r^2.
\end{align}
It is well defined because $r\in L^2(\R)$ implies $r_x\in H^{-1}(\R)$ and $r^2\in L^1(\R)\subset H^{-1}(\R)$. Section~\ref{sec:Gardner-lift} constructs an $L^2$ preimage of $q(t)$ at each time.

{\begin{remark}[The shifted Miura map is injective]
\label{rem:shifted-Miura-injective}
For every $\tau>0$, the map $M_\tau$ of \eqref{eq:Miura-map} is injective on $L^2(\R;\R)$. This agrees with \cite[Proposition~4.8(2)]{Klaus.Koch.Liu2023} and also follows directly. If $M_\tau(r_1)=M_\tau(r_2)$ and $d:=r_1-r_2$, then
\begin{align}
 \label{eq:shifted-Miura-ode} d_x+\bigl(2\tau+r_1+r_2\bigr)d=0 .
\end{align}
The coefficient is locally integrable, so $d$ has a locally absolutely continuous representative and, for any $x_0$,
\begin{align}
 \label{eq:shifted-Miura-solution} d(x)=d(x_0)\exp\left(-2\tau(x-x_0)-\int_{x_0}^x\bigl(r_1+r_2\bigr)(z)\,dz\right).
\end{align}
By Cauchy--Schwarz the integral is $O\bigl(\abs{x-x_0}^{\frac12}\bigr)$, whereas the linear term grows like $2\tau\abs{x-x_0}$ as $x\to-\infty$. Hence $\abs{d(x)}\to\infty$ as $x\to-\infty$ unless $d(x_0)=0$, and a nonzero $d$ cannot lie in $L^2(\R)$. Therefore $d\equiv0$. The positive shift gives the exponential factor $\ue^{-2\tau(x-x_0)}$. At $\tau=0$ this factor equals one, and \eqref{eq:one-Miura-insufficient-a2} gives a counterexample to injectivity of the unshifted map.
\end{remark}}

For smooth solutions $M_\tau$ carries \eqref{eq:Gardner} to \eqref{eq:KdV}. Indeed, differentiating the map in time gives $q_t=(\partial_x+2\tau+2r)r_t$, while a direct computation yields the operator identity
\begin{align}
 \label{eq:Miura-intertwining} \partial_x\bigl(-q_{xx}+3q^2\bigr) =(\partial_x+2\tau+2r)\,  \partial_x\bigl(-r_{xx}+6\tau r^2+2r^3\bigr),
\end{align}
so that \eqref{eq:Gardner} forces \eqref{eq:KdV}. The identity also certifies that the normalizations of \eqref{eq:KdV} and \eqref{eq:Gardner} match.

The correspondence also gives the regularity used in the present argument. The Gardner variable $r$ lies in $L^2(\R)$, a full derivative smoother than $q\in H^{-1}(\R)$. The quadratic nonlinearity $r^2$ of \eqref{eq:Gardner} is locally integrable. For a general $q\in H^{-1}(\R)$, $q^2$ has no canonical meaning. The cubic term asks for slightly more than $r\in L^2$: for $\chi\in C_c^\infty(\R)$ equal to one near the origin, the function $r=\chi(x)\abs{x}^{-2/5}$ lies in $L^2(\R)$ while $r^3\notin L^1_{\loc}(\R)$. What supplies it is the derivative bound carried by Definition~\ref{def:Gardner-solution}, which upgrades $r$ to $L^6_{\loc}$ through \eqref{eq:Gardner-local-L6}. The local smoothing estimate of Section~\ref{sec:Gardner-smoothing} therefore makes both Gardner nonlinearities locally integrable, so that \eqref{eq:Gardner} admits a distributional weak formulation. For KdV at $H^{-1}$ regularity, Definition~\ref{def:green} gives the formulation explained in Remark~\ref{rem:why-green}. Sections~\ref{sec:Gardner-lift}--\ref{sec:Gardner-dynamics} accordingly run the correspondence backwards: from a Green solution $q$ they construct $r$ and verify that it is a Gardner solution in the sense of Definition~\ref{def:Gardner-solution} below, after which the uniqueness theorem quoted as Lemma~\ref{lem:KKL-Gardner} applies to it.

{
\begin{remark}[Relation to mKdV]\label{rem:Gardner-is-mKdV}
The following change of variables relates \eqref{eq:Gardner} to the defocusing modified Korteweg--de Vries equation and explains the form of the shifted Miura map. Substituting
\begin{align}
 \label{eq:Gardner-mKdV-substitution} r=U-\tau
\end{align}
into \eqref{eq:Gardner} removes the quadratic term and leaves a constant-coefficient transport,
\begin{align}
 \label{eq:Gardner-as-boosted-mKdV} U_t=-U_{xxx}+6U^2U_x-6\tau^2U_x ,
\end{align}
which the Galilean change of frame $\widetilde U(t,x)=U(t,x+6\tau^2t)$ removes in turn. Thus $\widetilde U$ solves \eqref{eq:defocusing-mkdv} exactly: \eqref{eq:Gardner} is \eqref{eq:defocusing-mkdv} on a constant background, viewed in a moving frame, and at $\tau=0$ the two equations coincide. Under the same substitution $M_\tau$ becomes the unshifted Miura map,
\begin{align}
 \label{eq:shifted-Miura-is-ordinary} r_x+2\tau r+r^2=U^2+U_x-\tau^2 .
\end{align}
The right side uses the same plus-sign Miura expression as \eqref{eq:intro-miura-partners}, evaluated at $U$, and subtracts $\tau^2$. This is the origin of the word ``shifted'': the parameter $\tau$ is a shift of the modified variable, and $-\tau^2$ the corresponding shift of its Miura transform.

A Gardner solution has $r(t)\in L^2(\R)$, hence $U(t)\in\tau+L^2(\R)$: $U(t)$ is an $L^2$ perturbation of the constant background $\tau$. This condition allows disjoint smooth bumps of height one and summable widths, so it does not imply a pointwise limit at $\pm\infty$. Lemma~\ref{lem:KKL-Gardner} gives uniqueness for \eqref{eq:defocusing-mkdv} on this nonzero constant background under the transformed Kato conditions of $X_1$. Section~\ref{sec:mkdv} treats $u(t)\in L^2(\R)$ on the zero background and quotes the corresponding solution theory. It uses both Miura transforms $q_\pm$. The present section uses the single sign \eqref{eq:shifted-Miura-is-ordinary} on the background $\tau$.
\end{remark}
}

We use the parameter-dependent norm
\begin{align}
 \label{eq:H1-tau-norm} \norm{f}_{H^1_\tau}^2:=\norm{f_x}_2^2+\tau^2\norm{f}_2^2 .
\end{align}
By Plancherel,
\begin{align}
 \label{eq:Gardner-class-a1} \norm{f}_{H^1_\tau} =\norm[\big]{(\tau^2+\xi^2)^{1/2}\widehat f(\xi)}_{L^2_\xi},
\end{align}
so \eqref{eq:H1-tau-norm} equals the $s=1$ norm in \cite[Definition~5.1]{Klaus.Koch.Liu2023}.

{\begin{definition}
\label{def:Gardner-solution}
Let $\mathcal I\subset\R$ be open. A real-valued function $r$ on $\mathcal I\times\R$ is a \emph{Gardner solution} on $\mathcal I$ if
\begin{align}
 \label{eq:weak-Gardner-regularity} r\in L^\infty(\mathcal I;L^2(\R)), \qquad r_x\in L^2_{\loc}(\mathcal I\times\R),
\end{align}
and the first $\tau$-Gardner equation \eqref{eq:Gardner} holds in the sense of distributions on $\mathcal I\times\R$, that is,
\begin{align}
 \label{eq:weak-Gardner-formulation} \int_{\mathcal I}\int_\R \left[r\varphi_t+r\varphi_{xxx} -(6\tau r^2+2r^3)\varphi_x\right]\dd x\dd t=0 \quad\text{for every }\varphi\in C_c^\infty(\mathcal I\times\R).
\end{align}
A Gardner solution on $\mathcal I=\R$ is called \emph{global}.

For a compact time interval $J$, set
\begin{align}
 \label{eq:L2u-definition} \norm{F}_{\Lu(J\times\R)}^2 :=\sup\left\{\int_{J'}\int_{\cK'}\abs{F(t,x)}^2\dd x\dd t \ :\ J'\subset J\text{ and }\cK'\subset\R \text{ intervals of length at most }1\right\}.
\end{align}
A global Gardner solution belongs to the \emph{Kato class} $X_1$ if
\begin{align}
 \label{eq:X1-bounds} r\in L^\infty(\R_t;L^2_x), \qquad r_x\in \Lu(J\times\R_x) \quad\text{for every compact }J\subset\R,
\end{align}
and, for every $t_0\in\R$,
\begin{align}
 \label{eq:X1-tail} \lim_{y\to\pm\infty} \int_{t_0}^{t_0+1}\int_{y-1}^{y+1}\abs{r_x(t,x)}^2\dd x\dd t=0.
\end{align}
Gardner solutions are the case $N=1$ of the weak solutions of \cite[Definition~5.7]{Klaus.Koch.Liu2023}, and the class $X_1$ is the case $N=1$ of \cite[Definition~2.18, equation~(2.63)]{Klaus.Koch.Liu2023}.
\end{definition}}

{
\begin{remark}[The cells in \eqref{eq:L2u-definition} must be intervals]
\label{rem:interval-cells}
A supremum over arbitrary measurable sets of measure at most one would define a strictly stronger norm than \eqref{eq:L2u-definition}. Take $F(t,x)=\sum_{n\geq2}n\mathbf 1_{[3n,\,3n+n^{-3}]}(x)$, which does not depend on $t$, and a time interval of unit length. Consecutive bumps are separated by $3-n^{-3}\geq\frac{23}8$, so a spatial interval of length two meets at most one of them. If it meets the $n$th, the cell integral of $\abs{F}^2$ over it is at most $n^2\cdot n^{-3}=n^{-1}\leq\frac12$, and otherwise it is zero. As the centre $y$ of the cell $[y-1,y+1]$ in \eqref{eq:X1-tail} tends to $+\infty$, the index of any bump that cell can meet tends to infinity, so the integral tends to zero. As $y$ tends to $-\infty$ the cells are eventually empty. Thus $F$ has finite interval-cell norm, and the translated cell integrals vanish at both ends, which is \eqref{eq:X1-tail} with $F$ in place of $r_x$. The union $\cK'=\bigcup_{n\geq2}[3n,3n+n^{-3}]$ is not an interval, has measure $\zeta(3)-1<1$, and gives $\int_{\cK'}\abs{F}^2=\sum_{n\geq2}n^{-1}=\infty$.

The bounds below use interval cells, as does the class $X_1$ of \cite{Klaus.Koch.Liu2023}.
\end{remark}
}

Equation~\eqref{eq:weak-Gardner-formulation} follows by testing \eqref{eq:Gardner} against $\varphi$ and integrating by parts once in time, three times in space on $\partial_x(-r_{xx})$, and once in space on the polynomial terms. All derivatives then fall on the test function. The requirement $r_x\in L^2_{\loc}$ in \eqref{eq:weak-Gardner-regularity} makes the cubic term locally integrable through the estimate below.

For the global class $X_1$, replacing time intervals of length one by intervals of length at most one is equivalent by restriction and enlargement. Replacing unit spatial cells by centered cells of length two is equivalent by finite coverings. The regularity is sufficient for the cubic term: if $J\Subset\mathcal I$ and $\cK\Subset\cK^*$ are bounded intervals, then the localized Gagliardo--Nirenberg inequality gives
\begin{align}
 \label{eq:Gardner-local-L6} \norm{r}_{L^6(J\times\cK)}^6 \lesssim\norm{r}_{L^\infty(J;L^2)}^4 \left(\norm{r_x}_{L^2(J\times\cK^*)}^2 +\abs{J}\norm{r}_{L^\infty(J;L^2)}^2\right).
\end{align}
Thus $r^2,r^3\in L^2_{\loc}$ and every term in \eqref{eq:weak-Gardner-formulation} is defined.

The following is the uniqueness theorem of Klaus--Koch--Liu \cite{Klaus.Koch.Liu2023} in the form used to prove Theorem~\ref{thm:main}, with the solution class specified in Definition~\ref{def:Gardner-solution}.

{\begin{lemma}
\label{lem:KKL-Gardner}
Fix $\tau\geq2$. Let $r_1,r_2\in C(\R;L^2(\R;\R))$ be global Gardner solutions of class $X_1$, in the sense of Definition~\ref{def:Gardner-solution}, with $r_1(0)=r_2(0)$. Then $r_1(t)=r_2(t)$ in $L^2(\R)$ for every $t\in\R$. This is the uniqueness assertion of the case $N=1$ of \cite[Theorem~2.24]{Klaus.Koch.Liu2023}.
\end{lemma}}

\begin{proof}
The class $X_N$ of \cite[Theorem~2.24]{Klaus.Koch.Liu2023} is global in time, and the uniqueness argument there is written forward in time, from $t=0$. The reflection $(t,x)\mapsto(-t,-x)$ preserves \eqref{eq:Gardner}, Definition~\ref{def:Gardner-solution}, the class $X_1$ and continuity into $L^2$, and supplies the statement for $t\leq0$.
\end{proof}

{\begin{remark}[The input from Klaus--Koch--Liu \cite{Klaus.Koch.Liu2023}]
\label{rem:KKL-imports}
The input from \cite{Klaus.Koch.Liu2023} to the proof of Theorem~\ref{thm:main} is the uniqueness assertion of Lemma~\ref{lem:KKL-Gardner}, applied in Step~4. The two global Gardner solutions to which it is applied are constructed in Section~\ref{sec:uniqueness}, from the Green solutions extended by the canonical flow of \cite{Killip.Visan2019}, and every hypothesis of the lemma is verified in this paper: continuity into $L^2$ in \eqref{eq:v-r-basic-reg}, the uniform $L^2$ bound in Lemma~\ref{lem:Gardner-mass}, the Kato conditions in Proposition~\ref{prop:r-X1-local}, and the weak equation in Proposition~\ref{prop:green-Gardner}. Remark~\ref{rem:shifted-Miura-injective} proves injectivity of the Miura map, and Lemma~\ref{lem:Miura} proves the identity $M_\tau(r)=q$. Formulas of \cite{Klaus.Koch.Liu2023} reproduced elsewhere in this section, as in Remarks~\ref{rem:explicit-inverse-Miura} and~\ref{rem:good-variable-normalization}, are cited for attribution.

In the notation of \cite{Klaus.Koch.Liu2023}, the fixed parameter $\tau$ is $\tau_0$ and the Gardner variable $r$ is $w$. A Gardner solution in the sense of Definition~\ref{def:Gardner-solution} is a weak solution in the sense of \cite[Definition~5.7]{Klaus.Koch.Liu2023} at $N=1$, and the class $X_1$ is \cite[Definition~2.18]{Klaus.Koch.Liu2023} at $N=1$, with the tail condition \cite[equation~(2.63)]{Klaus.Koch.Liu2023}.
\end{remark}}

%%%%%%%%%%%%%%%%%%%%%%%%%%%%%%%%%%%%%%%%%%%%%%%%%%%%%%%%%%%%%%%%%%%%%%%%%%%%%%%%%%
\subsection{Construction of the Variable \texorpdfstring{$r$}{r}}
\label{sec:Gardner-lift}
%%%%%%%%%%%%%%%%%%%%%%%%%%%%%%%%%%%%%%%%%%%%%%%%%%%%%%%%%%%%%%%%%%%%%%%%%%%%%%%%%%

We now construct an $L^2$-valued curve $r$ from the given Green solution $q$ and prove $M_\tau(r)=q$. This supplies the variable required by Definition~\ref{def:Gardner-solution}.

The construction must depend only on $q(t)$ and the fixed parameter $\tau$, so that equal initial data give equal values of $r(0)$. We also need $r\in C(I;L^2)$, because the uniqueness statement of Lemma~\ref{lem:KKL-Gardner} concerns curves that are continuous into $L^2$.

The recovery identity \eqref{eq:q-from-h} suggests the formula below. Substitution into the shifted Miura map will recover its right-hand side. Lemma~\ref{lem:Miura} justifies this calculation at the available rough regularity.

Fix $\tau\geq K(q)$, with $K(q)$ as in \eqref{eq:single-K}. Remark~\ref{rem:common-lambda} supplies admissibility, the Green identity at $\tau$, and the restriction $\tau\geq2$. Define
{
\begin{align}
 \label{eq:v-r-def} v:=\frac{h_\tau}{\tau}-1, \qquad r:=\tau v-\frac12\,\partial_x\log(1+v) .
\end{align}
}
The first formula gives
\begin{align}
 \label{eq:v-h-relation} h_\tau=\tau(1+v),
\end{align}
Since $\partial_x\log(1+v)=h_{\tau,x}/h_\tau$, the second formula can also be written directly in terms of $h_\tau$:
\begin{align}
 \label{eq:r-from-h-closed} r=h_\tau-\tau-\frac12\,\frac{h_{\tau,x}}{h_\tau} =\frac{1}{2g_\tau}-\tau+\frac12\,\frac{g_{\tau,x}}{g_\tau}.
\end{align}

{
\begin{remark}[The formula \eqref{eq:v-r-def} is an explicit inverse Miura map]
\label{rem:explicit-inverse-Miura}
Lemma~\ref{lem:Miura} proves $M_\tau(r)=q$, and Remark~\ref{rem:shifted-Miura-injective} proves injectivity on $L^2(\R;\R)$. Thus \eqref{eq:r-from-h-closed} gives the inverse of the shifted Miura map on the admissible potentials considered here, expressed through the diagonal Green function of $L_q+\tau^2$.

The formula itself is due to Klaus--Koch--Liu, and so are its mapping properties. In \cite[\S2.3]{Klaus.Koch.Liu2023} the function $v$ is defined from the diagonal Green function $\beta$ by $v=(2\tau\beta)^{-1}-1$, and the Gardner variable by $w=\tau v-\frac12\partial_x\log(1+v)$, which is their equation~(2.29). Under the dictionary $\beta=g_\tau$, which is \eqref{eq:intro-resolvent}, and $w=r$, these are exactly \eqref{eq:v-h-relation} and \eqref{eq:v-r-def}, hence \eqref{eq:r-from-h-closed}. Their \cite[Proposition~4.8(2)]{Klaus.Koch.Liu2023} at $N=0$ makes $M_\tau$ a diffeomorphism of $L^2(\R;\R)$ onto its range in $H^{-1}(\R)$, and \cite[Theorem~4.13]{Klaus.Koch.Liu2023} at $s=0$ does the same for $v\mapsto r$. The formula does need $h_\tau$, hence a parameter $\tau$ that is \emph{admissible} for $q$ in the sense of Proposition~\ref{prop:static}, which follows from $\tau\geq K(q)$ throughout Section~\ref{sec:Gardner}.

Along a Green trajectory in the class of Definition~\ref{def:green}, Section~\ref{sec:Gardner} proves that \eqref{eq:r-from-h-closed} defines a Gardner solution in the sense of Definition~\ref{def:Gardner-solution}. Lemma~\ref{lem:internal-H2} obtains one more spatial derivative of $v$ from the Miura identity, Proposition~\ref{prop:r-X1-local} proves the Kato conditions locally in time, and Proposition~\ref{prop:green-Gardner} proves the weak Gardner equation.

For fixed $\tau$, the right side of \eqref{eq:r-from-h-closed} depends on $q$ through $h_\tau(q)$. Thus two competitors with the same initial data have the same initial value of $r$, as required in the proof of Theorem~\ref{thm:main}.
\end{remark}
}

By Proposition~\ref{prop:static}\,\textup{(ii)}, $h_\tau-\tau\in H^1(\R)$ at each time, with $t\mapsto h_\tau(t)-\tau$ continuous into $H^1$, and by \eqref{eq:g-quantitative} $h_\tau$ is bounded below by a positive constant, uniformly on $I\times\R$. Hence
\begin{align}
 \label{eq:v-r-basic-reg} v\in C(I;H^1), \qquad r\in C(I;L^2), \qquad 1+v\geq c_I>0 \quad\text{on }I\times\R .
\end{align}
The middle assertion uses both parts: $\tau v\in C(I;H^1)$ directly, while $\partial_x\log(1+v)=v_x/(1+v)$ is the product of $v_x\in C(I;L^2)$ with the bounded factor $(1+v)^{-1}$. The lower bound also follows from positivity and decay: at each fixed $t$, $h_\tau(t)-\tau\in C_0(\R)$ and $h_\tau(t,\cdot)>0$, so $h_\tau(t,x)\to\tau>0$ as $\abs{x}\to\infty$ and $\inf_xh_\tau(t,x)>0$. Continuity in $t$ and compactness of $I$ make it uniform.

{\begin{lemma}\label{lem:Miura}
Let $q\in C(I;H^{-1}(\R;\R))$ be a Green solution in the sense of Definition~\ref{def:green}, let $\tau\geq K(q)$ be the parameter fixed above, and define $v$ and $r$ from $q$ by \eqref{eq:v-r-def}. Then $r\in C(I;L^2(\R;\R))$ by \eqref{eq:v-r-basic-reg}, and $q$ is the image of $r$ under the shifted Miura map \eqref{eq:Miura-map}:
\begin{align}
 \label{eq:Miura} q=M_\tau(r)=r_x+2\tau r+r^2 \quad\text{in }C(I;H^{-1}(\R)).
\end{align}
\end{lemma}}

\begin{proof}
The proof has two parts: an algebraic identity valid for smooth functions, and a passage to the regularity \eqref{eq:v-r-basic-reg} actually available.

\proofstep{Step 1: the algebraic identity}
Let $v$ be smooth with $1+v$ bounded below by a positive constant, and let $r=\tau v-\frac12v_x/(1+v)$ as in \eqref{eq:v-r-def}. Expanding the three terms of $M_\tau(r)$ and collecting gives
\begin{align}
 \label{eq:composite-v} r_x+2\tau r+r^2 =-\frac{v_{xx}}{2(1+v)} +\frac{3v_x^2}{4(1+v)^2} +\tau^2v^2+2\tau^2v .
\end{align}
On the other hand, substituting \eqref{eq:v-h-relation} into the recovery identity \eqref{eq:q-from-h} at $\lambda=\tau$, and using $h_{\tau,x}=\tau v_x$ and $h_{\tau,xx}=\tau v_{xx}$, gives
\begin{align}
 \label{eq:recovery-in-v} q=-\frac{v_{xx}}{2(1+v)} +\frac{3v_x^2}{4(1+v)^2} +\tau^2(1+v)^2-\tau^2 ,
\end{align}
whose last two terms are $\tau^2v^2+2\tau^2v$. The right sides of \eqref{eq:composite-v} and \eqref{eq:recovery-in-v} agree term by term, which is \eqref{eq:Miura} for smooth $v$.

\proofstep{Step 2: the actual regularity}
For $v\in C(I;H^1)$ as in \eqref{eq:v-r-basic-reg}, we prove \eqref{eq:composite-v} by approximation. Equation~\eqref{eq:recovery-in-v} is the recovery identity \eqref{eq:q-from-h} at $\lambda=\tau$ with $h_\tau=\tau(1+v)$ substituted, each of its terms read through the multiplier calculus in \eqref{eq:affine-multiplier-definition}--\eqref{eq:multipliers-a4}, and it holds for the rough $v$ as it stands. Mollify $v(t)$ in the spatial variable with a nonnegative mollifier. Nonnegativity preserves the bound $1+v\geq c_I$, so every denominator in \eqref{eq:composite-v} stays away from zero along the approximation. Since the image of $v\colon I\to H^1$ is compact, the $H^1$ convergence is uniform in $t\in I$.

It remains to see that each term of \eqref{eq:composite-v} converges in $H^{-1}$. The reciprocal $(1+v)^{-1}$ lies in the affine algebra \eqref{eq:affine-Hone-algebra}, with the reciprocal and products justified by \eqref{eq:affine-multiplier-definition}--\eqref{eq:multipliers-a4}. Moreover $v\mapsto(1+v)^{-1}$ is continuous from $\{v\in H^1:1+v\geq c_I\}$ into $\mathcal A$, since
\begin{align}
 \label{eq:reciprocal-continuity} (1+v^\varepsilon)^{-1}-(1+v)^{-1}=\frac{v-v^\varepsilon}{(1+v^\varepsilon)(1+v)}
\end{align}
has $L^2$ norm at most $c_I^{-2}\norm{v-v^\varepsilon}_2$ and a derivative of $L^2$ norm at most $c_I^{-2}\norm{(v-v^\varepsilon)_x}_2+c_I^{-3}\norm{v-v^\varepsilon}_\infty\bigl(\norm{v^\varepsilon_x}_2+\norm{v_x}_2\bigr)$; together with the bilinear bound \eqref{eq:multipliers-a1} this is what makes the products below converge. The second-derivative term $v_{xx}/(1+v)$ is the action of the multiplier $(1+v)^{-1}$ on $v_{xx}\in H^{-1}$, whose continuity follows from \eqref{eq:multipliers-a1}. The squared first derivative $v_x^2$ lies in $L^1(\R)\subset H^{-1}(\R)$ and converges there, and $(1+v)^{-2}$ acts on it by the same calculus. The terms $\tau^2v^2+2\tau^2v$ converge in $H^1$. The left side converges as well: $v^\varepsilon\to v$ in $H^1$ and the uniform lower bound give $r^\varepsilon\to r$ in $L^2$, hence $r^\varepsilon_x\to r_x$ in $H^{-1}$ and $(r^\varepsilon)^2\to r^2$ in $L^1(\R)\subset H^{-1}(\R)$, so $M_\tau(r^\varepsilon)\to M_\tau(r)$ in $H^{-1}$. Every term of \eqref{eq:composite-v} is therefore continuous into $H^{-1}$, uniformly in $t$, and the identity passes to the limit in $C(I;H^{-1})$.
\end{proof}

{
\begin{remark}
Three points are recorded for later use. First, the lemma is proved at each fixed time. Its proof uses the admissibility of $\tau$ for $q(I)$ and the recovery identity of Proposition~\ref{prop:static}. Second, \eqref{eq:v-r-def} determines $r(t)$ from $q(t)$ and the fixed parameter $\tau$. In particular two Green solutions with $q_1(0)=q_2(0)$ satisfy $r_1(0)=r_2(0)$. Third, $r$ is the only curve in $C(I;L^2)$ satisfying \eqref{eq:Miura}, because $M_\tau$ is injective on $L^2(\R;\R)$ for $\tau>0$ by Remark~\ref{rem:shifted-Miura-injective}.
\end{remark}
}

The $L^2$ norm of the Gardner variable is the conserved mass of Section~\ref{sec:alpha}. This is what bounds $r$ in $L^2$ uniformly in time once the trajectory has been extended to all times in Section~\ref{sec:uniqueness}.

{\begin{lemma}\label{lem:Gardner-mass}
Let $q$ be a Green solution, let $\tau\geq K(q)$ be the parameter fixed above, and let $r$ be given by \eqref{eq:v-r-def}. Then
\begin{align}
 \label{eq:Gardner-mass} \norm{r(t)}_2^2=2\tau\,\alpha_\tau(q(t)) \qquad(t\in I).
\end{align}
In particular, $\norm{r(t)}_2=\norm{r(0)}_2$ for every $t\in I$.
\end{lemma}}

\begin{proof}
Fix $t\in I$ and suppress it from the notation. Since $v\in H^1$ and $1+v\geq c_I$ by \eqref{eq:v-r-basic-reg}, the function $\log(1+v)$ lies in $H^1$. By \eqref{eq:v-h-relation} and \eqref{eq:v-r-def},
\begin{align}
 \label{eq:Gardner-mass-a1} w_\tau=h_\tau-\tau=\tau v=r+\frac12\,\partial_x\log(1+v) .
\end{align}
Apply $2\tau\Rf_\tau$, which maps $H^{-1}$ into $H^1$ and commutes with $\partial_x$, to the Miura identity \eqref{eq:Miura}. By the definition \eqref{eq:R0}, $4\tau^2\Rf_\tau=\mathrm{Id}+\partial_x^2\Rf_\tau$, so that
\begin{align}
 \label{eq:Gardner-mass-a2} 2\tau\Rf_\tau q=2\tau\Rf_\tau(r^2)+r+\partial_x\bigl(2\tau\Rf_\tau r+\partial_x\Rf_\tau r\bigr).
\end{align}
Subtracting \eqref{eq:Gardner-mass-a1} and recalling \eqref{eq:def-rho} gives
\begin{align}
 \label{eq:Gardner-mass-a3} \rho_\tau(q)=2\tau\Rf_\tau(r^2)+\partial_xP, \qquad P:=2\tau\Rf_\tau r+\partial_x\Rf_\tau r-\frac12\log(1+v) .
\end{align}
Here $P\in H^1\subset H^{-1}$, because $\Rf_\tau$ maps $L^2$ into $H^2$. By \eqref{eq:R0-kernel}, $\Rf_\tau$ is convolution with the nonnegative kernel $\ue^{-2\tau\abs{x}}/(4\tau)$, whose integral is $(4\tau^2)^{-1}$. Since $r^2\in L^1$, the function $2\tau\Rf_\tau(r^2)$ is integrable, with integral $\norm{r}_2^2/(2\tau)$. As $\tau\geq K(q)\geq\Lambda(q(I))$, the density $\rho_\tau(q)$ is integrable by \eqref{eq:rho-nonneg}. Hence $\partial_xP\in L^1$, and Lemma~\ref{lem:zero-mass} gives $\int_\R\partial_xP\dd x=0$. Integrating \eqref{eq:Gardner-mass-a3} over $\R$ therefore gives $\alpha_\tau(q)=\norm{r}_2^2/(2\tau)$, which is \eqref{eq:Gardner-mass}. The last assertion follows from Proposition~\ref{prop:alpha-conservation} with $\kappa=\tau\geq K(q)$.
\end{proof}

Identity \eqref{eq:Gardner-mass} is the case $z=\ui\tau$ of \cite[equation~(3.12)]{Klaus.Koch.Liu2023}, which reads $\mathcal T^{\mathrm{KdV}}_{-1}(\ui\tau,q)=-\frac12\norm{r}_2^2$ in the present notation, once $\mathcal T^{\mathrm{KdV}}_{-1}(\ui\tau,q)=-\tau\alpha_\tau(q)$ is read off from \cite[equation~(2.7)]{Klaus.Koch.Liu2023} and \cite[equation~(1.5)]{Killip.Visan2019}.

%%%%%%%%%%%%%%%%%%%%%%%%%%%%%%%%%%%%%%%%%%%%%%%%%%%%%%%%%%%%%%%%%%%%%%%%%%%%%%%%%%
\subsection{Uniformly Local \texorpdfstring{$H^2$}{H2} Regularity of \texorpdfstring{$v$}{v} and Kato Conditions for \texorpdfstring{$r$}{r}}
\label{sec:Gardner-bootstrap}
\label{sec:Gardner-smoothing}
%%%%%%%%%%%%%%%%%%%%%%%%%%%%%%%%%%%%%%%%%%%%%%%%%%%%%%%%%%%%%%%%%%%%%%%%%%%%%%%%%%

Having constructed $r$, we now prove the derivative estimate needed for both the Kato conditions on $r_x$ and the Gardner equation. This subsection proves the estimate and verifies the Kato conditions.

The definition \eqref{eq:v-r-def} of $r$ is justified by the three assertions of \eqref{eq:v-r-basic-reg}. By the embedding $H^1(\R)\hookrightarrow L^\infty(\R)$ of Lemma~\ref{lem:products}\,\textup{(i)}, $v$ is bounded and continuous at each time. The lower bound $1+v\geq c_I>0$ keeps the denominator away from zero uniformly on $I\times\R$. $v_x\in L^2(\R)$. Hence the logarithmic derivative
\begin{align}
 \label{eq:mu-def} \mu:=\partial_x\log(1+v)=\frac{v_x}{1+v} ,
\end{align}
a quantity that recurs throughout this subsection and the next, is already an honest element of $C(I;L^2(\R))$, and $r=\tau v-\mu/2$ is the $C(I;L^2)$ curve constructed in Section~\ref{sec:Gardner-lift}.

The obstacle is the \emph{flux}. Section~\ref{sec:Gardner-dynamics} writes the equation for $v$ as $v_t=\partial_x\Psi$, with the flux $\Psi$ of \eqref{eq:Psi-def} below, and then turns it into the Gardner equation for $r$ by identifying $\Psi$ with an algebraic expression $\cW(v)$ in $v$ and its derivatives, given in \eqref{eq:algebraic-A} below. That identification is an identity between \emph{smooth} functions, extended to $v$ at the available regularity by mollifying $v$ and passing to the limit in the nonlinear expression $\cW(v^\varepsilon)$. Its leading term is $-(1+v)\mu_x$, where
\begin{align}
 \label{eq:mux-needs-H2} \mu_x=\frac{v_{xx}}{1+v}-\frac{v_x^2}{(1+v)^2} ,
\end{align}
and its next term is quadratic in $\mu$. To justify $\cW(v)\in L^2_{\loc}$ and convergence of its displayed terms under spatial mollification, we establish $\mu_x\in L^2_{\loc}$ and $\mu\in L^4_{\loc}$. Neither follows from $v\in C(I;H^1)$.

We obtain both properties from the Miura identity \eqref{eq:Miura} by first estimating $r_x$. By \eqref{eq:v-r-def} and \eqref{eq:mu-def},
\begin{align}
 \label{eq:mu-from-r} \mu=2(\tau v-r) ,
\end{align}
and $v$ is bounded with $v_x\in L^2(\R)$ at each time, so $\mu_x\in L^2_{\loc}$ and $\mu\in L^4_{\loc}$ are equivalent to $r_x\in L^2_{\loc}$ and $r\in L^4_{\loc}$. Written as $r_x+r^2=q-2\tau r$, the identity \eqref{eq:Miura} has a locally square integrable right side, by Proposition~\ref{prop:auto-smoothing}. Squaring it against a cutoff separates the two terms on the left, because the cross term $2r_xr^2=\frac23\partial_x(r^3)$ is a derivative and, after an integration by parts, costs only a lower-order term in $r$.

This subsection carries out that argument, in the uniformly local form the next subsection needs. Its output is Lemma~\ref{lem:internal-H2}, and then Corollary~\ref{cor:mu-regularity}, which records the two properties of $\mu$ named above. The next subsection consumes both: the lemma, together with the identity \eqref{eq:Miura-pointwise} established in its proof, in its Steps~1, 3 and~4, and the corollary in its Steps~4 and~5. The Kato conditions for $r_x$ are Proposition~\ref{prop:r-X1-local} at the end of the subsection. Its cell bound and its translated-cell tail are both read off from the cutoff estimate \eqref{eq:r-cell-estimate} below.

All uniformly local estimates below take the supremum after integration in time. For an integer $m\geq0$ and a compact time interval $J$, the condition is
\begin{align}
 \label{eq:Hm-unif-definition} \sup_{y\in\R}\sum_{\ell=0}^m\int_J\int_{y-1}^{y+1}\abs{\partial_x^\ell f(t,x)}^2\dd x\dd t<\infty.
\end{align}
This is equivalent to $\partial_x^\ell f\in \Lu(J\times\R)$ for $0\leq\ell\leq m$, with $\Lu$ as in \eqref{eq:L2u-definition}. To pass from these bounds to \eqref{eq:Hm-unif-definition}, partition $J$ into at most $\lceil\abs{J}\rceil+1$ intervals of length at most one and each spatial cell into two unit intervals. Conversely, each unit spacetime cell is contained in $J$ times a cell of radius one. The comparison constants depend only on $m$ and $\abs{J}$. Changing the spatial radius requires only another finite covering. Throughout the subsection $\cK_y:=(y-1,y+1)$ and $\cK_y^*:=(y-2,y+2)$ denote the cell centred at $y$ and its enlargement.

The zero-order bound for $r$ follows from \eqref{eq:v-r-basic-reg}: for every compact $J\subset I$,
\begin{align}
 \label{eq:r-zero-order-unif} \sup_{y\in\R}\int_J\int_{y-1}^{y+1}\abs{r(t,x)}^2\dd x\dd t \leq \abs{J}\sup_{t\in J}\norm{r(t)}_2^2<\infty .
\end{align}

{\begin{lemma}
\label{lem:internal-H2}
Let $q$ be a Green solution, let $\tau\geq K(q)$ be the parameter fixed in Section~\ref{sec:Gardner-lift}, and let $v$ be given by \eqref{eq:v-r-def}. Then
\begin{align}
 \label{eq:v-final-reg} v\in L^\infty(I;H^1(\R)),\qquad \partial_x^\ell v\in \Lu(J\times\R)\quad(0\leq\ell\leq2)
\end{align}
for every compact $J\subset I$.
\end{lemma}}

The two inputs are the Miura identity \eqref{eq:Miura} of Lemma~\ref{lem:Miura} and the uniform cell bound \eqref{eq:q-Kato-bound} of Proposition~\ref{prop:auto-smoothing}. Remark~\ref{rem:good-variable-normalization} explains the regularity required by the corresponding equivalence in \cite{Klaus.Koch.Liu2023}.

\begin{proof}
Throughout, \eqref{eq:v-r-basic-reg} and the embedding $H^1(\R)\hookrightarrow L^\infty(\R)$ give, uniformly in $t\in I$,
\begin{align}
 \label{eq:internal-H2-a1} \norm{v(t)}_{H^1}+\norm{v(t)}_{\infty}+\norm{r(t)}_2\leq C_I .
\end{align}
Fix a compact interval $J\subset I$. All constants below depend only on $C_I$ and $\tau$, never on the centre of a cell or on $t$.

\proofstep{Step 1: the Miura identity as an identity of functions}
Work at a time $t\in J$ of the common set of Remark~\ref{rem:measurability}, at which the distribution $q(t)$ is given by its local $L^2$ representative. By Lemma~\ref{lem:Miura}, at this time
\begin{align}
 \label{eq:Miura-pointwise} r_x+r^2=q-2\tau r \qquad\text{in }\mathcal D'(\R) .
\end{align}
The right side lies in $L^2_{\loc}(\R)$ and $r^2\in L^1(\R)$, so $r_x\in L^1_{\loc}(\R)$. Hence $r(t)$ is locally absolutely continuous, and in particular locally bounded, so that $r^2\in L^\infty_{\loc}(\R)$ and then $r_x\in L^2_{\loc}(\R)$ by \eqref{eq:Miura-pointwise}. Thus $r(t)\in H^1_{\loc}(\R)$, and \eqref{eq:Miura-pointwise} holds almost everywhere on $\R$ as an identity between functions in $L^2_{\loc}(\R)$.

\proofstep{Step 2: the cutoff estimate}
Keep $t$ as in Step~1, let $\chi\in C_c^\infty(\R)$ be real, multiply the square of \eqref{eq:Miura-pointwise} by $\chi^2$, and integrate over $\R$. By Step~1, $r^3$ is locally absolutely continuous with $\partial_x(r^3)=3r^2r_x$, so the cross term integrates by parts, $2\int_\R\chi^2r^2r_x\dd x=\frac23\int_\R\chi^2\partial_x(r^3)\dd x=-\frac43\int_\R\chi\chi'r^3\dd x$. Therefore
\begin{align}
 \label{eq:cutoff-exact} \int_\R\chi^2(q-2\tau r)^2\dd x=\int_\R\chi^2r_x^2\dd x+\int_\R\chi^2r^4\dd x-\frac43\int_\R\chi\chi'r^3\dd x .
\end{align}
Young's inequality gives $\frac43\abs{\chi r^2}\abs{\chi'r}\leq\frac12\chi^2r^4+\frac89(\chi')^2r^2$, and $(q-2\tau r)^2\leq2q^2+8\tau^2r^2$. Inserting both bounds in \eqref{eq:cutoff-exact} gives
\begin{align}
 \label{eq:cutoff-estimate} \int_\R\chi^2\Bigl(r_x^2+\frac12r^4\Bigr)\dd x\leq2\int_\R\chi^2q^2\dd x+\int_\R\Bigl(8\tau^2\chi^2+\frac89(\chi')^2\Bigr)r^2\dd x .
\end{align}

\proofstep{Step 3: integration over cells and in time}
Fix $\chi\in C_c^\infty(\R)$ with $0\leq\chi\leq1$, $\chi=1$ on $[-1,1]$ and $\supp\chi\subset[-2,2]$, and apply \eqref{eq:cutoff-estimate} to the translate $\chi(\cdot-y)$ at every time of the common set in $J$. Integrating over $J$ gives, for every $y\in\R$,
\begin{align}
 \label{eq:r-cell-estimate} \int_J\int_{\cK_y}\bigl(r_x^2+r^4\bigr)\dd x\dd t\leq4\int_J\int_{\cK_y^*}q^2\dd x\dd t+\Bigl(16\tau^2+\frac{16}9\norm{\chi'}_\infty^2\Bigr)\int_J\int_{\cK_y^*}r^2\dd x\dd t .
\end{align}
The integrands are jointly measurable on $J\times\R$, because \eqref{eq:Miura-pointwise} expresses $r_x$ through the representative of $q$ and through $r\in C(I;L^2)$, so the time integration is legitimate. Up to the single point $y$, the enlarged cell $\cK_y^*$ is the union of the two cells $\cK_{y-1}$ and $\cK_{y+1}$ of unit radius, so the first integral on the right of \eqref{eq:r-cell-estimate} is at most twice the supremum in \eqref{eq:q-Kato-bound}. The second is at most $\abs{J}C_I^2$, by \eqref{eq:internal-H2-a1}, as in \eqref{eq:r-zero-order-unif}. Hence
\begin{align}
 \label{eq:r-cell-bound} \sup_{y\in\R}\int_J\int_{\cK_y}\bigl(r_x^2+r^4\bigr)\dd x\dd t<\infty .
\end{align}

\proofstep{Step 4: from $r$ back to $v$}
At every time of the common set, \eqref{eq:mu-from-r} and Step~1 place $\mu$ in $H^1_{\loc}(\R)$, with $\mu_x=2(\tau v_x-r_x)$. By \eqref{eq:mu-def}, $v_x=(1+v)\mu$ is the product of two locally absolutely continuous functions, hence is itself locally absolutely continuous, and the product rule gives
\begin{align}
 \label{eq:vxx-from-mu} v_{xx}=(1+v)\bigl(\mu_x+\mu^2\bigr) .
\end{align}
By \eqref{eq:internal-H2-a1}, $\abs{1+v}\leq1+C_I$ and $v^4\leq C_I^2v^2$, while $\abs{\mu_x}\leq2\tau\abs{v_x}+2\abs{r_x}$ and $\mu^2\leq8\tau^2v^2+8r^2$ by \eqref{eq:mu-from-r}. Therefore, almost everywhere on $J\times\R$,
\begin{align}
 \label{eq:vxx-pointwise} \abs{v_{xx}}^2\leq C\bigl(v^2+v_x^2+r_x^2+r^4\bigr) ,
\end{align}
with $C$ depending only on $C_I$ and $\tau$. Integrating \eqref{eq:vxx-pointwise} over $J\times\cK_y$ and using \eqref{eq:r-cell-bound} for the last two terms and \eqref{eq:internal-H2-a1} for the first two bounds the second-order term of \eqref{eq:Hm-unif-definition} uniformly in $y$. The zero- and first-order terms are controlled by $v\in L^\infty(I;H^1)$, which is \eqref{eq:v-r-basic-reg}. This proves \eqref{eq:v-final-reg}, once the slicewise functions are identified with spacetime distributions. The right sides of \eqref{eq:Miura-pointwise} and \eqref{eq:vxx-from-mu} are jointly measurable on $J\times\R$, being built from the representative of $q$ and from $r$, $v$ and $v_x$, so the functions $(t,x)\mapsto r_x(t)(x)$ and $(t,x)\mapsto v_{xx}(t)(x)$ lie in $L^2_{\loc}(J\times\R)$ by the bounds just proved. Since $\pair*{v_{xx}(t)}{\varphi}=\pair*{v(t)}{\varphi''}$ for every $\varphi\in C_c^\infty(\R)$ at almost every $t$, Fubini applied to test functions on $J^\circ\times\R$ shows that the second of these functions is the distribution $\partial_x^2v$ on $J^\circ\times\R$. In the same way the first is $\partial_xr$, and $(t,x)\mapsto v_x(t)(x)$ is $\partial_xv$.
\end{proof}

The lemma is stated for $v$, because Section~\ref{sec:Gardner-dynamics} mollifies $v$ and uses \eqref{eq:v-final-reg} in that form. What that subsection also consumes is the corresponding statement for the logarithmic derivative $\mu$ of \eqref{eq:mu-def}: the two properties $\mu_x\in L^2$ and $\mu\in L^4$ on cells, uniformly in the centre, used below to justify the flux and its approximation. The proof has already produced them. We record them here, so that Section~\ref{sec:Gardner-dynamics} can quote a single statement.

{\begin{corollary}[Uniformly local $H^1$ regularity of $\mu$]
\label{cor:mu-regularity}
Under the hypotheses of Lemma~\ref{lem:internal-H2}, the function $\mu=\partial_x\log(1+v)$ of \eqref{eq:mu-def} satisfies
\begin{align}
 \label{eq:mu-regularity} \mu\in L^\infty(I;L^2(\R)),\qquad \mu,\mu_x\in \Lu(J\times\R)
\end{align}
for every compact $J\subset I$. In particular $\mu_x\in L^2(J\times\cK)$ and $\mu\in L^4(J\times\cK)$ for every bounded interval $\cK$, with bounds depending on $\cK$ only through its length.
\end{corollary}}

\begin{proof}
Equation~\eqref{eq:mu-def} gives $\mu\in C(I;L^2(\R))$, hence $\mu\in L^\infty(I;L^2)$. By \eqref{eq:mu-from-r}, $\mu_x=2(\tau v_x-r_x)$ and $\mu^4\leq128(\tau^4v^4+r^4)$. The cell bound for $\mu_x$ in \eqref{eq:mu-regularity} therefore follows from \eqref{eq:r-cell-bound} and $v_x\in L^\infty(I;L^2)$. The two consequences follow by covering $\cK$ with finitely many cells of radius one, whose number depends only on the length of $\cK$: the first directly from \eqref{eq:mu-regularity}, and the second from \eqref{eq:r-cell-bound} together with $v^4\leq C_I^2v^2$.
\end{proof}

The cutoff estimate \eqref{eq:r-cell-estimate} gives both the smoothing condition \eqref{eq:X1-bounds} and the translated-cell tail \eqref{eq:X1-tail} for $r_x$, on every compact subinterval of $I$.

{\begin{proposition}\label{prop:r-X1-local}
Let $q$ be a Green solution, let $\tau\geq K(q)$ be the parameter fixed in Section~\ref{sec:Gardner-lift}, and let $r$ be as in \eqref{eq:v-r-def}. Then, for every compact $J\subset I$, the derivative $r_x$ belongs to $L^2_{\loc}(J\times\R)$ and satisfies the uniform cell bound
\begin{align}
 \label{eq:r-local-smoothing} \sup_{y\in\R}\int_J\int_{y-1}^{y+1}\abs{r_x(t,x)}^2\dd x\dd t<\infty
\end{align}
together with the translated-cell tail
\begin{align}
 \label{eq:r-local-tail} \int_J\int_{y-1}^{y+1}\abs{r_x(t,x)}^2\dd x\dd t \rightarrow0\qquad(y\to\pm\infty).
\end{align}
Thus $r$ satisfies the requirements \eqref{eq:X1-bounds} and \eqref{eq:X1-tail} of the class $X_1$ on every compact subinterval of $I$.
\end{proposition}}

\begin{proof}
By the end of the proof of Lemma~\ref{lem:internal-H2}, the distribution $r_x$ is the function $q-2\tau r-r^2$ of \eqref{eq:Miura-pointwise}, and it lies in $L^2_{\loc}(J\times\R)$. The cell bound \eqref{eq:r-local-smoothing} is contained in \eqref{eq:r-cell-bound}. For the tail, the right side of \eqref{eq:r-cell-estimate} tends to zero as $y\to\pm\infty$: its first term by \eqref{eq:q-Kato-tail} applied to the two cells $\cK_{y-1}$ and $\cK_{y+1}$, and its second by dominated convergence, since for each $t\in J$ the integral of $r(t)^2$ over $\cK_y^*$ tends to zero as $y\to\pm\infty$, $r(t)$ being a fixed function in $L^2(\R)$, and is bounded by $C_I^2$.
\end{proof}

{\begin{remark}
\label{rem:good-variable-normalization}
The formulas for $v$ and $r$ in \eqref{eq:v-r-def}, together with the constraint $1+v\geq c>0$, use the normalization of \cite[equations~(2.29), (2.31) and (4.14)]{Klaus.Koch.Liu2023}. The reverse direction of the rough KdV--Gardner equivalence proved there assumes global $L^2$ control of the KdV variable. Definition~\ref{def:green} gives continuity into $H^{-1}$, and Proposition~\ref{prop:auto-smoothing} gives uniformly local $L^2$ control in spacetime. At this regularity, Lemma~\ref{lem:internal-H2} proves uniformly local $H^2$ regularity of $v$, and Proposition~\ref{prop:green-Gardner} derives the Gardner equation from \eqref{eq:v-r-def}.
\end{remark}}
%%%%%%%%%%%%%%%%%%%%%%%%%%%%%%%%%%%%%%%%%%%%%%%%%%%%%%%%%%%%%%%%%%%%%%%%%%%%%%%%%%
\subsection{Green Identity Becomes Gardner Equation}
\label{sec:Gardner-dynamics}
%%%%%%%%%%%%%%%%%%%%%%%%%%%%%%%%%%%%%%%%%%%%%%%%%%%%%%%%%%%%%%%%%%%%%%%%%%%%%%%%%%

Define the flux for the equation satisfied by $v$:
\begin{align}
 \label{eq:Psi-def} \Psi:=-4\tau^2v+2(1+v)q .
\end{align}
Dividing the Green identity \eqref{eq:integrated-w} at parameter $\tau$ by $\tau$ and using \eqref{eq:J-2hq} gives
\begin{align}
 \label{eq:v-integrated} v(t)-v(s) =\partial_x\int_s^t\Psi(t')\dd t'.
\end{align}
This is the normalized form of \cite[equation~(2.34)]{Klaus.Koch.Liu2023}, whose flux $2(1+v)(q-2\tau^2)$ differs from the one in \eqref{eq:v-integrated} by the spatial constant $-4\tau^2$ and therefore has the same spatial derivative; the flux of \eqref{eq:v-integrated} is exactly that of \cite[equation~(2.36)]{Klaus.Koch.Liu2023}, and, expressed through $v$ alone, it is the flux $-v_{xx}+\frac32v_x^2/(1+v)+6\tau^2v^2+2\tau^2v^3$ of the good-variable equation \cite[equation~(3.37)]{Klaus.Koch.Liu2023}, see \eqref{eq:algebraic-A}. The flux in \eqref{eq:v-integrated} is locally $L^2$ by Proposition~\ref{prop:auto-smoothing}.

{\begin{proposition}\label{prop:green-Gardner}
Let $q\in C(I;H^{-1}(\R;\R))$ be a Green solution on $I=[-T,T]$, let $\tau\geq K(q)$ be the parameter fixed in Section~\ref{sec:Gardner-lift}, and let $v$ and $r$ be given by \eqref{eq:v-r-def}. Then $r$ is a Gardner solution on the open interval $(-T,T)$ in the sense of Definition~\ref{def:Gardner-solution}. That is, $r\in L^\infty((-T,T);L^2(\R))$ and $r_x\in L^2_{\loc}((-T,T)\times\R)$, and the first $\tau$-Gardner equation \eqref{eq:Gardner} holds in the sense of distributions on $(-T,T)\times\R$.
\end{proposition}}

The proof uses two approximation arguments. Steps~1 and~2 establish the time chain rule for $\log(1+v)$, which enters $r$ through its spatial derivative. Steps~3 and~4 identify the resulting flux with $-r_{xx}+6\tau r^2+2r^3$. Step~5 collects the regularity required by Definition~\ref{def:Gardner-solution}.

\begin{proof}
Throughout, $v$ and $r$ are the functions \eqref{eq:v-r-def} attached to the given Green solution, and all regularity used is that of \eqref{eq:v-r-basic-reg} and Lemma~\ref{lem:internal-H2}.

\proofstep{Step 1: the flux in the equation for $v$}
The flux is $\Psi$ from \eqref{eq:Psi-def}. Here $(1+v)q$ initially denotes the affine $H^1$ multiplier product defined in \eqref{eq:affine-multiplier-definition}. For almost every $t$ in the common full-measure set of Remark~\ref{rem:measurability} and every $\phi\in C_c^\infty(\R)$,
\begin{align}
 \label{eq:hq-representative-agreement} \pair*{(1+v)q}{\phi}:=\pair* q{(1+v)\phi} =\int_\R q(t,x)\bigl(1+v(t,x)\bigr)\phi(x)\dd x,
\end{align}
where $q$ on the right is its local $L^2$ representative. The equality follows by approximating $(1+v)\phi\in H^1$ by compactly supported smooth functions in $H^1$. Thus the distributional multiplier product agrees almost everywhere with the ordinary product. Since $1+v$ is uniformly bounded and $q\in L^2_{\loc}$, this also proves $(1+v)q\in L^2_{\loc}$. Equation \eqref{eq:v-integrated} says
\begin{align}
 \label{eq:v-conservation} v_t=\partial_x\Psi
\end{align}
in spacetime distributions, and $\Psi\in L^2_{\loc}$ by the preceding paragraph together with $v\in C(I;H^1)$. By \eqref{eq:v-final-reg} and \eqref{eq:v-r-basic-reg} we may use throughout that
\begin{align}
 \label{eq:v-standing-regularity} \begin{aligned}
 &v\in L^\infty(I;H^1),\qquad \partial_x^\ell v\in \Lu(J\times\R)\quad(0\leq\ell\leq2),\\
 &1+v\geq c_I>0,
 \end{aligned}
\end{align}
for every compact $J\subset I$. We now combine \eqref{eq:v-conservation} with the relation between $v$ and $r$, expressed through $\mu$ in \eqref{eq:mu-def}.

\proofstep{Step 2: the chain rule for $\log(1+v)$}
The function $r$ is defined through $\partial_x\log(1+v)$, so $r_t$ requires differentiating a logarithm of $v$ in time. This is done by mollifying in space, where the ordinary chain rule applies, and passing to the limit. Let $\eta_\varepsilon$ be a nonnegative spatial mollifier of mass one and write
\begin{align}
 \label{eq:v-mollification} v^\varepsilon:=\eta_\varepsilon*v, \qquad \Psi_\varepsilon:=\eta_\varepsilon*\Psi .
\end{align}
Since $\eta_\varepsilon$ has mass one, mollifying $1+v$ produces $1+v^\varepsilon$. In particular $1+v^\varepsilon\geq c_I$ for every $\varepsilon$, and $v^\varepsilon_t=\partial_x\Psi_\varepsilon\in L^2_{\loc}$. Thus $v^\varepsilon$ is bounded and continuous with $v^\varepsilon,\partial_tv^\varepsilon\in L^2_{\loc}((-T,T)\times\R)$, so for almost every $x$ the slice $t\mapsto v^\varepsilon(t,x)$ is absolutely continuous with derivative $\partial_x\Psi_\varepsilon(\cdot,x)$; the chain rule for the logarithm on $[c_I,\infty)$, where its derivative is bounded, applies slicewise to $1+v^\varepsilon$, and Fubini turns it into a distributional identity:
\begin{align}
 \label{eq:green-Gardner-a1} \partial_t\log (1+v^\varepsilon) =\frac{\partial_x\Psi_\varepsilon}{1+v^\varepsilon} =\partial_x\!\left(\frac{\Psi_\varepsilon}{1+v^\varepsilon}\right) +\frac{\Psi_\varepsilon v_x^\varepsilon}{(1+v^\varepsilon)^2}.
\end{align}
Because $v\in C(I;H^1)$, the mollification $v^\varepsilon$ converges to $v$ uniformly in time in $H^1$, and hence uniformly in space. On compact spacetime sets, $\Psi_\varepsilon\to \Psi$ and $v_x^\varepsilon\to v_x$ strongly in $L^2$, while $(1+v^\varepsilon)^{\pm1}\to (1+v)^{\pm1}$ uniformly. The product term therefore converges in $L^1$, while $\log(1+v^\varepsilon)\to\log(1+v)$ uniformly, so that the left side of \eqref{eq:green-Gardner-a1} converges in distributions, and we obtain
\begin{align}
 \label{eq:log-chain} \partial_t\log (1+v)=\frac{\partial_x\Psi}{1+v} :=\partial_x\!\left(\frac{\Psi}{1+v}\right)+\frac{\Psi v_x}{(1+v)^2}
\end{align}
in distributions. Using $r=\tau v-\frac12\partial_x\log(1+v)$ from \eqref{eq:v-r-def}, together with \eqref{eq:v-conservation} and \eqref{eq:log-chain},
\begin{align}
 \label{eq:r-from-A} r_t=\partial_x\!\left[\tau \Psi-\frac{\Psi_x}{2(1+v)}\right].
\end{align}

\proofstep{Step 3: the flux identity for smooth functions}
It remains to identify the bracket in \eqref{eq:r-from-A} with the Gardner flux. This is first done algebraically, for smooth functions. Step~4 transports it to $v$ at the regularity already proved. Recall the logarithmic derivative $\mu$ of \eqref{eq:mu-def}, so that $r=\tau v-\mu/2$ by \eqref{eq:v-r-def}. In this notation the composite identity \eqref{eq:composite-v} reads
\begin{align}
 \label{eq:q-r-in-mu} q=-\frac12\mu_x+\frac14\mu^2+\tau^2\bigl(2v+v^2\bigr), \qquad r=\tau v-\frac{\mu}{2}.
\end{align}
Define the algebraic functional
\begin{align}
 \label{eq:algebraic-A} \cW(v) :=-(1+v)\mu_x+\frac{1+v}{2}\mu^2+2\tau^2(1+v)^3-6\tau^2(1+v)+4\tau^2.
\end{align}
At almost every $t$, in the common set of Remark~\ref{rem:measurability}, the identity \eqref{eq:q-r-in-mu} holds pointwise almost everywhere in $x$ as an identity of $L^2_{\loc}$ functions: it is \eqref{eq:Miura-pointwise}, which was established as an identity of functions in the proof of Lemma~\ref{lem:internal-H2}, with $r=\tau v-\mu/2$ and $v_x=(1+v)\mu$ substituted. At the same times $(1+v)q$ is the pointwise product, by \eqref{eq:hq-representative-agreement}. Substitution of \eqref{eq:q-r-in-mu} into \eqref{eq:Psi-def} then shows that $\cW(v)=\Psi$ in $L^2_{\loc}(I\times\R)$. For smooth $v$, the next two identities follow by differentiation and expansion:
\begin{align}
 \label{eq:Ax-h-p} \frac{\partial_x\cW(v)}{2(1+v)} =-\frac12\mu_{xx}+\frac14\mu^3+3\tau^2\mu\bigl(2v+v^2\bigr),
\end{align}
and direct expansion gives
\begin{align}
 \label{eq:smooth-flux-identity} &-r_{xx}+6\tau r^2+2r^3\\
 ={}&\frac12\mu_{xx}-\tau(1+v)\mu_x+\frac{\tau(1+v)}{2}\mu^2-\frac14\mu^3 -3\tau^2\mu\bigl(2v+v^2\bigr) +2\tau^3(1+v)^3-6\tau^3(1+v)+4\tau^3\notag \\
 ={}&\tau\cW(v)-\frac{\partial_x\cW(v)}{2(1+v)}.\notag
\end{align}

\proofstep{Step 4: passage to the available regularity}
With $v^\varepsilon$ as in \eqref{eq:v-mollification}, put
\begin{align}
 \label{eq:mu-W-mollified} \mu^\varepsilon:=\frac{v_x^\varepsilon}{1+v^\varepsilon}, \qquad \cW^\varepsilon:=\cW(v^\varepsilon).
\end{align}
The approximation $\cW^\varepsilon=\cW(v^\varepsilon)$ differs in general from $\Psi_\varepsilon=\eta_\varepsilon*\Psi$ because $\cW$ is nonlinear. The pair $(v^\varepsilon,\cW^\varepsilon)$ satisfies the algebraic identity \eqref{eq:smooth-flux-identity}, which holds for every smooth function. Step~2 used the mollified conservation law $v^\varepsilon_t=\partial_x\Psi_\varepsilon$ to prove the chain rule. Here we prove convergence of the recomputed flux $\cW^\varepsilon$ to $\Psi$ and pass to the limit in \eqref{eq:smooth-flux-identity}.

We use the regularity supplied by Corollary~\ref{cor:mu-regularity}:
\begin{align}
 \label{eq:mu-regularity-recalled} \begin{aligned}
 &\mu\in L^\infty(I;L^2),\qquad \mu,\mu_x\in \Lu(J\times\R),\\
 &\mu_x\in L^2(J\times\cK),\qquad \mu\in L^4(J\times\cK),
 \end{aligned}
\end{align}
for every compact $J\subset I$ and every bounded interval $\cK$.

\emph{The mollified function $v^\varepsilon$.} We spell out the approximation on compact spacetime cylinders. Fix a compact interval $J\subset I$ and bounded intervals
\begin{align}
 \label{eq:green-Gardner-a2} \cK\Subset\cK_1\Subset\cK_2.
\end{align}
Take the mollifier $\eta$ to be nonnegative and compactly supported. Since $t\mapsto v(t)$ maps the compact interval $I$ continuously into $H^1$, its image is compact in $H^1$. Uniform convergence of the approximate identity on compact subsets of $H^1$ therefore gives
\begin{align}
 \label{eq:v-mollifier-uniform-H1} \norm{v^\varepsilon-v}_{L^\infty(I;H^1)}\rightarrow0.
\end{align}
On the other hand, the local $H^2$ estimate \eqref{eq:v-final-reg} and the convergence of spatial mollification in $L^2(J\times\cK_2)$, applied to $\partial_x^\ell v\in L^2(J\times\cK_3)$ for $0\leq\ell\leq2$ and a bounded interval $\cK_3\supset\cK_2+\supp\eta_\varepsilon$, give
\begin{align}
 \label{eq:v-mollifier-local-H2} \norm{v^\varepsilon-v}_{L^2(J;H^2(\cK_2))}\rightarrow0.
\end{align}
Because $\eta_\varepsilon$ is nonnegative and has mass one, the lower bound for $1+v$ is preserved:
\begin{align}
 \label{eq:green-Gardner-a3} 1+v^\varepsilon(t,x)\geq c_I.
\end{align}
Moreover, \eqref{eq:v-mollifier-uniform-H1} and the one-dimensional Sobolev embedding imply
\begin{align}
 \label{eq:reciprocal-mollifier-convergence} \norm{v^\varepsilon-v}_{L^\infty(I\times\R)} +\norm*{\frac1{1+v^\varepsilon}-\frac1{1+v}}_{L^\infty(I\times\R)} \rightarrow0.
\end{align}

\emph{The quotient $\mu$ and its derivative.} From
\begin{align}
 \label{eq:mu-mollifier-difference} \mu^\varepsilon-\mu =\frac{v_x^\varepsilon-v_x}{1+v^\varepsilon} +v_x\,\frac{v-v^\varepsilon}{(1+v^\varepsilon)(1+v)},
\end{align}
the uniform reciprocal bound, \eqref{eq:v-mollifier-uniform-H1}, and \eqref{eq:reciprocal-mollifier-convergence} yield
\begin{align}
 \label{eq:mu-mollifier-Linfty-L2} \norm{\mu^\varepsilon-\mu}_{L^\infty(I;L^2)}\rightarrow0.
\end{align}

We next prove convergence of the derivative. Put $d^\varepsilon:=v^\varepsilon-v$. A localized one-dimensional Gagliardo--Nirenberg inequality gives
\begin{align}
 \label{eq:vx-mollifier-L4-estimate} \norm{d_x^\varepsilon}_{L^4(J\times\cK)}^4 \lesssim \norm{d_x^\varepsilon}_{L^\infty(J;L^2)}^3 \left( \abs{J}^{\frac12}\norm{d_{xx}^\varepsilon}_{L^2(J\times\cK_1)} +\abs{J}\norm{d_x^\varepsilon}_{L^\infty(J;L^2)} \right).
\end{align}
Indeed, apply the localized inequality $\norm{f}_4^4\lesssim\norm{f}_2^3 (\norm{f_x}_2+\norm{f}_2)$ to a cutoff equal to one on $\cK$, and then integrate in time. Equations \eqref{eq:v-mollifier-uniform-H1} and \eqref{eq:v-mollifier-local-H2} show that the right side of \eqref{eq:vx-mollifier-L4-estimate} tends to zero. Hence
\begin{align}
 \label{eq:vx-mollifier-L4} v_x^\varepsilon\rightarrow v_x \quad\text{in }L^4(J\times\cK).
\end{align}
The same localized inequality applied to $v_x$ shows that $v_x\in L^4(J\times\cK)$. Consequently,
\begin{align}
 \label{eq:vx-square-mollifier} (v_x^\varepsilon)^2\rightarrow v_x^2 \quad\text{in }L^2(J\times\cK).
\end{align}

Using
\begin{align}
 \label{eq:green-Gardner-a4} \mu_x^\varepsilon =\frac{v_{xx}^\varepsilon}{1+v^\varepsilon} -\frac{(v_x^\varepsilon)^2}{(1+v^\varepsilon)^2}, \qquad \mu_x =\frac{v_{xx}}{1+v}-\frac{v_x^2}{(1+v)^2},
\end{align}
we combine \eqref{eq:v-mollifier-local-H2}, \eqref{eq:reciprocal-mollifier-convergence}, and \eqref{eq:vx-square-mollifier} to obtain
\begin{align}
 \label{eq:mux-mollifier-convergence} \mu_x^\varepsilon\rightarrow \mu_x \quad\text{in }L^2(J\times\cK).
\end{align}
Thus, for every compact $J\subset I$,
\begin{align}
 \label{eq:mu-mollifier-convergence} \mu^\varepsilon\rightarrow \mu \quad\text{in } L^\infty(I;L^2)\cap L^2(J;H^1_{\loc}).
\end{align}

\emph{The quadratic term.} Apply the same localized interpolation estimate to $\mu^\varepsilon-\mu$. Equations \eqref{eq:mu-mollifier-Linfty-L2} and \eqref{eq:mux-mollifier-convergence} give
\begin{align}
 \label{eq:mu-mollifier-L4} \mu^\varepsilon\rightarrow \mu \quad\text{in }L^4(J\times\cK).
\end{align}
The limit $\mu$ lies in $L^4(J\times\cK)$ by \eqref{eq:mu-regularity-recalled}, and hence $(\mu^\varepsilon)^2\to \mu^2$ in $L^2(J\times\cK)$.

\emph{The algebraic flux.} Every term of \eqref{eq:algebraic-A} may now be passed to the limit. In the following two estimates, the $L^\infty$ norms are over $I\times\R$:
\begin{align}
 &\bnorm{(1+v^\varepsilon)\mu_x^\varepsilon-(1+v)\mu_x}_{L^2(J\times\cK)} \label{eq:green-Gardner-a5}\leq\norm{1+v^\varepsilon}_\infty\norm{\mu_x^\varepsilon-\mu_x}_{L^2(J\times\cK)}+\norm{v^\varepsilon-v}_\infty\norm{\mu_x}_{L^2(J\times\cK)}, 
\end{align}
\begin{align}
 &\bnorm{(1+v^\varepsilon)(\mu^\varepsilon)^2-(1+v)\mu^2}_{L^2(J\times\cK)} \label{eq:green-Gardner-square-estimate}\leq\norm{1+v^\varepsilon}_\infty\norm{(\mu^\varepsilon)^2-\mu^2}_{L^2(J\times\cK)}+\norm{v^\varepsilon-v}_\infty\norm{\mu^2}_{L^2(J\times\cK)}. 
\end{align}
Both right sides tend to zero. The polynomial terms converge in $L^2(J\times\cK)$ by the uniform convergence and uniform boundedness of $1+v^\varepsilon$. Therefore
\begin{align}
 \label{eq:Acal-convergence} \cW^\varepsilon\rightarrow\cW(v)=\Psi \quad\text{in }L^2_{\loc}(I\times\R).
\end{align}

Furthermore,
\begin{align}
 \label{eq:green-Gardner-a6} \frac{\partial_x\cW^\varepsilon}{1+v^\varepsilon} =\partial_x\!\left(\frac{\cW^\varepsilon}{1+v^\varepsilon}\right) +\frac{\cW^\varepsilon v_x^\varepsilon}{(1+v^\varepsilon)^2}.
\end{align}
The first term converges distributionally because $\cW^\varepsilon/(1+v^\varepsilon)\to \Psi/(1+v)$ in $L^2_{\loc}$. The second converges in $L^1_{\loc}$, since $\cW^\varepsilon\to \Psi$ and $v_x^\varepsilon\to v_x$ in $L^2_{\loc}$ while the reciprocals converge uniformly. Hence
\begin{align}
 \label{eq:Acal-derivative-convergence} \frac{\partial_x\cW^\varepsilon}{1+v^\varepsilon} \rightarrow \partial_x\!\left(\frac{\Psi}{1+v}\right)+\frac{\Psi v_x}{(1+v)^2} =\frac{\Psi_x}{1+v} \quad\text{in distributions},
\end{align}
where $\Psi_x/(1+v)$ is understood as in \eqref{eq:log-chain}. Halving it gives the coefficient $\Psi_x/(2(1+v))$ that appears in \eqref{eq:r-from-A} and \eqref{eq:smooth-flux-identity}.

\emph{The Gardner variable and its powers.} Finally, let
\begin{align}
 \label{eq:green-Gardner-a7} r^\varepsilon:=\tau v^\varepsilon-\frac12\mu^\varepsilon .
\end{align}
The preceding estimates imply
\begin{align}
 \label{eq:r-mollifier-H1} r^\varepsilon-r\rightarrow0 \quad\text{in } L^\infty(I;L^2)\cap L^2(J;H^1_{\loc}).
\end{align}
To carry the cubic flux in $L^2_{\loc}$, apply
\begin{align}
 \label{eq:green-Gardner-a8} \norm{f}_{L^6(\cK)}^6 \lesssim\norm{f}_{L^2(\cK_1)}^4 \bigl(\norm{f_x}_{L^2(\cK_1)}^2+\norm{f}_{L^2(\cK_1)}^2\bigr)
\end{align}
to $f=r^\varepsilon-r$ and integrate in time. This gives
\begin{align}
 \label{eq:green-Gardner-a9} \norm{r^\varepsilon-r}_{L^6(J\times\cK)}^6 \lesssim \norm{r^\varepsilon-r}_{L^\infty(J;L^2)}^4 \left( \norm{(r^\varepsilon-r)_x}_{L^2(J\times\cK_1)}^2 +\abs{J}\norm{r^\varepsilon-r}_{L^\infty(J;L^2)}^2 \right),
\end{align}
and hence
\begin{align}
 \label{eq:r-mollifier-L6} r^\varepsilon\rightarrow r \quad\text{in }L^6_{\loc}(I\times\R).
\end{align}
Consequently,
\begin{align}
 \label{eq:green-Gardner-a10} (r^\varepsilon)^2&\rightarrow r^2 \quad\text{and}\quad (r^\varepsilon)^3\rightarrow r^3 \qquad\text{in }L^2_{\loc}(I\times\R)
\end{align}
where the first convergence follows first in $L^3_{\loc}$ and then in $L^2_{\loc}$ on bounded cylinders. Also, $r_{xx}^\varepsilon\to r_{xx}$ in $L^2(J;H^{-1}_{\loc})$ by \eqref{eq:r-mollifier-H1}. We may therefore apply \eqref{eq:smooth-flux-identity} to $v^\varepsilon$ and pass to the limit, obtaining
\begin{align}
 \label{eq:rough-flux-identity} \tau \Psi-\frac{\Psi_x}{2(1+v)} =-r_{xx}+6\tau r^2+2r^3 \quad\text{in distributions}.
\end{align}
Inserting this identity in \eqref{eq:r-from-A} proves \eqref{eq:Gardner} in spacetime distributions.

\proofstep{Step 5: the regularity demanded by Definition~\ref{def:Gardner-solution}}
It remains to check that $r$ lies in the class in which \eqref{eq:Gardner} was formulated, so that the nonlinear flux is meaningful. We have $r\in L^\infty_tL^2_x$ by \eqref{eq:v-r-basic-reg} and $r_x=\tau v_x-\frac12\mu_x\in L^2_{\loc}$ by Corollary~\ref{cor:mu-regularity} (and again by Proposition~\ref{prop:r-X1-local}), which are exactly the hypotheses under which \eqref{eq:Gardner-local-L6} was derived. That estimate therefore gives $r^2,r^3\in L^2_{\loc}$, and $r_{xx}\in L^2_{\loc,t}H^{-1}_{\loc,x}$ follows from $r_x\in L^2_{\loc}$. This is the weak-solution regularity in Definition~\ref{def:Gardner-solution}. Integrating the distributional equation against a test function gives exactly \eqref{eq:weak-Gardner-formulation}, completing the transfer.
\end{proof}

%%%%%%%%%%%%%%%%%%%%%%%%%%%%%%%%%%%%%%%%%%%%%%%%%%%%%%%%%%%%%%%%%%%%%%%%%%%%%%%%%%
\subsection{Globalization and Proof of Main Theorem}
\label{sec:uniqueness}
%%%%%%%%%%%%%%%%%%%%%%%%%%%%%%%%%%%%%%%%%%%%%%%%%%%%%%%%%%%%%%%%%%%%%%%%%%%%%%%%%%

The final stage extends the time interval to $\R$. By Proposition~\ref{prop:r-X1-local} and Proposition~\ref{prop:green-Gardner} $r$ is a Gardner solution on $(-T,T)$ obeying the smoothing and tail conditions of $X_1$ on every compact subinterval, whereas $X_1$ itself, and with it the uniqueness assertion of Lemma~\ref{lem:KKL-Gardner}, is a global-in-time class. Each of the three conditions defining $X_1$ quantifies over the whole time axis: $r\in L^\infty(\R_t;L^2_x)$, the smoothing bound \eqref{eq:X1-bounds} on every compact $J\subset\R$, and the tail \eqref{eq:X1-tail} at every $t_0\in\R$. Applying Lemma~\ref{lem:KKL-Gardner} therefore requires an extension to all times satisfying these three conditions.

We extend the Green solution $q$ itself beyond $I$, before the lift is taken, by the canonical flow of Killip--Vi\c{s}an \cite{Killip.Visan2019}. The Green identities on the adjacent intervals combine at $t=\pm T$ by additivity of their time integrals. What has to be checked is that the thresholds of the extension on $[-L,L]$ do not grow with $L$, so that a single parameter $\tau$ serves all times. This is where it matters that the thresholds \eqref{eq:static-threshold} and \eqref{eq:high-energy-threshold} depend on a compact set only through the largest $H^{-1}$ norm of its elements. The lift of the extension is then bounded in $L^2$ uniformly in time by Lemma~\ref{lem:Gardner-mass}.

{\begin{lemma}[Extension by the canonical flow]\label{lem:global-extension}
Let $q\in C(I;H^{-1}(\R;\R))$ be a Green solution, let $\Phi_{\mathrm{KdV}}$ be the flow of Lemma~\ref{lem:canonical-is-green}, and define
\begin{align}
 \label{eq:global-extension} \widetilde q(t):=\begin{cases}
 \Phi_{\mathrm{KdV}}\bigl(t+T,q(-T)\bigr),&t<-T,\\
 q(t),&-T\leq t\leq T,\\
 \Phi_{\mathrm{KdV}}\bigl(t-T,q(T)\bigr),&t>T.
 \end{cases}
\end{align}
Then $\widetilde q\in C(\R;H^{-1}(\R;\R))$ and
\begin{align}
 \label{eq:global-extension-bound} \sup_{t\in\R}\,\norm{\widetilde q(t)}_{H^{-1}}<\infty .
\end{align}
Moreover, let $R$ be at least the supremum in \eqref{eq:global-extension-bound}, and let $\delta$ be the absolute constant of Proposition~\ref{prop:static}. Then, for every $L\geq T$, the restriction of $\widetilde q$ to $[-L,L]$ is a Green solution in the sense of Definition~\ref{def:green}, with the threshold
\begin{align}
 \label{eq:global-extension-threshold} \lambda_*\bigl(\widetilde q,[-L,L]\bigr):=\max\left\{\lambda_*(q,I),\ \frac{4(R+1)^2}{\delta^2}\right\},
\end{align}
which does not depend on $L$.
\end{lemma}}

\begin{proof}
Since $\Phi_{\mathrm{KdV}}$ is jointly continuous and $\Phi_{\mathrm{KdV}}(0,p)=p$ for every $p\in H^{-1}(\R;\R)$, the three pieces of \eqref{eq:global-extension} are continuous and agree at $t=\pm T$, so $\widetilde q\in C(\R;H^{-1})$. The supremum in \eqref{eq:global-extension-bound} is finite because $q(I)$ is compact and because every orbit of $\Phi_{\mathrm{KdV}}$ is bounded in $H^{-1}(\R)$, by \cite[Theorem~1.1]{Killip.Visan2019}.

\proofstep{Step 1: the identity beyond $T$}
Choose Schwartz functions $p_n\to q(T)$ in $H^{-1}$, and let $q_n(s):=\Phi_{\mathrm{KdV}}(s,p_n)$ be the corresponding Schwartz solutions of \eqref{eq:KdV}. Fix $S>0$. Joint continuity of $\Phi_{\mathrm{KdV}}$ on the compact set $[0,S]\times(\{p_n:n\geq1\}\cup\{q(T)\})$ gives $q_n\to\widetilde q(T+\cdot\,)$ in $C([0,S];H^{-1})$, so there is $n_S$ such that $\norm{q_n(s)}_{H^{-1}}\leq R+1$ for all $n\geq n_S$ and $s\in[0,S]$. Hence
\begin{align}
 \label{eq:global-extension-compact} \cQ:=\widetilde q([T,T+S])\cup\bigcup_{n\geq n_S}q_n([0,S])
\end{align}
is compact, being the image of a compact set under a continuous map, and all its elements have $H^{-1}$ norm at most $R+1$. By \eqref{eq:static-threshold} and $\delta\leq\frac12$, this gives $\lambda_0(\cQ)\leq4(R+1)^2/\delta^2$. Lemma~\ref{lem:microscopic}\,\textup{(ii)} applies to each $q_n$ on $[0,S]$, because the law \eqref{eq:microscopic-h} holds at each time along a Schwartz solution whose trajectory is admissible. The argument of Lemma~\ref{lem:canonical-is-green}, run on the time interval $[0,S]$ with \eqref{eq:global-extension-compact} in place of \eqref{eq:canonical-is-green-compact}, therefore gives
\begin{align}
 \label{eq:ext-right} h_\lambda(\widetilde q(t))-h_\lambda(q(T)) =\partial_x\int_T^t\J_\lambda(\widetilde q(t'))\dd t' \quad\text{in }H^{-2}(\R)
\end{align}
for every $t\in[T,T+S]$ and every $\lambda\geq4(R+1)^2/\delta^2$. Since $S$ is arbitrary, \eqref{eq:ext-right} holds for every $t\geq T$.

\proofstep{Step 2: the identity before $-T$, and concatenation}
The same argument, with $q(-T)$ in place of $q(T)$ and the time interval $[-S,0]$ in place of $[0,S]$, gives
\begin{align}
 \label{eq:ext-left} h_\lambda(\widetilde q(t))-h_\lambda(q(-T)) =\partial_x\int_{-T}^t\J_\lambda(\widetilde q(t'))\dd t' \quad\text{in }H^{-2}(\R)
\end{align}
for every $t\leq-T$ and every $\lambda\geq4(R+1)^2/\delta^2$. Now let $\lambda$ be at least the right side of \eqref{eq:global-extension-threshold}. For $t\in I$ the Green identity \eqref{eq:green-identity-recalled} of $q$ applies as it stands. For $t>T$, adding \eqref{eq:ext-right} to that identity at time $T$ gives
\begin{align}
 \label{eq:ext-concatenate} h_\lambda(\widetilde q(t))-h_\lambda(\widetilde q(0)) =\partial_x\int_0^t\J_\lambda(\widetilde q(t'))\dd t' \quad\text{in }H^{-2}(\R),
\end{align}
because Bochner integrals over adjacent time intervals add. For $t<-T$ the same follows from \eqref{eq:ext-left} and the identity of $q$ at time $-T$. Thus \eqref{eq:ext-concatenate} holds for every $t\in\R$. Finally, all values of $\widetilde q$ have $H^{-1}$ norm at most $R$, so \eqref{eq:static-threshold} gives $\lambda_0(\widetilde q([-L,L]))\leq4(R+1)^2/\delta^2$, and the threshold \eqref{eq:global-extension-threshold} meets the requirement of Definition~\ref{def:green} that it be at least $\lambda_0(\widetilde q([-L,L]))$. Hence the restriction of $\widetilde q$ to $[-L,L]$ is a Green solution with that threshold.
\end{proof}

\begin{proof}[Proof of Theorem~\ref{thm:main}]
Let $q_1,q_2\in C(I;H^{-1}(\R;\R))$ be Green solutions with $q_1(0)=q_2(0)$.

\proofstep{Step 1: extension and a common parameter}
Let $\widetilde q_1$ and $\widetilde q_2$ be the extensions of $q_1$ and $q_2$ given by Lemma~\ref{lem:global-extension}, and let $R$ be at least the supremum \eqref{eq:global-extension-bound} for both. For every $L\geq T$, the restrictions of $\widetilde q_1$ and $\widetilde q_2$ to $[-L,L]$ are Green solutions whose values have $H^{-1}$ norm at most $R$, with the thresholds \eqref{eq:global-extension-threshold}. Form the pair threshold $K_{12}$ of \eqref{eq:common-K} from these thresholds and this $R$. It does not depend on $L$, and we fix one $\tau\geq K_{12}$. By \eqref{eq:K-monotone}, for every $L\geq T$, every statement of Section~\ref{sec:Gardner} proved for a single Green solution at $\tau\geq K(q)$ is available for both restrictions to $[-L,L]$ at this same $\tau$. Since $K_{12}\geq2$, the parameter $\tau$ also satisfies the restriction of Section~\ref{sec:Gardner-class}.

\proofstep{Step 2: two global Gardner solutions}
Define $\widetilde v_j$ and $\widetilde r_j$ from $\widetilde q_j$ by \eqref{eq:v-r-def}, at every $t\in\R$. The value $\widetilde r_j(t)$ depends only on $\widetilde q_j(t)$ and $\tau$, so for every $L\geq T$ the restriction of $\widetilde r_j$ to $[-L,L]$ is the variable that Section~\ref{sec:Gardner-lift} attaches to the restriction of $\widetilde q_j$. The results of Sections~\ref{sec:Gardner-lift}--\ref{sec:Gardner-dynamics}, applied on $[-L,L]$ for every $L\geq T$, give the following. By \eqref{eq:v-r-basic-reg}, $\widetilde r_j\in C(\R;L^2(\R;\R))$. By Lemma~\ref{lem:Gardner-mass}, $\norm{\widetilde r_j(t)}_2=\norm{\widetilde r_j(0)}_2$ for every $t\in\R$, so that $\widetilde r_j\in L^\infty(\R;L^2)$. By Proposition~\ref{prop:green-Gardner}, $\widetilde r_j$ is a Gardner solution on $(-L,L)$ for every $L$, hence on $\R$, since every test function in \eqref{eq:weak-Gardner-formulation} is supported in $(-L,L)\times\R$ for some $L$. By Proposition~\ref{prop:r-X1-local}, applied with $L$ so large that a given compact interval $J\subset\R$ lies in $[-L,L]$, the derivative $\partial_x\widetilde r_j$ satisfies \eqref{eq:X1-bounds} on $J$ and \eqref{eq:X1-tail} on every unit window contained in $J$. Thus $\widetilde r_1$ and $\widetilde r_2$ are global Gardner solutions of class $X_1$, and they are continuous into $L^2$.

\proofstep{Step 3: equal data}
The two lifts agree at time zero. Indeed \eqref{eq:v-r-def} depends only on $q_j(0)$ and the parameter $\tau$, which is the same for both by Step~1, so $q_1(0)=q_2(0)$ forces $\widetilde r_1(0)=\widetilde r_2(0)$. (Injectivity of $M_\tau$ on $L^2$, Remark~\ref{rem:shifted-Miura-injective}, gives the same conclusion from $M_\tau(\widetilde r_1(0))=q_1(0)=q_2(0)=M_\tau(\widetilde r_2(0))$.)

\proofstep{Step 4: uniqueness}
By Steps~2 and~3, Lemma~\ref{lem:KKL-Gardner} applies to $\widetilde r_1$ and $\widetilde r_2$ and gives $\widetilde r_1(t)=\widetilde r_2(t)$ for every $t\in\R$. Since $\widetilde q_j=q_j$ on $I$, the restriction $r_j$ of $\widetilde r_j$ to $I$ is the lift of $q_j$ itself, and
\begin{align}
 \label{eq:r-equal} r_1=r_2\qquad\text{on }I .
\end{align}

\proofstep{Step 5: equality of the KdV solutions}
Apply $M_\tau$ to \eqref{eq:r-equal}. By Lemma~\ref{lem:Miura},
\begin{align}
 \label{eq:gluing-a5} q_1-q_2 =\partial_x(r_1-r_2)+2\tau(r_1-r_2)+(r_1^2-r_2^2)=0
\end{align}
in $H^{-1}(\R)$ at every $t\in I$. Hence $q_1\equiv q_2$ on $I$.
\end{proof}

\begin{proof}[Proof of Corollary~\ref{cor:all-s}]
For $s\geq-1$, the continuous embedding $H^s\hookrightarrow H^{-1}$ puts both competitors in the class of Theorem~\ref{thm:main}.
\end{proof}

\begin{proof}[Proof of Corollary~\ref{cor:canonical}]
By \cite[Theorem~1.1]{Killip.Visan2019}, the Schwartz KdV solution map extends uniquely to a jointly continuous global map
\begin{align}
 \label{eq:gluing-a6} \Phi_{\mathrm{KdV}}\colon\R\times H^{-1}(\R)\rightarrow H^{-1}(\R).
\end{align}
By Lemma~\ref{lem:canonical-is-green}, for every $q_0\in H^{-1}$ the canonical trajectory $t\mapsto\Phi_{\mathrm{KdV}}(t,q_0)$, restricted to $I$, satisfies Definition~\ref{def:green}. (The same inclusion is stated in the paragraph following equation~(7.17) of \cite{Killip.Murphy.Visan2020}; the equivalence of the two flux normalizations was verified in \eqref{eq:KMV-normalization-a1}.)

Now let $q$ be any Green solution with $q(0)=q_0$. Applying Theorem~\ref{thm:main} to $q$ and $\Phi_{\mathrm{KdV}}(\,\cdot\,,q_0)|_I$ gives
\begin{align}
 \label{eq:gluing-a7} q(t)=\Phi_{\mathrm{KdV}}(t,q_0) \qquad\text{for every }t\in I.
\end{align}
\end{proof}

% \newpage

%%%%%%%%%%%%%%%%%%%%%%%%%%%%%%%%%%%%%%%%%%%%%%%%%%%%%%%%%%%%%%%%%%%%%%%%%%%%%%%%%%
\section{Defocusing mKdV and the Two Miura Transforms}
\label{sec:mkdv}
%%%%%%%%%%%%%%%%%%%%%%%%%%%%%%%%%%%%%%%%%%%%%%%%%%%%%%%%%%%%%%%%%%%%%%%%%%%%%%%%%%

We now turn to the real defocusing modified Korteweg--de Vries equation
{
\begin{align}
 \label{eq:defocusing-mkdv} u_t=-u_{xxx}+6u^2u_x =-u_{xxx}+2\partial_x(u^3), \qquad u(0)=u_0\in L^2(\R;\R),
\end{align}
}
and prove Theorem~\ref{thm:mkdv-main} and Proposition~\ref{prop:mkdv-smoothing}. We transfer an mKdV trajectory to two KdV Green solutions in $H^{-1}$ and apply the uniqueness theorem proved in Part~I.

We use the Dirac-resolvent Green notion of Forlano--Killip--Vi\c{s}an \cite[Definition~3.11]{Forlano.Killip.Visan2024}, restricted to curves in $C(I;L^2)$ and stated in Definition~\ref{def:mkdv-green}. Its defining conditions are the integrated evolution identities for $g_{-,\kappa}$ and $\gamma_\kappa$ at every $\kappa\geq1$. We prove the implication
\begin{align}
 \label{eq:mkdv-task} \begin{aligned}[t]
 &\text{the two mKdV Green identities for }u\\
 &\qquad\Longrightarrow\quad \text{the KdV Green identity of Definition~\ref{def:green} for \emph{both} } q_\pm=u^2\pm u_x .
 \end{aligned}
\end{align}
Once \eqref{eq:mkdv-task} is available, Theorem~\ref{thm:main} applies once for each sign and the difference $q_+-q_-=2u_x$ recovers $u$. Section~\ref{sec:mkdv-notion} verifies the meaning of the defining laws and explains the multiplier obstruction in proving \eqref{eq:mkdv-task} at $L^2$ regularity. The three stages in which \eqref{eq:mkdv-task} is then proved are set out at the end of that subsection. 
% Figure~\ref{fig:dependence-mkdv} shows the resulting dependence graph.

%%%%%%%%%%%%%%%%%%%%%%%%%%%%%%%%%%%%%%%%%%%%%%%%%%%%%%%%%%%%%%%%%%%%%%%%%%%%%%%%%%
\subsection{Properties of the Miura Maps}
\label{sec:mkdv-static}
%%%%%%%%%%%%%%%%%%%%%%%%%%%%%%%%%%%%%%%%%%%%%%%%%%%%%%%%%%%%%%%%%%%%%%%%%%%%%%%%%%

We use the maps $\cM_\sigma$ and the potentials $q_\sigma$ defined in \eqref{eq:intro-miura-partners}. This subsection proves their continuity, factorization, and intertwining properties.

{
\begin{remark}
We use both signs for two reasons.

The first observation is that both transforms arise from the same first-order operator. Write $C$ for this operator and $C^*$ for its formal adjoint. The two Miura transforms occur as potentials when $C$ and $C^*$ are composed in the two possible orders:
\begin{align}
 \label{eq:Miura-supersymmetry} \begin{aligned}
 &C:=\partial_x+u , \qquad C^*=-\partial_x+u ,\\
 &L_{\cM_-(u)}=C^*C , \qquad L_{\cM_+(u)}=CC^* .
 \end{aligned}
\end{align}
Parts~\textup{(ii)} and~\textup{(iii)} of the lemma below are \eqref{eq:Miura-supersymmetry} made precise: the first reads it at the level of quadratic forms, which is what makes both operators nonnegative, and the second identifies the operator associated with the factorized form. The matrix operator behind \eqref{eq:Miura-supersymmetry}, from which the pair is usually read off, is the mKdV Lax operator \eqref{eq:intro-dirac} of Section~\ref{sec:mkdv-notion}.

The second reason is that the pair determines $u$: their difference is $\cM_+(u)-\cM_-(u)=2u_x$, and an $L^2$ function is determined by its distributional derivative. Remark~\ref{rem:one-Miura-insufficient} gives an example showing that $\cM_+$ is not injective on $L^2(\R)$.
\end{remark}
}

{\begin{lemma}
\label{lem:two-Miura-static}
Let $\sigma\in\{+1,-1\}$ and let $u,v\in L^2(\R;\R)$.
\begin{enumerate}[
label=\textup{(\roman*)}, leftmargin=*, itemsep=.7em, before={\let\fullwidthdisplay\relax} ]
\item \emph{Mapping property and Lipschitz bound.} The map \eqref{eq:intro-miura-partners} sends $L^2$ into $H^{-1}$,
\begin{align}
 \label{eq:two-Miura-static-a1} \cM_\sigma\colon L^2(\R;\R)\rightarrow H^{-1}(\R;\R),
\end{align}
and is Lipschitz on bounded subsets:
\begin{align}
 \label{eq:two-Miura-Lipschitz} \norm{\cM_\sigma(u)-\cM_\sigma(v)}_{H^{-1}} \lesssim \bigl(1+\norm{u}_2+\norm{v}_2\bigr)\norm{u-v}_2.
\end{align}

\item \emph{The Schr\"odinger form is a perfect square.} The quadratic form \eqref{eq:intro-form} of the potential $\cM_\sigma(u)$ factorizes: for every $\phi\in H^1(\R)$,
\begin{align}
 \label{eq:Miura-form-factorization} \mathfrak l_{\cM_\sigma(u)}[\phi] =\norm{\phi_x-\sigma u\phi}_2^2 .
\end{align}
In particular $L_{\cM_\sigma(u)}\geq0$, so that $L_{\cM_\sigma(u)}+\kappa^2$ is boundedly invertible for every $\kappa>0$.

\item \emph{Equality of the form domains.} The right side of \eqref{eq:Miura-form-factorization} is defined on the maximal first-order domain
\begin{align}
 \label{eq:Miura-maximal-domain} \Dom_\sigma(u) :=\bigl\{\phi\in L^2(\R)\cap H^1_{\loc}(\R): \phi_x-\sigma u\phi\in L^2(\R)\bigr\} \ \supseteq\ H^1(\R) ,
\end{align}
and this domain equals the form domain in \eqref{eq:intro-form}:
\begin{align}
 \label{eq:Miura-domain-equality} \Dom_\sigma(u)=H^1(\R) .
\end{align}
Hence the factorized description \eqref{eq:Miura-form-factorization} and the defining form \eqref{eq:intro-form} produce one and the same operator $L_{\cM_\sigma(u)}$, and either may be used for it.

\item \emph{Intertwining with \eqref{eq:KdV}.} The identities below are used in Lemma~\ref{lem:mkdv-smooth-identification}. For smooth spacetime functions $u$ put
\begin{align}
 \label{eq:two-Miura-static-a2} \mathcal E_{\mathrm{mKdV}}(u):=u_t+u_{xxx}-6u^2u_x ,
\end{align}
so that $\mathcal E_{\mathrm{mKdV}}(u)=0$ is \eqref{eq:defocusing-mkdv}. Then
\begin{align}
 \label{eq:classical-Miura-intertwining} \begin{aligned}
 &\bigl(\cM_+(u)\bigr)_t +\bigl(\cM_+(u)\bigr)_{xxx} -6\cM_+(u)\bigl(\cM_+(u)\bigr)_x =(\partial_x+2u)\mathcal E_{\mathrm{mKdV}}(u),\\
 &\bigl(\cM_-(u)\bigr)_t +\bigl(\cM_-(u)\bigr)_{xxx} -6\cM_-(u)\bigl(\cM_-(u)\bigr)_x =(-\partial_x+2u)\mathcal E_{\mathrm{mKdV}}(u).
 \end{aligned}
\end{align}
In particular both Miura maps carry smooth solutions of \eqref{eq:defocusing-mkdv} to smooth solutions of \eqref{eq:KdV}.
\end{enumerate}
\end{lemma}}

Parts \textup{(i)}--\textup{(iii)} place the Miura transforms in the framework of Section~\ref{sec:solution-notion}: each potential belongs to $H^{-1}$, its Schr\"odinger operator is nonnegative, and the factorized and form-defined operators agree.

Part~\textup{(iv)} gives the classical intertwining identity for smooth $u$. Step~2 of the proof of Lemma~\ref{lem:mkdv-smooth-identification} uses it to show that each Miura transform of a smooth mKdV solution solves \eqref{eq:KdV}. Sections~\ref{sec:mkdv-transfer} and~\ref{sec:mkdv-identification} establish the transfer for the rough trajectories of Theorem~\ref{thm:mkdv-main} through the Green identities.

Nonnegativity of $L_{\cM_\sigma(u)}$ makes $-\kappa^2$ a resolvent point at every $\kappa>0$. Admissibility in the sense of Proposition~\ref{prop:static} also requires the quantitative bound \eqref{eq:g-quantitative}, which can fail at a given $\kappa$ even for a Miura transform of an $L^2$ function, as shown in Remark~\ref{rem:positivity-not-admissibility}. Proposition~\ref{prop:mkdv-transfer} places $\cM_\sigma(u)$ in the class of Definition~\ref{def:green} above the threshold $K_u$ of \eqref{eq:mkdv-common-threshold}.

{
\begin{remark}
Part~\textup{(iv)} is a well-known fact for real defocusing mKdV, which is also the main driving force in this section. With uniqueness for KdV Green solutions established in the preceding sections, our goal is to prove that every real defocusing mKdV Green solution induces two KdV Green solutions. Their uniqueness then determines the mKdV trajectory through \eqref{eq:two-Miura-recovery}.
\end{remark}
}

\begin{proof}
The four parts are proved in turn. Part~\textup{(iii)} uses a cutoff estimate to identify the form domain.

\proofstep{Part \textup{(i)}}
The map $\partial_x\colon L^2\to H^{-1}$ is bounded, and $L^1(\R)\hookrightarrow H^{-1}(\R)$ by Lemma~\ref{lem:products}. Since $u^2-v^2=(u-v)(u+v)\in L^1$ for $u,v\in L^2$,
\begin{align}
 \label{eq:two-Miura-static-a3} \norm{\cM_\sigma(u)-\cM_\sigma(v)}_{H^{-1}} \lesssim \norm{u-v}_2+\norm{(u-v)(u+v)}_1 \leq\bigl(1+\norm{u}_2+\norm{v}_2\bigr)\norm{u-v}_2 ,
\end{align}
the last step by Cauchy--Schwarz. Taking $v=0$ gives \eqref{eq:two-Miura-static-a1}.

\proofstep{Part \textup{(ii)}}
Let $\phi\in H^1(\R)$. Then $u\phi\in L^2$, because $\phi\in L^\infty$ by the one-dimensional Sobolev embedding, and
\begin{align}
 \label{eq:two-Miura-static-a4} \pair*{u_x}{\abs{\phi}^2} =-2\operatorname{Re}\int_\R u\phi_x\overline\phi\dd x ,
\end{align}
the pairing being between $H^{-1}$ and $H^1$ and the right side finite for the same reason. Expanding the square,
\begin{align}
 \label{eq:two-Miura-square-expansion} \norm{\phi_x-\sigma u\phi}_2^2 =\norm{\phi_x}_2^2 -2\sigma\operatorname{Re}\int_\R u\phi_x\overline\phi\dd x +\int_\R u^2\abs{\phi}^2\dd x ,
\end{align}
and substituting \eqref{eq:two-Miura-static-a4} turns the middle term into $\sigma\pair*{u_x}{\abs{\phi}^2}$, the pairing being the one in \eqref{eq:intro-form}. The right side of \eqref{eq:two-Miura-square-expansion} is therefore $\norm{\phi_x}_2^2+\pair*{u^2+\sigma u_x}{\abs{\phi}^2}$, which is $\mathfrak l_{\cM_\sigma(u)}[\phi]$, the closed form constructed in Lemma~\ref{lem:form-bound}. Nonnegativity of the associated self-adjoint operator follows from \eqref{eq:Miura-form-factorization}.

\proofstep{Part \textup{(iii)}}
The inclusion $H^1\subset\Dom_\sigma(u)$ follows because $\phi\in H^1$ gives $\phi_x\in L^2$ and $u\phi\in L^2$ as in part~\textup{(ii)}. For the converse, take $\phi\in\Dom_\sigma(u)$ and set $\varrho:=\phi_x-\sigma u\phi\in L^2$. The point is to show $\phi_x\in L^2$, and the difficulty is that the term $u\phi$ is controlled only after $\phi$ is known to be bounded, which is what we are trying to prove. The following cutoff estimate bounds the derivative and the product together.

Choose cutoffs $\chi_R$ equal to one on $[-R,R]$, supported in $[-2R,2R]$, with $\norm{\chi_R'}_\infty\lesssim R^{-1}$, and put $Y_R:=\norm{(\chi_R\phi)_x}_2$ for $R\geq1$. From $(\chi_R\phi)_x=\chi_R\varrho+\sigma\chi_Ru\phi+\chi_R'\phi$ and the one-dimensional Sobolev inequality $\norm{\chi_R\phi}_\infty\lesssim\norm{\chi_R\phi}_2^{1/2}Y_R^{1/2}$,
\begin{align}
 \label{eq:two-Miura-static-a5} Y_R \leq\norm{\varrho}_2+C\norm{\phi}_2+\norm{u}_2\norm{\chi_R\phi}_\infty \leq\norm{\varrho}_2+C\norm{\phi}_2 +C\norm{u}_2\norm{\phi}_2^{\frac12}Y_R^{\frac12} .
\end{align}
This is a scalar quadratic inequality of the form $Y\leq A+BY^{1/2}$ with $A$ and $B$ independent of $R$, which forces $Y_R\leq 2A+B^2$ and so bounds $Y_R$ uniformly. Since $\chi_R=1$ on $[-R,R]$, monotone convergence gives $\phi_x\in L^2(\R)$, that is, $\phi\in H^1(\R)$.

\proofstep{Part \textup{(iv)}}
The first identity in \eqref{eq:classical-Miura-intertwining} is a direct differentiation, all terms being smooth. The second follows from the first on replacing $u$ by $-u$, which exchanges $\cM_+$ and $\cM_-$ by \eqref{eq:intro-miura-partners} and changes the sign of $\partial_x+2u$ accordingly. If $u$ solves \eqref{eq:defocusing-mkdv} then $\mathcal E_{\mathrm{mKdV}}(u)=0$, so both right sides vanish and both images solve \eqref{eq:KdV}.
\end{proof}

{\begin{remark}
\label{rem:one-Miura-insufficient}
The map $\cM_+$ is not injective on $L^2(\R)$, so a single Miura transform need not determine $u$. Indeed, set
\begin{align}
 \label{eq:one-Miura-insufficient-a1} u(x)=\frac{x+\frac12}{1+x^2}, \qquad v(x)=\frac{x-\frac12}{1+x^2}.
\end{align}
Then $u,v\in L^2$, $u\neq v$, and, with $\varpi=u-v=(1+x^2)^{-1}$,
\begin{align}
 \label{eq:one-Miura-insufficient-a2} \cM_+(u)-\cM_+(v) =\varpi_x+(u+v)\varpi=0.
\end{align}
The pair satisfies
\begin{align}
 \label{eq:two-Miura-recovery} \cM_+(u)-\cM_-(u)=2u_x, \qquad \cM_+(u)+\cM_-(u)=2u^2.
\end{align}
The first identity determines $u$ from the pair, since an $L^2$ function with zero distributional derivative vanishes.

Definition~\ref{def:mkdv-green} imposes an identity for each of the two diagonal Green functions of the Dirac resolvent \eqref{eq:intro-dirac}. Proposition~\ref{prop:mkdv-transfer} uses them to prove the KdV Green identity for $\cM_\sigma(u)$ at each sign $\sigma$, so that KdV uniqueness can be applied to both terms in \eqref{eq:two-Miura-recovery}.
\end{remark}}

Remark~\ref{rem:shifted-Miura-injective} proves that the shifted map $M_\tau$ is injective for $\tau>0$.

%%%%%%%%%%%%%%%%%%%%%%%%%%%%%%%%%%%%%%%%%%%%%%%%%%%%%%%%%%%%%%%%%%%%%%%%%%%%%%%%%%
\subsection{The mKdV Green Solution Notion}
\label{sec:mkdv-notion}
%%%%%%%%%%%%%%%%%%%%%%%%%%%%%%%%%%%%%%%%%%%%%%%%%%%%%%%%%%%%%%%%%%%%%%%%%%%%%%%%%%

We now return to Definition~\ref{def:mkdv-green}. We recall the diagonal Green functions and the two evolution laws, check the meaning of their terms, and explain how they will imply the scalar Green law for both Miura transforms.

For $u\in L^2(\R;\R)$, recall the operator and resolvent from \eqref{eq:intro-dirac}:
\begin{align}
 \label{eq:AKNS-operator-recalled} \bL_u=\begin{pmatrix}-\partial_x&u\\-u&\partial_x\end{pmatrix}, \qquad \bR_\kappa(u)=(\bL_u+\kappa)^{-1}, \qquad \kappa\geq1.
\end{align}

We take $\bL_u$ with its maximal domain
\begin{align}
 \label{eq:AKNS-maximal-domain} \Dom(\bL_u) :=\left\{ \binom{\phi}{\psi}\in L^2(\R;\mathbb C^2): \ \phi,\psi\in H^1_{\loc}(\R),\ u\psi-\phi_x,\ \psi_x-u\phi\in L^2(\R) \right\}.
\end{align}
With the domain \eqref{eq:AKNS-maximal-domain}, $\bL_u$ is anti-self-adjoint by \cite[Proposition~3.1]{Forlano.Killip.Visan2024}. Write $\langle x\rangle=(1+x^2)^{1/2}$. That proposition, like \cite[Proposition~3.4 and Corollary~7.2]{Forlano.Killip.Visan2024}, is stated under the weighted hypothesis $\langle x\rangle^{-b}u\in L^2(\R)$ for some $b>0$, which an unweighted $u\in L^2(\R)$ satisfies for every $b\geq0$ and in particular for every $b>0$. All of these results are therefore available here. What is taken from them below is qualitative and pointwise: the domain \eqref{eq:AKNS-maximal-domain}, the continuity of $\bG_\kappa-\bG_{\kappa,0}$, the identities and inequalities \cite[equations~(3.26)--(3.30)]{Forlano.Killip.Visan2024}, and the diagonal formula \cite[equation~(7.8)]{Forlano.Killip.Visan2024}. The norm bounds recorded in those results carry the weights $\langle x\rangle^{-3b}$ and $\langle x\rangle^{-4b}$. The unweighted Sobolev bounds used here, with locally Lipschitz dependence on $u\in L^2(\R)$, follow from \cite[Proposition~3.7]{Forlano.Killip.Visan2024} with $L=\infty$, as quoted in Lemma~\ref{lem:mkdv-static-reciprocal}\,\textup{(ii)} below. Definition~\ref{def:mkdv-green} is the case $b=0$ of a definition that \cite{Forlano.Killip.Visan2024} states for $b\geq0$.

Write $\bG_\kappa$ for the kernel of $\bR_\kappa(u)$ and $\bG_{\kappa,0}$ for the free kernel. The three diagonal Green functions from \eqref{eq:intro-dirac-fields}--\eqref{eq:intro-dirac-p} are
\begin{align}
 \gamma_\kappa&=[\bG_\kappa-\bG_{\kappa,0}]_{11}(x,x)+[\bG_\kappa-\bG_{\kappa,0}]_{22}(x,x), \label{eq:AKNS-gamma-recalled}\\
 g_{-,\kappa}&=[\bG_\kappa]_{21}(x,x)-[\bG_\kappa]_{12}(x,x), \label{eq:AKNS-m-recalled}\\
 g_{+,\kappa}&=[\bG_\kappa]_{21}(x,x)+[\bG_\kappa]_{12}(x,x). \label{eq:AKNS-p-recalled}
\end{align}
The difference $\bG_\kappa-\bG_{\kappa,0}$ is continuous, and the free off-diagonal entries vanish. Hence all three diagonal expressions are continuous, despite the jump of the full resolvent kernel.

For a real curve $u\in C(I;L^2)$, Definition~\ref{def:mkdv-green} requires the following two identities for every $\kappa\geq1$ and all $s,t\in I$:
\begin{align}
 g_{-,\kappa}(t)-g_{-,\kappa}(s)&=\int_s^t\bigl[-g_{-,\kappa,xxx}+6u^2g_{-,\kappa,x}\bigr](t')\dd t', \label{eq:mkdv-green-m-recalled}\\
 \gamma_\kappa(t)-\gamma_\kappa(s)&=\int_s^t\Bigl[-\gamma_{\kappa,xxx}+\partial_x\bigl\{6u^2(1+\gamma_\kappa)-12\kappa u g_{-,\kappa}-12\kappa^2\gamma_\kappa\bigr\}\Bigr](t')\dd t'. \label{eq:mkdv-green-gamma-recalled}
\end{align}
These are \eqref{eq:intro-mkdv-green-m}--\eqref{eq:intro-mkdv-green-gamma}, recalled for the transfer argument. Their Bochner integrals take values in $H^{-1}$ and $H^{-2}$, respectively.

The regularity quoted in Lemma~\ref{lem:mkdv-static-reciprocal}\,\textup{(ii)}, from Forlano--Killip--Vi\c{s}an \cite[Proposition~3.7]{Forlano.Killip.Visan2024}, gives $\gamma_\kappa,g_{+,\kappa}\in C(I;H^1)$ and $g_{-,\kappa}\in C(I;H^2)$. We use it to check the two integrands separately.

In \eqref{eq:mkdv-green-m-recalled}, the term $g_{-,\kappa,xxx}$ lies in $H^{-1}$. Also $g_{-,\kappa,x}\in H^1\hookrightarrow L^\infty$ and $u^2\in L^1$, so $u^2g_{-,\kappa,x}\in L^1\hookrightarrow H^{-1}$ by Lemma~\ref{lem:products}. The first integrand is therefore continuous into $H^{-1}$.

In \eqref{eq:mkdv-green-gamma-recalled}, the term $\gamma_{\kappa,xxx}$ lies in $H^{-2}$. Inside the braces, $u^2(1+\gamma_\kappa)\in L^1$, $ug_{-,\kappa}\in L^2$, and $\gamma_\kappa\in H^1$. The expression in braces thus lies in $H^{-1}$, and its spatial derivative lies in $H^{-2}$. These products depend continuously on time, so the second integrand is continuous into $H^{-2}$.

The different orders of the two laws determine which multipliers the transfer can use. The discussion below explains why the extra derivative of $g_{-,\kappa}$ is sufficient.

\textit{Required Green functions.}
For each Miura potential $q_\sigma=u^2+\sigma u_x$, recall the scalar Green function and its reciprocal from \eqref{eq:mkdv-target-coordinate}:
\begin{align}
 \label{eq:mkdv-target-recalled} \mg_{\sigma,\kappa}=g_\kappa(\,\cdot\,;q_\sigma), \qquad h_{\sigma,\kappa}=(2\mg_{\sigma,\kappa})^{-1}.
\end{align}
Definition~\ref{def:green} requires the reciprocal law \eqref{eq:green-identity-recalled} for $h_{\sigma,\kappa}$, with values in $H^{-2}$. Thus the task is to pass from the two laws \eqref{eq:mkdv-green-m-recalled}--\eqref{eq:mkdv-green-gamma-recalled} to that scalar law for each sign. Sections~\ref{sec:mkdv-transfer} and~\ref{sec:mkdv-identification} carry out this passage.

\textit{Main difficulty.}
A first approach is to use the scalar diagonal Green function. Forlano--Killip--Vi\c{s}an \cite{Forlano.Killip.Visan2024} work with $\mg_{\sigma,\kappa}$ and derive for it a local conservation law, recorded below as \eqref{eq:smooth-scalar-g-evolution}. At $L^2$ regularity of $u$ that law gives only
\begin{align}
 \label{eq:scalar-derivative-regularity} \partial_t\mg_{\sigma,\kappa}\in H^{-2}(\R),
\end{align}
and the passage from it to the reciprocal,
\begin{align}
 \label{eq:forbidden-reciprocal-derivative} \partial_th_{\sigma,\kappa} =-\frac{\partial_t\mg_{\sigma,\kappa}} {2\mg_{\sigma,\kappa}^2},
\end{align}
is not justified by the available multiplier estimates. By part~\textup{(iii)} of Lemma~\ref{lem:mkdv-static-reciprocal} below, which holds at every $\kappa\geq1$, the coefficient $\mg_{\sigma,\kappa}^{-2}$ is a constant plus an $H^1$ function, and $H^1$ multiplies $H^{-1}$ but not $H^{-2}$.

\textit{Resolution.}
The evolution law \eqref{eq:mkdv-green-m-recalled} governs $g_{-,\kappa}\in H^2(\R)$. This gives one more derivative than is available for $\mg_{\sigma,\kappa}$ and $h_{\sigma,\kappa}$ after their constant terms are subtracted. Multiplication by an $H^2$ function is bounded on $H^{-2}$. Lemma~\ref{lem:mkdv-static-reciprocal} gives the formula $h_{\sigma,\kappa}=\kappa(1+\gamma_\kappa-\sigma g_{+,\kappa})/(1-g_{-,\kappa}^2)$, which is the second expression in \eqref{eq:reciprocal-through-m}. In the resulting time derivative, every coefficient multiplying an $H^{-2}$ distribution depends only on $g_{-,\kappa}$.

The denominator in this formula has the regularity needed for the two evolution laws. Its reciprocal lies in $1+H^2(\R)$ and multiplies the $H^{-2}$-valued law inherited from $\gamma_\kappa$, using $(1+H^2)\cdot H^{-2}\subset H^{-2}$. The coefficient built from $\gamma_\kappa$ and $g_{+,\kappa}$ multiplies the $H^{-1}$-valued law for $g_{-,\kappa}$, using $H^1\cdot H^{-1}\subset H^{-1}\hookrightarrow H^{-2}$. Thus the proof uses both evolution identities imposed by Definition~\ref{def:mkdv-green}.

\textit{Further discussion.}
Two points remain in differentiating the quotient in time. First, its numerator and denominator satisfy integrated identities in negative Sobolev spaces of different orders. We need a chain rule compatible with these identities and the positive Sobolev regularity of the factors. Second, $\gamma_\kappa$, $g_{+,\kappa}$ and $g_{-,\kappa}$ vanish at spatial infinity, so the numerator $1+\gamma_\kappa-\sigma g_{+,\kappa}$, the denominator $1-g_{-,\kappa}^2$, its reciprocal, and the normalized quotient $h_{\sigma,\kappa}/\kappa$ all tend to $1$ there. Subtracting this constant places them in the required Sobolev spaces. The lemma in Section~\ref{sec:mkdv-chain-rule} treats both points.

\textit{The plan.}
The argument then runs in three stages.
\begin{enumerate}[
label=\textup{(\arabic*)}, leftmargin=*, itemsep=.5em, before={\let\fullwidthdisplay\relax} ]
\item \emph{A formula for $h_{\sigma,\kappa}$.} Lemma~\ref{lem:mkdv-static-reciprocal} expresses $h_{\sigma,\kappa}$ through $\gamma_\kappa,g_{-,\kappa},g_{+,\kappa}$ and records the regularity and the uniform gap for $g_{-,\kappa}$ needed to use the formula. These assertions hold for each fixed function $u$.
\item \emph{Dynamics.} Lemma~\ref{lem:mkdv-time-regularity} identifies, along an arbitrary mKdV Green solution, the time derivative of $h_{\sigma,\kappa}$ with an explicit expression $V_{h_\sigma}$, as an identity in $H^{-2}$.
\item \emph{Identification of the flux.} It remains to identify $V_{h_\sigma}$ with the spatial derivative of the reciprocal Green flux. Both
\begin{align}
 \label{eq:two-static-maps} u\mapsto V_{h_\sigma}(u) \qquad\text{and}\qquad u\mapsto \partial_x\J_\kappa(\cM_\sigma(u))
\end{align}
are maps of a fixed function $u$. Each is an algebraic expression in $u$ and the diagonal Green functions associated with $u$. Their equality is therefore a statement about a single function. Lemma~\ref{lem:mkdv-smooth-identification} proves it for smooth $u$, and Lemmas~\ref{lem:mkdv-AKNS-continuity} and~\ref{lem:mkdv-map-continuity} show that both maps in \eqref{eq:two-static-maps} are continuous from $L^2$ into $H^{-2}$, so that the identity extends by density.
\end{enumerate}
Combining stage~\textup{(2)} with stage~\textup{(3)} proves Proposition~\ref{prop:mkdv-transfer}. Stage~\textup{(3)} uses the classical mKdV equation for smooth data to identify the two maps, then approximates each fixed $L^2$ function by smooth functions inside this identity. The evolution of the given Green trajectory is supplied by stage~\textup{(2)}. Stages~\textup{(1)} and~\textup{(2)} occupy Section~\ref{sec:mkdv-transfer}, and stage~\textup{(3)} occupies Section~\ref{sec:mkdv-identification}.

%%%%%%%%%%%%%%%%%%%%%%%%%%%%%%%%%%%%%%%%%%%%%%%%%%%%%%%%%%%%%%%%%%%%%%%%%%%%%%%%%%
\subsection{A Preliminary Lemma}
\label{sec:mkdv-chain-rule}
%%%%%%%%%%%%%%%%%%%%%%%%%%%%%%%%%%%%%%%%%%%%%%%%%%%%%%%%%%%%%%%%%%%%%%%%%%%%%%%%%%

The following lemma treats integrated identities in negative Sobolev spaces for functions that differ from constants by Sobolev functions. It is stated for a generic pair $(m,P)$ and applied to the quotient \eqref{eq:reciprocal-through-m} in Step~4 of the proof of Lemma~\ref{lem:mkdv-time-regularity}.

{\begin{lemma}
\label{lem:affine-negative-chain}
Let $I$ be a compact interval. Suppose that $m$ is real valued, with
\begin{align}
 \label{eq:negative-chain-m-hypotheses} m&\in C(I;H^2(\R;\R)), \qquad m(t)-m(s)=\int_s^t V(t')\dd t' \quad\text{in }H^{-1}(\R),
\end{align}
where $V\in L^1(I;H^{-1})$, and assume that, for some $\delta>0$,
\begin{align}
 \label{eq:negative-chain-gap} \sup_{(t,x)\in I\times\R}\abs{m(t,x)}\leq1-\delta.
\end{align}
Set
\begin{align}
 \label{eq:affine-negative-chain-a1} D:=1-m^2, \qquad F:=D^{-1}.
\end{align}

\begin{enumerate}[
label=\textup{(\roman*)}, leftmargin=*, itemsep=.7em, before={\let\fullwidthdisplay\relax} ]
\item \emph{The reciprocal.} The reciprocal satisfies
\begin{align}
 \label{eq:negative-chain-first-conclusion} F-1\in C(I;H^2)\cap W^{1,1}(I;H^{-1}), \qquad \partial_tF=2mD^{-2}V \quad\text{in }L^1(I;H^{-1}).
\end{align}

\item \emph{The product.} Suppose in addition that
\begin{align}
 \label{eq:negative-chain-P-hypotheses} \begin{aligned}
 P&=1+\widetilde P, \qquad \widetilde P\in C(I;H^1(\R)),\\
 \widetilde P(t)-\widetilde P(s) &=\int_s^t Z(t')\dd t' \qquad\text{in }H^{-2}(\R),
 \end{aligned}
\end{align}
where $Z\in L^1(I;H^{-2})$. Then
\begin{align}
 \label{eq:negative-chain-product-conclusion} FP-1\in C(I;H^1)\cap W^{1,1}(I;H^{-2})
\end{align}
and
\begin{align}
 \label{eq:negative-chain-product-derivative} \partial_t(FP) =FZ+2PmD^{-2}V \quad\text{in }L^1(I;H^{-2}).
\end{align}
Equivalently, for every $s,t\in I$,
\begin{align}
 \label{eq:negative-chain-product-integrated} F(t)P(t)-F(s)P(s) =\int_s^t \bigl[FZ+2PmD^{-2}V\bigr](t')\dd t' \quad\text{in }H^{-2}.
\end{align}
\item \emph{Upgrade under continuity.} If $V\in C(I;H^{-1})$ and $Z\in C(I;H^{-2})$, then the $W^{1,1}$ regularity in \eqref{eq:negative-chain-first-conclusion} and \eqref{eq:negative-chain-product-conclusion} improves to $C^1(I;H^{-1})$ and $C^1(I;H^{-2})$ respectively.
\end{enumerate}
\end{lemma}}

The Sobolev conclusions concern the differences $F-1$ and $FP-1$, consistently with the assumption $P-1\in C(I;H^1)$ in \eqref{eq:negative-chain-P-hypotheses}. The gap hypothesis \eqref{eq:negative-chain-gap} keeps $D$ uniformly away from zero, so that $F$ is defined and the composition $s\mapsto(1-s)^{-1}$ may be applied in $H^2$.

\begin{proof}
The two parts are proved in the same way: an exact algebraic increment formula is written down, its coefficient is shown to be a good multiplier on the relevant negative space, and the difference quotient is then taken.

\proofstep{Step 1: regularity of $F$}
Three facts combine. First, the gap \eqref{eq:negative-chain-gap} bounds $\abs{m}$ by $1-\delta$ on the whole of $I\times\R$, so that
\begin{align}
 \label{eq:negative-chain-D-range} \delta(2-\delta)\leq D(t,x)=1-m(t,x)^2\leq1 \qquad\bigl((t,x)\in I\times\R\bigr) ,
\end{align}
and $D$ is bounded away from zero uniformly in time. In the application this uniform bound follows from part~\textup{(iii)} of Lemma~\ref{lem:mkdv-static-reciprocal}, where compactness of $m(I)$ in $H^2(\R)\hookrightarrow C_0(\R)$ turns pointwise positivity into \eqref{eq:AKNS-uniform-gap}. Second, $H^2(\R)$ is an algebra, so $D-1=-m^2$ lies in $C(I;H^2)$, and $\norm{D-1}_{H^2}$ is bounded on $I$ because $m\in C(I;H^2)$ and $I$ is compact. Third, the function
\begin{align}
 \label{eq:negative-chain-composition-function} G(s):=(1-s)^{-1}-1
\end{align}
is smooth on a neighbourhood of $[0,(1-\delta)^2]$, which contains the range of $m^2$ by \eqref{eq:negative-chain-D-range}, and satisfies $G(0)=0$. The one-dimensional Sobolev composition theorem therefore applies to $G\circ m^2$ and gives
\begin{align}
 \label{eq:affine-negative-chain-a2} F-1=(1-m^2)^{-1}-1\in C(I;H^2) ,
\end{align}
with $\norm{F-1}_{H^2}$ again bounded on $I$. The condition $G(0)=0$ places the composition $F-1$ in $H^2$.

The uniform bounds ensure that the multiplier constants below remain bounded as the pair of times varies over $I\times I$.

\proofstep{Step 2: the time derivative of $F$}
The increment of $F$ factors exactly through the increment of $m$: for $s,t\in I$,
\begin{align}
 \label{eq:negative-chain-exact-increment} F(t)-F(s) =C_{t,s}\,[m(t)-m(s)], \qquad C_{t,s}:=\frac{m(t)+m(s)}{D(t)D(s)}
\end{align}
holds pointwise, as one checks by putting the difference of the two fractions over a common denominator. The nonlinearity enters through the coefficient $C_{t,s}$, which is a fixed function at each pair of times. The increment of $m$ enters linearly and may therefore be replaced by its integrated law. The coefficients $C_{t,s}$ are uniformly bounded in $H^2$ by \eqref{eq:negative-chain-gap} and the $H^2$ algebra property. Since multiplication by an $H^2$ function is bounded on $H^{-1}$, \eqref{eq:negative-chain-m-hypotheses} gives
\begin{align}
 \label{eq:affine-negative-chain-a3} \norm{F(t)-F(s)}_{H^{-1}} \lesssim \int_s^t\norm{V(t')}_{H^{-1}}\dd t', \qquad s<t.
\end{align}
Thus $F-1$ is absolutely continuous as an $H^{-1}$-valued function, with the integrable modulus $\norm{V(\cdot)}_{H^{-1}}$. Since $H^{-1}(\R)$ is a Hilbert space, and so has the Radon--Nikod\'ym property, an absolutely continuous function into it is differentiable almost everywhere and is the Bochner integral of its derivative, which lies in $L^1(I;H^{-1})$. Therefore $F-1\in W^{1,1}(I;H^{-1})$, and its weak derivative is its almost-everywhere strong derivative, which we now compute.

By \eqref{eq:negative-chain-m-hypotheses} and the Lebesgue differentiation theorem for Bochner integrals, $m$ is differentiable as an $H^{-1}$-valued function, with derivative $V(t)$, at every Lebesgue point $t$ of $V$, hence at almost every $t\in I$. At such a time, divide \eqref{eq:negative-chain-exact-increment} by $t-s$. As $s\to t$, continuity of $m$ in $H^2$ implies
\begin{align}
 \label{eq:affine-negative-chain-a4} C_{t,s}\rightarrow2m(t)D(t)^{-2} \quad\text{in }H^2,
\end{align}
while the difference quotient of $m$ converges to $V(t)$ in $H^{-1}$. The multiplier estimate therefore yields
\begin{align}
 \label{eq:affine-negative-chain-a5} \partial_tF=2mD^{-2}V \quad\text{in }H^{-1}
\end{align}
for almost every $t$, proving \eqref{eq:negative-chain-first-conclusion}.

\proofstep{Step 3: regularity of the product}
Since $F-1\in C(I;H^2)$ and $P-1\in C(I;H^1)$, the $H^1$ algebra property gives
\begin{align}
 \label{eq:affine-negative-chain-a6} FP-1=(F-1)+(P-1)+(F-1)(P-1)\in C(I;H^1).
\end{align}
Moreover, the exact product increment
\begin{align}
 \label{eq:negative-chain-product-increment} F(t)P(t)-F(s)P(s) &=F(t)[P(t)-P(s)] +[F(t)-F(s)]P(s)
\end{align}
is meaningful in $H^{-2}$. Indeed, $1+H^2$ acts boundedly on $H^{-2}$, while $1+H^1$ acts boundedly on $H^{-1}\hookrightarrow H^{-2}$. The preceding estimates consequently give
\begin{align}
 \label{eq:affine-negative-chain-a7} \norm{F(t)P(t)-F(s)P(s)}_{H^{-2}} \lesssim \int_s^t \bigl(\norm{Z(t')}_{H^{-2}}+\norm{V(t')}_{H^{-1}}\bigr)\dd t', \qquad s<t.
\end{align}
Hence $FP-1$ is absolutely continuous into $H^{-2}$, and $FP-1\in W^{1,1}(I;H^{-2})$ by the Radon--Nikod\'ym argument of Step~2.

\proofstep{Step 4: the time derivative of the product, and the upgrade}
Let $t$ be a common Lebesgue point of $V$ and $Z$. Such points form a set of full measure. At $t$, Step~2 gives differentiability of $F$ in $H^{-1}$, and the Lebesgue differentiation theorem gives differentiability of $\widetilde P$ in $H^{-2}$ with derivative $Z(t)$. At such a time, divide \eqref{eq:negative-chain-product-increment} by $t-s$. Strong continuity of $F$ in $1+H^2$ and of $P$ in $1+H^1$, together with the two multiplier estimates, gives
\begin{align}
 \label{eq:affine-negative-chain-a8} \partial_t(FP) &=FZ+(\partial_tF)P=FZ+2PmD^{-2}V \quad\text{in }H^{-2}.
\end{align}
This proves \eqref{eq:negative-chain-product-derivative}. The fundamental theorem for Bochner--Sobolev functions gives \eqref{eq:negative-chain-product-integrated}.

If $V$ and $Z$ are continuous, the first derivative is continuous into $H^{-1}$ and the product derivative is continuous into $H^{-2}$. This gives the final $C^1$ assertion.
\end{proof}

%%%%%%%%%%%%%%%%%%%%%%%%%%%%%%%%%%%%%%%%%%%%%%%%%%%%%%%%%%%%%%%%%%%%%%%%%%%%%%%%%%
\subsection{Evolution of \texorpdfstring{$h_{\sigma,\kappa}$}{h}}
\label{sec:mkdv-transfer}
%%%%%%%%%%%%%%%%%%%%%%%%%%%%%%%%%%%%%%%%%%%%%%%%%%%%%%%%%%%%%%%%%%%%%%%%%%%%%%%%%%

Stages~\textup{(1)} and~\textup{(2)} of the plan are carried out here.

\medskip
Throughout this subsection and the next, $\kappa\geq1$ is fixed, and it is suppressed from every symbol that carries it: $\gamma$, $g_+$ and $g_-$ stand for the diagonal Green functions $\gamma_\kappa$, $g_{+,\kappa}$ and $g_{-,\kappa}$ of \eqref{eq:AKNS-gamma-recalled}--\eqref{eq:AKNS-p-recalled}. Likewise, $\mg_\sigma$ and $h_\sigma$ denote the functions $\mg_{\sigma,\kappa}$ and $h_{\sigma,\kappa}$ of \eqref{eq:mkdv-target-coordinate}. We write $\J$ for the flux $\J_\kappa$ of \eqref{eq:green-flux}. The parameter itself continues to appear as a coefficient, as in $g_{-,x}=-2\kappa g_+$. What is suppressed is only its appearance as a subscript. That subscript is restored in the statement of Proposition~\ref{prop:mkdv-transfer}, whose conclusion is quantified over $\kappa$, and dropped again in its proof. As in Lemma~\ref{lem:affine-negative-chain}, and here also with an approximating index in Lemma~\ref{lem:mkdv-AKNS-continuity}, we write
\begin{align}
 \label{eq:AKNS-D-def} D:=1-g_-^2 .
\end{align}

%%%%%%%%%%%%%%%%%%%%%%%%%%%%%%%%%%%%%%%%%%%%%%%%%%%%%%%%%%%%%%%%%%%%%%%%%%%%%%%%%%

Part~\textup{(i)} of the next lemma gives two formulas for $h_\sigma$. The first, $\kappa/(1+\gamma+\sigma g_+)$, follows from the quoted formula for the diagonal Green function. Its denominator is a constant plus an $H^1$ function, so differentiating its reciprocal in time would require the undefined product in \eqref{eq:forbidden-reciprocal-derivative}. The second formula writes the reciprocal using $D=1-g_-^2$, whose regularity is better because $g_-\in H^2$. Parts~\textup{(ii)} and~\textup{(iii)} give the regularity and lower bound needed to use $D^{-1}$ and $D^{-2}$ as multipliers on $H^{-2}$.

{\begin{lemma}
\label{lem:mkdv-static-reciprocal}
Let $\kappa\geq1$ and $\sigma\in\{+1,-1\}$, and let $D$ be as in \eqref{eq:AKNS-D-def}.
\begin{enumerate}[
label=\textup{(\roman*)}, leftmargin=*, itemsep=.7em, before={\let\fullwidthdisplay\relax} ]
\item \emph{The formula.} For every real $u\in L^2(\R)$, the operator $L_{q_\sigma}+\kappa^2$ is boundedly invertible by Lemma~\ref{lem:two-Miura-static}\,\textup{(ii)}, its inverse has a continuous integral kernel by \cite[Corollary~7.2]{Forlano.Killip.Visan2024}, and $\mg_\sigma$ in \eqref{eq:mkdv-target-coordinate} is the diagonal of that kernel, as in \eqref{eq:intro-resolvent}. These conclusions hold for every $\kappa\geq1$. The diagonal Green function of $L_{q_\sigma}+\kappa^2$ and its reciprocal are
\begin{align}
 \mg_\sigma &=\frac{1+\gamma+\sigma g_+}{2\kappa}, \label{eq:scalar-diagonal-from-AKNS}\\
 h_\sigma =\frac1{2\mg_\sigma} &=\frac{\kappa}{1+\gamma+\sigma g_+} =\kappa\,\frac{1+\gamma-\sigma g_+}{D} . \label{eq:reciprocal-through-m}
\end{align}

\item \emph{Regularity.} The functions satisfy $\gamma,\ g_+\in H^1(\R)$ and $g_-\in H^2(\R)$, with locally Lipschitz dependence on $u\in L^2$. If moreover $u\in C(I;L^2)$ on a compact interval $I$, then
\begin{align}
 \label{eq:AKNS-global-regularity} \gamma,\ g_+\in C(I;H^1), \qquad g_-\in C(I;H^2).
\end{align}

\item \emph{Uniform gap and $H^2$ multipliers.} If $u\in C(I;L^2)$, then
\begin{align}
 \label{eq:AKNS-uniform-gap} \sup_{(t,x)\in I\times\R}\abs{g_-(t,x)}<1,
\end{align}
and consequently, for $j=1,2$,
\begin{align}
 \label{eq:AKNS-reciprocal-Htwo} D^{-j}-1\in C(I;H^2), \qquad h_\sigma-\kappa\in C(I;H^1).
\end{align}

\item \emph{The quoted spatial identities.} For every real $u\in L^2(\R)$,
\begin{align}
 \label{eq:AKNS-static-identities} g_{-,x}=-2\kappa g_+, \qquad (1+\gamma)^2-g_+^2=D, \qquad \abs{g_+(x)}<1+\gamma(x), \qquad \abs{g_-(x)}<1 .
\end{align}
These are \cite[Proposition~3.4, equations~(3.28)--(3.30)]{Forlano.Killip.Visan2024}, rewritten in the present notation. Parts of \eqref{eq:AKNS-static-identities} are consumed as inputs by Lemmas~\ref{lem:mkdv-time-regularity}, \ref{lem:mkdv-smooth-identification} and \ref{lem:mkdv-AKNS-continuity}.
\end{enumerate}
\end{lemma}}

\begin{proof}
The four parts lean on \cite{Forlano.Killip.Visan2024} to different degrees: \textup{(ii)} and \textup{(iv)} are quotations, \textup{(i)} is a quotation combined with the factorization of Lemma~\ref{lem:two-Miura-static}, and \textup{(iii)} is the only one that requires an argument of its own. They are taken out of order, because \textup{(i)} and \textup{(iii)} both consume \textup{(iv)}.

\proofstep{Part \textup{(iv)}} The relations \eqref{eq:AKNS-static-identities} are \cite[Proposition~3.4, equations~(3.28)--(3.30)]{Forlano.Killip.Visan2024}. Only the second needs comment: the source states $\gamma+\frac12\gamma^2=\frac12(g_+^2-g_-^2)$, and multiplying by $2$, adding $1$, and subtracting $g_+^2$ from both sides turns the left side into $(1+\gamma)^2-g_+^2$ and the right side into $1-g_-^2=D$, which is the form used below.

\proofstep{Part \textup{(i)}} By the form factorization \eqref{eq:Miura-form-factorization} of Lemma~\ref{lem:two-Miura-static}, together with part~\textup{(iii)} of that lemma, which identifies the factorized operator with the form-defined one through equality of their form domains, the scalar Schr\"odinger operators $L_{q_\sigma}$ are exactly the nonnegative factorized operators appearing in \cite[Proposition~7.1]{Forlano.Killip.Visan2024}. The resolvent identity and diagonal formula \cite[Corollary~7.2, equation~(7.8)]{Forlano.Killip.Visan2024} therefore apply and give \eqref{eq:scalar-diagonal-from-AKNS}.

The reciprocal follows from part~\textup{(iv)}. Its third relation gives $1+\gamma+\sigma g_+>0$ pointwise, so \eqref{eq:scalar-diagonal-from-AKNS} may be inverted, which is the first equality in \eqref{eq:reciprocal-through-m}. Multiplying numerator and denominator by $1+\gamma-\sigma g_+$ and using the second identity of \eqref{eq:AKNS-static-identities} converts the result into $\kappa(1+\gamma-\sigma g_+)/D$, which is the second equality. The point of that second form is that the inversion has been moved off the sum $1+\gamma+\sigma g_+$, whose regularity is only $1+H^1$, and onto $D$, which is built from $g_-\in H^2$. The factors $D^{-1}$ and $D^{-2}$ have the $1+H^2$ regularity needed to multiply $H^{-2}$ distributions, as shown in part~\textup{(iii)}.

\proofstep{Part \textup{(ii)}} At a single potential the assertion is quoted. The regularity $\gamma,g_+\in H^1(\R)$ and $g_-\in H^2(\R)$, and the locally Lipschitz dependence of all three on $u\in L^2(\R;\R)$, are the unweighted case of \cite[Proposition~3.7, equations~(3.51)--(3.53)]{Forlano.Killip.Visan2024}, obtained by taking $L=\infty$ there. That proposition is formulated for $2L$-periodic potentials measured in the one-period $L^2$ norm, and the case $L=\infty$, that is $u\in L^2(\R)$, is included explicitly. Its norms $H^1_\kappa$ and $H^2_\kappa$ coincide with $H^1$ and $H^2$ up to constants depending on the fixed $\kappa$.

Along a trajectory the assertion is a composition. If $u\in C(I;L^2)$ with $I$ compact, then $u(I)$ is a compact, hence bounded, subset of $L^2$, so the quoted maps are Lipschitz on it with a single constant. Composing $t\mapsto u(t)$ with them gives the continuity in $t$ recorded in \eqref{eq:AKNS-global-regularity}. Compactness of $u(I)$ is what makes that constant uniform in $t$, and it is used again, for a different purpose, in part~\textup{(iii)}.

\proofstep{Part \textup{(iii)}} We first prove the gap \eqref{eq:AKNS-uniform-gap}. We combine the pointwise bound in \eqref{eq:AKNS-static-identities} with compactness to obtain a uniform gap. Since $H^2(\R)\hookrightarrow C_0(\R)$, each $g_-(t)$ vanishes at spatial infinity. Because $g_-(I)$ is compact in $H^2$, this vanishing is uniform in $t$, so there is $R>0$ with $\abs{g_-(t,x)}\leq\frac12$ for $\abs{x}\geq R$ and all $t\in I$. On the complementary set $I\times[-R,R]$, which is compact, the function $\abs{g_-}$ is continuous and everywhere $<1$, so its maximum there is $<1$ as well. Taking the larger of the two bounds proves \eqref{eq:AKNS-uniform-gap}, say with $\abs{g_-}\leq1-\delta$.

For \eqref{eq:AKNS-reciprocal-Htwo}, write
\begin{align}
 \label{eq:AKNS-Fj-def} D^{-j}-1=F_j(g_-), \qquad F_j(z):=(1-z^2)^{-j}-1 .
\end{align}
Each $F_j$ is smooth on the interval $\abs{z}\leq1-\delta$ supplied by the gap, and $F_j(0)=0$. Since $H^2(\R)$ is a Banach algebra, the one-dimensional Sobolev composition theorem applies and gives $F_j(g_-)\in C(I;H^2)$. Finally, $h_\sigma-\kappa=\kappa\bigl(D^{-1}(1+\gamma-\sigma g_+)-1\bigr)$ lies in $C(I;H^1)$, because $H^1$ is an algebra, $H^2\subset H^1$, and both $D^{-1}-1$ and $\gamma-\sigma g_+$ lie in $C(I;H^1)$.
\end{proof}

%%%%%%%%%%%%%%%%%%%%%%%%%%%%%%%%%%%%%%%%%%%%%%%%%%%%%%%%%%%%%%%%%%%%%%%%%%%%%%%%%%

The next lemma uses both identities in Definition~\ref{def:mkdv-green}, taken from \cite{Forlano.Killip.Visan2024}. The identity for $\gamma$ gives its contribution to the time derivative of the numerator in \eqref{eq:reciprocal-through-m}. The identity for $g_-$ gives the evolution of the denominator and, after one spatial derivative, the evolution of $g_+$. These identities supply the two time derivatives required by Lemma~\ref{lem:affine-negative-chain}.

{\begin{lemma}
\label{lem:mkdv-time-regularity}
Let $u\in C(I;L^2(\R;\R))$ be an mKdV Green solution in the sense of Definition~\ref{def:mkdv-green}, let $\kappa\geq1$, and let $\sigma\in\{+1,-1\}$. Define the expressions
\begin{align}
 V_-&:=-g_{-,xxx}+6u^2g_{-,x}, \label{eq:AKNS-vector-fields}\\
 V_\gamma&:=-\gamma_{xxx}+\partial_x\left\{6u^2(1+\gamma)-12\kappa ug_--12\kappa^2\gamma\right\}, \label{eq:AKNS-vector-gamma}\\
 V_+&:=-g_{+,xxx}+6\partial_x(u^2g_+). \label{eq:AKNS-vector-p}
\end{align}
Then
\begin{align}
 \label{eq:AKNS-field-spaces} V_-\in C(I;H^{-1}), \qquad V_\gamma,\ V_+\in C(I;H^{-2}),
\end{align}
the expression
\begin{align}
 \label{eq:AKNS-h-vector-field} V_{h_\sigma} :=\kappa D^{-1}(V_\gamma-\sigma V_+) +2\kappa\,(1+\gamma-\sigma g_+)\,g_-D^{-2}V_-
\end{align}
is well defined and lies in $C(I;H^{-2})$, and
\begin{align}
 \label{eq:AKNS-reciprocal-C1} h_\sigma-\kappa \in C(I;H^1)\cap C^1(I;H^{-2}), \qquad \partial_th_\sigma=V_{h_\sigma} \quad\text{in }H^{-2} .
\end{align}
Equivalently, for all $s,t\in I$,
\begin{align}
 \label{eq:weak-h-evolution} h_\sigma(t)-h_\sigma(s) =\int_s^tV_{h_\sigma}(t')\dd t' \quad\text{in }H^{-2}(\R).
\end{align}
\end{lemma}}

\begin{proof}
The plan is to place the increments of $g_-,\gamma,g_+$ in the spaces where they actually live, to convert the two defining identities of Definition~\ref{def:mkdv-green} into integrated laws for the factors $F$ and $P$ separately, and only then to apply the chain rule of Lemma~\ref{lem:affine-negative-chain} to their product.

\proofstep{Step 1: regularity of the right-hand sides} We verify \eqref{eq:AKNS-field-spaces} term by term, using the regularity \eqref{eq:AKNS-global-regularity} supplied by Lemma~\ref{lem:mkdv-static-reciprocal}\,\textup{(ii)}.

For $V_-$: since $g_-\in H^2$, we have $g_{-,xxx}\in H^{-1}$, and $g_{-,x}\in H^1\hookrightarrow L^\infty$ gives $u^2g_{-,x}\in L^1$, hence $u^2g_{-,x}\in H^{-1}$ by part~\textup{(iii)} of Lemma~\ref{lem:products}. Thus $V_-\in C(I;H^{-1})$.

For $V_\gamma$: the three terms inside the brace satisfy $u^2(1+\gamma)\in L^1\subset H^{-1}$, $ug_-\in L^2\subset H^{-1}$ and $\gamma\in H^1\subset H^{-1}$, so the brace lies in $H^{-1}$ and one spatial derivative places it in $H^{-2}$. Also $\gamma\in H^1$ gives $\gamma_{xxx}\in H^{-2}$. Thus $V_\gamma\in C(I;H^{-2})$.

For $V_+$: the same reasoning applies with $u^2g_+\in L^1\subset H^{-1}$ and $g_+\in H^1$.

\proofstep{Step 2: the three integrated identities} The identity \eqref{eq:mkdv-green-m-recalled} of Definition~\ref{def:mkdv-green} is exactly
\begin{align}
 \label{eq:AKNS-m-evolution} g_-(t)-g_-(s)=\int_s^tV_-(t')\dd t' ,
\end{align}
and \eqref{eq:mkdv-green-gamma-recalled} is the corresponding statement for $\gamma$ with integrand $V_\gamma$,
\begin{align}
 \label{eq:AKNS-gamma-evolution} \gamma(t)-\gamma(s)=\int_s^tV_\gamma(t')\dd t' .
\end{align}
Applying $-(2\kappa)^{-1}\partial_x$ to \eqref{eq:AKNS-m-evolution}, and using the first spatial identity $g_{-,x}=-2\kappa g_+$ of \eqref{eq:AKNS-static-identities} on the left and its consequence $-(2\kappa)^{-1}\partial_xV_-=V_+$ on the right, yields
\begin{align}
 \label{eq:AKNS-p-evolution} g_+(t)-g_+(s)=\int_s^tV_+(t')\dd t' .
\end{align}
Both defining identities of Definition~\ref{def:mkdv-green} have now been used, the $g_-$-identity twice: once as it stands, and once differentiated. Subtracting $\sigma$ times \eqref{eq:AKNS-p-evolution} from \eqref{eq:AKNS-gamma-evolution} gives
\begin{align}
 \label{eq:AKNS-P-evolution} \bigl[\gamma-\sigma g_+\bigr](t)-\bigl[\gamma-\sigma g_+\bigr](s) =\int_s^t\bigl[V_\gamma-\sigma V_+\bigr](t')\dd t' \quad\text{in }H^{-2} .
\end{align}
These identities are stated in $\mathcal D'(\R)$ in Definition~\ref{def:mkdv-green}. By Step~1 every term of each one lies in the Sobolev space indicated, so they are identities there.

\proofstep{Step 3: why $V_{h_\sigma}$ is well defined} Since $H^2(\R)$ is a Banach algebra, duality gives
\begin{align}
 \label{eq:Htwo-Hminustwo-product} \norm{aF}_{H^{-2}} \lesssim\norm{a}_{H^2}\norm{F}_{H^{-2}}, \qquad a\in H^2,\quad F\in H^{-2},
\end{align}
and the same estimate holds when $a$ is a constant plus an $H^2$ function. Both products in \eqref{eq:AKNS-h-vector-field} are therefore legitimate:
\begin{enumerate}[
label=\textup{(\alph*)}, leftmargin=*, itemsep=.3em, before={\let\fullwidthdisplay\relax} ]
\item the first is of the form $(\R+H^2)\cdot H^{-2}$, by \eqref{eq:AKNS-reciprocal-Htwo} and \eqref{eq:Htwo-Hminustwo-product}.
\item the second is of the form $H^1\cdot H^{-1}\subset H^{-1}\subset H^{-2}$, by part~\textup{(ii)} of Lemma~\ref{lem:products}, since $(1+\gamma-\sigma g_+)g_-D^{-2}\in C(I;H^1)$ by \eqref{eq:AKNS-global-regularity} and \eqref{eq:AKNS-reciprocal-Htwo}.
\end{enumerate}
The extra derivative on $g_-$ supplies the $H^2$ coefficient in \textup{(a)}. Hence $V_{h_\sigma}\in C(I;H^{-2})$.

\proofstep{Step 4: the chain rule} We apply Lemma~\ref{lem:affine-negative-chain} with $m=g_-$ and
\begin{align}
 \label{eq:chain-rule-substitution} V=V_-, \qquad Z=V_\gamma-\sigma V_+, \qquad F=D^{-1}, \qquad P=1+\gamma-\sigma g_+ .
\end{align}
We verify its hypotheses using \eqref{eq:AKNS-global-regularity}, which places $\gamma$ and $g_+$ in $C(I;H^1)$ and $g_-$ in $C(I;H^2)$. Thus \eqref{eq:negative-chain-m-hypotheses} holds with $V=V_-$, its increment being \eqref{eq:AKNS-m-evolution}, an identity in $H^{-1}$, and the uniform gap \eqref{eq:negative-chain-gap} is \eqref{eq:AKNS-uniform-gap}. Also, \eqref{eq:negative-chain-P-hypotheses} holds with $\widetilde P=\gamma-\sigma g_+$, which lies in $C(I;H^1)$ by the regularity of $\gamma$ and $g_+$, and whose increment is \eqref{eq:AKNS-P-evolution}, an identity in $H^{-2}$ with $Z=V_\gamma-\sigma V_+$. Both $V_-$ and $Z$ are continuous in time by \eqref{eq:AKNS-field-spaces}, so part~\textup{(iii)} applies as well.

By the second expression in \eqref{eq:reciprocal-through-m}, the product $FP$ of the lemma is $h_\sigma/\kappa$. Parts~\textup{(ii)} and~\textup{(iii)} therefore give \eqref{eq:AKNS-reciprocal-C1} and \eqref{eq:weak-h-evolution} after multiplication by the constant $\kappa$, the formula \eqref{eq:negative-chain-product-derivative} becoming
\begin{align}
 \label{eq:mkdv-time-regularity-a1} \kappa\bigl[FZ+2Pg_-D^{-2}V\bigr] =\kappa D^{-1}(V_\gamma-\sigma V_+) +2\kappa(1+\gamma-\sigma g_+)g_-D^{-2}V_- =V_{h_\sigma}.
\end{align}
This is \eqref{eq:AKNS-h-vector-field}, and the proof is complete.
\end{proof}

%%%%%%%%%%%%%%%%%%%%%%%%%%%%%%%%%%%%%%%%%%%%%%%%%%%%%%%%%%%%%%%%%%%%%%%%%%%%%%%%%%
\subsection{Identification of the Flux and Both Reciprocal KdV Laws}
\label{sec:mkdv-identification}
%%%%%%%%%%%%%%%%%%%%%%%%%%%%%%%%%%%%%%%%%%%%%%%%%%%%%%%%%%%%%%%%%%%%%%%%%%%%%%%%%%

At this point the time derivative of $h_\sigma$ has been identified with the expression $V_{h_\sigma}$. It remains to prove that this expression equals the spatial derivative of the reciprocal Green flux. By the discussion around \eqref{eq:two-static-maps}, this is an identity for each fixed function $u$. The next lemma proves it for smooth $u$. The notational conventions of Section~\ref{sec:mkdv-transfer} remain in force.

To verify the identity at a fixed smooth $u$, we embed $u$ in its smooth mKdV trajectory and compute both sides as time derivatives at time zero.

{\begin{lemma}
\label{lem:mkdv-smooth-identification}
Let $\kappa\geq1$ and $\sigma\in\{+1,-1\}$. For every real $u\in C_c^\infty(\R)$,
\begin{align}
 \label{eq:identify-AKNS-h-vector-field} V_{h_\sigma}(u)=\partial_x\J\bigl(\cM_\sigma(u)\bigr) \quad\text{in }H^{-2}(\R),
\end{align}
where $V_{h_\sigma}$ is the expression \eqref{eq:AKNS-h-vector-field} and $\J$ is the KdV Green flux \eqref{eq:green-flux}.
\end{lemma}}

\begin{proof}
Both sides of \eqref{eq:identify-AKNS-h-vector-field} depend only on the single function $u$. Embed $u$ as the time-zero trace of its defocusing mKdV evolution, which for Schwartz data exists globally, is smooth, and is rapidly decaying: this is the classical Schwartz-class flow, see Kato \cite{Kato1983}, whose extension to $L^2(\R)$ is the object of \cite[Theorem~1.1]{Harrop-Griffiths.Killip.Visan2020}. Along such a trajectory the evolution identities (3.72)--(3.74) of \cite{Forlano.Killip.Visan2024}, used in Step~1, hold classically. Every identity below is then evaluated at time zero, where the trace is the given $u$.

\proofstep{Step 1: the left side is a time derivative} Equations~(3.72)--(3.74) of \cite{Forlano.Killip.Visan2024} identify the time derivatives of $g_+$, $g_-$ and $\gamma$ along this smooth evolution with $V_+$, $V_-$ and $V_\gamma$ respectively. Differentiating \eqref{eq:reciprocal-through-m} by the ordinary chain rule, which is legitimate for smooth functions, therefore gives $\partial_th_\sigma=V_{h_\sigma}$ at time zero, with $V_{h_\sigma}$ the same expression \eqref{eq:AKNS-h-vector-field}. This gives the smooth instance of Lemma~\ref{lem:mkdv-time-regularity}.

\proofstep{Step 2: the scalar conservation law} Write
\begin{align}
 \label{eq:smooth-g-h-abbrev} g:=\mg_\sigma, \qquad h:=\frac1{2g}=h_\sigma.
\end{align}
By the intertwining identities \eqref{eq:classical-Miura-intertwining} of Lemma~\ref{lem:two-Miura-static}, the Miura transform $q_\sigma=\cM_\sigma(u)$ is itself a smooth solution of \eqref{eq:KdV}. For $\sigma=+1$, equations~(3.72)--(3.74) of \cite{Forlano.Killip.Visan2024}, together with the spatial identities (3.26)--(3.28) of that reference, of which (3.28) is recorded in \eqref{eq:AKNS-static-identities}, and the diagonal formula~(7.8) there, give the scalar conservation law
\begin{align}
 \label{eq:smooth-scalar-g-evolution} g_t =\partial_x\bigl(2g_{xx}-6q_\sigma g-12\kappa^2g\bigr),
\end{align}
the calculation recorded as equation~(7.20) there, in the proof of \cite[Theorem~7.6\,(c)]{Forlano.Killip.Visan2024}: the law for $g_-$ is differentiated in $x$ and rewritten by means of (3.26)--(3.28) as a law for $g_+$, equation~(7.21) there, the law for $\gamma$ is rewritten as equation~(7.22) there, and the two are added and combined with (7.8). That reference states (7.20) for every $\kappa>0$, which includes the range $\kappa\geq1$ used here. For $\sigma=-1$ the same identity follows by applying the case $\sigma=+1$ to $-u$, which is again a defocusing mKdV solution and satisfies $\cM_-(u)=\cM_+(-u)$.

\proofstep{Step 3: from the scalar law to the reciprocal law} We use the calculation in part~\textup{(ii)} of Lemma~\ref{lem:microscopic}. For the smooth Miura potentials at hand, the two decaying solutions in Step~1 of the proof of Lemma~\ref{lem:Drach-density} give \eqref{eq:Drach} whenever $L_{q_\sigma}+\kappa^2>0$. Lemma~\ref{lem:two-Miura-static} gives $L_{q_\sigma}\geq0$, so this condition holds for every $\kappa>0$, including the full range $\kappa\geq1$ asserted here. Substituting the recovery formula \eqref{eq:Drach} into \eqref{eq:smooth-scalar-g-evolution} expresses the scalar flux through $h$,
\begin{align}
 \mathcal B&:=\frac{h_{xx}}{2h^2}-\frac{h_x^2}{4h^3}-3h-\frac{3\kappa^2}{h}, \qquad g_t=\partial_x\mathcal B, \label{eq:smooth-B-and-flux}\\
 \J(q_\sigma)&=-h_{xx}+\frac{3h_x^2}{2h}+2h^3-6\kappa^2h+4\kappa^3, \qquad \partial_x\J(q_\sigma)=-2h^2\partial_x\mathcal B. \label{eq:smooth-J-and-derivative}
\end{align}
the last equality being the expansion \eqref{eq:microscopic-exact}. Since $h$ is smooth and bounded away from zero, the ordinary chain rule gives
\begin{align}
 \label{eq:smooth-h-from-g} h_t=-2h^2g_t=-2h^2\partial_x\mathcal B=\partial_x\J(q_\sigma).
\end{align}
Comparing \eqref{eq:smooth-h-from-g} with Step~1 proves \eqref{eq:identify-AKNS-h-vector-field} at time zero, that is, at the given $u$.
\end{proof}

%%%%%%%%%%%%%%%%%%%%%%%%%%%%%%%%%%%%%%%%%%%%%%%%%%%%%%%%%%%%%%%%%%%%%%%%%%%%%%%%%%

We now extend the identity between the two maps in \eqref{eq:two-static-maps} from smooth functions to all of $L^2$. This requires continuity with respect to a fixed spatial function $u$. We first quote the convergence of the associated diagonal Green functions, and then prove convergence of the two expressions built from them.

{\begin{lemma}
\label{lem:mkdv-AKNS-continuity}
Let $\kappa\geq1$, let $u\in L^2(\R;\R)$, and let $u_n\to u$ in $L^2(\R;\R)$. Attach a subscript $n$ to each diagonal Green function associated with $u_n$, and write $D_n:=1-g_{-,n}^2$. Then
\begin{align}
 \label{eq:AKNS-approximation-fields} \begin{aligned}
 \gamma_n&\rightarrow\gamma \quad\text{and}\quad g_{+,n}\rightarrow g_+ \qquad\text{in }H^1,\\
 g_{-,n}&\rightarrow g_- \qquad\text{in }H^2,
 \end{aligned}
\end{align}
and, for $j=1,2$,
\begin{align}
 \label{eq:AKNS-reciprocal-approximation} \norm{D_n^{-j}-D^{-j}}_{H^2}\rightarrow0 .
\end{align}
\end{lemma}}

\begin{proof}
The convergences \eqref{eq:AKNS-approximation-fields} follow from the local Lipschitz dependence on $u\in L^2$ recorded in part~\textup{(ii)} of Lemma~\ref{lem:mkdv-static-reciprocal}, applied to the convergent sequence $u_n\to u$.

To prove convergence of the reciprocal powers $D^{-j}$, we need a common gap from $\pm1$ for $g_-$ and all sufficiently large $g_{-,n}$. This places their ranges in a single compact subinterval of $(-1,1)$ on which the composition estimates have uniform constants. Step~1 constructs the interval, and Step~2 applies the estimate.

\proofstep{Step 1: a gap common to the approximating functions}
Since $g_-\in H^2\hookrightarrow C_0$ and $\abs{g_-(x)}<1$ for every $x$ by \eqref{eq:AKNS-static-identities}, the number
\begin{align}
 \label{eq:common-gap-def} \delta_-:=1-\norm{g_-}_{L^\infty}
\end{align}
is positive. By \eqref{eq:AKNS-approximation-fields} and $H^2\hookrightarrow L^\infty$, we have $\norm{g_{-,n}-g_-}_{L^\infty}\to0$, hence
\begin{align}
 \label{eq:common-gap-bound} \norm{g_{-,n}}_{L^\infty}\leq1-\frac{\delta_-}2
\end{align}
for all sufficiently large $n$. Thus $g_-$ and all these $g_{-,n}$ take values in the single compact interval $[-1+\delta_-/2,\,1-\delta_-/2]\subset(-1,1)$, on which the functions $F_j$ of \eqref{eq:AKNS-Fj-def} are smooth.

\proofstep{Step 2: the composition}
On that fixed compact interval the composition estimate in the $H^2$ algebra applies, and gives
\begin{align}
 \label{eq:Fj-composition-difference} \norm{D_n^{-j}-D^{-j}}_{H^2} =\norm{F_j(g_{-,n})-F_j(g_-)}_{H^2} \leq C_{j,\delta_-,R}\norm{g_{-,n}-g_-}_{H^2} \rightarrow0 ,
\end{align}
where
\begin{align}
 \label{eq:composition-radius} R:=\norm{g_-}_{H^2}+\sup_n\norm{g_{-,n}}_{H^2}<\infty ,
\end{align}
which is finite because $g_{-,n}\to g_-$ in $H^2$. The composition constant depends on $R$ because the $H^2$ difference estimate controls derivatives of $F_j(g_{-,n})$ and $F_j(g_-)$. The uniform bound on $R$ therefore proves \eqref{eq:AKNS-reciprocal-approximation}. Each $D_n^{-j}$ tends to one at spatial infinity, and the estimate applies to its difference from $D^{-j}$ in $H^2$.
\end{proof}

{\begin{lemma}
\label{lem:mkdv-map-continuity}
Let $\kappa\geq1$ and $\sigma\in\{+1,-1\}$. Then both maps in \eqref{eq:two-static-maps} are continuous from $L^2(\R;\R)$ into $H^{-2}(\R)$. Since the two agree on the dense subset $C_c^\infty(\R;\R)$ by Lemma~\ref{lem:mkdv-smooth-identification}, they agree everywhere: the identity \eqref{eq:identify-AKNS-h-vector-field} holds for every real $u\in L^2(\R)$.
\end{lemma}}

The proof occupies four steps. The first three establish the continuity of the two maps in turn, the first being the substantial one because $V_{h_\sigma}$ is assembled from the diagonal Green functions through several products in which one factor lies in a negative Sobolev space. The fourth applies density to the fixed spatial function $u$ inside the two continuous maps.

\begin{proof}
Let $u_n\to u$ in $L^2(\R;\R)$ and keep the notation of Lemma~\ref{lem:mkdv-AKNS-continuity}.

\proofstep{Step 1: convergence of $V_-,V_\gamma,V_+$} Set
\begin{align}
 \label{eq:epsilon-n-def} \varepsilon_n :=\norm{u_n^2-u^2}_{L^1} \leq(\norm{u_n}_2+\norm{u}_2)\norm{u_n-u}_2 \rightarrow0 .
\end{align}
Every nonlinear term of \eqref{eq:AKNS-vector-fields}--\eqref{eq:AKNS-vector-p} is a product of $u^2$ or $u$ with a diagonal Green function, so it is handled by splitting the difference into the two factors. Using \eqref{eq:AKNS-approximation-fields} and the embedding $H^1\hookrightarrow L^\infty$,
\begin{align}
 \norm{u_n^2g_{-,n,x}-u^2g_{-,x}}_{L^1} &\leq \varepsilon_n\norm{g_{-,n,x}}_\infty +\norm{u}_2^2\norm{g_{-,n,x}-g_{-,x}}_\infty \rightarrow0,\label{eq:Vm-product-limit}\\
 \norm{u_n^2(1+\gamma_n)-u^2(1+\gamma)}_{L^1} &\leq \varepsilon_n\norm{1+\gamma_n}_\infty +\norm{u}_2^2\norm{\gamma_n-\gamma}_\infty \rightarrow0,\label{eq:VA-first-product-limit}\\
 \norm{u_ng_{-,n}-ug_-}_2 &\leq \norm{u_n-u}_2\norm{g_{-,n}}_\infty +\norm{u}_2\norm{g_{-,n}-g_-}_\infty \rightarrow0,\label{eq:VA-second-product-limit}\\
 \norm{u_n^2g_{+,n}-u^2g_+}_{L^1} &\leq \varepsilon_n\norm{g_{+,n}}_\infty +\norm{u}_2^2\norm{g_{+,n}-g_+}_\infty \rightarrow0.\label{eq:Vp-product-limit}
\end{align}
The linear terms of \eqref{eq:AKNS-vector-fields}--\eqref{eq:AKNS-vector-p} converge because $\gamma_n\to\gamma$, $g_{+,n}\to g_+$ in $H^1$ and $g_{-,n}\to g_-$ in $H^2$. Feeding these into \eqref{eq:AKNS-vector-fields}--\eqref{eq:AKNS-vector-p} through $L^1\hookrightarrow H^{-1}$ and the continuous maps $\partial_x^3\colon H^2\to H^{-1}$, $\partial_x^3\colon H^1\to H^{-2}$ and $\partial_x\colon H^{-1}\to H^{-2}$, we obtain
\begin{align}
 \label{eq:AKNS-vector-field-approximation} \begin{aligned}
 V_{-,n}&\rightarrow V_- \qquad\text{in }H^{-1},\\
 V_{\gamma,n}&\rightarrow V_\gamma \quad\text{and}\quad V_{+,n}\rightarrow V_+ \qquad\text{in }H^{-2}.
 \end{aligned}
\end{align}

\proofstep{Step 2: the first map} It remains to pass through the two multipliers in \eqref{eq:AKNS-h-vector-field}. The coefficient of $V_-$ there is long and occurs repeatedly below, so we abbreviate
\begin{align}
 \label{eq:C-sigma-def} C_\sigma:=(1+\gamma-\sigma g_+)\,g_-D^{-2}, \qquad C_{\sigma,n}:=(1+\gamma_n-\sigma g_{+,n})\,g_{-,n}D_n^{-2} .
\end{align}
Since $\gamma_n\to\gamma$ and $g_{+,n}\to g_+$ in $H^1$, $g_{-,n}\to g_-$ in $H^2$, and \eqref{eq:AKNS-reciprocal-approximation} holds, the $H^1$ algebra estimate gives
\begin{align}
 \label{eq:AKNS-C-approximation} C_{\sigma,n}\rightarrow C_\sigma \quad\text{in }H^1 .
\end{align}
Now split the difference $V_{h_\sigma,n}-V_{h_\sigma}$ into three groups, according to whether the multiplier or one of $V_-,V_\gamma,V_+$ is being varied. Combining $H^2\cdot H^{-2}\subset H^{-2}$ from \eqref{eq:Htwo-Hminustwo-product}, $H^1\cdot H^{-1}\subset H^{-1}$ from Lemma~\ref{lem:products}\,\textup{(ii)}, and \eqref{eq:AKNS-vector-field-approximation} with \eqref{eq:AKNS-C-approximation}, we obtain
\begin{align}
 \label{eq:Vh-approximation} \norm{V_{h_\sigma,n}-V_{h_\sigma}}_{H^{-2}} &\lesssim \norm{D_n^{-1}-D^{-1}}_{H^2} \bigl(\norm{V_{\gamma,n}}_{H^{-2}}+\norm{V_{+,n}}_{H^{-2}}\bigr)\\
 &\quad+ \bigl(1+\norm{D^{-1}-1}_{H^2}\bigr) \bigl(\norm{V_{\gamma,n}-V_\gamma}_{H^{-2}} +\norm{V_{+,n}-V_+}_{H^{-2}}\bigr) \notag\\
 &\quad+ \norm{C_{\sigma,n}-C_\sigma}_{H^1}\norm{V_{-,n}}_{H^{-1}} +\norm{C_\sigma}_{H^1}\norm{V_{-,n}-V_-}_{H^{-1}} \rightarrow0,\notag
\end{align}
where the $H^{-1}$ products on the last line are viewed in $H^{-2}$ through the embedding $H^{-1}\hookrightarrow H^{-2}$, and where the first two groups are bounded because \eqref{eq:AKNS-vector-field-approximation} makes $\norm{V_{\gamma,n}}_{H^{-2}}$, $\norm{V_{+,n}}_{H^{-2}}$ and $\norm{V_{-,n}}_{H^{-1}}$ convergent, hence bounded. This is the continuity of the first map in \eqref{eq:two-static-maps}.

\proofstep{Step 3: the second map} By the Lipschitz bound \eqref{eq:two-Miura-Lipschitz} of Lemma~\ref{lem:two-Miura-static}, $q_{\sigma,n}\to q_\sigma$ in $H^{-1}$. From the second expression in \eqref{eq:reciprocal-through-m}, together with \eqref{eq:AKNS-approximation-fields} and \eqref{eq:AKNS-reciprocal-approximation}, we also have
\begin{align}
 \label{eq:h-sigma-approximation} h_{\sigma,n}-h_\sigma\rightarrow0 \quad\text{in }H^1 .
\end{align}
The affine $H^1$ multiplier estimate of Lemma~\ref{lem:products}\,\textup{(ii)} therefore gives
\begin{align}
 \label{eq:hQ-approximation} h_{\sigma,n}q_{\sigma,n} \rightarrow h_\sigma q_\sigma \quad\text{in }H^{-1} .
\end{align}
Inserting \eqref{eq:h-sigma-approximation} and \eqref{eq:hQ-approximation} into the grouped flux formula
\begin{align}
 \label{eq:mkdv-flux-continuity} \J(q_\sigma) =2h_\sigma q_\sigma -4\kappa^2(h_\sigma-\kappa),
\end{align}
which is \eqref{eq:J-2hq} at the parameter $\kappa$, gives $\J(q_{\sigma,n})\to\J(q_\sigma)$ in $H^{-1}$, and hence $\partial_x\J(q_{\sigma,n})\to\partial_x\J(q_\sigma)$ in $H^{-2}$. This is the continuity of the second map.

\proofstep{Step 4: the density argument}
Let $u\in L^2(\R;\R)$ be arbitrary and choose real $u_n\in C_c^\infty(\R)$ with $u_n\to u$ in $L^2$, which is possible because $C_c^\infty(\R;\R)$ is dense in $L^2(\R;\R)$. Steps~1--3 were proved for an arbitrary convergent sequence, so the notation is consistent and they may be applied to this particular one. Each $u_n$ is smooth and compactly supported, so Lemma~\ref{lem:mkdv-smooth-identification} applies to it and gives
\begin{align}
 \label{eq:identify-along-approximants} V_{h_\sigma}(u_n)=\partial_x\J\bigl(\cM_\sigma(u_n)\bigr) \quad\text{in }H^{-2}(\R)
\end{align}
for every $n$. Step~2 sends the left side of \eqref{eq:identify-along-approximants} to $V_{h_\sigma}(u)$ in $H^{-2}$, and Step~3 sends the right side to $\partial_x\J(\cM_\sigma(u))$ in the same space. A sequence in the Hausdorff space $H^{-2}(\R)$ has at most one limit, so the two limits agree, which is \eqref{eq:identify-AKNS-h-vector-field} at the given $u$.

Since $u$ was an arbitrary real element of $L^2(\R)$, the identity holds at every such function. It therefore applies to $u(t)$ at each fixed $t$ along a continuous $L^2$ trajectory, as required by Proposition~\ref{prop:mkdv-transfer}.
\end{proof}

%%%%%%%%%%%%%%%%%%%%%%%%%%%%%%%%%%%%%%%%%%%%%%%%%%%%%%%%%%%%%%%%%%%%%%%%%%%%%%%%%%

The lemmas of the last two subsections now combine.

{\begin{proposition}
\label{prop:mkdv-transfer}
Let $u\in C(I;L^2(\R;\R))$ be an mKdV Green solution in the sense of Definition~\ref{def:mkdv-green}, and let $q_\sigma=\cM_\sigma(u)$, $\sigma\in\{+1,-1\}$, be its two Miura transforms \eqref{eq:intro-miura-partners}. Then $q_\sigma\in C(I;H^{-1}(\R;\R))$, and both $q_+$ and $q_-$ are KdV Green solutions in the reciprocal sense of Definition~\ref{def:green}, with a common threshold $\lambda_*$ that depends only on the two compact trajectories $q_\pm(I)$ and is exhibited in \eqref{eq:mkdv-common-threshold} below. More precisely, for every $\kappa\geq1$ and all $s,t\in I$, the reciprocal evolution identity is
\begin{align}
 \label{eq:mkdv-to-reciprocal-law} h_{\sigma,\kappa}(t)-h_{\sigma,\kappa}(s) =\partial_x\int_s^t \J_\kappa(q_\sigma(t'))\dd t' \quad\text{in }H^{-2}(\R),
\end{align}
where $h_{\sigma,\kappa}$ is defined by \eqref{eq:mkdv-target-coordinate} and $\J_\kappa$ is the flux \eqref{eq:intro-green-flux}.
\end{proposition}}

\begin{proof}
That $q_\sigma\in C(I;H^{-1})$ follows by composing the Lipschitz bound \eqref{eq:two-Miura-Lipschitz} of Lemma~\ref{lem:two-Miura-static} with $u\in C(I;L^2)$. Fix $\kappa\geq1$ and $\sigma$.

Lemma~\ref{lem:mkdv-time-regularity} gives the integrated law \eqref{eq:weak-h-evolution} with integrand $V_{h_\sigma}$. Since $u(t')\in L^2$ for each $t'\in I$, Lemma~\ref{lem:mkdv-map-continuity} identifies $V_{h_\sigma}(u(t'))$ with $\partial_x\J(q_\sigma(t'))$ at each such $t'$. Substituting this identification into \eqref{eq:weak-h-evolution} gives
\begin{align}
 \label{eq:mkdv-law-before-swap} h_\sigma(t)-h_\sigma(s) =\int_s^t\partial_x\J(q_\sigma(t'))\dd t' \quad\text{in }H^{-2} .
\end{align}
The integrand of \eqref{eq:mkdv-law-before-swap} is continuous into $H^{-2}$, so its integral is a Bochner integral in that space. The flux $\J(q_\sigma(t'))$ is continuous into $H^{-1}$, as verified below. The bounded linear map $\partial_x\colon H^{-1}\to H^{-2}$ therefore commutes with its Bochner integral:
\begin{align}
 \label{eq:bochner-commutation} \underbrace{\int_s^t\partial_x\J(q_\sigma(t'))\dd t'}_{\text{Bochner integral in }H^{-2}} =\partial_x\underbrace{\int_s^t\J(q_\sigma(t'))\dd t'}_{\text{Bochner integral in }H^{-1}} .
\end{align}
Here $t'\mapsto\J(q_\sigma(t'))$ is continuous into $H^{-1}$ for every $\kappa\geq1$: by \eqref{eq:mkdv-flux-continuity} the flux is built from $h_\sigma-\kappa\in C(I;H^1)$ and $q_\sigma\in C(I;H^{-1})$ through the multiplication of Lemma~\ref{lem:products}\,\textup{(ii)}, and these are supplied for every $\kappa\geq1$ by Lemmas~\ref{lem:mkdv-static-reciprocal} and~\ref{lem:two-Miura-static}. Above a large threshold, part~\textup{(i)} of Lemma~\ref{lem:microscopic} gives the same conclusion directly. Combining \eqref{eq:mkdv-law-before-swap} with \eqref{eq:bochner-commutation} proves \eqref{eq:mkdv-to-reciprocal-law}.

It remains to place $q_\sigma$ in the class of Definition~\ref{def:green}. Nonnegativity of $L_{q_\sigma}$ puts $-\kappa^2$ in its resolvent set for every $\kappa\geq1$, so that $L_{q_\sigma}+\kappa^2$ is boundedly invertible by \eqref{eq:negative-spectral-point}. Admissibility in the sense of Proposition~\ref{prop:static} also requires the quantitative bounds \eqref{eq:g-quantitative}, as explained in Remark~\ref{rem:positivity-not-admissibility}. The images $q_+(I)$ and $q_-(I)$ are compact in $H^{-1}$, being continuous images of the compact interval $I$, so Proposition~\ref{prop:static} supplies finite thresholds $\lambda_0(q_+(I))$ and $\lambda_0(q_-(I))$. Choose
\begin{align}
 \label{eq:mkdv-common-threshold} K_u:=\max\bigl\{1,\ \lambda_0(q_+(I)),\ \lambda_0(q_-(I))\bigr\} .
\end{align}
For every $\kappa\geq K_u$ the identity \eqref{eq:mkdv-to-reciprocal-law} holds, as just proved, and every hypothesis and convention of Definition~\ref{def:green} is in force. Hence both $q_+$ and $q_-$ are KdV Green solutions, with threshold $\lambda_*=K_u$. The identity itself remains valid for every $\kappa\geq1$.
\end{proof}

{\begin{remark}
\label{rem:full-mkdv-law}
Forlano--Killip--Vi\c{s}an \cite[equation~(7.20)]{Forlano.Killip.Visan2024} use the scalar diagonal-$g$ law, and explain before their Section~7.1 its distinction from the Killip--Murphy--Vi\c{s}an \cite{Killip.Murphy.Visan2020} Green law. In the present notation, the latter is the reciprocal identity used in Definition~\ref{def:green}. The multiplier obstruction to \eqref{eq:forbidden-reciprocal-derivative} is described in Section~\ref{sec:mkdv-notion}. The proof above uses both mKdV identities and the $H^2$ function $g_-$ through \eqref{eq:reciprocal-through-m}. Proposition~\ref{prop:auto-smoothing} then gives smoothing for the two KdV Green solutions obtained in Proposition~\ref{prop:mkdv-transfer}.
\end{remark}}

%%%%%%%%%%%%%%%%%%%%%%%%%%%%%%%%%%%%%%%%%%%%%%%%%%%%%%%%%%%%%%%%%%%%%%%%%%%%%%%%%%
\subsection{Canonical Solutions, Uniqueness, and Automatic Smoothing}
\label{sec:mkdv-uniqueness}
%%%%%%%%%%%%%%%%%%%%%%%%%%%%%%%%%%%%%%%%%%%%%%%%%%%%%%%%%%%%%%%%%%%%%%%%%%%%%%%%%%

We first verify nonemptiness of the intrinsic class. Let $\Phi_{\mathrm{mKdV}}\colon\R\times L^2(\R)\to L^2(\R)$ be the jointly continuous extension of the Schwartz defocusing mKdV flow to $L^2(\R)$ supplied by \cite[Theorem~1.1]{Harrop-Griffiths.Killip.Visan2020}. Indeed, taking $s=0$ and the upper, defocusing sign in that theorem gives the unique jointly continuous extension of the Schwartz solution map $\R\times L^2(\R)\to L^2(\R)$. Restricting to real initial data yields the flow for \eqref{eq:defocusing-mkdv}: real Schwartz data have real classical solutions, and the real functions form a closed subspace of $L^2(\R)$, so the extension preserves realness. The comparison of arbitrary Green solutions with this continuous extension is proved below.

{\begin{lemma}
\label{lem:canonical-mkdv-green}
Let $u_0\in L^2(\R;\R)$ and let $u(t):=\Phi_{\mathrm{mKdV}}(t,u_0)$ be the canonical trajectory through it. Then $u$ is an mKdV Green solution in the sense of Definition~\ref{def:mkdv-green} on every compact time interval. This gives a Green solution through every $L^2$ datum and permits comparison with the canonical flow in the proof of Theorem~\ref{thm:mkdv-main}.
\end{lemma}}

\begin{proof}
The two identities of Definition~\ref{def:mkdv-green} are classical along Schwartz trajectories. The content of the lemma is that both sides are continuous in the datum, so that the identities survive the passage to a general $L^2$ datum. Choose real $u_{0,n}\in\cS(\R)$ with $u_{0,n}\to u_0$ in $L^2$, and set $u_n(t)=\Phi_{\mathrm{mKdV}}(t,u_{0,n})$. Joint continuity of the flow gives
\begin{align}
 \label{eq:canonical-mkdv-green-a1} u_n\rightarrow u\quad\text{in }C(I;L^2)
\end{align}
on every compact interval $I$. For each fixed $\kappa\geq1$, Forlano--Killip--Vi\c{s}an \cite[Proposition~3.7]{Forlano.Killip.Visan2024} then give, uniformly in $t\in I$,
\begin{align}
 \label{eq:canonical-mkdv-green-a2} \begin{aligned}
 \gamma_\kappa(u_n)&\rightarrow\gamma_\kappa(u) \quad\text{and}\quad g_{+,\kappa}(u_n)\rightarrow g_{+,\kappa}(u) \qquad\text{in }H^1,\\
 g_{-,\kappa}(u_n)&\rightarrow g_{-,\kappa}(u) \qquad\text{in }H^2.
 \end{aligned}
\end{align}
For the smooth trajectories, equations~(3.73) and~(3.74) of the same reference, which hold classically along Schwartz solutions as in Step~1 of the proof of Lemma~\ref{lem:mkdv-smooth-identification}, give the differential forms of \eqref{eq:mkdv-green-m-recalled} and \eqref{eq:mkdv-green-gamma-recalled}, which integrate to the stated identities for every $\kappa\geq1$. The product estimates of Step~1 in the proof of Lemma~\ref{lem:mkdv-map-continuity}, whose bounds involve only $\norm{u_n(t)-u(t)}_2$, $\norm{u(t)}_2$ and the differences in \eqref{eq:canonical-mkdv-green-a2}, show that their right sides converge respectively in $C(I;H^{-1})$ and $C(I;H^{-2})$, uniformly in $t\in I$ by \eqref{eq:canonical-mkdv-green-a1}--\eqref{eq:canonical-mkdv-green-a2}, while the left sides converge by \eqref{eq:canonical-mkdv-green-a2}. Passing to the limit in the integrated identities, for all $s,t\in I$, proves the claim.
\end{proof}

\begin{proof}[Proof of Theorem~\ref{thm:mkdv-main}]
The proof consumes Theorem~\ref{thm:main} twice, once for each Miura sign, and then recovers $u$ from its two images.

\proofstep{Step 1: uniqueness}
For $j\in\{1,2\}$ and $\sigma\in\{+1,-1\}$, set
\begin{align}
 \label{eq:canonical-mkdv-green-a3} q_{j,\sigma}:=\cM_\sigma(u_j).
\end{align}
Proposition~\ref{prop:mkdv-transfer} makes all four curves KdV Green solutions in $C(I;H^{-1})$, with thresholds $K_{u_1}$ and $K_{u_2}$. Step~1 in the proof of Theorem~\ref{thm:main} chooses a common parameter above both thresholds. Equality of the initial mKdV data implies $q_{1,\sigma}(0)=q_{2,\sigma}(0)$. Applying Theorem~\ref{thm:main} separately for the two signs yields
\begin{align}
 \label{eq:both-Miura-images-equal} q_{1,+}=q_{2,+}, \qquad q_{1,-}=q_{2,-} \quad\text{in }C(I;H^{-1}).
\end{align}
Subtracting these identities and using \eqref{eq:two-Miura-recovery} gives
\begin{align}
 \label{eq:canonical-mkdv-green-a4} \partial_x(u_1(t)-u_2(t))=0 \quad\text{in }\mathcal D'(\R)
\end{align}
for every $t\in I$. A distribution with zero derivative is spatially constant. The only such distribution represented by an $L^2(\R)$ function is zero. So $u_1(t)=u_2(t)$ for every $t\in I$.

\proofstep{Step 2: identification with the canonical flow}
Lemma~\ref{lem:canonical-mkdv-green} supplies a canonical Green solution at every initial datum. We compare it with the given solution using Step~1.

Let $u$ be an arbitrary mKdV Green solution on $I$ and set $u_0:=u(0)\in L^2(\R;\R)$. By that lemma the canonical trajectory
\begin{align}
 \label{eq:canonical-competitor} \widetilde u(t):=\Phi_{\mathrm{mKdV}}(t,u_0) , \qquad t\in I ,
\end{align}
is itself an mKdV Green solution on $I$, and it satisfies $\widetilde u(0)=u_0=u(0)$. The pair $(u,\widetilde u)$ therefore meets the hypothesis of Step~1, which gives $u=\widetilde u$ on $I$. As $u$ was arbitrary, every mKdV Green solution coincides with the canonical trajectory through its own datum. This is the asserted equality with $\Phi_{\mathrm{mKdV}}$.

This completes the proof of Theorem~\ref{thm:mkdv-main}.
\end{proof}

\begin{proof}[Proof of Proposition~\ref{prop:mkdv-smoothing}]
The smoothing conclusions follow from the two Miura transforms of a single Green solution. Let $u$ be an mKdV Green solution and fix a compact $J\subset I$. By Proposition~\ref{prop:mkdv-transfer} and Proposition~\ref{prop:auto-smoothing}, the two distributions $q_\pm=\cM_\pm(u)$ admit representatives $F_\pm\in L^2_{\loc}(J\times\R)$, in the sense of \eqref{eq:representative-def}, satisfying the uniform cell bounds and translated-cell tails in \eqref{eq:q-Kato-bound}--\eqref{eq:q-Kato-tail}. As spacetime distributions on $J\times\R$,
\begin{align}
 \label{eq:Miura-partner-recovery} F_+-F_-=2u_x, \qquad F_++F_-=2u^2,
\end{align}
because $F_\pm$ realize $q_\pm(t)$ against every spatial test function for almost every $t\in J$, by \eqref{eq:representative-def}, while $q_+(t)-q_-(t)=2u_x(t)$ and $q_+(t)+q_-(t)=2u(t)^2$ hold in $H^{-1}(\R)$ at every $t$ by \eqref{eq:two-Miura-recovery}, and $u\in C(I;L^2)$ is itself a function in $L^2_{\loc}(J\times\R)$ whose distributional $x$-derivative is the curve $t\mapsto u_x(t)$. Thus $u_x=\frac12(F_+-F_-)$ and $u^2=\frac12(F_++F_-)$ as functions, so $u_x,u^2\in L^2_{\loc}$ and $u\in L^4_{\loc}$. Almost everywhere,
\begin{align}
 \label{eq:Miura-parallelogram} \abs{u_x}^2+\abs{u}^4 =\frac14\bigl(\abs{F_+-F_-}^2+\abs{F_++F_-}^2\bigr) =\frac12\bigl(\abs{F_+}^2+\abs{F_-}^2\bigr).
\end{align}
Both assertions of \eqref{eq:intro-mkdv-smoothing} follow by adding the corresponding conclusions for $F_+$ and $F_-$.
\end{proof}

{\begin{remark}
\label{rem:mkdv-scope}
The original Forlano--Killip--Vi\c{s}an \cite{Forlano.Killip.Visan2024} definition requires only $\langle x\rangle^{-b}u\in C_tL^2_x$ for some $b\geq0$. The cases $b>0$ include infinite-volume Gibbs trajectories that need not lie in $L^2(\R)$. The present result treats the unweighted case $u\in C(I;L^2)$. That hypothesis is what delivers the global mapping properties and regularity of the diagonal Green functions used here, and with them the membership of both Miura transforms in $H^{-1}(\R)$ to which Theorem~\ref{thm:main} applies. For fixed functions, the Miura transforms can belong to $H^{-1}$ under weaker decay assumptions. For example, the function $u(x)=(1+x^2)^{-1/5}$ lies outside $L^2(\R)$, yet $u^2$ and $u_x$ both lie in $L^2(\R)$, so that $u^2\pm u_x\in L^2(\R)\subset H^{-1}(\R)$, while $\langle x\rangle^{-b}u\in L^2(\R)$ for every $b>\frac1{10}$. To extend the result to weighted trajectories, one would need the corresponding mapping estimates and a uniqueness theorem for nondecaying KdV Green solutions.
\end{remark}}

{\begin{remark}
\label{rem:defocusing-sign}
The real factorizations \eqref{eq:Miura-form-factorization} correspond to the sign $+6u^2u_x$ in \eqref{eq:defocusing-mkdv}. For focusing real mKdV, the analogous Miura transforms are complex-valued, so such a transfer would require a uniqueness theorem for complex KdV Green solutions.
\end{remark}}

% \newpage

\bibliographystyle{siam}
\bibliography{Reference}

% \begin{thebibliography}{99}

% \bibitem[CKV24]{Chapouto.Killip.Visan2024} Andreia Chapouto, Rowan Killip, and Monica Vi\c{s}an, \emph{Bounded solutions of {K}d{V}: uniqueness and the loss of almost periodicity}, Duke Math. J. \textbf{173} (2024), no.~7, 1227--1267.

% \bibitem[FKV24]{Forlano.Killip.Visan2024} Justin Forlano, Rowan Killip, and Monica Vi\c{s}an, \emph{Invariant measures for {mKdV} and {KdV} in infinite volume}, preprint, arXiv:2401.04292, 2024.

% \bibitem[HKV24]{Harrop-Griffiths.Killip.Visan2020} Benjamin Harrop-Griffiths, Rowan Killip, and Monica Vi\c{s}an, \emph{Sharp well-posedness for the cubic {NLS} and m{K}d{V} in {$H^s(\mathbb R)$}}, Forum Math. Pi \textbf{12} (2024), Paper No.~e6, 86.

% \bibitem[KMV20]{Killip.Murphy.Visan2020} Rowan Killip, Jason Murphy, and Monica Vi\c{s}an, \emph{Invariance of white noise for {K}d{V} on the line}, Invent. Math. \textbf{222} (2020), no.~1, 203--282.

% \bibitem[KV19]{Killip.Visan2019} Rowan Killip and Monica Vi\c{s}an, \emph{Kd{V} is well-posed in {$H^{-1}$}}, Ann. of Math. (2) \textbf{190} (2019), no.~1, 249--305.

% \bibitem[KKL26]{Klaus.Koch.Liu2023} Friedrich Klaus, Herbert Koch, and Baoping Liu, \emph{Wellposedness for the {K}d{V} hierarchy}, preprint, arXiv:2309.12773v2, 2026.

% \bibitem[RS75]{Reed.Simon1975} Michael Reed and Barry Simon, \emph{Methods of modern mathematical physics. {II}. {F}ourier analysis, self-adjointness}, Academic Press, New York, 1975.

% \end{thebibliography}

\end{document}